\documentclass[11pt,a4paper]{article} 
\usepackage{a4wide}
\usepackage{amsmath,amsfonts}

\newtheorem{theorem}{Theorem}
\newtheorem{definition}{Definition}
\newtheorem{corollary}[theorem]{Corollary}
\newtheorem{proposition}{Proposition}
\newtheorem{lemma}{Lemma} 
\newtheorem{claim}{Claim}
\numberwithin{theorem}{section}
\numberwithin{claim}{section}
\numberwithin{equation}{section}
\numberwithin{lemma}{section}
\numberwithin{proposition}{section}
\def\R{\mathbb{R}}

\begin{document}
 
\title{Description of two soliton collision \\ for the quartic gKdV equation
\footnote{
This research was supported in part by the Agence Nationale de la Recherche (ANR ONDENONLIN).
}}
\author{Yvan Martel$^{(1)}$   and  Frank Merle$^{(2)}$ }
\date{(1) Universit\'e de Versailles Saint-Quentin-en-Yvelines,
 Math\'ematiques, \\
 45, av. des Etats-Unis,
 78035 Versailles cedex, France\\  
 martel@math.uvsq.fr\\
\quad \\ 
(2)
Universit\'e de Cergy-Pontoise, IHES and CNRS,
Math\'ematiques   \\
2, av. Adolphe Chauvin,
95302 Cergy-Pontoise cedex, France \\
 Frank.Merle@math.u-cergy.fr\\
}
\maketitle
 
\begin{abstract}
This paper concerns the subcritical gKdV equations
\begin{equation}\label{kdvabs}
    \partial_t u + \partial_x (\partial_x^2 u + u^p)=0
\end{equation}
for $p=2,$ $3$ and $4$. We mainly focus on the nonintegrable case $p=4$.
In \cite{MMcol2}, we will extend the main results to more general nonlinearities with
stable traveling waves. 

Equation \eqref{kdvabs} is known to have special solutions of the type
$u(t,x)=Q_{c_0}(x-x_0-c_0 t)$, called solitons. 
The general problem is the following: one knows the existence of solutions of the equation
which behave as $t\to -\infty$ like
\begin{equation}\label{absun}
  u(t,x)=Q_{c_1}(x{-}x_1{-}c_1 t) + Q_{c_2}(x{-}x_2{-}c_2 t) + \eta(t,x),
\end{equation}
where $c_1>c_2$ and  $\eta(t)$ is a dispersion term small in the energy space $H^1$ with respect to $Q_{c_1}$, $Q_{c_2}$ (see \cite{MMT}, \cite{Martel}).
From the Physics point of view, the two solitons $Q_{c_1}$ and $Q_{c_2}$ have collide at some time $t_0$.
Can one understand the collision and determine what happens after the collision?
In nonlinear analysis, except for some completely integrable equations,
these questions are completely open. 

In this paper, we introduce a new framework to understand these problems 
for \eqref{kdvabs} in the case
$c_2 \ll c_1$ (or equivalently, $\|Q_{c_2}\|_{H^1}\ll \|Q_{c_1}\|_{H^1}$) and  $\|\eta(t)\|_{H^1} \ll \|Q_{c_2}\|_{H^1}$, for $t$ close to $-\infty$.
The understanding of the collision region is based on explicit computations, in particular  
on the introduction of a new nonlinear ``basis'' which allows us
to  write and compute the solution up to any order of size.
After the collision, i.e. for $t\to +\infty$, computation in this basis is not valid anymore and 
we rely on analysis in the original space variable using refined asymptotic techniques from \cite{MMas2},
\cite{MM1} and \cite{MMnonlinearity}. 

First, this approach allows us to describe for all time solutions satisfying
\eqref{absun} for $t$ close to $-\infty$. In particular, we prove that the two solitons
survive the collision up to a correction of lower order, i.e. for all $t$ large, we have
\begin{equation}\label{absde}
  u(t,x)=Q_{\tilde c_1}(x{-}y_1(t)) + Q_{\tilde c_2}(x{-}y_2(t)) + \tilde \eta(t,x),
\end{equation}
where $\tilde c_1\sim c_1$, $\tilde c_2\sim c_2$  and $\|\tilde \eta(t)\|_{H^1}
\ll \|Q_{c_2}\|_{H^1}$.
From the explicit decomposition in the interaction region, we can describe precisely
the collision, in particular, we are able to compute explicitly the main orders of the resulting shifts on the solitons. 

For $p=2,$ $3$, we check that
our approach matches at the main orders classical results based on the inverse scattering transform.

For $p=4$, this description is completely new and we
point out the following surprizing points:
(a) the slower soliton survives the collision and is not destroyed; (b) the shifts $\Delta_1$ and $\Delta_2$ on $Q_{c_1}$ and ${Q_{c_2}}$ are explicit. In particular, the shift on $Q_{c_1}$ is negative and tends to $-\infty$ as $c_2/c_1\to 0$, which is in contrast with the integrable cases.

Second,  our analysis in the nonintegrable case $p=4$
proves that for a solution which is asymptotically a pure $2$-soliton solution at $-\infty$, i.e. $\|\eta(t)\|_{H^1}\to 0$ as $t\to -\infty$
in \eqref{absun},
a nonzero part of the energy transfers into dispersion during the collison, 
which means that  in \eqref{absde}, $\|\tilde \eta(t)\|_{H^1}\sim \delta_c>0$, for $t$ large.
Therefore   no pure $2$-soliton solution exists in this regime ($c_2\ll c_1$).
This is clearly in contrast with the integrable case for which explicit multi-soliton solutions exist. 

Nevertheless, we are able to exhibit new exceptional solutions  for $p=4$:
we prove that for all $c_1,$ $c_2>0$, $c_2\ll c_1$, for all $y_1,$ $y_2$, there exists a
solution $\varphi(t)$ such that
\begin{equation*}
		\varphi(t,x)=Q_{c_1}(x{-}y_1{-}c_1 t+\tfrac 12 \Delta_1) 
		+ Q_{c_2}(x{-}y_2{-}c_2 t+\tfrac 12 \Delta_2) + \eta(t,x), 
		\text{ for $t\ll -1$,}
\end{equation*}
\begin{equation*}
		\varphi(t,x)=Q_{c_1}(x{-}y_1{-}c_1 t-\tfrac 12 \Delta_1)
		+ Q_{c_2}(x{-}y_2{-}c_2 t-\tfrac 12 \Delta_2) + \eta(t,x), 
		\text{ for $t\gg 1$,}
\end{equation*}
where $\eta(t)$ converges to $0$ around the solitons as $t\to \pm\infty$.
These solutions are
natural extensions in the nonintegrable case of  the multi-solitons for the integrable case. 	
\end{abstract}

\section{Introduction}
We consider the generalized Korteweg-de Vries (gKdV) equations:
\begin{equation}\label{gkdv}
    \partial_t u + \partial_x (\partial_x^2 u + u^p)=0, \quad 
    x,t\in \mathbb{R},
\end{equation}
in the subcritical case, i.e. for $p=2$, $3$ or $4$. 
Our main results concern the nonintegrable case $p=4$.
An extension to the case of a general nonlinearity $f(u)$ for which
the traveling waves are stable is considered in \cite{MMcol2}.

It is well-known that the Cauchy problem for equation \eqref{gkdv} is globally well-posed in the energy space $H^1(\mathbb{R})$
(see Kenig, Ponce and Vega \cite{KPV}):  for any $u_0\in H^1(\mathbb{R})$, there exists
 a unique   solution $u(t)\in C(\mathbb{R},H^1(\mathbb{R}))$ of (\ref{gkdv}) 
with $u(0)=u_0$, uniformly bounded in $H^1(\mathbb{R})$.
Moreover, the following quantities are conserved (if they are well-defined):
\begin{equation}\label{masse}
  \int u(t)=\int u(0),\quad \int u^2(t)=\int u^2(0),
\end{equation}
\begin{equation}\label{energie}
  E(u(t))={\frac 1 2} \int u^2_x(t) -{\frac 1 {p+1}} \int u^{p+1}(t)
         ={\frac1 2} \int u_{x}^2(0) -{\frac 1 {p+1}} \int u^{p+1}(0).
\end{equation}
Recall that for $p=2,3,4$, global well posedness follows from local well posedness,  
(\ref{masse})--(\ref{energie})  and the 
Gagliardo--Nirenberg inequality:
$
\forall v\in H^1$, $
\int |v|^{p+1}\le C  \left(\int v^2\right)^{\frac {p+3} 4}\left(\int v_x^2\right)^{\frac {p-1} 4}.
$

Recall also that there exist explicit traveling wave solutions of \eqref{gkdv}.
Denote by $Q$ the unique even solution of
\begin{equation}\label{defQ}
    Q>0, \quad Q''+Q^p=Q, \quad Q\in H^1(\mathbb{R}) \quad \hbox{i.e.}\quad 
    Q(x)=\bigg( \frac{p+1}  {2\cosh^2\big(\frac {p-1} 2 \,x\big)}\bigg)^{\frac1 {p-1}},
\end{equation}
and, for any $c>0$, let
\begin{equation}
    Q_c(x)=c^{\frac 1{p-1}} Q(\sqrt{c} \,x)
    \quad \text{be solution of $Q_{c}''+Q_{c}^p=c\,Q_{c}$.}
\end{equation}
Then, for any $x_{0}\in \mathbb{R}$, $c>0$, the functions
$
R_{c,x_{0}}(t,x)=Q_{c}(x-x_{0}-ct)
$
are solutions of \eqref{gkdv}, called \emph{solitons}.
These solutions have been intensively studied especially in the integrable cases, i.e. $p=2$
and $p=3$ in equation \eqref{gkdv}.

\subsection{Known results on soliton and multi-soliton solutions}

\noindent\emph{a. Integrable case $p=2,3$: N-solitons for the KdV and mKdV equations.}

Pioneering works of Fermi, Pasta and Ulam \cite{FPU} and Zabusky and Kruskal  \cite{KZ} have exhibited from the numerical point of view
remarkable phenomena related to soliton collision. Then, Lax (\cite{LAX1}) has developed a mathematical framework to study these problems,
known now as complete integrability. Many other developements
appeared, such as the inverse scattering transform (for a review on this theory, we refer  for example to
Miura \cite{Miura}).

This nonlinear transformation led to one of the most striking property of the KdV and mKdV equations 
which is the existence of pure
$N$-soliton solutions (Hirota \cite{HIROTA}). Namely, let $p=2$ or $p=3$,
and let $c_1>\ldots>c_N>0$, $\delta_1,\ldots,\delta_N\in \mathbb{R}$. There exists an explicit multi-soliton solution $U(t,x)$ of \eqref{gkdv} that satisfies
	\[ 	\biggr\|U(t,x)- \sum_{j=1}^N Q_{c_j}(.-c_jt-\delta_j)\biggr\|_{H^1} \mathop{\longrightarrow}_{t\to -\infty} 0,\quad
	\biggl\|U(t,x)- \sum_{j=1}^N Q_{c_j}(.-c_jt-\delta_j')\biggr\|_{H^1}
	\mathop{\longrightarrow}_{t\to +\infty} 0,
\]
for some $\delta_j'$ such that the shifts $\Delta_j=\delta_j'-\delta_j$ depends on the $(c_k)$.
Recall that  explicit formulas for such solutions were derived using the inverse scattering transform.
For example, the following function $U_{1,c}$, solution of \eqref{gkdv} with $p=2$,
is a 2-soliton solution ($0<c<1$):
	\begin{equation}\label{NSOL}
	 U_{1,c}(t,x)=6 \frac {
	\partial^2}{
	\partial x^2} \log\bigl(1+e^{x-t}+e^{\sqrt{c}(x-ct)}+\alpha e^{x-t}e^{\sqrt{c}(x-ct)}\bigr) \quad \textrm{with}\quad  \alpha=\left(\frac {1-\sqrt{c}}{1+\sqrt{c}}\right)^2.\end{equation}
The $N$-solitons are fundamental in studying the properties of general solutions of the KdV equation because of the following 
(Kruskal \cite{KRUSKAL}, Eckhaus and Schuur \cite{EcSc}, \cite{Schuur},
Cohen \cite{Cohen}):

\medskip
 
\noindent\textbf{Decomposition property (\cite{EcSc}, \cite{Schuur},
 \cite{Cohen})}
\textit{Let $u(t)$ be a solution of 
\eqref{gkdv} with $p=2$. Suppose that $u(0)\in C^4(\mathbb{R})$ satisfies for 
$k\in \{0,...,4 \}$, $\forall x\in \mathbb{R}$,
$\left|( {\partial^k u}/{\partial x^k} )(0,x)  \right|\leq    C/(1+{|x|^{10}}).$}
\textit{Then, there exist  $N\in \mathbb{N}$, $x_1,\ldots,x_N$ and $c_1>\ldots>c_N>0$ such that for all $x>0$,}
\[ u(t,x)-\sum_{j=1}^N Q_{c_j}(x-x_j-c_jt) \to 0 \quad \text{as $t\to +\infty$}. \]

\medskip

This result means that the asymptotic behavior in large time of any sufficiently regular and decaying solution is governed by a finite number of solitons.

\medskip

\noindent\emph{b. PDE results for the subcritical generalized  KdV equations $(p=2,\,3,\,4)$.}

First, we recall the following well-known orbital stability result.	
	
\medskip
\noindent \textbf{Stability of  soliton for the gKdV equation (\cite{Be}, \cite{Bo}, \cite{CL}, \cite{We2})} \textit{ Let $1<p<5$. Let $u(t)$ be an $H^1$ solution of the gKdV equation \eqref{gkdv}. For all $\epsilon>0$, there exists $\delta>0,$ such that if $\|u(0)-Q\|_{H^1}\le \delta$, then for all $t\in \R$, there exists $\rho(t)$, such that} 
$
		\|u(t)-Q(.-\rho(t))\|_{H^1(\mathbb{R})}\leq \epsilon. 
$

\medskip
		
By invariance by scaling and translation of the gKdV equation, the result is the same for $Q_{c_0}(x-x_{0})$, for any $c_0>0$, $x_{0}\in \mathbb{R}$. The proof of this result only relies on the  conservation laws \eqref{masse}--\eqref{energie}, and  the variational characterization of $Q(x)$ (see \cite{CL}, \cite{We2}).

The family of solitons $(R_{c,x_{0}}(t,x))$  is actually asymptotically stable, for equation \eqref{gkdv} in the subcritical case: $p=2$, $3$ or $4$.

\medbreak

\noindent\textbf{Asymptotic stability for the gKdV equation (\cite{MM1}, \cite{MMnonlinearity})}
\emph{Let $u(t)$ be an $H^1$ solution of (\ref{gkdv}).
There exists $\alpha  >0$ such that if $\|u(0)-Q\|_{H^1}\le \alpha $, then
there exist $c^+$ with $|c^+-1|=O(\alpha)$ and a $C^1$ function $\rho:[0,+\infty)\rightarrow \R$ such that 
		\begin{equation}
			\label{convv}w(t,x)=u(t,x)-Q_{c^+}(x-\rho(t))\quad \hbox{satisfies}\quad \lim_{t\rightarrow +\infty}\|w(t)\|_{H^1(x> \frac 1{10} t)}=0.
		\end{equation}
		Moreover, $\lim_{t\rightarrow +\infty} \frac{d\rho}{dt}(t)=c^+$. }

\medskip
		
	This result means that by taking $\alpha$ small enough, we know the   behavior of $u(t)$ as $t\to +\infty$ in the space time region $x>\frac 1 {10}t$ (in fact the result can be extended to the region $x>\beta t$, for any $\beta>0$, for $\alpha$ small). 
 This region of convergence  is in some some sharp since there exist solutions which behave asymptotically as $t\rightarrow +\infty$ as $Q(x-t)+Q_c(x-ct)$, where $c>0$ is arbitrarily small.   Note that the above theorem, proved only for 
$p=2,$ $3$ and $4$ in  \cite{MM1}, \cite{MMnonlinearity} also holds for \eqref{gkdv} with a general nonlinearity $f(u)$, see \cite{yvanSIAM}, \cite{MMas1}. The stability and asymptotic stability results above can be extended to the sum of $N$ solitons (and then to multi-solitons), when the various solitons are decoupled, see \cite{MMT}. Moreover,
assuming $\int_{x>0} x^2 u^2<+\infty$ implies that $\lim_{t\to +\infty} (\rho(t)-c^+ t)$
exists (see \cite{MMas2} and Section 4.2).

\medskip

Let us introduce the notion of asymptotic $N$-soliton solutions and pure $N$-soliton solution.
\begin{definition}\label{DEF1}
1. A  solution $u(t)$ of \eqref{gkdv} is an
         \emph{asymptotic $N$-soliton solution at $-\infty$} if
          there exist $c_1^->\ldots>c_N^->0$ and  $x_1^-(t),\ldots,x_N^-(t)$
          such that
          \begin{equation}\label{limit-}
             \lim_{t\to -\infty}
             \Bigl\|u(t)-\sum_{j=1}^N Q_{c_j^-}(.-x_j^-(t))\Bigr\|_{H^1(\mathbb{R})}=0.
          \end{equation}
2. A  solution $u(t)$ of \eqref{gkdv} is an
        \emph{asymptotic $N$-soliton solution at $+\infty$} 
         if there exist         
         $c_1^+>\ldots>c_N^+>0$ and $x_1^+(t),\ldots,x_N^+(t)$
         such that
         $
            \lim_{t\to +\infty}
\Bigl\|u(t)-\sum_{j=1}^N Q_{c_j^+}(.-x_j^+(t))\Bigr\|_{H^1(\mathbb{R})}=0.
         $\\
3. An $H^1$ solution $u(t)$ of \eqref{gkdv} is a 
          \emph{pure $N$-soliton solution} if $u(t)$ is an asymptotic $N$-soliton
          solution at both $+\infty$ and $-\infty$.
\end{definition}

We recall the following existence result.

\medskip 

\noindent\textbf{Asymptotic $N$-soliton solutions for the gKdV equation (\cite{Martel})}
		\emph{Let $p=2,$ $3$ or $4$. Let $N\geq 1$, $c_1>\ldots>c_N>0$, and $x_1,\ldots,x_N\in \R$. There exists a unique $H^1$ solution $U$  of (\ref{gkdv}) such that }
		\begin{equation}\label{limit0+}
            \lim_{t\to -\infty}
            \Bigl\|U(t)-\sum_{j=1}^N Q_{c_j }(.-x_j -c_j t)\Bigr\|_{H^1(\mathbb{R})}=0.
         \end{equation}

\medskip

See Proposition \ref{2SOL} for more properties on $U(t)$.
A similar statement holds true as $t\to +\infty$, since equation \eqref{gkdv} is invariant under the transformation $x\to -x$, $t\to -t$.

This result means that there exist asymptotic $N$-soliton solutions at $-\infty$ for $p=4$,
 similarly as in the   integrable cases $p=2,3$. However, for $p=4$, no information
is known concerning the collision phenomenon 
or the behavior as $t\to +\infty$ for such solutions.

Recent works have completed the above asymptotic results.
C\^ote \cite{Cote1}, \cite{Cote2} has proved, for $p=4$, $5$,
the existence of solutions satisfying
\begin{equation*}
      \lim_{t\to +\infty}
      \Bigl\|u(t,x)-\sum_{j=1}^N Q_{c_j^+}(x-c_j t -x_j ) -\mathcal{W}(t) v_0\Bigr\|_{H^1 }=0,
\end{equation*}
where $\mathcal{W}(t)$ is  the linear Airy group  and $v_0$ is a given function with suitable properties. 

Tao \cite{Tao} has established a well-posedness and scattering result (small data)
for \eqref{gkdv} with $p=4$ in the critical space $\dot H^{-1/6}(\mathbb{R})$. As a corollary of the estimates in \cite{Tao} and of the asymptotic stability result above, it is proved that if
 $u_0$ is close to $Q$ in $\dot H^{-1/6}\cap H^1$, then there exists $v_0\in \dot H^{-1/6}\cap H^1$
   such that
\begin{equation*}
      \lim_{t\to +\infty}
      \Bigl\|u(t,x)-Q_{c^+}(x-x(t)) -\mathcal{W}(t) v_0\Bigr\|_{H^1 }=0.
\end{equation*}

\subsection{Motivation of the problem}
We consider in this paper the problem of collision of two solitary waves. 

In the integrable case $p=2,3$, the explicit $2$-soliton solutions give a precise description of this phenomenon, and allow calculations of the interaction effects such as the resulting shifts on the solitons after the collision.

In the nonintegrable situation, the collision problem is an open question since the 70's.
Recall that  PDE theory was first related to existence and stability properties of
solitary waves. More recently, the focus has been put on trying to understand 
interaction between solitary waves and dispersion.
Finally, the results presented above are asymptotic results for large time, without description of the collision of solitons.

\smallskip

The problem of describing the collison of two traveling waves is  a general problem for nonlinear PDEs, which is completely open, except in some integrable situations where explicit formulas are known.
It is the simplest case of interaction between two nonlinear dynamics.
If one conjectures that any general solution (under suitable assumptions) decomposes as time goes to $+\infty$ as a sum of decoupled solitons and a dispersive part, as in the integrable case,  a   natural question is to try to relate the decomposition as $t\to +\infty$ to the one as 
$t\to -\infty$ by understanding the interaction of the different parts of the solution.

\smallskip

Apart from integrability theory,
this kind of problems have been studied since the 60's from both experimental and numerical points of view.

First, Fermi, Pasta and Ulam \cite{FPU}, Zabusky and Kruskal \cite{KZ} and Zabusky \cite{Z} have
introduced nonlinear systems and computed interaction of nonlinear objects by numerics.
Later, the theory of integrability justified these numerics as explained above.
Since then,  many other systems have been studied numerically.

There is also an extensive literature devoted to experiments on water tanks.
A key question  is whether or not the collision between two solitary waves is elastic
(equivalently, whether   the collision is pure or generates dispersion).
From experiments related to wave propagation in shalow water (see Weidman and Maxworthy \cite{WM},
Hammack et al. \cite{HHGY}, Craig et al. \cite{Craig}), it seems that   collisions are inelastic
but very close to be elastic, for solitary waves of different amplitude.

We recall some numerical works for equations of gKdV type.   Bona et al. \cite{BPS}, and
Kalisch and Bona \cite{KB}, studied numerically the problem of collision of two solitary waves for the Benjamin   and the BBM equations.
Shih \cite{SHIH} studied the case of the gKdV equation \eqref{gkdv} with some 
half-integer values of $p$.  Li and Sattinger \cite{LS} investigated the collision problem for the
case of the
Ion Acoustic Plasma equation, and Craig et al. \cite{Craig} report on numerics for the Euler equation with free surface. In all these works, the numerics match the experiments and show that, unlike for
the pure solitons of the integrable case, the collision of two solitary waves fails to be elastic
by a very small dispersion (difficult to see numerically).

\smallskip

Let us now review some more recent mathematical results related to these problems. 
First, Haragus and Sattinger \cite{HS} have studied perturbation of the KdV equation around the explicit
$N$-soliton solutions, in particular the invertibility of the linearized operator around these solutions.
Second, Mizumachi \cite{Mizu} for equation \eqref{gkdv} with $p=4$, has treated the case of two solitons
with close sizes, in a situation of repulsive interaction without collision 
(using scattering techniques).
Finally, 
the multi-soliton solutions
of the NLS (nonlinear Schr\"odinger) model, with special nonlinearity and under spectral assumptions (ruling out the existence of small solitary waves) have been studied by Perelman \cite{P} and
Rodnianski, Schlag and Soffer \cite{RSS}. Using Galilean invariance, speeds and sizes are independent
(in particular, high speed is possible for size one solitary waves). Thus, one can consider the case
where the collision has a negligeable effect on the solitary waves due to a very small time of interaction. In all these works, the interaction of two nonlinear objects in a non perturbative
case is not considered. In addition, up to now, no example of 
inelastic collision is known rigorously.

\subsection{Main results}
Our main results in this paper concern the problem of collision of two solitons for \eqref{gkdv}
in the (nonintegrable) case  $p=4$.
We consider the situation where one soliton, $Q_{c_1}$ is supposed to be large with respect to the other one, $Q_{c_2}$; thus we assume $c=c_2/c_1\ll 1$. 
This is not a perturbative setting, related to the integrable case or to a linearized equation.
In addition, the techniques of this paper can be applied in a general context for \eqref{gkdv}, $p=2,3,4$ or with a general nonlinearity $f(u)$ (see \cite{MMcol2}).
In this situation, we are able to compute the interaction term during the collision up to any order of $c$, which allows us to describe very precisely the collision phenomenon.

First, this approach allows us to prove that for $p=4$, there does not exist pure $2$-soliton solutions in the regime $c_2\ll c_1$:
 an asymptotic $2$-soliton solution at $-\infty$ cannot be 
  an asymptotic $2$-soliton solution  at $+\infty$.

\begin{theorem}[Non existence of a pure $2$-soliton solution for $p=4$]\label{NONEX4}
     Let $c_1>c_2>0$.
    There exists $\epsilon_{0}>0$ such that if $c=\frac{c_2}{c_1}<\epsilon_{0}$, then there exists no 
    pure $2$-soliton solution with speeds $c_1, \, c_2$ at $-\infty$. 

More precisely,
let $x_1^-,\ x_2^- \in \mathbb{R}$, and let $u(t)$ be the unique $H^1$ solution
    of \eqref{gkdv} such that
    \begin{equation*}
    		\lim_{t\to -\infty} \|u(t)-Q_{c_1}(.-x_1^- -c_1 t)-Q_{c_2}(.-x_2^- -c_2 t)\|_{H^1}=0.
    \end{equation*}
    Then, there exist $x_1^+,x_2^+$, $c_1^+>c_2^+>0$ and  $T_0,K>0$ such that
    \begin{equation*}
    w^+(t,x)=u(t,x)-Q_{c_1^+}(x-x_1^+-c_1^+ t)-Q_{c_2^+}(x-x_2^+-c_2^+ t) \quad
        \end{equation*}
        satisfies
       \begin{equation}\label{THUN1}
        \lim_{t\to +\infty} \|w^+(t)\|_{H^1(x>\frac 1{10} c_2 t)}=0,
   \end{equation}
        \begin{equation}\label{THUN2}
            \tfrac 1 K  c^{\frac {17} 6}\leq \frac {c_1^+}{c_1} - 1 \leq K  c^{\frac {11}{12}},\quad
   \tfrac 1K  c^{\frac {8} 3}  \leq 1- \frac {c_2^+}{c_2}  \leq  K c^{\frac 13},
      \end{equation}
        \begin{equation}\label{THUN3} 
     \tfrac 1 K    c_1^{\frac 7 {12}}  c^{\frac {17} {12}}\leq
    \|\partial_x w^+(t)\|_{L^2}+\sqrt{c_1c }\,  \|w^+(t)\|_{L^2}
    \leq K  c_1^{\frac 7 {12}} c^{\frac {11}{12}},\quad \text{for $t\geq T_0$}.
   \end{equation}
\end{theorem}
   Theorem \ref{NONEX4} confirms the common belief that the existence of pure $2$-soliton solutions,
    in particular, the elactic collision between two solitons, is a property which is specific to integrable models.
However, we observe that the $2$-soliton structure persists, in the sense that the slow soliton
is not destroyed by the collision, and remains approximately of the same size as $t\to +\infty$
(see also Remark 2  after the statement of Theorem \ref{EXIST4}).

Note that the size of $w^+(t)$ 
  measures the distance of the solution at $+\infty$ to a pure $2$-soliton solution. 
The bound below in \eqref{THUN3} is thus a qualitative version of nonexistence of a pure $2$-soliton solution. As a corollary of the proof, asymptotically in time, 
the minimal distance of any solution to a pure $2$ soliton solution at $+\infty$ or
at $-\infty$ is $K c^{\frac {17}{12}}$ at $\pm \infty$, in the
same sense as in \eqref{THUN3}.
We also see from \eqref{THUN2} 
how the speeds and the sizes of $Q_{c_1}$ and $Q_{c_2}$ are altered through the collision, 
the fast soliton accerelates while the slow soliton slows down.

Note that this is the first rigorous result describing a property of inelastic (but almost elastic)
collision, and thus a first illustration of nonintegrability of the equation from the dynamics of the solitary waves.

\medskip

\noindent\textbf{Remark 1.} Using the invariant $\int u(t)$ of equation \eqref{gkdv} in the framework of Theorem \ref{NONEX4}, one   proves that
$w^+(t)$ has to contain some dispersive part as $t\to +\infty$, in the sense that it does  not converge to a pure sum of small solitons, i.e.
$u(t)$ is not an asymptotic $N$-soliton solution at $+\infty$, for any $N\geq 1$
(see end of section 5.1). See also Remark 2 (4).

\medskip

In spite of the nonexistence result above, 
we prove  for $p=4$
the existence of exceptional solutions related to the $2$-soliton structure. These solutions
are the illustration of the persistence of the two solitons structure through the collision,
and provide a sharp description of the collision (conservation of the speeds and 
explicit shifts). This is a surprizing result for a nonintegrable equation.

\begin{theorem}[Existence of $2$-soliton like solutions for $p=4$]\label{EXIST4}
       Let $c_1>c_2>0$.
    There exists $\epsilon_0>0$ such that  
    if $c=\frac {c_2}{c_1}<\epsilon_0$, then there exist an $H^1$ solution
    $\varphi(t)$ of \eqref{gkdv}, $\Delta_1, \, \Delta_2\in \mathbb{R}$, satisfying, for all $t,x\in \mathbb{R}$,
    \begin{equation}\label{TH1Z}
   		\varphi(-t,-x)=\varphi(t ,x ),
		\end{equation}
and such that the following holds for $w^{\pm}(t):$
\begin{equation*}
w^{-}(t,x)=\varphi(t,x) - Q_{c_1}(x-c_1 t + \tfrac 12 \Delta_1) - Q_{c_2}(x-c_2 t + \tfrac 12 \Delta_2),
\end{equation*}
\begin{equation*}
w^{+}(t,x)=\varphi(t,x) - Q_{c_1}(x-c_1 t - \tfrac 12 \Delta_1) - Q_{c_2}(x-c_2 t - \tfrac 12 \Delta_2),
\end{equation*}
    \begin{enumerate}
        \item Asymptotic behavior at $\pm\infty$
                \begin{equation}\label{TH1A}
        \begin{split}
           & \lim_{t\to -\infty}
            \|w^-(t)\|
            _{H^1(x< \frac{c_2 t}{10})} =0,\quad  \lim_{t\to +\infty}
            \|w^+(t)\|
            _{H^1(x>\frac{c_2 t}{10})} =0,
        \end{split}
        \end{equation}
		where the shifts $\Delta_1$, $\Delta_2$ satisfy
		$\Delta_1<0$, $\Delta_2<0$ and 
		\begin{equation}\label{TH1B} 
			 \left|c_1^{\frac 12} \Delta_1 -    c^{-\frac 16}  \left({-} 2 \frac{\left(\int Q\right)^2}{\int Q^2} \right)\right|  +
             \left|c_1^{ \frac 12}\Delta_2-    \left(  \frac 13 \frac {\left(\int Q\right)^2}{\int Q^2} - \int Q^3\right) \right|\leq K {c^{\frac {1}{12}}}.
        \end{equation}\item Distance to the sum of two solitons : there exists $T_0>0$
        such that,
   \begin{equation}\label{TH1E}\begin{split}
                    				&		\tfrac 1K c_1^{\frac 7 {12}} c^{\frac{17}{12}}\leq 
		  \|\partial_x w^+(t)\|_{L^2}
		+\sqrt{c_1c}\,\|w^+(t)\|_{L^2} 
            \leq  K c_1^{\frac 7 {12}}   c^{\frac {17}{12}}, \quad\text{for all $t\geq T_0$.}
            \end{split}
\end{equation}
    \end{enumerate}
\end{theorem}

\noindent\textbf{Remark 2.}\quad 
1. From the stability result of one soliton (variational argument), it follows immediately that the  soliton $Q_{c_1}$ is preserved up to a certain order through the collision by a slow soliton $Q_{c_2}$. What is quite surprizing, and very similar to the integrable situation, is the fact that the second soliton, which is small, is also preserved by the collision (dynamical argument).  One could have expected the small soliton to be destroyed by such a collision.

Moreover, the solutions constructed in Theorem \ref{EXIST4} describe very precisely the effect of the collision on the two solitons,
since the speeds at $\pm \infty$ are the  same and since we have explicit formulas 
for the main order of the shifts on $Q_{c_1}$ and $Q_{c_2}$.
From the proof of Theorem~\ref{EXIST4}, the shifts are a consequence of the
collision and are observed in the relatively short period of time around $t=0$.

Concerning the  shifts, we point out two main differences with the integrable cases:

\smallskip

- The shift  $\Delta_1$ on  $Q_{c_1}$ and the shift $\Delta_2$ on $Q_{c_2}$
are both negative.

\smallskip

- The shift $\Delta_1\to -\infty$  as $c=c_2/c_1 \to 0$, which means that the effect of the soliton $Q_{c_2}$ on the trajectory of $Q_{c_1}$ becomes larger when $c_2/c_1$ is smaller
(note also that in this case the period of interaction is larger since the support of $Q_{c_2}$ becomes larger).

\smallskip

Both are new remarkable properties of the collision of two solitons for $p=4$.

\medskip

2.
First note that by the symmetry property of $\varphi(t)$ (see \eqref{TH1Z}), a statement similar
to \eqref{TH1E} for $w^-$ holds as $t\to -\infty$.
Now, let $w_{\rho_1,\rho_2}(t)=\varphi(t) - Q_{c_1}(.-\rho_1 ) - Q_{c_2}(.-\rho_2 )$, then from the proof of Theorem \ref{EXIST4}, we also have
        \begin{equation}\label{TH1Ebis}
                    						\tfrac 1K c_1^{\frac 7 {12}} c^{\frac{17}{12}}\leq 
						\inf_{\rho_1,\rho_2\in \R} 
        	 \left\{\|\partial_x w_{\rho_1,\rho_2}(t)\|_{L^2}
		+\sqrt{c_1c}\|w_{\rho_1,\rho_2}(t)\|_{L^2}\right\}
            \leq  K c_1^{\frac 7 {12}}   c^{\frac {17}{12}},
            \quad \text{for $|t|$ large,}
\end{equation}
            \begin{equation}\label{TH1F}
				\inf_{\rho_1,\rho_2\in \R}
						\left\{\|\partial_x w_{\rho_1,\rho_2}(t)\|_{L^2}
		+\sqrt{c_1}\|w_{\rho_1,\rho_2}(t)\|_{L^2}\right\}            \leq K   c_1^{\frac 7 {12}}   c^{\frac {1}{3}},\quad 			 \text{for all $t\in \mathbb{R}$.}           \end{equation}
Estimate   \eqref{TH1F} is sharp, indeed, at $t=0$, 
we have $$\inf_{\rho_1,\rho_2\in \R}
						\left\{\|\partial_x w_{\rho_1,\rho_2}(0)\|_{L^2}
		+\sqrt{c_1}\|w_{\rho_1,\rho_2}(0)\|_{L^2}\right\}\geq K_1 c_1^{\frac 7 {12}}   c^{\frac {1}{3}}.  $$  
		
For $p=4$, $\|Q_c\|_{L^2}=c^{\frac 1{12}} \|Q\|_{L^2}$,
 this is to compare with   \eqref{TH1Ebis}-\eqref{TH1E} giving sharp
estimates  of the distance of $\varphi(t)$ to
the sum of two solitons.

\medskip

3. By time and translation invariances, for all $x_1,x_2\in \mathbb{R}$, one 
derives from Theorem \ref{EXIST4} the existence of a solution $\varphi_{x_1,x_2}$ 
such that 
\begin{equation*}
        \begin{split}
           & \lim_{t\to -\infty}
            \|\varphi_{x_1,x_2}(t) - Q_{c_1}(.-c_1 t-x_1+\tfrac 12\Delta_1) - Q_{c_2}(.-c_2 t-x_2 +\tfrac 12\Delta_2)\|
            _{H^1(x< \frac{c_2 t}{10})} =0,\\
            & \lim_{t\to +\infty}
            \|\varphi_{x_1,x_2}(t) - Q_{c_1}(.-c_1 t-x_1 -\tfrac 12\Delta_1) - Q_{c_2}(.-c_2 t-x_2-\tfrac 12 \Delta_2)\|
            _{H^1(x>\frac{c_2 t}{10})} =0.
        \end{split}
        \end{equation*}
Moreover, there exist  an infinite number of solutions $\varphi(t)$ satisfying the conclusions of Theorem \ref{EXIST4} for given $c_1>c_2>0$, $x_1$, $x_2$. Indeed, it is enough to perturb the initial data
$\varphi(0)$ in a suitable way to obtain a solution with similar properties (see 
proof of Theorem \ref{EXIST4}).
Note finally that the solution $\varphi(t)$ which we have constructed belongs to 
$H^s$ for all $s\geq 0$.

\medskip

4. Remark 1 also applies to the solution $\varphi(t)$ constructed in Theorem \ref{EXIST4}, 
i.e. $\varphi(t)$ has some dispersive part as $t\to \pm \infty$.
Using Tao \cite{Tao} (specific for $p=4$), we should obtain some more information on the solution
since $\varphi(0)\in L^{\frac 32}(\mathbb{R})$. Indeed, 
$\varphi(t)$ is conjectured to satisfy, for some $v_0\in H^1 $,
\begin{equation}\label{RK1A}
        \begin{split}
           & \lim_{t\to -\infty}
            \|\varphi(t) - Q_{c_1}(.-c_1 t+\tfrac 12 \Delta_1 ) - Q_{c_2}(.-c_2 t
            +\tfrac 12 \Delta_2 )-\mathcal{W}(t) v_0\|
            _{H^1 } =0,\\
            & \lim_{t\to +\infty}
            \|\varphi(t) - Q_{c_1}(.-c_1 t -\tfrac 12\Delta_1) - Q_{c_2}(.-c_2 t-\tfrac 12\Delta_2)-\mathcal{W}(t)v_0\|
            _{H^1 } =0,
        \end{split}
        \end{equation}
where
$K_1   c_1^{\frac 7 {12}}  c^{\frac {17} {12}}\leq
    \|\partial_x v_0\|_{L^2}+\sqrt{c_1 c }\,  \|v_0\|_{L^2}
    \leq K_2 c_1^{\frac 7 {12}} c^{\frac {17}{12}}.
$

We conjecture that there exists a universal $v_0$, minimizer of a certain functional  related to energy quantities (for example $\int (\partial_x v_0)^2 + c_1\, c \int v_0^2$). 
This function $v_0$ should have additional special properties, such as smoothness and
exponential decay in space.

\medskip

5. Precise information concerning the solution $\varphi(t)$ at $t=0$ can be obtained from the proof of Theorem \ref{EXIST4}. See in particular Theorem \ref{APSOL4}.

\medskip

Finally, the behavior of such solutions  
is proved to be stable in $H^1$, which means that if
a solution $u(t)$ of \eqref{gkdv} is close at $t=0$ to the solution $\varphi$
constructed above, then $u(t)$ has a $2$-soliton structure for all time.

\begin{theorem}[Stability of the $2$-soliton structure for $p=4$]\label{STAB4}
   \quad Let $c_1>c_2>0$.
    Assume that $c=\frac{c_2}{c_1}<\epsilon_0$ is small enough and let $\varphi(t)$ be constructed
    in Theorem \ref{EXIST4}.
  Let   $u(t)$ be an $H^1$ solution of \eqref{gkdv} such that 
 for some $\delta>0$,    
    \begin{equation*} 
    \begin{split} 
	\|\partial_x u(0) - \partial_x  \varphi(0)\|_{L^2}+
	 \sqrt{c_1}\|u(0) -\varphi(0)\|_{L^2}\leq c_1^{\frac {7}{12}}   c^{\delta+\frac {7}{12}}.
	 \end{split}
    \end{equation*}
    Then,  there exist  $  \rho_1(t),\, \rho_2(t) \in \mathbb{R}$ and
    $c_1^\pm,\,c_2^\pm >0$ such that
    \begin{enumerate}
        \item Global in time stability,   $w(t,x)=u(t,x) - Q_{c_1}(x-  \rho_1(t))
- Q_{c_2}(x- \rho_2(t))$ satisfies
        \begin{align*}
				&  
            \|\partial_x w(t)\|_{L^2}+\sqrt{c_1} \|w(t)\|_{L^2}
\leq K c_1^{\frac 7{12}} ( c^{\delta+\frac 1{12}}+c^{\frac {1}{3}} ),
\quad \text{for all $t\in \mathbb{R},$}                     
        \end{align*}
        \item Asymptotic stability 
        \begin{equation*} 
         \lim_{t\to -\infty}
            \|u(t) - Q_{c_1^-}(. -  \rho_1(t)) - Q_{c_2^-}(.- \rho_2(t))\|
            _{H^1(x {<} \frac{c_1 t}{10})} =0,
        \end{equation*}
        \begin{equation*} 
         \lim_{t\to +\infty}
            \|u(t) - Q_{c_1^+}(. -  \rho_1(t)) - Q_{c_2+}(.- \rho_2(t))\|
            _{H^1(x > \frac{c_1 t}{10})} =0,
        \end{equation*}
        \begin{equation*}
			\left|\frac {c_1^\pm}{c_1} -1\right| \leq K c^{\frac 7{12}} (c^\delta + c^{\frac 13}),\quad
			\left|\frac {c_2^\pm}{c_2} -1\right| \leq K  (c^{\delta } + c^{\frac 13}).
		\end{equation*}
		               \end{enumerate}
\end{theorem}

\noindent\textbf{Remark 3.}\quad 
  Theorem \ref{STAB4} shows that the various properties exhibited in  Theorem \ref{EXIST4} 
are stable by perturbation of the initial data (during and after the collision).
This constructs in particular a large set of initial data having globally in time
a $2$-soliton structure (as for the integrable case). The stability property can also be proved
assuming $u(T_0)$ close to $\varphi(T_0)$ for some $T_0$
(see proof of Theorem \ref{STAB4}).

\medskip

The paper is organized as follows.
Section 2 is devoted to the collision region. We introduce a new method allowing to compute
a function $v(t)$ describing the collision which is a solution up to any order of $\frac {c_2}{c_1}$.
This part is mostly algebraic.

Section 3 is concerned with the recomposition of $v(t)$ after the collision.
We mainly focus on the case $p=4$.
For $p=2$, we only compare at the main orders the function $v(t)$ to the explicit $2$-soliton solutions.
The same could be done for $p=3$, but this case is omitted.

In Section 4 we present some analysis tool  necessary to relate the  function $v(t)$
to exact solutions of the equation. This part  mainly recalls asymptotic results from 
  \cite{MMas2}, which are refinements and generalizations of results in 
\cite{MM1}, \cite{MMT} and \cite{MMnonlinearity}.

Section 5 is devoted to the proofs of the main results, i.e.
Theorems \ref{NONEX4},  \ref{EXIST4} and \ref{STAB4}.

\section{Construction of an approximate $2$-soliton  of the gKdV eq.}

In the proof of the main results (Theorems 1, 2 and 3), we restrict ourselves to the case
$c_1=1$, $c_2=c$ small by a scaling argument. Therefore, in this section, we concentrate
on this case.

Let   $p=2,$ $3$, $4$ and define
\begin{equation}\label{defTc}
    T_c=c^{-\frac 12 -\frac 1{100}} \quad \text{and}\quad 
    q=\frac 1{p-1}-\frac 14.
\end{equation}
In this section, for any $n_0\in \mathbb{N}$, for $0<c<c_0$ small enough,
we construct a function $v_{n}(t,x)=v(t,x)$ which satisfies the following two properties.
\begin{itemize}
	\item $v(t,x-t)$ is a solution of the gKdV equation \eqref{gkdv} on $[-T_c,T_c]$ up to an error term of polynomial order $c^n$, 
	\begin{equation*}
       \forall t\in [-T_c,T_c],\quad
       \|\partial_{t} v + \partial_{x} (\partial_{x}^2 v - v + v^p)\|_{H^1(\mathbb{R})}\leq K(n) c^n.
    \end{equation*}
    \item $v(-T_c)$ and $v(+T_c)$ are at the main order the sum to two solitons $Q$ and $Q_c$
    respectively before and after their collision.
    \end{itemize}
The function $v(t)$  is the new fundamental object of this paper. Its existence and properties will lead to the main results stated in the Introduction.

\medskip

Our approach is to consider $c$ as a small parameter and look for such a function $v$ in terms of  expansions in powers of $c$, both in the functions and the space variables.
More precisely, the construction of the function $v(t,x)$ is related to the method of separation of variables: the variable $y$ of the large soliton $Q(y)$ is separated from the variable $y_c$ of the small soliton $Q_c(y_c)$.

First, we set
\begin{equation*}
    y_{c}=x+(1-c)t  \quad \text{and} \quad R_c(t,x) = Q_c(y_c);
\end{equation*}
note that $R_c(t)$ is then solution of 
$\partial_t R_c + \partial_x (\partial_x^2 R_c - R_c + R_c^p) =0$.

We look for a function $v(t,x)$ having the structure
\begin{equation}\label{defv} 
    v(t,x)=Q(y)+Q_{c}(y_{c})+W(t,x).
\end{equation}
We choose the function $W$ and the variable $y$ under the form of series.
Let $k_0 \geq 1$, $\ell_0\geq 0$ and define
$$
\Sigma_0=\{(k,\ell),~1\leq k\leq k_0,~0\leq \ell\leq \ell_0\}.
$$
For real unknown parameters $(a_{k,\ell})_{(k,\ell)\in \Sigma_0}$, we consider the variable $y$ of the form
\begin{equation*}
    y=x-\alpha(y_{c})=x-\alpha\bigl(x+(1-c)t\bigr) \quad \text{and}\quad     R(t,x)=Q(y),
\end{equation*}
where
\begin{equation}\label{defALPHA} 
    \alpha(s)=\int_{0}^s \beta(s') ds',\quad \text{and}\quad    
    \beta(s)=\sum_{(k,\ell)\in \Sigma_0} a_{k,\ell} \, c^\ell Q_{c}^k(s).
\end{equation}
The form of $W$ is
\begin{equation}\label{defW}
    W(t,x)=\sum_{(k,\ell)\in \Sigma_0} 
        c^\ell\left(Q_{c}^k(y_{c}) A_{k,\ell}(y)+(Q_{c}^k)'(y_{c}) B_{k,\ell}(y)\right),
\end{equation}
where the functions $A_{k,\ell}$, $B_{k,\ell}$, as the parameters $(a_{k,\ell})$, are yet undetermined.
Note that the functions $c^\ell Q_c^k$ and $c^\ell (Q_c^k)'$ used to define the series play the role of a set of nonlinear eigenfunctions for the interaction problem. Thus, the structure of $W$ will allow us to compute the interaction terms at any order of power of $c$.
Moreover, choosing the variable $y$ as above will allow us to understand the effect of the soliton $Q_c$ on the position of $Q$, that is, the shift phenomenon which appears through the interaction of two solitons.

\begin{theorem}[Construction of an approximate solution of the gKdV equation]\label{APSOL4} 
\ \\
Let $p=2$, $3$ or $4$. For all $k\geq 1$, $\ell\geq 0$, there exist
$a_{k,\ell}\in \mathbb{R}$ and $C^\infty$ functions  $A_{k,\ell}$, $B_{k,\ell}:\mathbb{R}\to \mathbb{R}$ such that,
for any $0<c<1$, for any $k_0\geq 1$ and for any $\ell_0\geq 0$, the function $v(t)$ defined by
\begin{equation}\label{defvth} 
    v(t,x)=Q(y)+Q_{c}(y_{c})+\sum_{(k,\ell)\in \Sigma_0} 
        c^\ell\left(Q_{c}^k(y_{c}) A_{k,\ell}(y)+(Q_{c}^k)'(y_{c}) B_{k,\ell}(y)\right)
\end{equation}
where
$ 
    y_c=x+(1{-}c)t, $ $  y=x-\alpha(y_{c}) $  and $ 
    \alpha(s)=\sum_{(k,\ell)\in \Sigma_0} a_{k,\ell} \, c^\ell \int_{0}^s    Q_{c}^k(s')ds',
$ 
satisfies
\begin{enumerate}
	 \item The function $v(t,x-t)$ is an approximate solution:  $S(t)$ defined by
    \begin{equation}\label{defS}
        S(t,x)=\partial_{t} v + \partial_{x} (\partial_{x}^2 v - v + v^p)
    \end{equation}
    satisfies, for all $j\geq 0$,
    \begin{equation}\label{almost}
        \forall t\in [-T_c,T_c],\quad    \|\partial_x^{(j)}S(t)\|_{L^2(\mathbb{R})}\leq K c^{n_0},
    \end{equation}
    where $n_0=(\tfrac 12 -\tfrac 1{100} ) \min\big(\frac{k_0}{p-1},1+\ell_0\big)$ and $K=K(j,k_0,\ell_0)>0$.

  	\item The function $v(t)$  belongs to $H^1(\mathbb{R})$ for all $t\in \mathbb{R}$ and satisfies for $K=K(k_0,\ell_0)>0$
    \begin{equation}\label{vclose}
        \forall t\in [-T_c,T_c],\quad \|v(t)-R(t)-R_c(t)\|_{H^1(\mathbb{R})}\leq K c^{\frac 1{p-1}}.
    \end{equation}
   \end{enumerate}
\end{theorem}

\noindent\emph{Remark.} (a) Size comparison in \eqref{vclose}.
First, note that
\begin{equation}\label{cScale}
    \|Q_c\|_{L^2}=c^{q} \|Q\|_{L^2},\quad \|Q_c'\|_{L^2}= c^{q+\frac 12} \|Q'\|_{L^2} \quad
    \text{and} \quad \|Q_c\|_{L^\infty}=c^{q+\frac 14} \|Q\|_{L^\infty}.
\end{equation}
Since $\frac 1{p-1}=q+\frac 14$,    \eqref{vclose} says that $v(t)-R(t)-R_c(t)$ is smaller in $H^1$ norm than $R_c(t)$
by a factor $c^{1/4}$. Thus, in $H^1$, $v(t)=R(t)+R_c(t)+$ smaller order terms in $c$. 

Remark that the $L^\infty$ norm is not adequate in this framework, indeed, we also have $\|v(t)-R(t)-R_c(t)\|_{L^\infty}\leq K
\|v(t)-R(t)-R_c(t)\|_{H^1}\leq K c^{\frac 1 {p-1}}\leq K\|Q_c\|_{L^\infty}$.
Moreover, from \eqref{defvth} and from the fact $A_{1,0}\neq 0$ (see proofs),
we have for $t\sim 0$, $\|v(t)-R(t)-R_c(t)\|_{L^\infty}\sim \|Q_c\|_{L^\infty}$.
Observe also that $\|Q_c'\|_{L^2}$ is   smaller than $\|Q_c\|_{L^2}$ for $c$ small. In all this paper, the norm that really matters in the various estimates is the  $L^2$ norm.

Note that \eqref{vclose} is only a first estimate concerning the relation between $v$ and the sum of two solitons.
This estimate does not take into account the shift of the soliton $Q_c$, and thus cannot be sharp. 
In sections 3 and 4, by recompositing $v$ at $t=\pm T_c$, we will prove a better estimate for $v(t)-Q(y)-Q_c(y_c\pm \Delta_c)$, for some $\Delta_c$ and for $t=\pm T_c$ (see Proposition \ref{AVRIL4}).
Estimate \eqref{almost} is also not optimal, especially for small $k_0$ and $\ell_0$ (but $n_0\to +\infty$ as $k_0,\ell_0\to +\infty$).

Note also that $k_0\geq 5$ and $\ell_0\geq 1$ in Theorem \ref{APSOL4} would be enough to prove
the main results of this paper. Nevertheless, the result as stated for all $k_0$, $\ell_0$ clearly indicates that there is no algebraic obstruction to the complete understanding of the interaction process
and we expect it to be useful in future works.

\smallskip

(b) The time interval $[-T_c,T_c]$ contains the interaction region.
Since for $t=-T_c$, $y\ll y_c$ and for $t=T_c$, $y\gg y_c$, the interaction of the two solitons $Q$ and $Q_c$  takes place in the time interval 
$[-T_c,T_c]$.
Moreover, since $y_c=y+\alpha(y_c)+(1-c)t$ we have$|y_c|\geq (1-c) |t| - |\alpha(y_c)| - |y|$.
Thus, if $\sqrt{c}<2,$ we obtain $\sqrt{c} |y_c| \geq (1-c)\sqrt{c} |t| - \sqrt{c}|\alpha(y_c)| - \tfrac 12 |y|$,
and 
by neglecting $\sqrt{c}|\alpha(y_c)|$, we obtain for $|t|\geq T_c$,
\begin{equation*}
    0\leq R(t)R_c(t)
    \leq K c^{\frac 1 {p-1}} e^{-\frac {|y|}2 } \frac 1{e^{\frac 12c^{-{\frac 1{100}}}}},
\end{equation*}
which is an exponentially small term when $c$ is small, which says that the interaction between
$Q$ and $Q_c$ is very weak for such $t$.

\smallskip

(c) Decomposition of $W$. The function constructed in Theorem \ref{APSOL4} is not unique. For given $k_0$ and $\ell_0$ there exist in fact several such functions $v$ corresponding to the fact that the decomposition at $t=0$, for example, is not unique.

We refer to Proposition \ref{SYSkl4} for more properties of the functions $A_{k,\ell}$ and $B_{k,\ell}$ introduced in Theorem \ref{APSOL4}.

Note that choosing $k_0=\ell_0=+\infty$ in this expression of $v$ would formally give an exact solution of the gKdV equation at least for $t\in [-T_c,T_c]$.
However, one has to verify that the resulting series in \eqref{defW} converges in some appropriate sense, which is an open problem.

\medskip

\noindent We give a first interpretation of the function $v$ constructed in Theorem \ref{APSOL4}.
\begin{description}
    \item[Integrable case] ($p=2$ and $3$) In this case, one of the function $v$ constructed in Theorem~\ref{APSOL4} coincides at the main orders to the explicit 2-soliton solution.
    \item[Nonintegrable case] ($p=4$) In this case, explicit $2$-soliton solutions are not known and will be proved not to exist later in this paper. The function $v$ is a \textit{completely new object}. Note that this object, up to the order $c^{n_0}$, plays the same role as a
    $2$-soliton solution in the collision region.
     This will allow us to prove the main results of this paper.
\end{description}

\noindent The proof of Theorem \ref{APSOL4}  is organized as follows:

\smallskip

In Section 2.1, we claim 
that the decomposition of $v(t)$ is preserved by gKdV equation, see Proposition \ref{SYSTEME}). The main part of the proof of Proposition \ref{SYSTEME} is  given in Appendix A. 

\smallskip

In Section 2.2, we derive the systems $(\Omega_{k,\ell})$ to be solved at each rank $(k,\ell)$. Next, we solve a model system $(\Omega)$ related to $(\Omega_{k,\ell})$. In particular, we choose a special structure for the functions $A_{k,\ell}$ and $B_{k,\ell}$ which follows from the resolution of the model system.

Then we solve by induction on $(k,\ell)$ all the systems $(\Omega_{k,\ell})$, for $1\leq k\leq k_0$, $0\leq \ell\leq \ell_0$.
This determines $(a_{k,\ell})$, $(A_{k,\ell})$ and
$(B_{k,\ell})$ for all  $1\leq k\leq k_0$, $0\leq \ell\leq \ell_0$ in the expression of $v$. Thus, at this point the function $v(t)$ is fixed.

Finally in Section 2.3, we prove some properties of $v(t)$ and to estimate  the size of $S(t)$ in terms of powers of $c$.

\medskip

 For $k,$ $k'$, $\ell$, $\ell'\in \mathbb{N}$, we denote 
\begin{equation*}
    (k',\ell') \prec (k,\ell) \quad \text{if
$k'<k$ and $\ell'\leq \ell$ or if $k'\leq k$ and $\ell'<\ell$.}
\end{equation*}
We denote by $\mathcal{Y}$ the set of functions $f\in C^\infty(\mathbb{R})$ such that
\begin{equation*}
    \forall j\in \mathbb{N},\ \exists K_j,\, r_j>0,\ \forall x\in \mathbb{R},\quad |f^{(j)}(x)|\leq K_j (1+|x|)^{r_j} e^{-|x|}. 
\end{equation*}
Note that the set $\mathcal{Y}$ is stable by sum, multiplication and differentiation.

\subsection{Preservation of the decomposition \eqref{defvth} by the equation}
The motivation for choosing $W$ of the form \eqref{defW} is the stability of the family of functions
\begin{equation}\label{ily}
    \left\{c^\ell Q_c^k,\ c^\ell (Q_c^k)', \ k\geq 1,\ \ell \geq 0 \right\}
\end{equation}
by  multiplication and differentiation (see Lemma \ref{surQc}). 
A consequence is that the term $S(t,x)=\partial_{t} v + \partial_{x} (\partial_{x}^2 v - v + v^p)$ has the same decomposition as the function $v$ in terms of functions \eqref{ily}.
Let $\mathcal{L} w = -\partial_x^2 w + w -p Q^{p-1} w$.

\begin{proposition}[Decomposition of $S(t,x)$]\label{SYSTEME}
    Let $K_0=(p+1)k_0+12$ and $L_0=(p+1)\ell_0+4$. Then,
    \begin{align*}
            S(t,x) & = 
        \sum_{(k,\ell)\in \Sigma_0}
        c^\ell Q_c^k(y_c)    \Big[a_{k,\ell} (-3 Q+2 Q^p)'(y)    -(\mathcal{L} A_{k,\ell})'(y)\Big]
        \\& \quad
        + \sum_{(k,\ell)\in \Sigma_0}
        c^\ell (Q_c^k)'(y_c)    
        \Big[a_{k,\ell} (-3 Q'')(y) + \left(3A_{k,\ell}'' +pQ^{p-1} A_{k,\ell}\right)(y)    - (\mathcal{L} B_{k,\ell})'(y)\Big]
        \\& \quad
        + \sum_{\substack{1\leq k\leq K_0 \\ 0\leq \ell \leq L_0}}
        c^\ell\left( Q_c^k(y_c)  F_{k,\ell}(y)   +  (Q_c^k)'(y_c) G_{k,\ell}(y)\right),
    \end{align*}
    where $F_{k,\ell}$ and $G_{k,\ell}$ are functions defined on $\mathbb{R}$ satisfying:
    \begin{itemize}
        \item[{\rm (i)}] Dependence property of $F$ and $G$:
        For any $k,\ell$, the expressions of $F_{k,\ell}$ and $G_{k,\ell}$ depend only on $(a_{k',\ell'})$ and $(A_{k',\ell'})$, $(B_{k',\ell'})$ for $k',$ $\ell'$ such that $(k',\ell') \prec (k,\ell)$.
        \item[{\rm (ii)}] Parity property of $F$ and $G$: Let $k\in \{1,\ldots,K_0\}$, $\ell\in \{0,\ldots,L_0\}$. 
        Assume that for any $(k',\ell')$ such that $(k',\ell')\prec (k,\ell)$
         $A_{k',\ell'}$ is even and    $B_{k',\ell'}$ is odd, then 
        $F_{k,\ell} $ is odd and  $G_{k,\ell}$ is  even.
    \end{itemize}
    Moreover,\\
    $\bullet$ If $p=2$, then
        \begin{align*}
             F_{1,0} &= 2 Q'     ,\quad G_{1,0}= 2Q ,    \\
            F_{2,0} &= (-A_{1,0}+A_{1,0}^2)'- (3 B_{1,0}''+ 2Q B_{1,0}) - a_{1,0}(Q+3 A_{1,0}''+ 2 Q A_{1,0})'
                + 3a_{1,0}^2 Q^{(3)}, \\
            G_{2,0} &= A_{1,0}+A_{1,0}^2+(-2 B_{1,0}+A_{1,0}B_{1,0})' 
            - \frac{a_{1,0}}2    (9 A_{1,0}' + 3B_{1,0}'' + 2Q B_{1,0})'
                +\frac 32 a_{1,0}^2 Q''.\quad
        \end{align*}
     $\bullet$ If $p=4$, then
        \begin{align*}
                 F_{1,0}&= (4Q^3)'  ,\quad G_{1,0}=4Q^3,    \\
            F_{2,0}&=    (6 Q^2(1+A_{1,0})^2)'- a_{1,0}(4 Q^3+3 A_{1,0}''+4 Q^3 A_{1,0})'+ 3 a_{1,0}^2 Q^{(3)}, \\
            G_{2,0}&= 6 Q^2(1+A_{1,0})^2+ (6Q^2 B_{1,0}(1+A_{1,0}))'
             -\frac {a_{1,0}}2 (9 A_{1,0}'+3B_{1,0}''+ 4Q^3 B_{1,0})' +\frac 32 a_{1,0}^2 Q''.
        \end{align*}
 \end{proposition}

See Proposition \ref{SYSkl4}, Lemma \ref{STRUCT4} and Claim \ref{new24} for additional properties of $F_{k,\ell}$ and $G_{k,\ell}$.

\medskip

\noindent\emph{Proof of Proposition \ref{SYSTEME}.} 
A large part of the proof of Proposition \ref{SYSTEME} is given in Appendix A. We present here some preliminary results.

We begin by proving that the family of functions
\eqref{ily} is stable by    multiplication and differentiation.

\begin{lemma}[Properties of $Q$ and $Q_c$]\label{surQc} 
1. The function $Q$ is even and belongs to $\mathcal{Y}$.\\
2. For any $k\in \mathbb{N}^*$,
\begin{align*}
& Q_c''  =cQ_c -Q_c^p, \quad (Q'_c)^2=cQ_c^2-\frac 2{p+1}{Q_c}^{p+1},\\ 
& (Q_c^k)''=ck^2Q_c^k-\frac {k(2k+p-1)}{p+1} Q_c^{k+p-1},\quad (Q_c^k)^{(3)}=ck^2(Q_c^k)'-\frac {k(2k+p-1)}{p+1} (Q_c^{k+p-1})',\\
& (Q_c^k)^{(4)}=c^2k^4 Q_c^k - c \frac {k(2k+p-1)}{p+1}(k^2+(k+p-1)^2) Q_c^{k+p-1}\\
& \qquad\qquad + k(k+p-1)\frac{(2k+p-1)(2k+3p-3)}{(p+1)^2} Q_c^{k+2p-2}.
\end{align*}
3. For any $k_1,$ $k_2\in \mathbb{N}^*$,
\begin{equation*}
(Q_c^{k_1})' Q_c^{k_2} = \frac {k_1}{k_1+k_2} (Q_c^{k_1+k_2})',\quad
(Q_c^{k_1})'(Q_c^{k_2})' = c k_1 k_2 Q_c^{k_1+k_2} - \frac {2 k_1 k_2}{p+1} Q_c^{k_1+k_2+p-1}.
\end{equation*}
\end{lemma}

\noindent\emph{Proof of Lemma \ref{surQc}.}
It is clear from \eqref{defQ}
that $Q$ is even and belongs to $\mathcal{Y}$.
From the equation of $Q_c$, i.e. $Q_c''  =cQ_c -Q_c^p$, we easily get the second equation by multiplying by $Q_c'$ and integrating over $(-\infty,x)$.

Next, we have
\begin{align*}
(Q_c^k)'' & = k (Q_c^{k-1} Q_c')' = k\left( (k-1) Q_c^{k-2} (Q_c')^2 + Q_c^{k-1} Q_c''\right)\\
& = c k^2 Q_c^k - k \left(\frac {2(k-1)}{p+1} +1\right) Q_c^{k+p-1} =
ck^2Q_c^k-\frac {k(2k+p-1)}{p+1} Q_c^{k+p-1}.
\end{align*}
From this we immediately obtain the expression of $(Q_c^k)^{(3)}$. Next,
we have
\begin{align*} (Q_c^k)^{(4)} &= c k^2 (Q_c^k)'' - \frac {k(2k+p-1)}{p+1} (Q_c^{k+p-1})''\\
& = ck^2 \left(c k^2 Q_c^k - \frac {k(2k{+}p{-}1)}{p+1} Q_c^{k+p-1}\right)\\
& \quad - \frac {k(2k{+}p{-}1)}{p+1}\left( c(k{+}p{-}1)^2 Q_c^{k+p-1} - (k{+}p{-}1)\frac{2k{+}3p{-}3}{p+1} Q_c^{k+2p-2}\right).
\end{align*}
The rest of the proof follows.

\medskip

Let now us give a preliminary decomposition of $S(t)$.
We insert $v=R+R_c+W$ into $S(t,x)$, and rearrange terms:
\begin{align*}
    S(t,x)&= \partial_t v + \partial_{x} (\partial_{x}^2 v - v + v^p)\\
    &= \partial_t (R+R_c+W) + \partial_x \left( \partial_x^2 (R+R_c+W) - (R+R_c+W) + (R+R_c+W)^p\right)\\
    &= \partial_t R + \partial_x(\partial_x^2R -R+ R^p)+\partial_t R_c + \partial_x(\partial_x^2 R_c-R_c+R_c^p)\\
    & \quad  +\partial_x((R+R_c)^p -R^p -R_c^p) \\
    & \quad + \partial_t W +\partial_x(\partial_x^2 W -W +(R+R_c+W)^p-(R+R_c)^p).
\end{align*}
By the equation of $Q_c$ ($Q_c''=cQ_c-Q_c^p$) and $y_c=x+(1-c)t$, it is straightforward that
\begin{equation}\label{eqRc}
    \partial_t R_c + \partial_x (\partial_x^2 R_c - R_c + R_c^p) = ((1-c) Q_c+ Q_c''-Q_c+Q_c^p)'(y_c)=0.
\end{equation}
Set 
\begin{equation}\label{defLy}
    \mathcal{L} w = -\partial_x^2 w + w -p Q^{p-1} w,\quad \overline{\mathcal{L}} w = -\partial_x^2  w + w - p R^{p-1} w.
\end{equation}
We decompose $S(t,x)$ as follows:
\begin{equation}\label{decS}
    S(t,x)= \mathbf{I} +\mathbf{II} +\mathbf{III}+{\mathbf{IV}},
\end{equation}
where $\mathbf{I}$, $\mathbf{II}$, $\mathbf{III}$ and $\mathbf{IV}$ are respectively:

- Contribution of terms containing only $R$: 
$\mathbf{I}=\partial_t R + \partial_x(\partial_x^2 R -R + R^p);$

- Nonlinear interaction terms between $R$ and $R_c$:
$ \mathbf{II}=\partial_x((R+R_c)^p-R^p-R_c^p);$

- Linear terms in $W$:
$ \mathbf{III}= \partial_t W - \partial_x (\overline{\mathcal{L}} W);$

- Higher order terms in $W$:
$ {\mathbf{IV}}= \partial_x((R+R_c+W)^p-(R+R_c)^p-pR^{p-1} W).$

\medskip

The expansion of $\mathbf{I}$, $\mathbf{II}$, $\mathbf{III}$ and $\mathbf{IV}$ is given in Appendix A, and allows to finish the proof of Proposition \ref{SYSTEME}.

\subsection{Resolution of the systems $(\Omega_{k,\ell})$}

From Proposition \ref{SYSTEME}, we observe that if for any $0\leq k\leq k_0$, $0\leq \ell\leq \ell_0$, $(a_{k,\ell},A_{k,\ell},B_{k,\ell})$ satisfies the following system
\begin{equation*}
    (\Omega_{k,\ell})\quad 
    \left\{
    \begin{array}{l}
         (\mathcal{L} A_{k,\ell})' + a_{k,\ell} (3Q -2 Q^p)' = F_{k,\ell}\\
         (\mathcal{L} B_{k,\ell})' + a_{k,\ell} (3Q'') - 3 A_{k,\ell}'' -p Q^{p-1} A_{k,\ell}=G_{k,\ell},
    \end{array}
    \right.
\end{equation*}
then $S(t,x)$  contains only terms of the form $c^{\ell} Q_c^k$ or
$c^{\ell}(Q_c^{k})'$ with $k\geq k_0+1$ or $\ell\geq \ell_0+1$.

\medskip

This observation leads us to consider the model system
\begin{equation*}
    (\Omega)\quad 
    \left\{
    \begin{array}{l}
         (\mathcal{L} A)' + a(3Q -2 Q^p)' = F\\
         (\mathcal{L} B)' + a (3Q'') - 3 A'' -p Q^{p-1} A=G,
    \end{array}
    \right.
\end{equation*}
where $F(x)$ and $G(x)$ are given functions (with a specific structure, see Proposition \ref{SYS4}) and $(a,A(x),B(x))$ is to be determined.
We study existence of solutions of the system $(\Omega)$. 
Before stating and proving the existence result for the model system $(\Omega)$, we introduce some notation and we recall well-known results concerning the operator $\mathcal{L}$.

First, let $\varphi:\mathbb{R}\rightarrow \mathbb{R}$ be defined by 
\begin{equation*}
    \forall x\in \mathbb{R}, \quad \varphi(x)=-\frac {Q'(x)}{Q(x)}.
\end{equation*}

\begin{claim}\label{surphi}
    The function $\varphi$ is odd and    satisfies the following properties.
    \begin{itemize}
        \item[{\rm (a)}] $\lim_{x\rightarrow -\infty} \varphi(x)=-1$; $\lim_{x\rightarrow +\infty} \varphi(x)=1$;
        \item[{\rm (b)}] $\forall x\in \mathbb{R}$, $|\varphi'(x)|+|\varphi''(x)|+|\varphi^{(3)}(x)|\leq C e^{-|x|}$.
        \item[{\rm (c)}] $\varphi'\in \mathcal{Y}$,  $(1-\varphi^2) \in \mathcal{Y}$.
        \item[{\rm (d)}] 
        {For $p=2$,} $ (\mathcal{L} \varphi)' = 2 Q -\frac 53 Q^2.$
        {For $p=4$,} $(\mathcal{L} \varphi)' = \frac {36}5 Q^3 - \frac {99}{25} Q^6.$    \end{itemize} 
\end{claim}
\noindent\emph{Proof of Claim \ref{surphi}.}
From the explicit formula
$Q(x)=\Big( \frac {p+1}  {2\cosh^2\bigl(\frac {p-1} 2 \,x\bigr)}\Big)^{\frac1 {p-1}},$
we have
\begin{equation*}
    Q'(x)=-\tanh\bigl(\tfrac {p-1} 2 \,x\bigr) Q(x),
\end{equation*}
and so $\varphi(x)=\tanh\bigl(\tfrac {p-1} 2 \,x\bigr)$. From $\tanh'=1-\tanh^2=\frac {1}{\cosh^2}$, we obtain (a), (b) and (c).

 By $Q''=Q-Q^p$ and $(Q')^2=Q^2-\frac 2{p+1} Q^{p+1}$, we have
\begin{equation*}
    \varphi' = - \frac 1 {Q^2} (Q'' Q- (Q')^2)=\frac {p-1}{p+1} Q^{p-1}, \quad \text{and} \quad
    \varphi''=\frac {(p-1)^2}{p+1} Q' Q^{p-2}.
\end{equation*}
Thus,
$
    -\varphi''-pQ^{p-1}\varphi=\bigl(-\frac {(p-1)^2}{p+1}+p\bigr) Q' Q^{p-2}= \frac {3p-1}{p+1} Q' Q^{p-2},
$
and
\begin{align*}
    (\mathcal{L} \varphi)' & = \frac {3p-1}{p+1} Q'' Q^{p-2} + \frac {(3p-1)(p-2)}{p+1} (Q')^2 Q^{p-3} + \frac {p-1}{p+1} Q^{p-1}\\
    & =\frac {3 p (p-1)}{p+1} Q^{p-1} - \frac {3 (3p-1)(p-1)}{(p+1)^2} Q^{2(p-1)}.
\end{align*}

\begin{lemma}[Properties of $\mathcal{L}$]\label{surL} Let $p\geq 2$. The operator $\mathcal{L}$ defined in $L^2(\mathbb{R})$ by
    \begin{equation*}
        \mathcal{L} f= -f''+f-pQ^{p-1} f
    \end{equation*}
    is self-adjoint and satisfies the following properties:
    \begin{itemize}
        \item[{\rm (i)}] First eigenfunction : $\mathcal{L} Q^{\frac {p+1} 2} = - \frac 14 (p{-}1)(p{+}3) Q^{\frac {p+1} 2}$;
        \item[{\rm (ii)}] Second eigenfunction : $\mathcal{L} Q'=0$; the kernel of $\mathcal{L}$ is 
        $\{\lambda Q', \lambda \in \mathbb{R}\}$;
        \item[{\rm (iii)}] For any   function $h \in L^2(\mathbb{R})$ orthogonal to $Q'$ for the $L^2$ scalar product, 
        there exists a unique function $f \in H^2(\mathbb{R})$ orthogonal to $Q'$ such that $\mathcal{L} f=h$; moreover,
        if $h$ is even (respectively, odd), then $f$ is even (respectively, odd).
        \item[{\rm (iv)}] Suppose that $f\in H^2(\mathbb{R})$ is such that $\mathcal{L} f \in \mathcal{Y}$.
            Then, $f\in \mathcal{Y}$.
    \end{itemize}
\end{lemma}

\medskip

\noindent\emph{Proof of Lemma \ref{surL}.} From 
$Q''=Q-Q^p$ and $(Q')^2=Q^2 -\frac 2{p + 1} Q^{p+1},$
\begin{align*}
     \frac {d^2}{d x^2} Q^{\frac {p+1} 2} &= \tfrac {p + 1} 2 \left[ \tfrac {p - 1}2 Q'^2 Q^{\frac {p-3}2} + Q'' Q^{\frac {p-1}2}\right]
    = \left(\tfrac {p + 1} 2 \right)^2 Q^{\frac {p+1} 2} -p Q^{p-1} Q^{\frac {p+1} 2},
\end{align*}
and so $\mathcal{L} Q^{\frac {p+1} 2} =  - \Bigl[\bigl(\frac{p+1}2\bigr)^2 - 1\Bigr] Q^{\frac {p+1} 2}=- \frac 14 (p{-}1)(p{+}3) Q^{\frac {p+1} 2}.$

The property $\mathcal{L} Q'=0$ is easily checked. Moreover, the fact that the spectrum of $\mathcal{L}$ is restricted to $\{\lambda Q', \lambda \in \mathbb{R}\}$ was proved by ordinary differential equations techniques (see Weinstein \cite{We2}, Proposition 2.8 (b)). The third property is a direct consequence of the structure of $\mathcal{L}$, and Lax-Milgram theorem. 

Property {\rm (iv)} is also a consequence of standard arguments of ordinary differential equations theory. First, we claim the following.

\begin{claim}\label{claimzz}
    Suppose that $f\in H^2(\mathbb{R})$ satisfies for $K>0$ and $r>0$,
    \begin{equation}\label{decayL}
        \forall x\in \mathbb{R},\quad |(f''-f)(x)|\leq K (1+|x|^r) e^{-|x|}.
    \end{equation}  
    Then, there exists $K'>0$ such that
    \begin{equation}\label{decayf}
        \forall x\in \mathbb{R},\quad |f(x)|\leq K' (1+|x|^{r+1}) e^{-|x|}.
    \end{equation}
\end{claim}

\noindent\emph{Proof of Claim \ref{claimzz}.}\quad 
We  set $g(x)=e^{-x}(f'+f)$. Then $g'=e^{-x}(f''-f)$, 
\begin{equation*}
    \forall x>0, \quad |g'(x)|\leq K (1+|x|^r) e^{-2 x}, \quad \text{and}
\end{equation*}
\begin{equation*}
    |g(x)|\leq K \int_x^{+\infty} (1+s^r) e^{- 2 s} ds \leq K' (1+x^r) e^{-2 x}.
\end{equation*}
Set $h=e^x f$. Then $|h'|=|e^{2x}g|\leq K (1+|x|^r)$. 
By integration between $0$ and $x$, we obtain
$    \forall x>0,\quad e^x |f(x)|=|h(x)| \leq K'' (1+|x|^{r+1}).$
The same property is true for $x<0$, by changing $x$ in $-x$. 

\medskip

We now finish the proof of {\rm (iv)}. 
Let $f\in H^2(\mathbb{R})$ be such that $\mathcal{L} f\in \mathcal{Y}$. Since
$f''=(- \mathcal{L} f +f - p Q^{p-1} f)$, 
by induction on $j$ and $Q\in \mathcal{Y}$, 
it is clear that $f\in C^j(\mathbb{R})$, for all $j\in \mathbb{N}$.
Since $(f^{(j)})''-f^{(j)}=-(\mathcal{L} f +p Q^{p-1} f)^{(j)}$ and  $\mathcal{L} f$, $Q\in \mathcal{Y}$, using Claim \ref{claimzz} we prove    by an induction argument on $j$ that for all
$j$ and all $x$, $|f^{(j)}(x)|    \leq K_j    (1+|x|^{r_j}) e^{-|x|}$.
Thus, $f\in \mathcal{Y}$.

\medskip

The next result of this section concerns the existence of solutions of system $(\Omega)$.

\begin{proposition}[Existence for the model problem $(\Omega)$]\label{SYS4}
    Let $F(x)$ and $G(x)$ be  such that
    \begin{equation*}
        F    = \overline{F}  + \widetilde{F}  + \varphi  \widehat{F} ,\quad
        G    = \overline{G}  + \widetilde{G}  + \varphi  \widehat{G} ,
    \end{equation*}
       \begin{itemize}
        \item $\overline{F}$, $\overline{G}\in \mathcal{Y}$;
          $\overline{F}$ is odd and $\overline{G}$ is even;
        \item $\widetilde{F}$ and $\widehat{G}$ are odd polynomial functions;
         $\widehat{F}$ and $\widetilde{G}$ are even polynomial functions.
    \end{itemize}
    Then, there exist $a\in \mathbb{R}$ and two functions $A(x) $, $B(x)$
    \begin{equation*}
        A    = \overline{A}  + \widetilde{A}  + \varphi  \widehat{A} ,\quad
        B    = \overline{B}  + \widetilde{B} + \varphi  \widehat{B} ,
    \end{equation*}
        \begin{itemize}
        \item $\overline{A}$, $\overline{B}\in \mathcal{Y}$;
          $\overline A$ is even and  $\overline B$ is odd;
        \item $\widetilde{A}$ and $\widehat{B}$ are even polynomial functions;
        $\widehat{A}$ and $\widetilde{B}$ are odd polynomial functions;
    \end{itemize}
     satisfying
    \begin{equation*}
        (\Omega)\quad 
        \left\{
        \begin{array}{ll}
             (\mathcal{L} A)' + a (3Q -2 Q^p)' = F                           & (\Omega_A)\\
             (\mathcal{L} B)' + a (3Q'') - 3 A'' -p Q^{p-1} A=G. \qquad         & (\Omega_B)
        \end{array}
        \right.
    \end{equation*}
    The degrees of the polynomial functions $\widetilde{A}$, $\widehat{A}$,  $\widetilde{B}$ and $\widehat{B}$ are related to the degrees of
    $\widetilde{F}$, $\widehat{F}$, $\widetilde{G}$  and $\widehat{G}$ as follows:
    \begin{align}
        &\deg \widetilde{A} \leq 1+ \deg \widetilde{F},\quad \deg \widetilde{B} \leq \max(1+\deg \widetilde{G},\deg \widetilde{F}),\label{degres1} \\
        &\deg \widehat{A} \leq 1+ \deg \widehat{F},\quad \deg \widehat{B} \leq \max(1+\deg \widehat{G},\deg \widehat{F}).\label{degres2}
    \end{align}
    Moreover,
    \begin{align}
        &\text{if $\widetilde F=0$ (respectively, $\widehat F=0$) then $\widetilde A=0$ (respectively, $\widehat A=0$);}
                 \label{DEGA1}\\
        &\text{if $\widetilde A''=0$ and $\widetilde G=0$ then $\widetilde B=0$;}\label{DEGB1}\\
        &\text{if $\widehat A''=0$ and $\widehat G=0$ then $\deg \widehat B=0$.}\label{DEGB2}
    \end{align}
\end{proposition}

\noindent\emph{Remark.}\quad Observe that the conclusions of \eqref{DEGB2} and \eqref{DEGA1}-\eqref{DEGB1} are different.
In \eqref{DEGB2}, only $\deg \widehat B=0$ which allows the possibility that $\widehat B=b$, a nonzero constant, even if no polynomial is present in $F$ and $G$. 
Without this freedom, the system cannot be solved in general.
This remark is    essential for    two reasons: 
\begin{enumerate}
    \item The fact that possibly $\widehat B\neq 0$ whereas $\widetilde F$, $\widehat F$, $\widetilde G$ and $\widehat G$ are zero, is responsible for the apparition of polynomial growths in $A_{k,\ell}$ and $B_{k,\ell}$ when solving the systems $(\Omega_
{k,\ell})$. Indeed, from the structure of the systems $(\Omega_{k,\ell})$, one cannot find solutions $A_{k,\ell}$,
    $B_{k,\ell}$ all in $\mathcal{Y}$.
    It is the reason why we need to allow polynomial growth in the functions $A,$ $B$, $F$ and $G$ as in Proposition \ref{SYS4}.
    \item In the next section, we will see that the shift on the soliton $Q_c$ resulting from the interaction with the soliton $Q$ is obtained from $\widehat B_{1,0}\neq 0$.
\end{enumerate}

\medskip

\noindent\emph{Remark.} In Proposition \ref{SYS4}, we find one solution of the system $(\Omega)$.
We refer to Corollary \ref{UNISYS4} for the uniqueness question.

\medskip

\noindent\emph{Proof of Proposition \ref{SYS4}.} \quad We first reduce the proof to the case where there is no polynomial functions in $F$ and $G$. Then we solve the problem using Lemma \ref{surL} and choosing the free parameter $a$.

\medskip

\noindent\emph{Step 1. Reduction to the case without polynomial functions.} \quad

Let $F$ and $G$ be two functions satisfying the assumptions of Proposition \ref{SYS4}.
First, we consider $\widetilde{A}$ and $\widehat{A}$ the two (unique) polynomial functions satisfying
\begin{align*}
    &    -\widetilde{A}''(x)+\widetilde{A}(x) = \int_0^x \widetilde{F}(z) dz
    \quad \text{and}\quad    -\widehat{A}''(x)+\widehat{A}(x) = \int_0^x \widehat{F}(z) dz
\end{align*}
(obtained by resolution of a system in the basis $\{x^r\}_{r\geq 0}$.
Observe that $\widetilde{A}$ is even and $\widehat{A}$ is odd; moreover
\begin{itemize}
    \item if $\widetilde{F}=0$ (respectively, $\widehat{F}=0$) then $\widetilde{A}=0$ (respectively, $\widehat{A}=0$);
    \item if $\widetilde{F}\not=0$ (respectively, $\widehat{F}\not=0$) then $\deg \widetilde{A}=1 + \deg \widetilde{F}$ (respectively, $\deg \widehat{A}=1 + \deg \widehat{F}$).
\end{itemize}
We have
\begin{equation*}
    (\mathcal{L} \widetilde{A})'=(-\widetilde{A}''+\widetilde{A}-pQ^{p-1}\widetilde{A})' = \widetilde{F} - p (Q^{p-1}\widetilde{A})',
\end{equation*}
\begin{align*}
    (\mathcal{L} (\varphi \widehat{A}))' &= \left(- \varphi \widehat{A}'' - 2 \varphi' \widehat{A}' - \varphi'' \widehat{A} + \varphi \widehat{A} - p Q^{p-1} \varphi \widehat{A} \right)'\\
    &=\varphi \widehat{F} + \varphi' \int_0^x \widehat{F} + \left(-2 \varphi' \widehat{A}'-\varphi'' \widehat{A} - p Q^{p-1}\varphi \widehat{A}\right)'.
\end{align*}
For $\overline{A}$ to be chosen later, let
$    A=\overline{A} + \widetilde{A} + \varphi \widehat{A}.$
Then,    $A$ solves $(\Omega_A)$ if and only if
\begin{equation*}
    (\mathcal{L} \overline{A})'+ (\mathcal{L} \widetilde{A})' + (\mathcal{L} (\varphi \widehat{A}))'+ a(3 Q - 2Q^p)= \overline{F} + \widetilde{F} + \varphi \widehat{F},
\end{equation*}
or equivalently by the previous calculations
$    (\mathcal{L} \overline{A})'+ a(3 Q - 2Q^p)= \mathcal{F},$
where 
\begin{align}
    \mathcal{F} & = \overline{F} + \widetilde{F} - (\mathcal{L} \widetilde{A})' + \varphi \widehat{F} - (\mathcal{L} (\varphi \widehat{A}))'\nonumber \\
     & = \overline{F} - \varphi' \int_0^x \widehat{F} + \left(2 \varphi' \widehat{A}'+\varphi'' \widehat{A}+p Q^{p-1}(\widetilde{A}+\varphi \widehat{A})\right)'. \label{defcF}
\end{align}
Since $\overline F$, $\varphi'$, $Q\in \mathcal{Y}$, and $\widetilde A$, $\widehat A$ and $\widehat F$ are polynomial functions, we get $\mathcal{F}\in \mathcal{Y}$.
Moreover, we observe that $\mathcal{F}$ is odd.

\smallskip

We proceed in a similar way for $B(x)$ except for the need of an additional parameter $b\in \mathbb{R}$ and the term $(-3A'')$ in equation $(\Omega_B)$.
Let $\widetilde{B}$ and $\widehat{B}^*$ be the two (unique) polynomial functions satisfying
\begin{align*}
    &    -\widetilde{B}''(x)+\widetilde{B}(x) = \int_0^x \left(\widetilde{G}(z)  + 3 \widetilde{A}''(z)\right)dz,\\
    &    -(\widehat{B}^*)''(x)+\widehat{B}^*(x) = \int_0^x \left(\widehat{G}(z)  + 3 \widehat{A}''(z)\right)dz.
\end{align*}
Observe that $\widetilde{B}$ is odd and $\widehat{B}^*$ is even; moreover  
\begin{itemize}
    \item if $\widetilde{A}''=0$ and $\widetilde{G}=0$  then $\widetilde{B}=0$;
    \item if $\widetilde{A}''\not=0$ or $\widetilde{G}\not=0$ then 
    $
        \deg \widetilde{B}=1 + \max(\deg \widetilde{G}, \deg \widetilde{A}''),
    $
    \item if $\widehat{A}''=0$ and $\widehat{G}=0$  then $\widehat{B}^*=0$;
    \item if $\widehat{A}''\not=0$ or $\widehat{G}\not=0$ then 
    $    \deg \widehat{B}^*=1 + \max(\deg \widehat{G}, \deg \widehat{A}''),
    $
\end{itemize}
In all cases, we have
\begin{equation}\label{degBA}
    \deg \widetilde{B}\leq \max(1+ \deg \widetilde{G}, \deg \widetilde{F} ),\quad       \deg \widehat{B}^*\leq \max(1+ \deg \widehat{G}, \deg \widehat{F}).
\end{equation}

We have
\begin{equation*}
    (\mathcal{L} \widetilde{B})'=(-\widetilde{B}''+\widetilde{B}-pQ^{p-1}\widetilde{B})' = \widetilde{G} + 3 \widetilde{A}''- p (Q^{p-1}\widetilde{B})',
\end{equation*}
\begin{align*}
    (\mathcal{L} (\varphi \widehat{B}^*))' &= \left(- \varphi (\widehat{B}^*)'' - 2 \varphi' (\widehat{B}^*)' - \varphi'' \widehat{B}^* + \varphi \widehat{B}^* - p Q^{p-1} \varphi \widehat{B}^* \right)' \\
    &=\varphi (\widehat{G}+ 3 \widehat{A}'') +  \varphi' \int_0^x (\widehat{G}(z) + 3 \widehat{A}''(z))dz  - \left(2 \varphi' (\widehat{B}^*)'+\varphi'' \widehat{B}^* + p Q^{p-1}\varphi \widehat{B}^* \right)' .
\end{align*}
For $\overline{B}$ and $b$ to be chosen later, let
\begin{equation*}
    B=\overline{B} + \widetilde{B} + \varphi \widehat{B}, \quad \text{with}\quad \widehat{B}= \widehat{B}^* + b,
\end{equation*}
Then,    $B$ solves $(\Omega_B)$ if and only if
\begin{align*}
    &(\mathcal{L} \overline{B})'+ (\mathcal{L} \widetilde{B})' + (\mathcal{L} (\varphi \widehat{B}))'+ 3 a Q'' 
    - 3 \overline{A}'' - pQ^{p-1} \overline{A} \\
    &\quad - 3 \widetilde{A}''- p Q^{p-1} \widetilde{A} - 3(\varphi \widehat{A})'' - p Q^{p-1} (\varphi \widehat{A})= \overline{G} + \widetilde{G} + \varphi \widehat{G},
\end{align*}
or equivalently by the previous calculations
\begin{equation*}
    (\mathcal{L} \overline{B})'+ 3 a Q'' - 3 \overline{A}'' - pQ^{p-1} \overline{A} = \mathcal{G}+ b (\mathcal{L} \varphi)',
\end{equation*}
where the function $\mathcal{G}$ is defined by
\begin{align}
    \mathcal{G} & = \overline{G} + \widetilde{G} + 3 \widetilde{A}''- (\mathcal{L} \widetilde{B})' 
        + \varphi \widehat{G} + 3 (\varphi \widehat{A})''- (\mathcal{L} (\varphi \widehat{B}^*))'
        + pQ^{p-1}(\widetilde{A} + \varphi \widehat{A})\nonumber\\
    & = \overline{G} + 6 \varphi' \widehat{A}'+3 \varphi'' \widehat{A}
        - \varphi' \int_0^x (\widehat{G}(z)+3 \widehat{A}''(z))dz\nonumber \\& \quad
        + \left(2 \varphi' (\widehat{B}^*)'+\varphi'' \widehat{B}^*+p Q^{p-1}(\widetilde{B}+\varphi \widehat{B}^*)\right)'.\label{defcG}
\end{align}
Since $\overline G$, $\varphi'$, $Q\in \mathcal{Y}$, and $\widetilde A$, $\widetilde B$, $\widehat B^*$ and $\widehat G$ are polynomial functions, $\mathcal{G}\in \mathcal{Y}$ is even.

\smallskip

Thus, in conclusion, the system $(\Omega)$ is equivalent to the following system in $(a,b,\overline{A},\overline{B})$:
\begin{equation*}
    \left\{
    \begin{array}{ll}
         (\mathcal{L} \overline{A})' + a (3Q -2 Q^p)' = \mathcal{F}                 \\
         (\mathcal{L} \overline{B})' + a (3Q'') - 3 \overline{A}'' -p Q^{p-1} \overline{A}=\mathcal{G} + b (\mathcal{L} \varphi)', 
    \end{array}
    \right.
\end{equation*}
where $\mathcal{F}\in \mathcal{Y}$ is odd, given by \eqref{defcF}, $\mathcal{G}\in \mathcal{Y}$ is even, given by \eqref{defcG}.    Note that $\mathcal{F}$ and $\mathcal{G}$ do not depend on the  parameters $a$ and $b$.

\medskip

\noindent\emph{Step 2. Existence of a solution of system $({\Omega})$.} \quad
We set
$   \mathcal{H}(x)=\int_{-\infty}^x \mathcal{F}(z) dz.$
Since $\mathcal{F}$ is odd, $\int_{\mathbb{R}} \mathcal{F}=0$ and so
$\mathcal{H}\in \mathcal{Y}$ is even.

To find a solution $(a,b,\overline A,\overline B)$ of $(\overline{\Omega})$, it is sufficient to solve
\begin{equation*}
    (\overline{\Omega})\quad 
    \left\{
    \begin{array}{ll}
         \mathcal{L} \overline{A} + a (3Q -2 Q^p)= \mathcal{H}  \\
         (\mathcal{L} \overline{B})' + a (3Q'') - 3 \overline A'' -p Q^{p-1} \overline A=\mathcal{G} + b (\mathcal{L} \varphi)'.
    \end{array}
    \right.
\end{equation*}
Since $\int \mathcal{H}Q'=0$ (by parity) and $\mathcal{H}\in \mathcal{Y}$, it follows from Lemma \ref{surL} (iii)-(iv) that there exists $\overline{H} \in \mathcal{Y}$, even, such that 
\begin{equation*}
    \mathcal{L}\overline{H}=\mathcal{H}.
\end{equation*}
By Lemma \ref{surL}, there also exists $V_0 \in \mathcal{Y}$, even, such that 
$    \mathcal{L} V_0=3 Q -2 Q^p.$
It follows that, for all $a$,
\begin{equation}
    \overline{A} = \overline{H} - a V_0
\end{equation}
is solution of $\mathcal{L} \overline{A} + a(3 Q - 2 Q^p)= \mathcal{H}$, 
moreover, $\overline{A}$ is even and $\overline A\in \mathcal{Y}$.  
Note that at this point $(a,b)$ are still free, they will be chosen when solving the second equation. 

Now, replacing $\overline{A}$ by $\overline{H} - a V_0$ in this equation, we only need to find $\overline{B}$ such that
\begin{equation}\label{interB} 
    (\mathcal{L} \overline{B})' =- a Z_0 + D + b (\mathcal{L} \varphi)',
\end{equation}
where
\begin{equation*}
    D= 3 \overline{H}''+pQ^{p-1} \overline{H}+\mathcal{G},
\quad
    Z_0=3 Q'' + 3 V_0'' + p Q^{p-1} V_0.
\end{equation*}
It follows from the properties of $Q$, $V_0$, $\mathcal{G}$ and $\overline{H}$ that $D$ and $Z_0$ are even and satisfy
$Z_0$, $D\in \mathcal{Y}$.
To solve \eqref{interB}, it suffices to find $\overline{B}\in \mathcal{Y}$ such that
\begin{equation}\label{fm}
    \mathcal{L} \overline{B} = E \quad \text{where} \quad E=\int_{0}^x (D-a Z_0)(z) dz + b \mathcal{L}\varphi.
\end{equation}
We can choose $(a,b)$ such that the function $E$ is orthogonal to $Q'$ and has decay at $+\infty$.

\begin{claim}\label{surE} 
    \begin{itemize}
        \item[{\rm (i)}] Nondegeneracy:
        \begin{equation}\label{NONDEGE}
    \int Z_0 Q = \frac {p-5}{4(p-1)} \int Q^2.
\end{equation}
             \item[{\rm (ii)}] Let
        $
            a= \frac {\int DQ}{\int Z_0 Q}$  and $b= - \int_0^{+\infty} (D-aZ_0)(z) dz.
        $
        Then, $E$ defined by \eqref{fm} satisfies
        \begin{equation}\label{decE}
            E\in \mathcal{Y}, \quad E \text{ is odd,} \quad \int E Q'=0 .
        \end{equation}
    \end{itemize}
\end{claim}

Assuming Claim \ref{surE}, we finish the proof of Proposition \ref{SYS4}.
We fix $(a,b)$ as in Claim \ref{surE}. Then, from \eqref{decE} and Lemma \ref{surL}, it follows that there exists
 $\overline{B}\in \mathcal{Y}$ such that
$\mathcal{L} \overline{B} = E .$
Setting
\begin{equation*}
    A=\overline A+\widetilde A+\widehat A, 
     \quad     B=\overline B+\widetilde B+\widehat B,
\end{equation*}
we have constructed a solution of system $(\Omega)$ with the structure described in
Proposition \ref{SYS4}.

\medskip 

Now, we only have to prove Claim \ref{surE}.

\medskip

\noindent\emph{Proof of Claim \ref{surE}.} \quad Proof of {\rm (i)}.
First, we check that 
\begin{equation}\label{Vzero}
    V_0=-\tfrac 1{p-1} Q-\tfrac 32 xQ'.
\end{equation} 
Indeed, $\mathcal{L}Q=-Q''+Q-p Q^{p-1}Q=-(p-1)Q^p$
and $\mathcal{L} (xQ')=-2 Q'' + x \mathcal{L}Q'=-2 Q + 2Q^p$. Thus,
\begin{equation*}
    \mathcal{L}(-\tfrac 1{p-1} Q-\tfrac 32 xQ') =-\frac 1{p-1} \mathcal{L}Q-\tfrac 32 \mathcal{L}(xQ')=-Q^p +3 Q -3 Q^p=3Q-2Q^p. 
\end{equation*}

Second, we compute $\int Z_0 Q$, where $Z_0=3 Q'' + 3 V_0'' + p Q^{p-1} V_0.$ By $Q''=Q-Q^p$, we get:
\begin{align*}
    \int Z_0 Q & = \int \left(3 Q'' + 3 V_0'' + p Q^{p-1} V_0\right) Q=3 \int Q^2 - 3 \int Q^{p+1} + \int V_0 (3 Q''+pQ^p)\\
    &=3 \int Q^2 - 3 \int Q^{p+1} + \int V_0 (3 Q+(p-3) Q^p).
\end{align*}
We compute the last term, integrating    by parts:
\begin{align*}
    &\int V_0 (3 Q+(p-3) Q^p) = - \int \left(\tfrac 1{p-1} Q+\tfrac 32 xQ'\right) (3 Q+(p-3) Q^p)\\
    &= -3 \left(\tfrac 1{p-1} - \tfrac 34\right)\int Q^2 + (p-3) \left(\tfrac 1{p-1}-\tfrac {3}{2(p+1)}\right) \int Q^{p+1}\\
    &=  \frac {3(3p-7)}{4(p-1)} \int Q^2 + \frac{(p-5)(p-3)}{2(p-1)(p+1)} \int Q^{p+1}.
\end{align*}
Finally, using Claim \ref{LemmaA1} in Appendix C,
\begin{align*}
    \int Z_0 Q & = \frac{3(7p-11)}{4(p-1)}\int Q^2 - \frac{(5p-7)(p+3)}{2(p-1)(p+1)} \int Q^{p+1}
     = \frac {p-5}{4(p-1)} \int Q^2.
\end{align*}

\medskip

Proof of {\rm (ii)}. Let $a$ and $b$ be defined as in Claim \ref{surE}. The function $E$ is odd
by its definition in \eqref{fm}.
By integration by parts, and decay properties of $Q$, we have
\begin{equation*}
    \int E Q'=-\int (D-a Z_0)Q+ b \int (\mathcal{L} \varphi) Q'= -\int DQ + a \int Z_0 Q
    +b \int \varphi (\mathcal{L} Q')=0,
\end{equation*}
by the definition of $a,b$ and $\mathcal{L}Q'=0$.
By Claim \ref{surphi} and the definition of $a,b$, we have
\begin{equation*}
    \lim_{+\infty}E=\int_{0}^{+\infty} \left(D-a Z_0 \right) dz+ b \lim_{+\infty} (\mathcal{L}\varphi)=0
\quad \text{and so $E\in \mathcal{Y}$.} \qquad \Box
\end{equation*}

\medskip

\noindent\emph{Resolution of the systems $(\Omega_{k,\ell})$}

\smallskip

Using Propositions \ref{SYSTEME} and \ref{SYS4}, we solve the systems $(\Omega_{k,\ell})$ 
by induction on $(k,\ell)$.  We check that at given $(k,\ell)$, the  systems $(\Omega_{k',\ell'})$    being solved for all $(k',\ell')\prec (k,\ell)$, we can apply Proposition \ref{SYS4} to 
 $(\Omega_{k,\ell})$.
The induction argument can be for example:
1) Initialization: $k=1,\ \ell=0$, 2) For $\ell=0$, all $k\geq 1$, by induction on $k$, 3) 
By induction on $\ell\geq 0$, all $k\geq 1$ similarly as in 2).

For future use in the proof of Theorem \ref{APSOL4}, we estimate in the next section
 the degrees of the polynomials $\widetilde A_{k,\ell}$, $\widehat A_{k,\ell}$, $\widetilde B_{k,\ell}$ and $\widehat B_{k,\ell}$ with respect to $k$ and $\ell$ (see Lemma \ref{DEGG}).

\begin{proposition}[Resolution of $(\Omega_{k,\ell})$ by induction on $(k,\ell)$]\label{SYSkl4}
    For all $k\in \{1,\ldots,k_0\},$ $\ell\in \{0,\ldots,\ell_0\}$, there exists $(a_{k,\ell},A_{k,\ell},B_{k,\ell})$ of the form
    \begin{equation}\label{ABkl1}
    \begin{split}
        & A_{k,\ell}(x) = \overline{A}_{k,\ell}(x) + \widetilde{A}_{k,\ell}(x) + \varphi(x) \widehat{A}_{k,\ell}(x),\\
        &  B_{k,\ell}(x) = \overline{B}_{k,\ell}(x) + \widetilde{B}_{k,\ell}(x) + \varphi(x) \widehat{B}_{k,\ell}(x), \quad\text{where}  \\
         &    \text{$\overline{A}_{k,\ell}$, $\overline{B}_{k,\ell}\in \mathcal{Y}$;
          $\overline{A}_{k,\ell}$ is even and  $\overline{B}_{k,\ell}$ is odd;}\\
         &   \text{$\widetilde{A}_{k,\ell}$ and $\widehat{B}_{k,\ell}$ are even polynomials;
        $\widehat{A}_{k,\ell}$ and $\widetilde{B}_{k,\ell}$ are odd polynomials;}
    \end{split}\end{equation}
    satisfying
    \begin{equation*}
        (\Omega_{k,\ell})\quad 
        \left\{
        \begin{array}{ll}
             (\mathcal{L} A_{k,\ell})' + a_{k,\ell} (3Q -2 Q^p)' = F_{k,\ell}                             &  \\
             (\mathcal{L} B_{k,\ell})' + a_{k,\ell} (3Q'') - 3 A_{k,\ell}'' -p Q^{p-1} A_{k,\ell}=G_{k,\ell}, \qquad         &    
        \end{array}
        \right.
    \end{equation*}
    where  $F_{k,\ell}$, $G_{k,\ell}$ are defined in Proposition \ref{SYSTEME}.
    Moreover, 
    for all $1\leq k\leq p-1$, $\ell=0$,
    \begin{equation}\label{debut}
            \widetilde A_{k,0}=\widehat A_{k,0}=\widetilde B_{k,0}=0,\quad \widehat B_{k,0}=b_{k,0}\in \mathbb{R}.
    \end{equation}
\end{proposition}

\noindent\emph{Remark.} 
\quad 
(i) The parity condition on $A_{k,\ell}$, $B_{k,\ell}$ are related to the resolution of the systems $(\Omega_{k',\ell'}$) for $(k,\ell)\prec(k',\ell')$. The use of the function $\varphi$ is related to the asymmetry of the gKdV equation.

(ii) The resolution of $(\Omega_{k,\ell})$ at each step $(k,\ell)$ does not give a unique solution.
Indeed, from Corollary \ref{UNISYS4}, if $(a_{k,\ell},A_{k,\ell},B_{k,\ell})$ is solution, then for any $(\gamma_{k,\ell},\delta_{k,\ell})\in \mathbb{R}^2$,
\begin{equation}\label{exct}
(a_{k,\ell}+\gamma_{k,\ell}a_{1,0},A_{k,\ell}+\gamma_{k,\ell}(1+A_{1,0}),
B_{k,\ell}+\gamma_{k,\ell}B_{1,0}+\delta_{k,\ell}Q')
\end{equation}
is also solution, which gives two degrees of freedom at each step.
From Corollary \ref{UNISYS4}, \eqref{exct} is exactly the set of solutions  in the class \eqref{ABkl1}. Note that for $p=4$, it seems that one cannot use the parameters to avoid polynomial growth. For $p=2$, there is a choice of parameters giving no polynomial growth 
  corresponding to the explicit $2$-soliton solutions. See Section~5.2.

\medskip

\noindent\emph{Proof of Proposition \ref{SYSkl4}.} \quad The proof of Proposition \ref{SYSkl4} is based on Proposition \ref{SYS4} and on the structure of $F_{k,\ell}$ and $G_{k,\ell}$
(see Lemma \ref{STRUCT4}).
Using the induction argument described above, it is enough to check that if 
$(a_{k',\ell'},A_{k',\ell'},B_{k',\ell'})$ satisfies \eqref{ABkl1} for all 
$(k',\ell')\prec (k,\ell)$, then we can find $(a_{k,\ell},A_{k,\ell},B_{k,\ell})$ as in
\eqref{ABkl1} and solving $(\Omega_{k,\ell})$.
This will follow from Proposition \ref{SYS4} and the following Claim.

\begin{claim}\label{new24}
Let $(k,\ell)$ be such that $(k,\ell)\in \Sigma_0$, with $(k,\ell)\neq (1,0)$.
    Assume that for all $(k',\ell')\prec (k,\ell)$,
    the functions $A_{k',\ell'}$ and $B_{k',\ell'}$ verify \eqref{ABkl1}.
    Then,
    \begin{equation}\label{FGkl1}
    \begin{split}
       &  F_{k,\ell}(x)  = \overline{F}_{k,\ell}(x) + \widetilde{F}_{k,\ell}(x) + \varphi(x) \widehat{F}_{k,\ell}(x),\\
       & G_{k,\ell}(x) = \overline{G}_{k,\ell}(x) + \widetilde{G}_{k,\ell}(x) + \varphi(x) \widehat{G}_{k,\ell}(x),\quad \text{where}\\
       &\text{$\overline{F}_{k,\ell}$, $\overline{G}_{k,\ell}\in \mathcal{Y}$;
          $\overline{F}_{k,\ell}$ is odd and  $\overline{G}_{k,\ell}$ is even;}\\
       & \text{$\widetilde{F}_{k,\ell}$ and $\widehat{G}_{k,\ell}$ are odd polynomials;
        $\widehat{F}_{k,\ell}$ and $\widetilde{G}_{k,\ell}$ are even polynomials.}
    \end{split}
    \end{equation}
    Moreover, let $2\leq k\leq p-1$,  if for any $1\leq k' <k$,
        \begin{equation*}
            \deg \widetilde A_{k',0}=\deg \widehat A_{k',0}=\deg \widetilde B_{k',0}=\deg \widehat B_{k',0}=0 
            \quad \text{then}\quad F_{k,0}, \ G_{k,0}\in \mathcal{Y}.
        \end{equation*}
\end{claim}
Claim \ref{new24} is a consequence of the more detailed Lemma \ref{STRUCT4} proved in Appendix B.

\medskip
\emph{- Case    $k=1$, $\ell=0$.} \quad The system $(\Omega_{1,0})$ is explicit from Proposition \ref{SYSTEME}, indeed,
$F_{1,0}=p(Q^{p-1})'$ and $G_{1,0}=pQ^{p-1}$. Thus 
\begin{equation*}
    \text{$F_{1,0}=\overline F_{1,0}\in \mathcal{Y}$, $G_{1,0}=\overline G_{1,0}\in \mathcal{Y}$ and $\widetilde F_{1,0}=\widehat F_{1,0}=\widetilde G_{1,0}=\widehat G_{1,0}=0.$}
\end{equation*} 
It follows from Proposition \ref{SYS4} that $(\Omega_{1,0})$ has a solution
$a_{1,0}$, $A_{1,0}$, $B_{1,0}$ with the desired properties.
Moreover, from \eqref{DEGA1}--\eqref{DEGB2}, we obtain
$
    \widetilde A_{1,0}=\widehat A_{1,0}=\widetilde B_{1,0}=0 \quad \text{and}\quad
    \widehat B_{1,0}=b_{1,0},
$
where $b_{1,0}\in \mathbb{R}$ is a constant. Whether or not $b_{1,0}$ is zero will be determined in Section 3 for each case 
$p=2,$ $3$ and $4$.

\medskip

\emph{- Case    $2\leq k\leq p-1$, $\ell=0$.} \quad
By induction on $1\leq k\leq p-1$, we solve $(\Omega_{k,0})$ and we prove:
\begin{equation}\label{induc}
    \widetilde A_{k,0}=\widehat A_{k,0}=\widetilde B_{k,0}=0,\quad \widehat B_{k,0}=b_{k,0}\in \mathbb{R}.
\end{equation}
Indeed, let $2\leq k\leq p-1$ and assume that \eqref{induc} is true for all $1\leq k'< k$. Then, it follows
from   Claim \ref{new24}    that $F_{k,0}, G_{k,0}\in \mathcal{Y}$, which means that $\widetilde F_{k,0}=\widehat F_{k,0}=
\widetilde G_{k,0}=\widehat G_{k,0}=0$.
Therefore, from Proposition \ref{SYS4}, we solve $(\Omega_{k,0})$ with property \eqref{induc} at    rank $k$, which completes the induction argument.
Thus \eqref{debut} is proved. 

\medskip

\emph{- Case $k\geq p$, $\ell\geq 0$ or $k\geq 1$, $\ell\geq 1$.} \quad By induction on $(k,\ell)$, we prove that
$(\Omega_{k,\ell})$ has a solution $(a_{k,\ell},A_{k,\ell},B_{k,\ell})$ satisfying
 \eqref{ABkl1}.
First, note that \eqref{ABkl1} holds for $1\leq k\leq p-1$, $\ell=0$ by \eqref{induc}.
 From Claim \ref{new24}, we know that
$F_{k,\ell}$ and $G_{k,\ell}$ have the required structure to apply Proposition \ref{SYS4}, thus
we obtain a solution $(a_{k,\ell},A_{k,\ell}, B_{k,\ell})$ with the structure \eqref{ABkl1}. Thus, the induction argument is complete and the system $(\Omega_{k,\ell})$ is solved up to $(k_0,\ell_0)$.

\subsection{Estimate on $S(t,x)$ and proof of Theorem \ref{AVRIL4}}
We consider the function $v(t)$ constructed in \eqref{defv}--\eqref{defW}
where $(a_{k,\ell})$, $(A_{k,\ell})$ and $(B_{k,\ell})$ are defined in Proposition \ref{SYSkl4}.
For this choice, we have
\begin{equation}\label{page23} 
    S(t,x)   = 
     \sum_{\substack{ 1\leq k\leq 5k_{0}+12 \\ 0\leq \ell \leq 5\ell_{0}+4 \\ k> k_0 ~\rm{or}~ \ell> \ell_0}}
    c^\ell\left( Q_c^k(y_c)  F_{k,\ell}(y)   +  (Q_c^k)'(y_c) G_{k,\ell}(y)\right) . 
 \end{equation}
Recall that
$    q=\tfrac 1{p-1}-\tfrac 14,$
 and
$    T_c=c^{-\frac 12 -{\frac 1{100}}}.
$

\begin{proposition}[Estimates on $W$ and $S$]\label{VandS}
    Let $k_0\geq 1$, $\ell_0\geq 0$. There exists $K$ such that, for any  $0<c<1$, 
    for any $t\in   [-T_c,T_c]$, $W(t)$, $S(t)$ belong to $H^s(\mathbb{R})$ for all $s\geq 1$ and satisfy
    \begin{equation}\label{estV}
        \|W(t)\|_{H^1}=\|v(t)-R(t)-R_c(t)\|_{H^1}\leq K c^{\frac 1{p-1}},
    \end{equation}
    \begin{equation}\label{estS} j=0,1,2,\quad
        \|\partial_x^{(j)}S(t)\|_{L^2}\leq K c^{n_0}, 
    \end{equation}
    where $n_0=(\tfrac 12-\tfrac 1{100} ) \min(\frac {k_0}{p-1}, \ell_0+1)$.
\end{proposition}
Before proving Proposition \ref{VandS}, we claim several preliminary results.  The first result
concerns the degrees of the polynomials in the decomposition of $W(t)$.

\begin{lemma}\label{DEGG}
    a) For all $1\leq k \leq k_0$ and $0\leq \ell \leq \ell_0$
        such that $\frac k{p-1} + \ell \leq 2$,
    \begin{equation}\label{suite1}
            \deg \widetilde A_{k,\ell}=\deg \widehat A_{k,\ell}=0,
    \end{equation}
     b)    For all     $1\leq k\leq k_0$, $0\leq \ell \leq \ell_0$,
    \begin{equation}\label{DEGkl2}
d_{AB}(k,\ell)   = \max\left(\deg \widetilde{A}_{k,\ell},\deg \widehat{A}_{k,\ell},\deg \widehat{B}_{k,\ell}, \deg \widetilde{B}_{k,\ell}\right) \leq
        \frac {k-1}{p-1}+\ell .
    \end{equation}
\end{lemma}
\noindent\emph{Proof of Lemma \ref{DEGG}.}
The proof proceeds by induction  on $(k,\ell)$. 

\emph{- Case $k\geq p$, $\ell\geq 0$ or $k\geq 1$, $\ell\geq 1$.} \quad By induction on $(k,\ell)$, we prove:
\eqref{DEGkl2} holds.
First, note that \eqref{DEGkl2} holds for $1\leq k\leq p-1$, $\ell=0$ by \eqref{induc}.
Let
\begin{equation}\label{defxi}
    \xi(k,\ell)=\frac{k-1}{p-1}+\ell.
\end{equation}
Assume 
\begin{equation}\label{hypo}
    \text{for all $(k',\ell')\prec (k,\ell)$,}\quad  d_{AB}(k',\ell')\leq \xi(k',\ell') \text{ holds true.}
\end{equation} 
From Lemma \ref{STRUCT4}, we know that
$F_{k,\ell}$ and $G_{k,\ell}$ satisfy:
\begin{equation}\label{moreover}
    d_{FG}(k,\ell)  \leq \max \left( 
    d_{AB}(k{-}1, \ell){-}1, d_{AB}(k {-}p+{1}, \ell), d_{AB}(k, \ell{-}1), d_{\mathbf{N}}(k,\ell)\right).  
\end{equation}
We claim 
\begin{align*}
    & \text{if $k\geq p$}, \quad d_{AB}(k-p+1,\ell)\leq \xi(k,\ell)-1,\\
    & \text{if $\ell\geq 1$}, \quad d_{AB}(k,\ell-1)\leq \xi(k,\ell)-1,\quad \\
    & \text{if $k\geq p$}, \quad d_{\mathbf{N}}(k,\ell)\leq \xi(k,\ell)-1
    \quad (\text{$d_{\mathbf{N}}$ is defined in Lemma \ref{STRUCT4}}).
\end{align*}
Indeed,  assume $k\geq p$, then by \eqref{hypo},
\begin{equation}
    d_{AB}(k-p+1,\ell)\leq \xi(k-p+1,\ell)=\xi(k,\ell)-1.
\end{equation}
Similarly, if $\ell\geq 1$, by \eqref{hypo}
\begin{equation}
    d_{AB}(k,\ell-1)\leq \xi(k,\ell-1)=\xi(k,\ell)-1.
\end{equation}
Finally, if $k\geq p$ and if $k_j,\ell_j$ satisfy $\sum_{j=1}^p k_j\leq k$ and $\sum_{j=1}^p \ell_j\leq \ell$, then  \eqref{hypo} implies
\begin{align*}
    \sum_{j=1}^p d_{AB}(k_j,\ell_j)\leq \sum_{j=1}^p \xi(k_j,\ell_j)=\frac {k-p}{p-1}+\ell=\frac {k-1}{p-1}+\ell-1=\xi(k,\ell)-1.
\end{align*}
Thus, $d_{\mathbf{N}}(k,\ell)\leq \xi(k,\ell)-1.$

By \eqref{moreover}, we obtain
$
    d_{FG}(k,\ell)  \leq \xi(k,\ell)-1.
$
Using Proposition \ref{SYS4},   $(a_{k,\ell},A_{k,\ell}, B_{k,\ell})$ satisfies
\eqref{degres1}--\eqref{degres2}, i.e.
\begin{equation*}
    \deg_{AB}(k,\ell)\leq \deg_{FG}(k,\ell)+1\leq \xi (k,\ell).
\end{equation*}
Thus, the induction argument is complete and the system $(\Omega_{k,\ell})$ is solved up to $(k_0,\ell_0)$.

\medskip

We now prove \eqref{suite1} to finish the proof of Proposition \ref{SYSkl4}.

\smallskip

- \emph{Case $p\leq k \leq 2(p-1)$, $\ell=0$.}
We prove \eqref{suite1} for the case $\ell=0$ by induction on $k$ starting at $k=p$. For $k=p$ and $\ell=0$ we know that
for all $k'< p$, $\widetilde A_{k',0}=\widehat A_{k',0}=\widetilde B_{k',0}=0$ and $\deg    \widehat B_{k',0}=0$. Thus, by Lemma \ref{STRUCT4} (b), we have $F_{p,0}\in \mathcal{Y}$, and thus by Proposition \ref{SYS4}, $A_{p,0}\in \mathcal{Y}$, which means
$\widetilde A_{p,0}=\widehat A_{p,0}=0$. In the statement of Proposition \ref{SYSkl4}, we give a weaker statement $\deg \widetilde A_{p,0}=\deg \widehat A_{p,0}=0$, since we want that the rest of the estimate to be compatible with nonzero (constant) $\widetilde A_{p,0}$ (see Section 3). The induction argument from $p$ to $2(p-1)$ is done in the same way and we omit it.

\smallskip

- \emph{Case $1\leq k\leq p-1$, $\ell=1$.}
We also omit this case, since it is similar.

\begin{claim}[Estimate on $\alpha(s)$]\label{surALPHA}
    Let $\alpha(s)$ be the function defined in \eqref{defALPHA}. Then
    \begin{equation*}
        \forall s\in \mathbb{R},\quad |\alpha(s)|\leq K c^{\frac {1}{p-1}-\frac 12},\quad
        |\alpha'(s)|\leq K c^{\frac {1}{p-1}}.
    \end{equation*}  
    In particular, there exists $c_0>0$ so that for all $0<c<c_0$, for all $s\in \mathbb{R}$,
    $|\alpha'(s)|\leq \frac 12$.
\end{claim}

\noindent\emph{Remark.}\quad
From now on, we choose $c>0$ small enough so that $1+\alpha'>1/2$ for all $s\in \mathbb{R}$. 
\medskip

\noindent\emph{Proof of Claim \ref{surALPHA}.}\quad
We have
\begin{equation*}
    |\alpha(s)|\leq\sum_{(k,\ell)\in \Sigma_0}\Biggl| a_{k,\ell}   \int_0^s Q_{c}^k(s') ds'\Biggr|\leq
    \max_{(k,\ell)\in \Sigma_0} |a_{k,\ell}| \ 
    \times \sum_{(k,\ell)\in \Sigma_0}    c^\ell \int  Q_{c}^k \leq K \  \int  Q_{c}  .
\end{equation*}
Since $Q_c(s')=c^{\frac 1{p-1}} Q(\sqrt{c} \,s')$, 
$  |\alpha(s)|\leq   K  \int  Q_{c}  = K c^{\frac{1}{p-1}-\frac 12} \int  Q 
$  and similarly
$    |\alpha'(s)|\leq K c^{\frac {1}{p-1}}$.
 
 \begin{claim}[$H^1$-estimates]\label{SIX} 
    Let $0<c<1/2$.
    Let $f\in \mathcal{Y}$ and let $P$ be a polynomial function of degree $d$. Then,    for all $k\geq 1$, $\ell\geq 0$, for all $t\in  [-T_c,T_c]$,   
    \begin{align*}   &
         \| c^{\ell} Q_c^k(y_c) f(y) \|_{H^1} 
            +c^{-\frac 12}
             \| c^{\ell} (Q_c^k)'(y_c) f(y) \|_{H^1}\leq K c^{\frac k{p-1}+\ell } e^{-(1-c)\sqrt{c}\, |t|}.
   \\   &      \|c^{\ell} Q_c^k(y_c) P(y) \|_{H^1}\leq
         K c^{\xi(k,\ell) + q -\frac d2 (1+{\frac 1{50}})}+c^{-\frac 12}
             \|c^{\ell} (Q_c^k)'(y_c) P(y) \|_{H^1}\leq K c^{\xi(k,\ell) + q   -\frac d2 (1+{\frac 1{50}})}.
    \end{align*}
\end{claim}

\noindent\emph{Proof of Claim \ref{SIX}.}\quad
Let $f\in \mathcal{Y}$, so that $|f(y)|\leq K |y|^r e^{-|y|}$ on $\mathbb{R}$.
By $Q_c(x)=c^{\frac 1{p-1}}Q(\sqrt{c}x)\leq K c^{\frac 1{p-1}} e^{-\sqrt{c}|x|}$, we have
\begin{equation*}
    | c^\ell Q_c^k(y_c) f(y)|^2\leq K c^{\frac {2k}{p-1}+2\ell} e^{-2k\sqrt{c}|y_c|} |y|^{2r} e^{-2|y|} 
    \leq K c^{\frac {2k}{p-1}+2\ell} e^{-2\sqrt{c}|y_c|} |y|^{2r} e^{-2|y|}. 
\end{equation*}
Since $y_c=x+(1-c)t$ and $y=x+\alpha(y_c)$, we have
$y_c=y+(1-c)t-\alpha(y_c)$, and so by   Claim \ref{surALPHA},
\begin{equation*}
   \sqrt{c} |y_c|\geq \sqrt{c}((1-c)|t|-|y|-|\alpha(y_c)|)
   \geq    (1-c)\sqrt{c} |t|-\sqrt{c}|y| -K.
\end{equation*}
Thus,
\begin{equation*}
    | c^\ell Q_c^k(y_c) f(y)|^2
    \leq K c^{\frac {2k}{p-1}+2\ell} e^{-2(1-c)\sqrt{c} |t|}    |y|^{2r} e^{-2(1-\sqrt{c})|y|}
    \leq K c^{\frac {2k}{p-1}+2\ell} e^{-|y|}.
\end{equation*}
By changing the variable, $y=x+\alpha(x+(1-c)t)$, and using Claim \ref{surALPHA}, we have\begin{equation*}
    \int  e^{- |y|} dx=\int  e^{- |y|} \frac {dy}{1+\alpha'(y_c)}
    \leq 2 \int  e^{- |y|} dy \leq K.
\end{equation*}
Thus,
$
    \| c^\ell Q_c^k(y_c) f(y)\|_{L^2}\leq K c^{\frac k{p-1}+\ell} e^{-(1-c)\sqrt{c}\, |t|}.
$

Since $|Q_c'|\leq \sqrt{c}\, Q_c$ (recall $(Q'_c)^2=cQ_c^2 - \frac {2}{p+1}Q_c^{p+1}$), we also get
\begin{equation*}
    \| c^\ell (Q_c^k)'(y_c) f(y)\|_{L^2}\leq K c^{\frac k{p-1}+\ell+\frac 12} e^{-(1-c)\sqrt{c}\, |t|}.
\end{equation*}
Since $\partial_x(c^\ell Q_c^k(y_c) f(y))= c^\ell (Q_c^k)'(y_c) f(y)+(1+\alpha'(y_c)) c^\ell Q_c^k(y_c) f'(y)$, and $f'\in \mathcal{Y}$, the above estimates and Claim \ref{surALPHA} give the $H^1$ estimate on $ c^\ell Q_c^k(y_c) f(y)$.
The proof of the estimates for $c^\ell (Q_c^k)'(y_c) f(y)$ is similar.
\medskip

Now, we consider a monomial function $P(y)=y^d$.
For all $t\in [-T_c,T_c]$, and by Claim \ref{surALPHA},
\begin{equation*}
    y=y_c-(1-c)t + \alpha(y_c) \quad \text{and so} \quad |y|\leq |y_c|+ T_c + K c^{\frac 1{p-1}-\frac 12}
    \leq |y_c|+ K c^{-\frac 12 (1+{\frac 1{50}})}.
\end{equation*}
Therefore,
\begin{equation*}
    |c^{\ell} Q_c^k(y_c) P(y) |=c^{\ell} Q_c^k(y_c) |y|^d\leq K\left(c^\ell |y_c|^d Q_c^k(y_c)+c^{\ell-\frac d2 (1+{\frac 1{50}})}  Q_c^k(y_c) \right).
\end{equation*}
By $Q_c(x)=c^{\frac 1{p-1}}Q(\sqrt{c}x)$, 
\begin{equation*}
    \|c^\ell |y_c|^d Q_c^k(y_c)\|_{L^2}=c^{\ell+\frac k{p-1} - \frac d2 -\frac 14} \||x'|^d Q^k(x)\|_{L^2},
\end{equation*}
\begin{equation*}
    \|c^{\ell-\frac d2 (1+{\frac 1{50}})}  Q_c^k(y_c)\|_{L^2} =c^{\ell+\frac k{p-1} - \frac d2(1+{\frac 1{50}}) -\frac 14} \|Q^k\|_{L^2}.
\end{equation*}
Thus, 
$    \|c^{\ell} Q_c^k(y_c) P(y) \|_{L^2}\leq K c^{\ell+\frac k{p-1} - \frac d2(1+{\frac 1{50}}) -\frac 14} \| Q^k\|_{L^2}.
$
The other estimates are obtained in a similar way.

\medskip

\noindent\emph{Proof of Proposition \ref{VandS}.}\quad
From Claim \ref{SIX}, we claim sharp estimates on the terms in $W(t)$ and $S(t)$. These estimates are applied to prove Proposition \ref{VandS} and will be used again in the rest of this paper.

\begin{claim}[Estimates for terms in $W(t)$]\label{SEPT}
     For all $t\in [-T_c,T_c]$,
    \begin{enumerate}
        \item[{\rm (a)}] For all $1\leq k\leq p-1$, $\ell=0$,
        \begin{align}
            & \|Q_c^k(y_c)A_{k,0}(y)\|_{H^1}\leq K c^{\frac k{p-1}} e^{-(1-c)\sqrt{c}|t|}\label{debutSEPT} \\
            &  \|c^{\ell} (Q_c^k)'(y_c) B_{k,0}(y) \|_{L^2}
            + \tfrac 1 {\sqrt{c}}\|c^{\ell} \partial_x((Q_c^k)'(y_c) B_{k,0}(y)) \|_{L^2} 
            \leq K c^{\frac {k}{p-1}+\frac 14}.\label{debutSEPT2}
        \end{align}
        \item[{\rm (b)}] For all $1\leq k \leq k_0$ and $0\leq \ell \leq \ell_0$    such that $\frac k{p-1} + \ell \leq 2$,
        \begin{align}
            & \|c^{\ell} Q_c^k(y_c) A_{k,\ell}(y) \|_{L^2}
            + \tfrac 1 {\sqrt{c}}\|c^{\ell} \partial_x(Q_c^k(y_c) A_{k,\ell}(y)) \|_{L^2}
            \leq K c^{\xi(k,\ell)+q} \label{suite1SEPT} \\
            & \|c^{\ell} (Q_c^k)'(y_c) B_{k,\ell}(y) \|_{L^2}
            + \tfrac 1 {\sqrt{c}}\|c^{\ell} \partial_x((Q_c^k)'(y_c) B_{k,\ell}(y) )\|_{L^2}
            \leq K c^{\frac 12 (1-{\frac 1{50}})\xi(k,\ell)+q+\frac 12}.\label{suite1SEPT2}
        \end{align}
            \item[{\rm (c)}] For all     $1\leq k\leq k_0$, $0\leq \ell \leq \ell_0$,
        \begin{align}
            & \|c^{\ell} Q_c^k(y_c) A_{k,\ell}(y) \|_{L^2}
            + \tfrac 1 {\sqrt{c}}\|c^{\ell}
            \partial_x(c^{\ell} Q_c^k(y_c) A_{k,\ell}(y) )\|_{L^2}
            \leq K c^{\frac 12 (1-{\frac 1{50}})\xi(k,\ell)+q}\label{DEGkl2SEPT} \\
            & \|c^{\ell} (Q_c^k)'(y_c) B_{k,\ell}(y) \|_{L^2}
            + \tfrac 1 {\sqrt{c}}\|c^{\ell}
            \partial_x(c^{\ell} (Q_c^k)'(y_c) B_{k,\ell}(y)) \|_{L^2}
            \leq K c^{\frac 12 (1-{\frac 1{50}})\xi(k,\ell)+q+\frac 12}.\label{DEGkl2SEPT2}
        \end{align}
    \end{enumerate}    
\end{claim}

\noindent\emph{Proof of Claim \ref{SEPT}.}
By Proposition \ref{SYSkl4}, we have
\begin{equation*}
    A_{k,\ell}  = \overline{A}_{k,\ell}  + \widetilde{A}_{k,\ell}  + \varphi  \widehat{A}_{k,\ell} , \quad
    B_{k,\ell}  = \overline{B}_{k,\ell}  + \widetilde{B}_{k,\ell}  + \varphi  \widehat{B}_{k,\ell} , 
\end{equation*}
where 
 $\overline{A}_{k,\ell}$, $\overline{B}_{k,\ell}\in \mathcal{Y}$ and 
    $\widetilde{A}_{k,\ell}$, $\widehat{A}_{k,\ell}$, $\widehat{B}_{k,\ell}$ $\widetilde{B}_{k,\ell}$
    are polynomial functions satisfying (Proposition \ref{SYSkl4} (c)):
\begin{equation}\label{216}
    \max\left(\deg \widetilde{A}_{k,\ell},\deg \widehat{A}_{k,\ell},\deg \widehat{B}_{k,\ell}, \deg \widetilde{B}_{k,\ell}\right) \leq
    \frac {k-1}{p-1}+\ell=\xi(k,\ell) .
\end{equation}
By Claim \ref{SIX}, $c^{\ell} Q_c^k(y_c) A_{k,\ell}(y)$ and $c^{\ell} (Q_c^k)'(y_c) B_{k,\ell}(y)$ belong to $H^1$.

- Proof of (c).
From the estimates of Claim \ref{SIX} applied to $A_{k,\ell}$ and $B_{k,\ell}$,  we obtain from \eqref{216}, for all $t\in [-T_c,T_c]$,
\begin{equation*}
    \|c^{\ell} Q_c^k(y_c) A_{k,\ell}(y) \|_{L^2}\leq
    K c^{\xi(k,\ell)+q-\frac {1} 2(1+{\frac 1{50}} )\xi(k,\ell)}\leq
    K c^{\frac 12(1-{\frac 1{50}}) \xi(k,\ell)+q},
\end{equation*}
and
\begin{equation*}
    \|c^{\ell} (Q_c^k)'(y_c) B_{k,\ell}(y) \|_{L^2}\leq K c^{\frac 12(1-{\frac 1{50}}) \xi(k,\ell)+\frac 12+q}.
\end{equation*}
The proof for $\|\partial_x(c^{\ell} Q_c^k(y_c) A_{k,\ell}(y)) \|_{L^2}$
is the same, except that since 
$\partial_x(c^{\ell} Q_c^k(y_c) A_{k,\ell}(y))=
c^{\ell} (Q_c^k)'(y_c) A_{k,\ell}(y)+c^{\ell} Q_c^k(y_c) A_{k,\ell}'(y)$
there is a gain of $\sqrt{c}$ due to derivation of $Q_c$ , and of $c^{-\frac 12(1+
{\frac 1{50}})}$ due to derivation of polynomial terms in  
$A_{k,\ell}$ (see Claim \ref{SIX}).

- Proof of (a).
Note that by Proposition \ref{SYSkl4} (a) for all $1\leq k\leq p-1$, $\widetilde A_{k,0}=\widehat A_{k,0}=0$, which means $A_{k,0}\in \mathcal{Y}$ and thus for such $k$,
by Claim \ref{SIX}, for all $t\in \mathbb{R}$,
\begin{equation*}
    \|Q_c^k(y_c)A_{k,0}(y)\|_{H^1}\leq K c^{\frac k{p-1}} e^{-(1-c)\sqrt{c}|t|}.
\end{equation*}
For such $k$ and $\ell$, $\widetilde B_{k,0}=0$ and $\deg \widehat B_{k,0}=0$, and thus, for all $t\in [-T_c,T_c]$
\begin{equation*}
    \|c^{\ell} (Q_c^k)'(y_c) B_{k,0}(y) \|_{L^2}\leq K c^{\frac {k}{p-1}+\frac 14}.
\end{equation*}

- Proof of (b). From Proposition \ref{SYSkl4} (b)  and Claim \ref{SIX}, in the case $\frac k{p-1}+\ell \leq 2$ we obtain:
\begin{align*}
    \forall t\in [-T_c,T_c], \quad 
    & \|c^{\ell} Q_c^k(y_c) A_{k,\ell}(y) \|_{L^2}\leq K c^{\xi(k,\ell)+q}  \\
    & \|c^{\ell} (Q_c^k)'(y_c) B_{k,0}(y) \|_{L^2}\leq K c^{\frac 12 (1-{\frac 1{50}})\xi(k,\ell)+q+\frac 12}. 
\end{align*}

\begin{claim}[Estimates for terms in $S(t)$]\label{HUIT}
    For all $(k,\ell)$ satisfying $k_0+1\leq k\leq K_0$ or $\ell_0+1\leq \ell\leq L_0$,
    for all $t\in [-T_c,T_c]$,
        \begin{equation}
            \|c^{\ell} Q_c^k(y_c) F_{k,\ell}(y) \|_{H^1}+
            \|c^{\ell} (Q_c^k)'(y_c) G_{k,\ell}(y) \|_{H^1}
            \leq K c^{n_0}, \label{DEGkl2HUIT2}
        \end{equation}
	for $n_0=\frac 12 (1{-}{\frac 1{50}}) \min\big(\frac{k_0}{p-1},1+\ell_0\big)$.
\end{claim}
\noindent\emph{Proof of Claim \ref{HUIT}.}
Assume, for example, that $k\geq k_0+1$.
By Claim \ref{SIX}, for all $t\in [-T_c,T_c]$,
\begin{align*}
    &\|c^{\ell} Q_c^k(y_c) F_{k,\ell}(y) \|_{H^1}+\|c^{\ell} (Q_c^k)'(y_c) G_{k,\ell}(y) \|_{H^1} \\ & \leq
    K c^{\xi(k,\ell)+q-\frac {1} 2(1+{\frac 1{50}} )\xi(k_0+1,\ell)+\frac {1} 2(1+{\frac 1{50}} )}
   \leq
    K c^{\frac 12(1-{\frac 1{50}}) \xi(k,\ell) }\leq c^{\frac 12(1-{\frac 1{50}}) \frac {k_0}{p-1}}.
\qquad \Box
\end{align*}

\medskip

Recall $
    W(t,x)=\sum_{(k,\ell)\in \Sigma_0} 
        c^\ell\left(Q_{c}^k(y_{c}) A_{k,\ell}(y)+(Q_{c}^k)'(y_{c}) B_{k,\ell}(y)\right).
$
We apply the estimates of Claim \ref{SEPT} to each term of $W(t)$, for all $t\in [-T_c,T_c]$,

- For $1\leq k\leq p-1$ and $\ell=0$, we have
\begin{equation*}
    \|Q_c^k(y_c)A_{k,0}(y)\|_{H^1}\leq K c^{\frac k{p-1}} e^{-(1-c)\sqrt{c}|t|}\leq K c^{\frac 1{p-1}},
\end{equation*}

- For $k\geq p$ and $\ell\geq 0$, or $k\geq 1$ and $\ell\geq 1$, we have $\xi(k,\ell)\geq 1$; \begin{equation*}
     \|c^\ell Q_c^k(y_c)A_{k,\ell}(y)\|_{H^1}\leq K c^{\frac 12(1-{\frac 1{50}}) \xi(k,\ell)+q}\leq K c^{\frac 14 +q}= K c^{\frac 1 {p-1}},
\end{equation*}
and similarly for $ \|c^\ell (Q_c^k)'(y_c)B_{k,\ell}(y)\|$.

Thus, for all $t\in [-T_c,T_c]$,
$\|W(t)\|_{H^1}\leq K c^{\frac 1 {p-1}}$.

\medskip

By \eqref{page23}, for a given $k_0\geq 1$ and $\ell_0\geq 0$, the rest $S(t,x)$ contains only terms for $k\geq k_0+1$ or terms for $\ell\geq \ell_0+1$.
Thus, from Claim \ref{SEPT}, for all $t\in [-T_c,T_c]$,
$  \|S(t)\|_{H^1}\leq K c^{n_0},
$
where $n_0=\frac 12 (1{-}{\frac 1{50}}) \min(\frac {k_0}{p-1}, \ell_0+1)$.
The proof for $\|\partial_{x}^j S(t)\|_{H^1}$, for   $j=1,2$, is the same.

\section{Recomposition of the approximate solution at $\pm T_c$}

In this section, we consider the function $v$ defined in Theorem \ref{APSOL4}. We prove further properties of $v$ by solving explicitely the first two systems $(\Omega_{1,0})$ and $(\Omega_{2,0})$. Detailled properties depend on the specific value of $p=2$  or $4$. 

\subsection{Explicit resolution of the first systems}

\noindent\emph{1. Resolution of the systems $(\Omega_{1,0})$, $(\Omega_{2,0})$ for $p=2,  4$.}
We begin with two technical results.
 
\begin{claim}[Expression of $V_1$]\label{ClaimVun}
Let $V_1\in \mathcal{Y}$, even, be solution of 
$
    \mathcal{L} V_1= p Q^{p-1}.
$ Then  $ V_1=-2 Q -xQ'$ for $p=2$  and
 $V_1=\frac 13\Bigl(Q'{\left(\textstyle \int_0^x Q^2\right)} - 2 Q^3\Bigr)$ for $p=4$.
\end{claim}
\noindent\emph{Proof.}\quad For $p=2$,  set $V_1=-2 Q -xQ'$. Then, using the equation of $Q$
\begin{equation*}
    \mathcal{L} V_1 =-V_1'' + V_1 -2 Q V_1 =(2 Q''+ 2 Q''+ x Q^{(3)}) -2 Q - x Q' + 4 Q^2 + 2 x Q Q' = 2 Q.
\end{equation*}
Now, let $p=4$. By $\mathcal{L} (fg) = g (\mathcal{L} f) -2 f' g'    - f g'' $, we have:
\begin{equation*}
    \mathcal{L} (Q'{\left(\textstyle \int_0^x Q^2\right)})= {\left(\textstyle \int_0^x Q^2\right)}   \mathcal{L} Q' - 2 Q'' Q^2 - 2(Q')^2 Q, 
\end{equation*}
but from Lemma \ref{surL}, $\mathcal{L} Q'=0$, so that by $Q''=Q-Q^4$ and $(Q')^2=Q^2-\frac 25 Q^5$,
\begin{equation*}
    \mathcal{L} (Q'{\left(\textstyle \int_0^x Q^2\right)})= - 2 Q^3 + 2 Q^6 - 2 Q^3 + \frac 45 Q^6 = -4 Q^3 + \frac {14} 5 Q^6. 
\end{equation*}
We also have:
$    \mathcal{L} Q^3   = - 3Q'' Q^2 - 6 (Q')^2 Q + Q^3 - 4 Q^6=  - 8 Q^3 + \frac 75 Q^6.
$
Thus, by combining these two calculations, we get
$
   \frac 13 \mathcal{L} ( Q' {\left(\textstyle \int_0^x Q^2\right)} - 2 Q^3)= 4 Q^3.
$

\begin{claim}[Computation of $\int Z_1 Q$]\label{surZ1Q}
Let $Z_1=3 V_1'' +p Q^{p-1} V_1 + pQ^{p-1}$. Then
    \begin{equation*}
        \int Z_1 Q= \frac{p-3}{2(p-1)} \int Q.
    \end{equation*}
\end{claim}
\noindent\emph{Proof of Claim \ref{surZ1Q}.}\quad 
\begin{align*}
    \int Z_1 Q   &= \int (3 V_1'' +p Q^{p-1} V_1 + pQ^{p-1}) Q
    = \int V_1 (3 Q'' + p Q^p) + p \int Q^p.
\end{align*}
Since $\mathcal{L} Q =  -(p-1) Q^p$ and $\mathcal{L} (\frac 2{p-1} Q + x Q')=-2 Q$, we have
\begin{align*}
    3 Q''+p Q^p     &=3 Q +(p-3) Q^p=\mathcal{L}(-\tfrac 32 (\tfrac 2{p-1} Q + x Q') - \tfrac {p-3}{p-1} Q) = - \mathcal{L}(\tfrac{p}{p-1} Q + \tfrac 32 x Q').
\end{align*}
Thus,
\begin{align*}
    \int Z_1 Q   &= - \int V_1 \mathcal{L}(\tfrac{p}{p-1} Q +\tfrac 32 x Q')+ p \int Q^p= - \int (\mathcal{L} V_1) (\tfrac{p}{p-1} Q +\tfrac 32 x Q')  + p\int Q^p\\
                &= - p \int Q^{p-1} (\tfrac{p}{p-1} Q +\tfrac 32 x Q')+p\int Q^p =\frac {p-3}{2(p-1)}\int Q^p,
\end{align*}
by integration by parts. Since $Q=Q^p+Q''$, we have $\int Q^p=\int Q$, 
and the claim follows. 

\begin{lemma}[Resolution of the first systems for $p=2, 4$]\label{RESUME}
$\bullet$ For $p=2$,
    \begin{align*}
        &a_{1,0}=\frac 23,\quad b_{1,0}=-2,\quad A_{1,0}=-\frac 43 Q,\quad B_{1,0}=-2 \varphi.\\
        & 
    a_{2,0}=-\frac 49,\quad a_{1,1}=\frac 23,\quad A_{2,0}=-2 + \frac 43 Q,\quad A_{1,1}=2 - \frac 23 Q -\frac 13 xQ',
        \quad b_{2,0}=\frac 43.\end{align*}
$\bullet$ For $p=4$, 
        \begin{align*}
         &a_{1,0}=-2 \frac {\int Q}{\int Q^2},\quad b_{1,0}=-\frac 12 \int Q^3 +\frac 16 \frac{\bigl(\int Q\bigr)^2}{\int Q^2}<0,\\
        & A_{1,0}=\frac 13 (Q'{\left(\textstyle \int_0^x Q^2\right)} - 2 Q^3) +2 \frac {\int Q}{\int Q^2} (-\tfrac 13 Q -\tfrac 32 x Q'),\\
        &    b_{2,0}=-\frac 1{18} \int Q^2 - \frac 34 \frac {\bigl(\int Q\bigr)\bigl(\int Q^3\bigr)}{\int Q^2} 
            -\frac 1{18} \frac {\bigl(\int Q\bigr)^3}{\bigl(\int Q^2\bigr)^2}
<0.
    \end{align*}    
\end{lemma}
From Corollary \ref{UNISYS4}, there are several solutions. The choice of the solution for $p=2$
above is related to the exact $2$-soliton solutions.

We only solve $(\Omega_{1,0})$ in this paper. The resolution of the next systems is done in \cite{MMb2}.

\medskip

\noindent\emph{Proof of Lemma \ref{RESUME}.} 
From Proposition \ref{SYSTEME} and Proposition \ref{SYSkl4}, the system $(\Omega_{1,0})$ writes, for $p=2$, $3$ and $4$:
\begin{equation*}
    (\Omega_{1,0})\quad  \left\{
    \begin{array}{l}
         \mathcal{L} A_{1,0} + a_{1,0} (3Q -2 Q^p) = pQ^{p-1} \\
         (\mathcal{L} B_{1,0})' + a_{1,0} (3Q'') - 3 A_{1,0}'' -p Q^{p-1} A_{1,0}=pQ^{p-1}.
    \end{array}
    \right.
\end{equation*}
- Computation of $A_{1,0}$.
Recall from Claim \ref{surE}
$
    V_0=-\frac 1{p-1}Q-\frac 32 x Q',$ {and} $\mathcal{L} V_0=3Q - 2 Q^p.
$
Thus, the function $A_{1,0}=V_1- a_{1,0} V_0$
solves the first line of $(\Omega_{1,0})$, independently of the value of $a_{1,0}$.
By replacing $A_{1,0}$ in the second line of the system $(\Omega_{1,0})$,
\begin{equation*}
    (\mathcal{L} B_{1,0})'+ a_{1,0} Z_0=Z_{1},
\end{equation*}
where
\begin{equation}\label{defZs}
    Z_0=3 Q''+ 3 V_0''+pQ^{p-1} V_0, \quad 
    Z_1=3 V_1'' +p Q^{p-1} V_1 + pQ^{p-1}.
\end{equation}
- Computation of $a_{1,0}$.
Since $\mathcal{L} Q'=0$, we have $\int (\mathcal L B_{1,0})' Q=0$ and so
$a_{1,0} \int Z_0 Q= \int Z_1 Q$.
Recall that by Claim \ref{surE},
$
    \int Z_0 Q = \frac {p-5} {4(p-1)} \int Q^2.
$
Assuming this, we obtain
\begin{equation}\label{fa1}
    a_{1,0}=2 \,\frac {p-3}{p-5}\, \frac {\int Q}{\int Q^2}.
\end{equation}

$\bullet\ p=2.$
Since $Q=Q^2+Q''$, we have $\int Q=\int Q^2$ and so \eqref{fa1} gives $a_{1,0}=\frac 23$.
Next, $A_{1,0}=V_1-\frac 23 V_0=-2 Q-xQ'-\frac 23 (-Q-\frac 32 xQ')=-\frac 43 Q$.

By the second line of the system $(\Omega_{1,0})$, we get
\begin{equation*}
    (\mathcal{L} B_{1,0})'=2 Q -(3 a_{1,0} Q'' - 3 A_{1,0}''- 2 Q A_{1,0})=
    2Q -(2 Q''+4 Q''+\frac 83 Q^2)=- 4 Q    +\frac {10} 3 Q^2.
\end{equation*}
From Claim \ref{surphi}, we have $(\mathcal{L} \varphi)' = 2 Q -\frac 53 Q^2$, thus,
$B_{1,0}=-2 \varphi$ is solution.

From Proposition \ref{SYSTEME}, we write the following two systems for $p=2$.
\begin{equation*}
    (\Omega_{2,0})\quad \left\{
    \begin{aligned}
        (\mathcal{L} A_{2,0})'+a_{2,0}(3 Q- 2 Q^2)'&     = (-A_{1,0}+A_{1,0}^2)'- (3 B_{1,0}''+ 2Q B_{1,0}) 
        \\ &    - a_{1,0}(Q+3 A_{1,0}''+ 2 Q A_{1,0})'+ 3a_{1,0}^2 Q^{(3)}, \\
         (\mathcal{L} B_{2,0})'+3 a_{2,0} Q'' - 3 A_{2,0}'' -2 Q A_{2,0}& = A_{1,0}+A_{1,0}^2+(-2 B_{1,0}+A_{1,0}B_{1,0})'
        \\ &        - \tfrac 12 {a_{1,0}}    (9 A_{1,0}' + 3B_{1,0}'' + 2Q B_{1,0})'+\tfrac 32 a_{1,0}^2 Q''.
    \end{aligned}\right.
\end{equation*}
\begin{equation*}
    (\Omega_{1,1}) \quad  \left\{
    \begin{aligned}
                &(\mathcal{L} A_{1,1})'+a_{1,1}(3 Q- 2 Q^2)' = 3A_{1,0}'+3 B_{1,0}''+ 2Q B_{1,0} \\
                &(\mathcal{L} B_{1,1})'+3 a_{1,1} Q'' - 3 A_{1,1}'' -2 Q A_{1,1} = 3 B_{1,0}'.
    \end{aligned}\right.
\end{equation*}
The resolution of these two systems is done in \cite{MMb2}.

\medskip

$\bullet\ p=4$.  From \eqref{fa1}, we obtain the expression of $a_{1,0}$, and from $A_{1,0}=V_1-a_{1,0} V_0$, and the expressions of $V_1$ and $V_0$, we obtain $A_{1,0}.$
Here $a_{1,0}<0$ which will have a surprizing consequence on the shift of $Q$ after collision (see Proposition \ref{AVRIL4}).

Next, 
$B_{1,0}$ is of the form $B_{1,0}=\overline B_{1,0} + \varphi b_{1,0}$, where $\overline B_{1,0}\in \mathcal{Y}$ and $b_{1,0}\in \mathbb{R}$ from Proposition \ref{SYS4}.
We do not compute $B_{1,0}$ in this case. Thus we only need to compute $b_{1,0}$.
By Claim \ref{surphi},
$2 b_{1,0}=\lim_{+\infty} B_{1,0} -\lim_{-\infty} B_{1,0}=\lim_{+\infty} \mathcal{L} B_{1,0}-\lim_{-\infty} \mathcal{L} B_{1,0}$.
Recall the equation of $B_{1,0}$:
$
    (\mathcal{L}B_{1,0})' =Z_1 - a_{1,0} Z_0,
$
where $Z_0=3 Q''+3 V_0'' + 4Q^3 V_0$ and $Z_1=4Q^3+3 V_1''+ 4Q^3 V_1$.
It follows that $2 b_{1,0}= \int Z_1 - a_{1,0} \int Z_0$.
By integration by parts,
\begin{equation*}
    \int Z_0=4 \int Q^3 V_0= - 4 \int Q^3 (\tfrac 13 Q + \tfrac 32 xQ')=\tfrac 16 \int Q^4 =\tfrac 16 \int Q,
\end{equation*}
\begin{equation*}
    \int Z_1=    4 \int Q^3 + \tfrac 4 3 \int (Q^3Q' {\left(\textstyle \int_0^x Q^2\right)} - 2 Q^6)   = 4 \int Q^3 - 3 \int Q^6= -\int Q^3
\end{equation*}
since $3 \int Q^6=5 \int Q^3$ (from the equation of $Q$).
Thus, $2 b_{1,0}=-\int Q^3 - \frac {a_{1,0}}{6}\int Q$, which gives the desired formula.

We justify that $b_{1,0}<0$. By Cauchy--Schwarz' inequality and Claim \ref{LemmaA1}, we have
\begin{equation*}
    \int Q= \int Q^4\leq \bigg(\int Q^2\bigg)^{1/2}\bigg(\int Q^6\bigg)^{1/2}=\sqrt{\frac 53} \bigg(\int Q^2\bigg)^{1/2}\bigg(\int Q^3\bigg)^{1/2}.
\end{equation*}
Thus, $
\frac 16 \frac {\bigl(\int Q\bigr)^2}{\int Q^2}\leq \frac {5}{18} \int Q^3
$ and so $b_{1,0}\leq-\frac 29 \int Q^3.
$
Numerically $b_{1,0}\sim -0.9.$

System $(\Omega_{2,0})$  for $p=4$ writes:
\begin{equation*}
    \left\{
    \begin{aligned}
        (\mathcal{L} A_{2,0})'+a_{2,0}(3 Q- 2 Q^4)'  
        &= (6 Q^2(1+A_{1,0})^2)'\\ 
        &  - a_{1,0}(4 Q^3+3 A_{1,0}''
        +4 Q^3 A_{1,0})'+ 3 a_{1,0}^2 Q^{(3)}, \\
        (\mathcal{L} B_{2,0})'+3 a_{2,0} Q'' - 3 A_{2,0}'' - 4 Q^3 A_{2,0}  
        &= 6 Q^2(1+A_{1,0})^2+ (6Q^2 B_{1,0}(1+A_{1,0}))'
        \\  &           -\tfrac 12 {a_{1,0}} (9 A_{1,0}'+3B_{1,0}''+ 4Q^3 B_{1,0})' +\tfrac 32 a_{1,0}^2 Q''.
    \end{aligned}\right.
\end{equation*}
The fact that $b_{2,0}\not =0$ can easily be checked by solving
$(\Omega_{2,0})$ numerically. However, we were able to give an explicit expression of $b_{2,0}$, by long but elementary calculations, see \cite{MMb2}.

\medskip

\noindent \emph{2. Determination of all solutions of  $(\Omega)$.}
Now, let us justify the remark following Proposition \ref{SYSkl4} concerning the existence of several solutions of system $(\Omega_{k,\ell})$. At each step of resolution, the number of solutions of $(\Omega_{k,\ell})$ is related to the existence of nontrivial solutions of the system $(\Omega_0)$
\begin{equation*}
    (\Omega_0)\quad 
    \left\{
    \begin{array}{l}
         (\mathcal{L} A_0)' + a_0 (3Q -2 Q^p)' = 0 \\
         (\mathcal{L} B_0)' + a_0 (3Q'') - 3 A_0'' -p Q^{p-1} A_0= 0.
    \end{array}
    \right.
\end{equation*}

\begin{corollary}\label{UNISYS4}
    Assume that  $(a_0,A_0,B_0)$ solves the system $(\Omega_0)$, where $A_0$ is a $C^\infty$ even function, with at most polynomial growth at $\infty$, and $B_0$ is a $C^\infty$ odd function, with at most polynomial growth at $\infty$. Then, there exists $\gamma\in \mathbb{R}$ and $\delta \in \mathbb{R}$ such that
    \begin{equation}\label{allsol}
        (a_0,A_0,B_0)=(\gamma a_{1,0},\gamma (1+A_{1,0}), \gamma B_{1,0} +\delta Q').
    \end{equation}
    Conversely,  for any $\gamma, \delta\in \mathbb{R}$, \eqref{allsol} defines a solution of $(\Omega_0)$.
\end{corollary}

\noindent\emph{Proof of Corollary \ref{UNISYS4}.}\quad
The first line of $(\Omega_0)$ is equivalent to 
\begin{equation*}
    \mathcal{L} A_0 + a_0 (3Q -2 Q^p)=\gamma,
\end{equation*}
where $\gamma$ is a constant.
Since $\mathcal{L} 1= 1-p Q^{p-1}$, we have
$\mathcal{L} (A_0-\gamma + a_0V _0-\gamma V_1)=0.$
Claim  \ref{claimzz} implies that if $\mathcal{L} f=0$ where $f$ is a function with at most polynomial growth, then $f\in L^2(\mathbb{R})$, and so $f=\overline \delta Q'$.
Since $A_0$ is even, and has at most polynomial growth,  we obtain 
\begin{equation*}
    A_0= \gamma V_1 - a_0 V_0+\gamma.
\end{equation*}
The resolution of the second line of the system is similar to the previous calculations on $(\Omega_{1,0})$:
\begin{equation*}
    (\mathcal{L} B_0)'= \gamma Z_1 - a_0 Z_0,
\end{equation*}
which gives a relation between $a_0$ and $\gamma$:
$a_0 \int Z_0 Q = \gamma \int Z_1 Q$, which means that $a_0 = \gamma a_{1,0}$, and so $A_0=\gamma (1+A_{1,0})$. Thus,
$(\mathcal{L} B_0)'= \gamma (Z_0-a_{1,0}Z_1)=\gamma (\mathcal{L} B_1)'$,
and so $\mathcal{L} (B_0-\gamma B_{1,0})=0$ by parity. Therefore, $B_0=\gamma B_{1,0}+\delta Q'$.
 
\subsection{Asymptotics of the approximate solution at $\pm T_c$}

So far, we have searched an approximate solution $v$ on $[-T_c,T_c]$ with a structure
adapted to the interaction problem. For $t\in [-T_c,T_c]$, 
$v(t)=Q(y)+Q(y_c)+W(t)$, with
$\|W(\pm T_c)\|_{H^1}\leq K c^{\frac 1{p-1}}\sim K c^{\frac 14}
\|Q_c\|_{H^1}.$
Nevertheless, since the functions $A_{k,\ell}$, $B_{k,\ell}$ may contain polynomial functions of degree larger than $1$,   the previous decomposition is not adapted for  $t>T_c$. 

At $t=T_c$, we note that
$y_c\sim x +T_c$ and $y\sim x-\tfrac \Delta 2$, where $|\Delta|/{T_c}\ll 1$.
Thus $v(T_c)$ is close to the sum of two exponentially decoupled solitons, and for $t>T_c$,
one can use asymptotic techniques (see Section 4) close to $2$-soliton solutions, or equivalently
close to the sum of two solitons for the proofs. This set of $2$-soliton solutions have several parameters, as the size and the position of each soliton. In this section, we understand what is the optimal
choice for these parameters. In fact, at the formal level, from the decomposition, the size parameters will not be changed, we will concentrate on the position parameters.

\medskip

First, we   point out that the function $v(t,x)$ is, as the (gKdV) equation, invariant by the transformation $x\to -x$, $t\to -t$.
Indeed, $y_c(-t,-x)=-y_c(t,x)$, $\alpha(-s)=-\alpha(s)$, $y(-t,-x)=-y(t,x)$ and 
\begin{equation*}
v(-t,-x)=Q(-y)+Q_c(-y_c)+\sum_{(k,\ell)\in \Sigma_0} 
        c^\ell\left(Q_{c}^k(-y_{c}) A_{k,\ell}(-y){+}(Q_{c}^k)'(-y_{c}) B_{k,\ell}(-y)\right)=v(t,x),
\end{equation*}
by the parity properties of the functions $Q$, $Q_c$, $A_{k,\ell}$ and $B_{k,\ell}$.
Thus it suffices to study the properties of the function $v$ for $t\geq 0$.

\medskip

Let us present formal computations to recompose $v(T_c)$ in terms of the asymptotic
$2$-soliton family at $t\to +\infty$. 
We first observe that $Q$ and $Q_c$ are well-ordered and located far away in the original space variable $x$ at $t=T_c$.
Indeed, if $x>-T_c/2$, then $y_c=x+(1-c)t>T_c/4$ (say, $0<c<1/4$), thus the soliton $Q_c$ is at the left of $x=-T_c/2$. Conversely, if $x<-T_c/2$, then $y=x-\alpha(y_c)<-T_c/4$ for $c$ small and thus the soliton $Q$ is at the right of $x=-T_c/2$.

\medskip

\noindent\emph{1. Position of $Q$ at $t=T_c$ (for $p=2 ,4$).}

\smallskip

We determine the position of $Q(y)$, and thus we consider $x>-T_c/2$. For such $x$, $\sqrt{c} \,y_c\geq
\sqrt{c} T_c/4\gg 1$, and so
$\alpha(y_c)= \int_{0}^{y_c} \beta(s) ds \sim \int_{0}^{+\infty} \beta(s) ds$.
 Since
$    \int_{0}^{\infty}Q_{c}^k(s) ds=\frac 12 c^{\frac k{p-1}-\frac 12} \int Q^k $  we obtain
    \begin{equation*}
    \alpha(y_c)\sim     \sum_{(k,\ell)\in \Sigma_0} a_{k,\ell} \, c^\ell \int_{0}^{\infty}Q_{c}^k(s) ds=\frac 12 \sum_{(k,\ell)\in \Sigma_0}
    a_{k,\ell} \, c^{\frac k{p-1}+\ell-\frac 12} \int Q^k.
\end{equation*}
This means that at $t=T_c$, the  soliton $Q$ is located at $x=\tfrac \Delta 2$, where
\begin{equation}\label{defDelta}
    \Delta= \sum_{(k,\ell)\in \Sigma_0} a_{k,\ell} \, c^{\frac k{p-1}+\ell-\frac 12} \int Q^k.
\end{equation}
By symmetry, at $t=-T_c$, the soliton $Q$ is located at $x=-\tfrac \Delta 2$. Thus, as a consequence of the interaction with the small soliton $Q_c$, the large soliton $Q$ is shifted by $\Delta$ defined by \eqref{defDelta}.

\medskip

\noindent\emph{2. Position of $Q_c$ at $t=T_c$ (for $p=2 ,4$).}

\smallskip

For the soliton $Q$, we have introduced the variable $y$ depending on $x$ and $t$ which follows
the trajectory of $Q$ and in particular the shift phenomenon. On $Q_c$, the variable $y_c=x+(1-c) t$ does not catch any shift of the trajectory of $Q_c$. 
However, in the integrable cases, it is known that the small soliton is also shifted by through the interaction. 
In fact, the shift on $Q_c$ is to be determined by examinating the rest of the expansion of $v$.
Since we want to locate the soliton $Q_c$ at $t=T_c$, we consider $x<-T_c/2$. In particular, 
$y=x-\alpha(y_c) <-T_c/4$, for $c$ small. 
Recall from Proposition \ref{SYSkl4} that $A_{1,0}, A_{2,0}\in \mathcal{Y}$, and at $t=T_c$, $B_{1,0}\sim \overline B_{1,0}-b_{1,0}$,
where $\overline B_{1,0}\in \mathcal{Y}$. Thus 
\begin{equation*}
\begin{split}
 Q_c(y_c)+W(T_c) &\sim (1+A_{1,0}(y)) Q_c(y_c)+A_{2,0}(y) Q_c^2(y_c)+B_{1,0}(y)Q_c'(y_c)
 \\ & \sim
 Q_c(y_c)-b_{1,0} Q_c'(y_c)\sim Q_c(y_c-b_{1,0}).
 \end{split}
\end{equation*}
Thus,
\begin{equation}\label{devv1}
    v(T_c,x)\sim Q(x-\tfrac \Delta 2)+ Q_c(y_c-b_{1,0}).
\end{equation}
By the symmetry $x\to -x$, $t\to -t$,
the value $-2 b_{1,0}$ can be interpreted as the first order of the shift $\Delta_c$ on the soliton $Q_c$. Thus, we can set
\begin{equation*}
\Delta_c=2 b_{1,0}.
\end{equation*}

\noindent\emph{3. The integrable cases $p=2,3$.}

$\bullet$    $p=2$. 
In this case, we consider the explicit $2$-soliton solution with speeds $1$ and $0<c<1$
 defined in \eqref{NSOL}.
It is classical to observe that   
for $t$ large, for $x\in \mathbb{R}$,
\begin{equation*}
U_2(t,x)\sim Q(x-t-\Delta')+Q_c(x-ct),\quad U_2(-t,x)\sim Q(x+t)+Q_c(x+ct+\Delta_c'),
\end{equation*}
where $\Delta'=-\log \alpha(c)>0$ and $\Delta_c'=- \frac 1{\sqrt{c}} \Delta'$.

Let us check that the function $v$ can be chosen to match the explicit $2$-soliton
at the main orders at $T_c$.
We are not able to check all the   relations up to any $k_0$, $\ell_0$ by an algebraic argument. However, one can expect that there exists a function 
$v$ matching precisely at any order the explicit $2$-soliton solution.

\noindent First, let us check that the shifts are matching
$
    \Delta'\sim 4 \sqrt{c} + \frac 43 c \sqrt{c} , $ $ 
    \Delta_c'= -\frac 1 {\sqrt{c}} \Delta'.
$
From \eqref{defDelta},
$$
\Delta\sim (a_{1,0} \int Q) \sqrt{c} + \Big(a_{2,0} \int Q^2+a_{1,1}\int Q\Big) c^{3/2},
\quad 
\Delta_c \sim 2 b_{1,0}.
$$
From $\int Q=\int Q^2=6$, $a_{1,0}=\frac 23$, $a_{2,0}+a_{1,1}=\frac 29$ (Lemma \ref{RESUME})
and $b_{1,0}=-2$ (Lemma \ref{RESUME}), $\Delta'$ and $\Delta$
(and $\Delta'_c$ and $\Delta_c$) math at the first main order.

Now, we check that $v(T_c)$ matches the $2$-soliton solution at the main orders.
Note that 
\begin{equation*}\begin{split}
    Q_c(y_c-b_{1,0})& \sim Q_c(y_c)- b_{1,0} Q_c'(y_c) +\tfrac 12 b_{1,0}^2 Q_c''(y_c)-\tfrac 16 b_{1,0}^3 Q_c^{(3)}(y_c)\\&
\sim Q_c(y_c)- b_{1,0} Q_c'(y_c) +\tfrac 12 b_{1,0}^2 (c Q_c -Q_c^2) (y_c)-\tfrac 16 b_{1,0}^3 (cQ_c - Q_c^2)'(y_c),
\end{split}\end{equation*}
since  $Q_c''(y_c)=cQ_c-Q_c^2$, and $Q^{(3)}=(cQ_c-Q_c^2)'$.
From the decomposition of $v$ at $t=T_c$, we have on the other hand
$$
v(T_c)-Q(y)\sim Q_c(y_c) -b_{1,0} Q_c'(y_c) +\tilde A_{2,0} Q_c^2(y_c)-
b_{2,0} (Q_c^2)'(y_c)+\tilde A_{1,1} cQ_c(y_c)-
b_{1,1} c(Q_c)'(y_c).
$$
From Lemma \ref{RESUME}, we have
$
    \widetilde A_{2,0}=-\tfrac 12 b_{1,0}^2=-2,$
$\widetilde A_{1,1}=\tfrac 12 b_{1,0}^2=2,$ $ b_{2,0}=b_{1,1}=\tfrac 16 b_{1,0}^3 = \frac {4}3,
$ 
and again the two functions match at $t=T_c$ at this order.
 
 \medskip

\noindent\emph{4. Case $p=4$.}
In this case, we recall that no explicit $2$-soliton solution is known, nor was any approximate solution. In the next sections, by analytical methods, we will use the function $v$ to describe any solution close in large time to the sum of two solitons $Q$, $Q_c$ for $c$ small. Therefore, the function $v$ really describes the interaction between a soliton $Q$ and a soliton $Q_c$ for $p=4$. In particular, from Eq. \eqref{defDelta} and Lemma  \ref{RESUME},
\begin{equation}\label{D1quatre}
    (p=4)\qquad \Delta\sim - 2\, \frac 1{c^{1/6 }} \,\frac {\big(\int Q\big)^2}{\int Q^2}.
\end{equation}
Eq. \eqref{D1quatre} is surprizing for two reasons.
First, the value of the shift is negative. This means that for $p=4$, the large soliton $Q$ is shifted to the left by interaction with the small soliton $Q_c$. This is in contrast with the two previously known situations $p=2$ and $p=3$, where the shift is positive.

The second surprize is that the shift becomes infinite as $c\to 0$.  Therefore, the smaller is $c$, the larger is the influence of $Q_c$ on the trajectory
 of $Q$. 
To obtain the next order of the shift $\Delta$ for $p=4$, it is sufficient to compute $a_{2,0}$ from Lemma \ref{RESUME}.
However, note that the next order is $c^{1/6}$ ($k=2$ and $\ell=0$) and thus it corresponds to a small perturbation of $\Delta$ as $c$ is small.

 The function $v$ also allows us to determine for the first time the shift    $\Delta_c$ on the small soliton.
From Lemma \ref{RESUME}, it is at the first order $\Delta_c= 2 b_{1,0}<0$. Thus, the small soliton is also shifted to the left through the interaction, for $c$ sufficiently small, as for $p=2,3$.

As for the case $p=2$, we note that $Q_c(y_c-b_{1,0})=Q(y_c)-b_{1,0}Q_c'(y_c)+\frac 12 b_{1,0}^2 Q_c''(y_c)+...$
Since $Q_c''=cQ_c-Q_c^4$,
we obtain at $t=T_c$,
\begin{equation*}
	W(T_c)\sim Q(x-\tfrac \Delta 2)+Q(y_c)-b_{1,0}Q_c'(y_c)+\frac 12 b_{1,0}^2 (cQ_c-Q_c^4).
\end{equation*}
From Lemma \ref{RESUME}, we prove that 
$b_{2,0}<0$, thus $B_{2,0}\not\in \mathcal{Y}$. Thus, the approximate solution $v$ does not match an exact $2$-soliton solution at $T_c$ by  a term exactly of order
$\|(Q_c^2)'\|_{H^1}\sim K c^{11/12}$. This fact and perturbative analytic arguments
around $2$-soliton solutions, allow us to prove (Section 5)  that there is no pure  $2$-soliton solution for the nonintegrable case $p=4$ and to estimate above and below
 the size of the nonzero error term created by the interaction.
Thus, there is no choice of parameters nor any other decomposition which gives an exact $2$-soliton solution.

\medskip

Now, we give a precise statement concerning $v$ at $\pm T_c$ for $p=4$ and then for $p=2$.
We prove it only for $p=4$, the proof for $p=2$ being similar.

\begin{proposition}\label{AVRIL4}
    Let $p=4$. Let $k_0\geq 5$ and $\ell_0\geq 1$.
    There exists a function $v$ as in Theorem \ref{APSOL4} and Proposition \ref{SYSkl4} satisfying, for $c$
    sufficiently small,
    \begin{enumerate}
        \item Approximate solution on $[-T_c,T_c]$: for all $j\geq 1$, there exists $K=K(j)>0$ such that
        \begin{equation}\label{AV0}
            \forall t\in [-T_c,T_c],\quad    \|\partial_x^{j}
            \big(\partial_{t} v + \partial_{x} (\partial_{x}^2 v - v + v^p)\big)\|_{L^2(\mathbb{R})}
                \leq K c^{n_0},
        \end{equation}
        where $n_0=\frac {11}{24}\min(\frac {k_0}{3},\ell_0+1)$.
        \item Closeness to the sum of two solitons for $t=\pm T_c$,
\begin{equation}\label{3-16}\begin{split}
 &   \|v(T_c)-Q(.-\tfrac \Delta 2)-Q_c(.+(1-c)T_c-\Delta_c/2)-b_{2,0} (Q_c^2)'(y_c)\|_{H^1}\leq K c,\\&
    \|v(-T_c)-Q(.+\tfrac \Delta 2)-Q_c(.-(1-c)T_c+\Delta_c/2)+b_{2,0} (Q_c^2)'(y_c)\|_{H^1}\leq K c,
\end{split}\end{equation}
\begin{equation}\label{3-15}
	\|\partial_x (v(T_c) - Q(.-\tfrac \Delta 2) - Q_c(.+(1-c)T_c -\Delta_c/2))\|_{L^2}
            \leq K c^{\frac {17}{12}},
\end{equation}
                where
        \begin{equation}\label{defDeltabis}
            \Delta = \sum_{(k,\ell)\in \Sigma_0} 
                    a_{k,\ell} \, c^{\frac k{p-1}+\ell-\frac 12} \int Q^k,
        \quad 
            \Delta_c = 2b_{1,0} = -\int Q^3 +\frac 13 \frac {\left(\int Q\right)^2}{\int Q^2} <0.
        \end{equation}
        \item Decay on the right:
        \begin{align}\label{AV4}
            & \|v(T_c) - Q(.-\tfrac \Delta 2)\|_{H^1(x>-T_c/2)}\leq K  \exp(-\tfrac 14 c^{-\frac 1{ 100}}).  
        \end{align}
        \begin{equation}\label{AVV100}
        \forall x\geq 0,\quad 
        |v(0,x)|\leq  C \exp(-\tfrac 12 \sqrt{c} x). 
        \end{equation}
    \end{enumerate}
\end{proposition}

\noindent\emph{Remark.} Recall that for $p=4$, $\|Q_c\|_{L^2}=c^{1/12}\|Q\|_{L^2}$.
By \eqref{defDeltabis} and Lemma \ref{RESUME}, we  have
\begin{equation}\label{AV2}
    \Bigg| \Delta -\biggl(- 2 \frac{\left(\int Q\right)^2}{\int Q^2} \frac 1{c^{1/6}} + d_1 \, c^{1/6}+
        d_2 \, c^{1/2}+ d_3 \, c^{5/6}\bigg) \Bigg| \leq K c^{7/6},
\end{equation}
where $d_1,$ $d_2$ and $d_3$ are universal constants.
The other terms in the sum \eqref{defDeltabis} are of higher order than $c^{7/6}$, in particular, these terms are not relevant in our estimate. We will not compute $d_1$, $d_2$ and $d_3$
and will just keep the first order  term to state the main results
(see Theorem \ref{EXIST4}).

Since $   \frac 1K c^{11/12}\leq \| (Q_c^2)'(y_c)\|_{H^1}\leq  K c^{11/12}$ and $b_{2,0}\neq 0$,
estimates \eqref{3-16} imply that
\begin{align}
            & \frac 1 K c^{\frac {11}{12}}\leq \|v(T_c) - Q(.-\tfrac \Delta 2) - Q_c(.+(1-c)T_c -\Delta_c/2)\|_{H^1}
            \leq K c^{\frac {11}{12}},\label{AV1} \\
    	        & \frac 1K c^{\frac {11}{12}}\leq \|v(-T_c)-Q(.+\tfrac \Delta 2)-Q_c(.-(1-c)T_c+\Delta_c/2)\|_{H^1}
	        \leq K c^{\frac {11}{12}}.
	        \label{AV1sym}
        \end{align} 

\noindent\emph{Proof of Proposition \ref{AVRIL4}.} We consider the function $v$ constructed in Theorem \ref{APSOL4},
for $k_0\geq 5$ and $\ell_0\geq 1$. Since $p=4$, we have $q= 1/12$. Thus estimate \eqref{AV0} is a consequence of Theorem \ref{APSOL4} \eqref{almost}.
Estimate \eqref{vclose} still holds for $v$ on $[-T_c,T_c]$, but our objective is to prove \eqref{AV1}-\eqref{AV1sym}, which is a much sharper estimate for $t=\pm T_c$. 
We consider only $t=T_c$ by symmetry.
We justify the formal approach above.

\medskip

\noindent 1. Estimates on the remaining terms in $W(t,x)$ using Claim \ref{SEPT} (a)-(b)-(c). 
We claim, at $t=T_c$,
\begin{equation}\label{estWW}
    \|v(T_c)-Q(y)-Q_c(y_c)-b_{1,0} Q_c'(y_c)-b_{2,0} (Q_c^2)'(y_c)\|_{H^1}\leq K c.
\end{equation}

From \eqref{debutSEPT}--\eqref{debutSEPT2}, for $k=1,2,3$, $\ell=0$, at $t=T_c$
\begin{align*}
    & \|Q_c(y_c) A_{1,0}(y)\|_{H^1}+\|Q_c^2(y_c) A_{2,0}(y)\|_{H^1} + \|Q_c^3(y_c) A_{3,0}(y)\|_{H^1}\leq K   e^{-(1-c)\sqrt{c}T_c}\leq K c  e^{-c^{-q/2}},\\
    & \|(Q_c^3)'(y_c) B_{3,0}(y)\|_{L^2}
    +\tfrac 1{\sqrt{c}} \|\partial_x((Q_c^3)'(y_c) B_{3,0}(y)\|_{L^2} \leq K c^{5/4}.
\end{align*}
By similar estimates, since $\overline B_{1,0}, \overline B_{2,0}\in \mathcal{Y}$, we have at $t=T_c$,
\begin{equation*}\begin{split}
    &\|(Q_c)'(y_c) \overline B_{1,0}(y)\|_{L^2}
    +\tfrac 1{\sqrt{c}} \|\partial_x((Q_c)'(y_c) \overline B_{1,0}(y))\|_{L^2}\\&
    +\|(Q_c^2)'(y_c) \overline B_{2,0}(y)\|_{L^2}
    +\tfrac 1{\sqrt{c}} \|\partial_x((Q_c^2)'(y_c) \overline B_{2,0}(y))\|_{L^2}
    \leq K  e^{-c^{-q/2}}.
\end{split}\end{equation*}
We also check using Claim \ref{SEPT}, \eqref{suite1SEPT}--\eqref{suite1SEPT2}, that for $4\leq k\leq 6=2(p-1)$, $\ell=0$
and for $1\leq k\leq 3=(p-1)$, $\ell=1$, at $t=T_c$,
\begin{equation*}\begin{split}
    &\|c^{\ell} Q_c^k(y_c) A_{k,\ell}(y)\|_{L^2}
    +\tfrac 1{\sqrt{c}} \|\partial_x(c^{\ell} Q_c^k(y_c) A_{k,\ell}(y))\|_{L^2}\\&
    + \|c^{\ell} (Q_c^k)'(y_c) B_{k,\ell}(y)\|_{L^2}
    +\tfrac 1{\sqrt{c}} \|\partial_x( c^{\ell} (Q_c^k)'(y_c) B_{k,\ell}(y))\|_{L^2}   \leq K c^{25/24}.
\end{split}\end{equation*}
Finally, by Claim \ref{SEPT}, \eqref{DEGkl2SEPT}--\eqref{DEGkl2SEPT2}, we check that for $(k,\ell)$ such that 
$\xi(k,\ell)\geq 2$,
\begin{equation*}\begin{split}
    &\|c^{\ell} Q_c^k(y_c) A_{k,\ell}(y)\|_{L^2}
    +\tfrac 1{\sqrt{c}} \|\partial_x(c^{\ell} Q_c^k(y_c) A_{k,\ell}(y))\|_{L^2}\\&
+ \|c^{\ell} (Q_c^k)'(y_c) B_{k,\ell}(y)\|_{L^2}
    +\tfrac 1{\sqrt{c}} \|\partial_x(c^{\ell} (Q_c^k)'(y_c) B_{k,\ell}(y))\|_{L^2}
    \leq K c.
\end{split}\end{equation*}
Thus \eqref{estWW} is proved.

\medskip

\noindent 2.  Position of the soliton $Q$ at $t=T_c$. We claim
    \begin{enumerate}
        \item[{\rm (a)}] For $x\geq -T_c/2$ and $t=T_c$,
        \begin{equation}\label{posQ1}
            |\alpha(y_c)-\tfrac \Delta 2|\leq K    e^{-\frac 14 c^{-q/2}}. 
        \end{equation}
        \item[{\rm (b)}] For $t=T_c$,
        \begin{equation}\label{posQ2}
            \|Q(y)-Q(.-\tfrac \Delta 2)\|_{H^1}\leq K   e^{-\frac 12 c^{-q/2}}. 
        \end{equation}
    \end{enumerate}

  We have
    $|\alpha(y_c)-\tfrac \Delta 2|\leq  K    \int_{y_c}^{+\infty}Q_{c} (s) ds,$
and, for any $k\geq 1$, for any $y_c>0$,
\begin{equation*}
    0\leq \int_{y_c}^{\infty} Q_{c} (s) ds
    \leq K c^{1/3}\int_{y_c}^{\infty} e^{-\sqrt{c}\, s } ds=K c^{-1/6} e^{-\sqrt{c}\, y_c},
\end{equation*}
we obtain
\begin{equation*}
    \left|\alpha(y_c)-\tfrac \Delta 2\right| \leq K c^{-1/6} e^{-\sqrt{c}\,y_c}.
\end{equation*}
For $x\geq -T_c/2$ and $t=T_c$, we have $y_c=x+(1-c)T_c$$\geq (\frac 12 -c)T_c$, thus $\sqrt{c}\, y_c \geq \frac 12 c^{-q/2}-1$, and so we obtain (a).

Proof of (b). 
For $x\geq -T_c/2$, using (a), we have
\begin{equation*}
    \|Q(y)-Q(.-\tfrac \Delta 2)\|_{H^1(x>-T_c/2)}\leq  K    e^{-\frac 14 c^{-\frac 1{100}}}.
\end{equation*}

For $x<-T_c/2$, since $y=x-\alpha(y_c)$, and $|\alpha(y_c)|\leq K c^{-1/6}$, we have $y<-T_c/4$.
Thus,
\begin{align*}
&    \|Q(y)-Q(.-\tfrac \Delta 2)\|_{H^1(x<-T_c/2)}\\ &\leq  \|Q(y)\|_{H^1(x<-T_c/2)}+\|Q(.-\tfrac \Delta 2)\|_{H^1(x<-T_c/2)}\leq K    e^{-\frac 12 c^{-\frac 1{100}}}.
\end{align*}

\medskip

\noindent 3. Position of the soliton $Q_c$ at $t=T_c$. We claim

    \begin{equation}\label{posQc1}
        \|Q_c(y_c)-b_{1,0}Q_c'(y_c)-Q_c(.-b_{1,0})\|_{H^1}\leq K c^{13/12}.
    \end{equation}
For example, for the $L^2$-norm, we have
\begin{align*}
    \|Q_c-b_{1,0}Q_c'-Q_c(.-b_{1,0})\|_{L^2}& =c^q \|Q-\sqrt{c}\,b_{1,0}Q'-Q(.-\sqrt{c}\,b_{1,0})\|_{L^2}  \\
    &\leq K c^{q}(\sqrt{c}\, b_{1,0})^2 =K c^{1+q}.
\end{align*}
Therefore,  we obtain \eqref{3-16}.

\medskip

\noindent 4. Estimate on the right. 
Finally, we prove \eqref{AV4}. It is sufficient to prove that at $t=T_c$,
\begin{equation}\label{aprouver}
    \|Q_c(y_c)+W(t,x)\|_{H^1(x>-T_c/2)}\leq K  e^{-\frac 12 c^{-\frac 1{100}}}.
\end{equation}
For $x>-T_c/2$ and $t=T_c$, we have $y_c=x+(1-c)T_c>(1/2-c)T_c$ and so $\sqrt{c}\, y_c \geq \frac 12 c^{-\frac 1{100}} -1$.
Thus, it is clear   that
$
    \|Q_c(y_c)\|_{H^1(x>-T_c/2)}\leq K   e^{-\frac 12 c^{-\frac 1{100}}}.
$
All the other terms in $W(t,x)$ are checked to satisfy the same estimate, using the control on the degrees of the polynomial functions $\widetilde A_{k,\ell}$, $\widehat A_{k,\ell}$ and $\widetilde B_{k,\ell}$, $\widehat B_{k,\ell}$
as in the proof of Claim \ref{SEPT}.

The pointwise estimate \eqref{AVV100} for $x>0$ is clear from the decay properties of $Q$ and $Q_c$.
Thus Proposition \ref{AVRIL4} is proved.

\medskip

Finally, we present without proof a similar result for $p=2$.

\begin{proposition}\label{AVRIL4p2}
    Let $p=2$. Let $k_0\geq 2$ and $\ell_0\geq 1$. There exist $K>0$ and
     a function $v$ as in Theorem \ref{APSOL4} and Proposition \ref{SYSkl4} satisfying, for $c$
    sufficiently small,
    \begin{enumerate}
        \item Approximate solution on $[-T_c,T_c]$: for all $j=0,1,2$ such that
        \begin{equation}\label{AV0p2}
            \forall t\in [-T_c,T_c],\quad    \|\partial_x^{j}
            \big(\partial_{t} v + \partial_{x} (\partial_{x}^2 v - v + v^p)\big)\|_{L^2(\mathbb{R})}
                \leq K c^{2}.
        \end{equation}
        \item Closeness to the sum of two solitons for $t=\pm T_c$,
        \begin{align}
            & \|v(T_c) - Q(.-2\sqrt{c}) - Q_c(.+(1-c)T_c +2)\|_{H^1(\mathbb{R})}
            \leq K c^{3/2},\label{AV1p2} \\
    	        &\|v(-T_c)-Q(.+2\sqrt{c})-Q_c(.-(1-c)T_c-2)\|_{H^1(\mathbb{R})}
	        \leq K c^{3/2}.
	        \label{AV1symp2}
        \end{align} 
    \end{enumerate}
\end{proposition}

\section{Preliminary results for stability of the $2$-soliton structure}

In this section, we gather several stability results (essentially refinements of tools 
developed in \cite{MM1}, \cite{MMT} and \cite{MMnonlinearity}). Section 4.1 concerns the
stability of $v(t)$ by the gKdV equation during the interaction. Sections 4.2 and 4.3 
concern the large time behavior (after interaction).

\subsection{Dynamic stability in the interaction region}
For any $c$ small enough, we consider a function $v(t)$ of the form
\begin{equation}\label{defvthINT}
    v(t,x)=Q(y)+Q_{c}(y_{c})+\sum_{(k,\ell)\in \Sigma_0} 
        c^\ell\left(Q_{c}^k(y_{c}) A_{k,\ell}(y)+(Q_{c}^k)'(y_{c}) B_{k,\ell}(y)\right)
\end{equation}
where
$
    y_c=x+(1{-}c)t, \quad y=x-\alpha(y_{c})$ and 
$ 
    \alpha(s)=\sum_{(k,\ell)\in \Sigma_0} a_{k,\ell} \, c^\ell \int_{0}^s    Q_{c}^k(s')ds',
$
and $(a_{k,\ell})$, $(A_{k,\ell})$, $(B_{k,\ell})$ satisfy the properties of Proposition \ref{SYSkl4}.
Set
$
    S(t)=\partial_{t} v + \partial_{x} (\partial_{x}^2 v - v + v^p).
$

\begin{proposition}[Exact solution close to the approximate solution $v$]\label{INTERACT4}
    Let $p=2$, $3$ or $4$.  Let $\theta>\frac 1{p-1}$. There exists $c_5>0$ such that the following holds for any $0<c<c_5$.
    Suppose that 
    \begin{equation}\label{INTkl}
      \text{for $j=1,2,3$},\
        \forall t\in [-T_c,T_c],\quad 
        \left\|\partial_{x}^j S(t)\right\|_{L^2(\mathbb{R})}\leq K \frac {c^\theta}{T_c} = K c^{\theta+\frac 12+\frac 1{100}},
    \end{equation}
    and that for some $T_0\in [-T_c,T_c]$,
    \begin{equation}\label{hypINT}
        \| u(T_0) - v(T_0) \|_{H^1(\mathbb{R})}\leq K c^{\theta},
    \end{equation}
where  $u(t)$ is an $H^1$ solution of the (gKdV) equation \eqref{gkdv}. 
    Then, there exist $K_0=K_0(\theta,K)$ and a function $\rho:[-T_c,T_c]\rightarrow \mathbb{R}$ such that, for all $t\in [-T_c,T_c]$,
    \begin{equation}\label{INT41}
        \|u(t)-v(t,.-\rho(t)) \|_{H^1} \leq K_0 c^{\theta},\quad |\rho'(t)-1|\leq K_0 c^{\theta}.
    \end{equation}
\end{proposition}

\noindent\emph{Remark.} By usual techniques related to the resolution of the Cauchy problem, one obtains for approximate solutions a divergence of order $e^{T_c}$ for a time interval
$[0,T_c]$. Here, such an estimate would not be sufficient since $T_c= c^{-\frac 12 
+{\frac 1{100}}}\gg
c^{-1/2}$. In this proof, we use the Hamiltonian properties of the gKdV equation. More precisely, the proof is based on
the fact that $v$ is close to $Q$ ($c$ is small), and
on refined stability analysis around $Q$ 
(on the one hand  standard arguments of long time stability (see Weinstein \cite{We2})
and on the other hand  some algebraic cancellations in the energy functional).
This leads us to  a simple ODE
estimate in time on the error term.

Note that $\theta>\frac 1{p-1}$ is arbitrary in Proposition \ref{INTERACT4}. Moreover, from the algebraic
argment (Theorem \ref{APSOL4}), there exists $v$ such that \eqref{INTkl} holds for any $\theta$ large.
This implies that  if (for example) $u(0)=v(0)$, then $\|u(T_c)-v(T_c)\|_{H^1}\leq K(\theta) c^\theta$,
for any $\theta$ large. Theorefore, the approximate function $v$ and 
its properties (for example the shift properties) are sharp up to any order $c^\theta$, and provide
a sharper description of the collision problem as $\theta\to+\infty$.

\medskip

\noindent\emph{Proof of Proposition \ref{INTERACT4}.}
We  prove the result on $[T_0,T_c]$. By using the transformation $x\to -x$, $t\to -t$, the proof is the same on $[-T_c,T_0]$.
Let $K^*>1$ be a constant to be fixed later. Since $\|u(T_0)-v(T_0)\|_{H^1}\leq c^\theta$, by continuity in time in $H^1(\mathbb{R})$, there exists $T^*>T_0$ such that
\begin{equation*}
    T^*=\sup\left\{T\in [T_0,T_c] \text{ s.t. $\forall t\in [T_0,T]$, $\exists r(t)\in \mathbb{R}$ with }
    \|u(t){-}v(t,.{-}r(t))\|_{H^1}\leq K^* c^{\theta} \right\}.
\end{equation*}
Note that the translation direction is degenerate and without the freedom in the translation parameter, the result would not be correct.
The objective is to prove that $T^*=T_c$ for $K^*$ large. For this, we argue by contradiction, assuming that $T^*<T_c$ and reaching a    contradiction with the definition of $T^*$ by proving independent estimates on $\|u(t)-v(t,.-r)\|_{H^1}$ on $[T_0,T^*]$.

First, we claim some estimates related to $v$.

\begin{claim}[Preliminary estimates]\label{TRANScl}
    The following hold.
    \begin{equation}\label{Tcl1}
        \|\partial_t v(t)\|_{L^\infty}  \leq K c^{\frac 1{p-1}},
    \end{equation}
    \begin{equation}\label{Tcl2}
        \|\partial_t v(t) + \alpha'(y_c) Q'(y)\|_{L^2}\leq  K c^{q+\frac 12},\quad
        \|\partial_t v(t) + \alpha'(y_c) Q'(y)\|_{L^\infty}\leq  K c^{m_0},
    \end{equation}
    \begin{equation}\label{Tcl3}
        \|v^{p-2}- Q^{p-2}(y)\|_{L^\infty}\leq K c^{\frac 1{p-1}},
    \end{equation}
    \begin{equation}\label{Tcl4}
        \|\partial_x v - Q'(y)\|_{L^2}\leq K c^{\frac 1{p-1}},
    \end{equation}
    \begin{equation}\label{Tcl5}
        \|\alpha''(y_c)\|_{L^\infty}+\frac 1c \|\alpha^{(4)}(y_c)\|_{L^\infty}\leq K c^{\frac 12 + \frac 1{p-1}},
    \end{equation}
    where $m_0=\min\left(\frac 2{p-1},\frac 1{p-1}+\frac 12\right)$.
\end{claim}

\noindent\emph{Proof of Claim \ref{TRANScl}.}
 \eqref{Tcl1}---\eqref{Tcl2}:  We differentiate formula \eqref{defvthINT} with respect to $t$: 
\begin{align*} 
    \partial_t v(t) & =-(1-c)\alpha'(y_c)Q'(y)+(1-c)Q_c'(y_c)\\
        &\quad +\sum_{(k,\ell)\in \Sigma_0} 
        c^\ell\left((1-c) (Q_{c}^k)'(y_{c}) A_{k,\ell}(y) -(1-c)\alpha'(y_c) Q_{c}^k(y_{c}) A_{k,\ell}'(y)\right)\\
        &\quad +\sum_{(k,\ell)\in \Sigma_0}
        c^\ell\left((1-c)(Q_{c}^k)''(y_{c}) B_{k,\ell}(y)-(1-c)\alpha'(y_c)(Q_{c}^k)'(y_{c}) B_{k,\ell}'(y) \right).
\end{align*}
By the same estimates as in the proofs of Proposition \ref{VandS} and Claim \ref{SEPT}, and by $|\alpha'(y_c)|\leq K c^{\frac 1{p-1}}$ (see Claim \ref{surALPHA}), we have
$
    \|\partial_t v(t)\|_{L^\infty}\leq K \|Q_c\|_{L^\infty}\leq K c^{\frac 1 {p-1}},
$
and \eqref{Tcl2}.

 From the expression of $v$ and estimates as in the proof of Proposition \ref{VandS}, we obtain
\eqref{Tcl3}.

\eqref{Tcl4}: Differentiating \eqref{defvthINT} with respect to $x$:
\begin{align*} 
    \partial_x v(t) & = Q'(y) - \alpha'(y_c) Q'(y) + Q_c'(y)+o(c^{\frac 1{p-1}}).
\end{align*}
\eqref{Tcl5}:
$
 |\alpha''(s)|\leq K \sum_{1\leq k\leq k_{0}  }   |(Q_{c}^k)'(s)|\leq K \|Q_c'\|_{L^\infty}\leq K c^{\frac 12 + \frac 1{p-1}}.$

\medskip

\emph{Step 1.} Choice of the translation parameter and control of the $Q'$ direction.

\begin{lemma}[Modulation]\label{TRANS}
    There exists a $C^1$ function $\rho:[T_0,T^*]\rightarrow \mathbb{R}$ such that, for all $t\in [T_0,T^*]$,
    the function $z(t)$ defined by
    $
        z(t)=u(t,x+\rho(t))-v(t,x)
    $
    satisfies, $ \forall t\in [T_0,T^*],$
    $             \int z(t) Q'(y) dx =0,\label{TRANS1}
     $ and for $K$ independent of $K^*$,
        \begin{equation}
           \|z(t)\|_{H^1}\leq  2 K^* c^{\theta},\
            |\rho(T_0)|+\|z(T_0)\|_{H^1}\leq K c^{\theta},\
 |\rho'(t)-1|\leq K \|z(t)\|_{H^1}+K \|S(t)\|_{H^1}.\label{TRANS3}
        \end{equation}
 \end{lemma}
\emph{Proof of Lemma \ref{TRANS}.}\quad
The existence of $\rho(t)$ is obtained at fixed time $t\in [T_0,T^*]$.
Let (recall $y=x-\alpha(y_c)$)
\begin{equation*}
    \zeta(U,r)=\int (U(x+r)-v(t,x)) Q'(y)dx.
\end{equation*}
Then $
    \frac{\partial \zeta}{\partial r}(U,r)= \int U'(x+r) Q'(y) dx, 
$
so that
from Claim \ref{TRANScl},  for $c$ small enough,
\begin{equation*}
    \frac{\partial \zeta}{\partial r}(v,0)= \int (\partial_x v)(t,x) Q'(y) dx >  \int (Q'(y))^2dx - K c^{\frac 1{p-1}} > \frac 12 \int (Q')^2
\end{equation*}
(note that $\int (Q'(y))^2 dx=\int (Q'(y))^2 \frac {dy}{1-\alpha_c'(y_c)}>\frac 34\int (Q'(y))^2dy $). 
Since    $\zeta(v,0)=0$, for $U$ close to $v(t)$ in $L^2$ norm, the existence    of a unique 
$\rho(U)$ satisfying $\zeta(U(x-\rho(U)),\rho(U))=0$ is  a consequence of the Implicit Function Theorem.

\smallskip

Let us prove that
\begin{equation}\label{estxxX}
    |\rho'(t)-1|\leq K \|z(t)\|_{H^1} + K \|S(t)\|_{H^1}.
\end{equation}
From the definition of $T^*$, it follows that there exists $\rho(t)=\rho(u(t))$, such that 
$\zeta(u(x-\rho(t)),\rho(t))=0$.
We set
\begin{equation}\label{defz}
    z(t,x)=u(t,x+\rho(t))-v(t,x),
\end{equation}
then $\int z(t)Q'=0$ follows from the definition of $\rho(t)$ and \eqref{TRANS3} from the Implicit
Function Theorem and the definition of $K^*$.
Moreover, since $\|u(T_0)-v(T_0)\|\leq c^{\theta}$, we have $
    |\rho(T_0)|+\|z(T_0)\|_{H^1}\leq K c^{\theta}$,
where $K$ is independent of $K^*$.

From the definition of $z(t)$, $u(t)$ being a solution of the (gKdV) equation, we obtain
\begin{align}
    \partial_t z + \partial_x (\partial_{x}^2 z - z + (z+v)^p -v^p) 
    & =\partial_{t} v + \partial_{x} (\partial_{x}^2 v - v + v^p) + (\rho'(t)-1) \partial_x u \nonumber\\
    & = -S(t) + (\rho'(t)-1) \partial_x(v+z).
    \label{eqz}
\end{align}
Since
$
    \int z(t,x) Q'(y) dx = 0,
$
by $y=x-\alpha(y_c)$ and $y_c=x+(1-c)t$, we have
\begin{equation*}
0 = \frac d{dt} \int z Q'(y) dx = \int \partial_t z Q'(y) - (1-c)  \int\alpha'(y_c) z Q''(y).
\end{equation*}
Thus, integrating by parts,
\begin{equation}\label{exprxxX}
    \begin{split}
        (\rho'(t)-1) \int (v+z) \partial_x (Q'(y)) & = 
         \int z (\partial_x^3 - \partial_x) (Q'(y)) 
            + \int \left((z+v)^p - v^p\right) \partial_x (Q'(y)) 
        \\& \qquad - \int S(t) Q'(y) - (1-c)  \int\alpha'(y_c) z Q''(y). 
    \end{split}
\end{equation}
Thus
$
\left| (\rho'(t)-1) \int (v+z) \partial_x (Q'(y)) \right|  \leq K \|z(t)\|_{H^1} + K \|S(t)\|_{H^1}.
$
The term $ - \int (v+z) \partial_x Q'(y)$ has a positive lower bound:
\begin{align*}
    \int (v+z) \partial_x (Q'(y)) & = \int (1-\alpha'(y_c)) (v+z) Q''(y)\\ 
    & = \int Q(y)Q''(y) +    \int (v-Q(y)+z) Q''(y) -  \int \alpha'(y_c) (v+z) Q''(y).
\end{align*}
Since $-\int Q(y)Q''(y)dx >\frac 34  \int (Q'(y))^2 dy>0$ and since the other terms are small for $c$ small, we have
$
    - \int (v+z) \partial_x Q'(y) > \frac 12 \int (Q')^2.
$
Note that $\rho(t)$ is $C^1$ since $Q(y)$ and $v$ are $C^\infty$ and $z(t)$ is continuous in $H^1(\mathbb{R})$.
\eqref{estxxX} is proved.
\quad 

\medskip

\emph{Step 2.} $L^2$ norm conservation and control of the direction $\int zQ(y)$. The use of the
$L^2$ norm conservation replaces a modulation argument in the scaling parameter.

\begin{lemma}[Control of the $Q$ direction]\label{Qdir}
    For all $t\in [T_0,T^*]$,
    \begin{equation}\label{Qdir1}
        \left|\int  z(t)Q(y)\right|\leq  K c^\theta + K c^{q} \|z(t)\|_{L^2}+ \|z(t)\|_{L^2}^2.
    \end{equation}
\end{lemma}

\noindent\emph{Proof of Lemma \ref{Qdir}.}
Remark that since $v(t)$ is an approximate solution of \eqref{gkdv}, its $L^2$ norm has a small variation.
Indeed, by multiplying the equation $S(t)= \partial_t v + \partial_x(\partial_x^2 v -v + v^p)$ by $v$ and integrating,
we obtain
$
    \left|\frac 12 \frac d {dt} \int v^2 \right|=\left| \int S(t,x) v(t,x) dx\right|
    \leq K \|S(t)\|_{L^2}.
$
Thus,    
\begin{equation}\label{vL2}
\forall t\in [T_0,T^*],\quad
    \left|\int v^2(t) - \int v^2(T_0)\right| \leq K T_c \sup_{t\in [-T_c,T_c]} \|S(t)\|_{H^1}
    \leq K c^\theta. 
\end{equation}
Since $u(t)$ is a solution of the (gKdV) equation, we have
\begin{equation}\label{uL2}
    \int u^2(t)=\int \left(v(t)+z(t)\right)^2=\int u^2(T_0)=\int \left(v(T_0)+z(T_0)\right)^2.
\end{equation}
By expanding \eqref{uL2} and using \eqref{vL2} and \eqref{TRANS3}, we obtain:
\begin{equation*}
    2 \left|\int v(t)z(t)\right|\leq K c^\theta + 2 \left|\int v(T_0)z(T_0)\right|
    + \|z(T_0)\|_{L^2}^2+ \|z(t)\|_{L^2}^2\leq K c^\theta+ \|z(t)\|_{L^2}^2.
\end{equation*}
 Using this  and $\|v(t)-Q(y)\|_{L^2}\leq K c^q$, we obtain:
\begin{equation*}
     \left|\int  z(t)Q(y)\right|\leq 
 \left|\int  z(t)v\right|+ \left|\int  z(t)(v-Q(y))\right|
\leq  K c^\theta + K c^{q} \|z(t)\|_{L^2}+ \|z(t)\|_{L^2}^2.
\end{equation*}
\quad 

\medskip

\emph{Step 3.} Introduction of a  energy functional for  $\|z(t)\|_{H^1}$.
We set
\begin{equation*}
    \mathcal{F}(t)=\frac 12 \int \left((\partial_x z)^2 + (1+\alpha'(y_c)) z^2)\right)
    -\frac 1{p+1} \int \left((v+z)^{p+1}- v^{p+1}-(p+1) v^p z\right).
\end{equation*}
The above definition is similar to  a linearized energy 
$\frac 12 \left(     \int    ((\partial_x z)^2    + z^2)
    - p \int Q^{p-1}z^2\right).
$

However, the terms $\int  \alpha'(y_c)    z^2 $ and
the nonlinear terms were added to remove some diverging terms in $\mathcal{F}'$.
This is the new ingredient of the proof of Proposition \ref{INTERACT4}.

We first claim that the functional $\mathcal{F}(t)$ indeed controls the size of $z(t)$ in $H^1$ up to the direction $Q(y)$, extending the similar classical result 
for the linearized energy.

\begin{claim}[Coercivity of $\mathcal{F}$]\label{POSf}
    There exists $\kappa_0>0$ such that
    \begin{equation}\label{posf1}
       \|z(t)\|_{H^1}^2\leq  \kappa_0 \mathcal{F}(t) + {\kappa_0} \left|\int z(t) Q(y)\right|^2.
    \end{equation}
\end{claim}
The proof of Claim \ref{POSf} is given in Appendix D.1.

Next, we claim the following control of the variation of $\mathcal{F}(t)$ through time.

\begin{lemma}[Control of the variation of the energy fonctional]\label{varF}
    \begin{equation}\label{varF1}
        \mathcal{F}(T^*)-\mathcal{F}(T_0)\leq  K c^{2 \theta}\left((K^*)^2(1+K^*)c^{q/2} + K^*\right).
    \end{equation}
    where $K$ is independent of $c$ and $K^*$.
\end{lemma}

\noindent\emph{Proof of Lemma \ref{varF}.}
We have
\begin{multline*}
    \mathcal{F}'(t)     =       \int \partial_t z \left(-\partial_x^2 z +  z - ((v+z)^p-v^p)\right) 
 +\int                        \partial_t z  \alpha'(y_c) z\\
+\int\left\{ \frac 12 (1-c)  \alpha''(y_c) z^2                        -  \partial_t v \left((v+z)^p-v^p-pv^{p-1}z\right)\right\}
 = \mathbf{F}_1 +\mathbf{F}_2 + \mathbf{F}_3.
\end{multline*}
Now, we claim:
\begin{equation}\label{f6}
    \left| \mathbf{F}_1 + (\rho'(t)-1) \int \alpha'(y_c) Q'(y) z \right| \leq K   c^{q+\frac 12} \|z(t)\|_{H^1}^2
    + K \|z(t)\|_{L^2} \left(\|\partial_x^2 S(t)\|_{L^2} + \| S(t)\|_{L^2}\right).
\end{equation}
\begin{equation}\label{f7}
\begin{split}
&    \left|\mathbf{F}_2 - (\rho'(t)-1) \int \alpha'(y_c) Q'(y) z + \frac {p(p-1)}2 \int \alpha'(y_c) z^2 Q'(y) Q^{p-2}(y)\right| 
\\ & \leq K \|z(t)\|_{H^1}^2 \left(c^{m_0}+c^{\frac 1{p-1}} \|z(t)\|_{H^1}\right)+K \|z(t)\|_{H^1}\left(\|\partial_x^2 S(t)\|_{L^2}+\|S(t)\|_{L^2}\right). 
\end{split}\end{equation}
\begin{equation}\label{f3}
    \left|\mathbf{F}_3 - \frac {p(p-1)}2 \int \alpha'(y_c) Q'(y) Q^{p-2}(y)\, z^2\right|\leq
     K c^{m_0} \|z\|_{H^1}^2 + K c^{\frac 1 {p-1}} \|z(t)\|_{H^1}^3.
\end{equation}
Assuming \eqref{f6}--\eqref{f3}, we conclude the proof of the lemma.

Note that $q+\frac 12=\frac 1{p-1}+\frac 14\leq m_0.$
From the cancellations of the main terms of $\mathbf{F}_1$, $\mathbf{F}_2$ and $\mathbf{F}_3$, and then from  \eqref{TRANS3}, \eqref{INTkl},   we  get
\begin{equation*}\begin{split}
|\mathcal{F}'(t)| &
\leq           K \|z(t)\|_{H^1}^2 \left(c^{q+\frac 12}+c^{\frac 1{p-1}} \|z(t)\|_{H^1}\right)+K \|z(t)\|_{H^1}\left(\|\partial_x^2 S(t)\|_{L^2}+\|S(t)\|_{L^2}\right)\\
& \leq  K\left[ (K^*)^2 c^{2\theta} ( c^{q+\frac 12}+ K^*c^{\frac 1{p-1}+\theta})+ K^* c^{ \frac 12(1+{\frac 1{50}})+2 \theta }\right].
\end{split}
\end{equation*}
Now, $q\geq {\frac 1{50}}$ and   $\theta+\frac 1 {p-1}\geq \frac 2 {p-1}\geq q+\frac 12\geq \frac 12 (1+{\frac 1{50}}) +\frac q2 $ imply
\begin{equation*}
|\mathcal{F}'(t)|\leq           K c^{\frac 12 (1+{\frac 1{50}})+ 2 \theta}\left((K^*)^2c^{q/2}(1+K^*)+ K^*\right). 
\end{equation*}

Integrating on the time interval $[T_0,T^*]$ where $T^*-T_0 \leq  2T_c= 2 c^{\frac 12 (1+{\frac 1{50}})}$, we obtain
\begin{equation*}
    |\mathcal{F}(T^*)-\mathcal{F}(T_0)|\leq   K c^{2 \theta}\left((K^*)^2(1+K^*) c^{q/2} + K^*\right).
\end{equation*}

\medskip

Proof of \eqref{f6}. We replace $\partial_t z$ by its expression:
\begin{align*}
    \mathbf{F}_1 & = -\int S(t) \left(-\partial_x^2 z + z - ((v+z)^p - v^p)\right)\\
                &    \quad + (\rho'(t)-1) \int \partial_x (v+z)  \left(-\partial_x^2 z + z - ((v+z)^p - v^p)\right)=
                \mathbf{g}_1+\mathbf{g}_2.
\end{align*}
By integration by parts, the Cauchy-Schwarz' inequality, we have
\begin{equation*}
    |\mathbf{g}_1|\leq K\|z(t)\|_{L^2} \left(\|\partial_x^2 S(t)\|_{L^2} + \| S(t)\|_{L^2}\right).
\end{equation*}
Since $\int \partial_x(v+z) (v+z)^p=0$, and by the definition of $S(t)$
\begin{align*}
    \mathbf{g}_2 &  =(\rho'(t)-1)\int \partial_x(v+z) (-\partial_x^2 z + z + v^p) =
 (\rho'(t)-1) \int (\partial_x v (-\partial_x^2 z + z ) + \partial_x z v^p)\\& = (\rho'(t)-1) \int z \partial_x(-\partial_x^2 v + v - v^p)
    = (\rho'(t)-1) \int z(\partial_t v - S(t)).
\end{align*}
By \eqref{Tcl2} and \eqref{TRANS3}, we obtain:
\begin{align*}
    \left| \mathbf{g}_2 + (\rho'(t)-1) \int \alpha'(y_c) Q'(y) z \right| 
        & \leq K |\rho'(t)-1| \|z(t)\|_{L^2} \left( \|\partial_t v - \alpha'(y_c) Q'(y)\|_{L^2} + \|S(t)\|_{L^2}\right)     \\
        & \leq K \|z(t)\|_{L^2}(\|z(t)\|_{H^1} + \|S(t)\|_{L^2})(c^{q+\frac 12}+ \|S(t)\|_{L^2}).
\end{align*}

\medskip

Proof of \eqref{f7}.
Note that the term $\mathbf{F}_2$ was introduced on purpose in the expression of $\mathcal{F}$ to cancel the main terms in 
$\mathbf{F}_1$ and $\mathbf{F}_3$.
\begin{align*}
    \mathbf{F}_2 & = \int \alpha'(y_c) z \partial_x(-\partial_x^2 z + z  - ((z+v)^p - v^p))
   \\ & \quad  - \int \alpha'(y_c) z S(t)
    + (\rho'(t)-1) \int \alpha'(y_c) \partial_x(v+z) z = \mathbf{g}_{3}+\mathbf{g}_{4}.
\end{align*}
First,
\begin{equation*}
    \mathbf{g}_{4}= - \int \alpha'(y_c) z S(t)+ (\rho'(t)-1) \int \alpha'(y_c) \partial_x v\,    z - \frac 12 (\rho'(t)-1)\int z^2 \alpha''(y_c).
\end{equation*}
By \eqref{Tcl4}--\eqref{Tcl5}, and \eqref{TRANS3}, we have
\begin{equation*}
    \left|\mathbf{g}_{4} - (\rho'(t)-1) \int \alpha'(y_c) Q'(y) z \right|
        \leq K c^{m_0} \|z(t)\|_{H^1}(\|z(t)\|_{H^1} + \|S(t)\|_{H^1}).
\end{equation*}

Second, for the term $\mathbf{g}_{3}$, we integrate by parts, to obtain:
\begin{equation}\label{ff10}
    \mathbf{g}_{3}  = - \int \alpha''(y_c) (\tfrac 32 (\partial_x z)^2 + \tfrac 12 z^2) + \int \alpha^{(4)} (\tfrac 12 z^2)
    - \int \alpha'(y_c) z \partial_x((z+v)^p-v^p). 
\end{equation}
Using the estimate on $\alpha''(y_c)$ and $\alpha^{(4)}(y_c)$ in Claim \ref{TRANScl}, we obtain
\begin{equation*}
    \left|-\int \alpha''(y_c) (\tfrac 32 (\partial_x z)^2 + \tfrac 12 z^2) + \int \alpha^{(4)} (\tfrac 12 z^2)\right|
    \leq K c^{\frac 1 {p-1} + \frac 12} \|z(t)\|_{H^1}^2.
\end{equation*}
In the last term of \eqref{ff10}, cubic and higher order terms are controlled by $K c^{\frac 1{p-1}} \|z(t)\|_{H^1}^3.$
The quadratic term is
\begin{equation*}
    \int \alpha'(y_c) z \partial_x(-p v^{p-1} z) 
    = \frac p2 \int \alpha''(y_c) z^2    v^{p-1}
    -\frac p 2 \int \alpha'(y_c) z^2 \partial_x(v^{p-1})=\mathbf{g}_{5}+\mathbf{g}_{6}.
\end{equation*}
As before, $|\mathbf{g}_{5}|\leq K c^{\frac 12+\frac 1{p-1}} \|z(t)\|_{H^1}^2$.
Finally, by \eqref{Tcl3}--\eqref{Tcl4},
\begin{equation*}
    \left|\mathbf{g}_{6} + \frac {p(p-1)}2 \int \alpha'(y_c) z^2 Q'(y) Q^{p-2}(y)\right| \leq K c^{\frac 2{p-1}} \|z(t)\|_{H^1}^2.
\end{equation*}

\medskip

Proof of \eqref{f3}. First note
$
    \left| \frac 12 (1-c)  \int \alpha''(y_c) z^2 \right| \leq K c^{\frac 12 + \frac 1 {p-1}}\|z(t)\|_{L^2}^2 .
$
We now estimate
\begin{equation*}
     - \int \partial_t v \left((v+z)^p-v^p-pv^{p-1}z-\tfrac {p(p-1)} 2 v^{p-2} z^2\right) - \frac {p(p-1)} 2 \int \partial_t v \, v^{p-2} z^2
=\mathbf{g}_7+\mathbf{g}_8
\end{equation*}
By \eqref{Tcl1}, we have
$    |\mathbf{g}_7|\leq K c^{\frac 1 {p-1}} \|z(t)\|_{H^1}^3.$
By \eqref{Tcl2}, \eqref{Tcl3}, and $|\alpha'(y_c)|\leq K c^{\frac 1{p-1}}$, we have
\begin{equation*}
    \left|\mathbf{g}_8 - \frac {p(p-1)}2 \int \alpha'(y_c) Q'(y) Q^{p-2} z^2\right|\leq
     K c^{m_0} \|z\|_{H^1}^2,
\end{equation*}

\emph{Step 4.} Conclusion of the proof.
By Claim \ref{Qdir}, Lemmas \ref{TRANS}--\ref{Qdir}, we have
\begin{equation*}
    \left|\int  z(T^*) Q(y) \right|\leq  K c^\theta + K c^{q} \|z(T^*)\|_{L^2}
    + \|z(T^*)\|_{L^2}^2,
\end{equation*}
and thus by Claim \ref{POSf},
$
    \|z(T^*)\|_{H^1}^2\leq K \mathcal{F}(T^*) + K (c^\theta +   c^{q} \|z(T^*)\|_{L^2}+ \|z(
    T^*)\|_{L^2}^2)^2.
$
It follows that for $c$ small enough,
$
    \|z(T^*)\|_{H^1}^2\leq (K+1) \mathcal{F}(T^*) + K c^{2 \theta}.
$

Next, by Lemma \ref{varF} and $|\mathcal{F}(T_0)|\leq K c^{2 \theta}$, we obtain
\begin{equation*}
    \|z(T^*)\|_{H^1}^2 
\leq (K+1) (\mathcal{F}(T^*)-\mathcal{F}(T_0))+K c^{2\theta}
\leq K_1 c^{2 \theta}\left((K^*)^2(1+K^*) c^{q/2} + K^*+1\right),
\end{equation*}
where $K_1$ is independent of $c$ and $K^*$.
Choose $c^*=c^*(K^*)$ such that  
\begin{equation*}
    (K^*)^2(1+K^*) (c^*)^{q/2}<1.
\end{equation*}
Then, for $0<c<c^*$,
\begin{equation*}
    \|z(T^*)\|_{H^1}^2 \leq K_1 c^{2 \theta}\left(2+ K^*\right).
\end{equation*}
Next, fix $K^*$ such that $K_1 (2+K^*)<\frac 12 (K^*)^2$. Then 
$
    \|z(T^*)\|_{H^1}^2 \leq \frac 12 (K^*)^2 c^{2 \theta}.
$
This contradict the definition of $T^*$, thus proving that $T^*=T_c$.
Thus estimate \eqref{INT41} is proved on $[T_0,T_c]$.

\subsection{Stability and asymptotic stability for large time}
In this section, we consider the stability of the $2$-soliton structure after the collision. These questions
have been considered in   \cite{MMas2}. See also \cite{MM1}, \cite{MMT}, \cite{Martel}.
Denote for $v\in H^1(\mathbb{R})$,
$
\|v\|_{H^1_c}=\left(\int_{\mathbb{R}} \left( (v'(x))^2 + c v^2(x) \right)dx \right)^{\frac 12},$
which corresponds to the natural norm to study the stability of $Q_c$.

\begin{proposition}[Stability of two decoupled solitons,  \cite{MMas2}]\label{ASYMPTOTIC}
There exists $K>0$, $\alpha_0>0$, $c_0>0$  such that for any $0<c<c_0$, $0<\alpha<\alpha_0$, the following holds.
Let $u(t)$ be an $H^1$ solution of \eqref{gkdv} such that
for some $t_1\in \mathbb{R}$ and $X_0\geq \frac 12 T_c$,
\begin{equation}
\label{D25}
\|u(t_1)-Q-Q_c(.+X_0)\|_{H^1}\le \alpha c^{q+\frac 12}.
\end{equation}
Then there exist   $C^1$  functions $\rho_1(t)$, $\rho_2(t)$ defined on $[t_1,+\infty)$  such that
\begin{enumerate}
\item Stability.
\begin{equation}\label{huit}
\sup_{t\ge t_1} \| u(t)-(Q(.-\rho_1(t))+Q_c(.-\rho_2(t)))  \|_{H^1_c}
\le K \alpha  c^{q+\frac 12}+K\exp (-  {c^{-\frac  1{400}}}  ),
\end{equation}
\begin{equation}\label{suppl}\begin{split}
& 
\forall t\ge t_1,\ \tfrac 12\leq \rho_1'(t)-\rho_2'(t) \leq \tfrac 32,
\\ &|\rho_1(t_1)|\leq K \alpha c^{q+\frac 12},\quad
|\rho_2(t_1)-X_0|\leq K \alpha.
\end{split}\end{equation}
\item Convergence of $u(t)$.
There exist $c_1^+, c_2^+ >0$  such that
\begin{equation}\label{neuf}
\lim_{t\rightarrow +\infty}\|u(t)-Q_{c_1^+}(x-\rho_1(t))-Q_{c_2^+}(x-\rho_2(t))\|
_{H^1(x> ct/10)}=0.
\end{equation}
\begin{equation}\label{sept}
|c_1^+-1|\leq   K  \alpha c^{q+\frac 12 }+ K \exp(-c^{- {\frac 1{400}}}),\quad
 \left|\frac {c_2^+}{c} - 1\right|\le K  \alpha  + K \exp(-c^{-   {{\frac 1{400}}}}).
\end{equation}
\item 
Assume further that  $\int_{x>0} x^2 \, u^2(t_1,x) dx <K_0$. Then, there exist  $x_1^+$ and 
$x_2^+$ such that
\begin{equation}
\lim_{t\to +\infty} \rho_1(t)-c_1^+ t=x_1^+, \quad
\lim_{t\to +\infty} \rho_2(t)-c_2^+ t=x_2^+.
\end{equation}
For $p=4$, if  for  $\kappa>0$,
\begin{equation}\label{decayfort}
\alpha< \kappa c^{\frac 13} \quad \text{and}\quad \int_{x>\frac {11} {12} |\ln c|} x^2 \, u^2(t_1,x) dx <  \kappa c^{\frac 54}
\end{equation}
then
\begin{equation}
|x_1^+-\rho_1(t_1)|\leq K c^{\frac 58}, \quad |x_2^+-\rho_2(t_1)|\leq K c^{\frac 1{12}}.
\end{equation}
\end{enumerate}
\end{proposition}
The proof of 
Proposition \ref{ASYMPTOTIC} is based on energy arguments, monotonicity results on local quantities,
and a Virial argument, see \cite{MMas1}, \cite{MMas2}.  

\medskip

\noindent\emph{Remark.} To obtain the convergence of the translation parameters, one has to add
an extra assumption on the initial data such as \eqref{decayfort}. Indeed, in the energy space,
one can construct an explicit example where convergence does not hold (see \cite{MMnonlinearity}).

\subsection{Decomposition and monotonicity result}
We recall a more precise   stability result related to the usual decomposition of the solution $u(t)$. See proof of Proposition 4 in \cite{MMas2}.
Define 
\begin{equation}\label{surphiff} 
  \psi(x)=\tfrac {2}{\pi}\arctan(\exp(-\tfrac x 4)),
\text{ so that } {\rm lim}_{+\infty} \psi=0,\ {\rm lim}_{-\infty} \psi=1.
\end{equation}
 
\begin{claim}[\cite{MMas2}]\label{pourstab}
Under the assumptions of Proposition \ref{ASYMPTOTIC}, 
there exist $C^1$  functions $\rho_1(t)$, $\rho_2(t)$, $c_1(t),$ $c_2(t)$,
defined on $[t_1,+\infty)$,  such that  $\eta(t,x)$ and $g(t)$ defined by
\begin{equation}\label{decompo}
\eta(t,x)=u(t,x)-R_1(t,x)-R_2(t,x),
  \text{ where for $ j=1,2,$
$R_j(t,x)=Q_{c_j(t)}(x-\rho_j(t))$,}
\end{equation}
\begin{equation}\label{defg} 
 g(t)=\int \left(\eta_x^2(t,x) +  \left(c+ \psi(x-m(t))\right) \eta^2(t,x)\right)dx,
\end{equation}
satisfy for all $t\in [t_1,+\infty),$
$
\int R_j(t)\eta(t) =\int  (x-\rho_j(t))R_{j}(t)\eta(t)=0,$ $j=1,2,$
\begin{equation}\label{MMnonlinearity}
\|\eta(t)\|_{H^1_c}^2\le g(t)\le K g(t_1) + K \exp(-   c^{-\frac 1{400}})
\leq K \alpha^2 c^{2q+1}  + K \exp(- c^{-\frac 1{400}}).
\end{equation}
\begin{equation}\label{mm5}
\left|\frac {c_1(t)}{c_1(t_1)}-1 \right|+c^{2q+1}  \left|\frac {c_2(t)}{c_2(t_1)}-1 \right|
 \leq K  g(t_1) +  K \exp(-   c^{-\frac 1{400}}),
\end{equation}
\begin{equation}\label{f4-54}
|c_1(t)-1|+c^{q+\frac 12} \left| \frac {c_2(t)}{c} -1 \right|
\leq K \alpha c^{q+\frac 12}+  K \exp(-   c^{-\frac 1{400}}).
\end{equation}
\end{claim}

Now, we recall   monotonicity results for quantities defined in $\eta(t)$, to be used in the proof of Theorem \ref{EXIST4}.
For  $0\leq t_0\leq t$, $x_0\geq 0$, $j=1,2$, let
\begin{equation*}
\begin{split}
& \mathcal{M}_j(t)=\int \eta^2 \psi_j ,\\
&\mathcal{E}_j(t)=\int  \left[\frac 12 \eta_x^2
-\frac 1{p{+}1}\left( (R_1{+}R_2{+}\eta)^{p+1} {-}(p{+}1) R_1^p \eta {-} (p{+}1) R_2^p \eta {-}(R_1{+}R_2)^{p+1}\right)\right] \psi_j,
\end{split}
\end{equation*}
where $\psi_1(t,x)=\psi(\tilde x)$,
$\tilde x= x-\rho_1(t)+x_0+\frac 12 (t-t_0)$,
and $\psi_2(t,x)=\psi(\sqrt{c}\tilde x_c)$, $\tilde x_c= x-\rho_2(t)+x_0+ \frac c2 (t-t_0) $.
  
\begin{claim}[\cite{MMas2}]\label{NEW}
Let $x_0>0$, $t_0>0$.
For all $t\geq t_0 $,
\begin{equation*}
\begin{split}
	& \frac d{dt}\left( c_1^{2q}(t) \int Q^2+ \mathcal{M}_1(t) \right) 
	   \leq K e^{-\frac 1{16}(t-t_0+x_0)} g_1(t) + Ke^{-\frac 1 {32} \sqrt{c} (t+T_c)},\\
	& \frac d{dt}\left(-\frac {2q}{2q+1} c_1^{2q+1}(t) \int Q^2+ 2 \mathcal{E}_1(t) 
	+\frac 1{100} \left( c_1^{2q}(t) \int Q^2+ \mathcal{M}_1(t) \right)
	\right) 
	   \\ &\qquad \leq K e^{-\frac 1{16}(t-t_0+x_0)} g_1(t)+ Ke^{-\frac 1 {32} \sqrt{c} (t+T_c)}.\\
	   & \frac d{dt}\left( \left(c_1^{2q}(t)+c_2^{2q}(t)\right) \int Q^2+ \mathcal{M}_2(t) \right) 
	   \leq K e^{-\frac {c\sqrt{c}}{16}    (t-t_0)} e^{-\frac {\sqrt{c} }{16}  x_0} \sqrt{c}\, g_2(t) + Ke^{-\frac 1 {32} \sqrt{c} (t+T_c)},\\
	& \frac d{dt}\left(-\frac {2q}{2q{+}1} \left(c_1^{2q+1}(t){+}c_2^{2q+1}(t)\right) \int Q^2+ 2 \mathcal{E}_2(t) 
	+\frac c{100} \left( \left(c_1^{2q}(t){+}c_2^{2q}(t)\right) \int Q^2+ \mathcal{M}_2(t) \right)
	\right) \\ &\qquad
	   \leq K e^{-\frac {c\sqrt{c}}{16}   (t-t_0)  } e^{-\frac {\sqrt{c} }{16}  x_0} c^{\frac 32} g_2(t)+ Ke^{-\frac 1 {32} \sqrt{c} (t+T_c)}.
\end{split}
\end{equation*}
\end{claim}

\section{Proofs of the main results ($p=4$)}

First, let us remark that for Theorems \ref{NONEX4}, \ref{EXIST4} and \ref{STAB4} (concerning the
case $p=4$) by considering $\tilde u(t,x)= \lambda^{\frac 13} u(\lambda^{\frac 32} t, \lambda^{\frac 12} x)$ with
$\lambda = \frac 1{c_1}$ instead of $u(t)$, we can restrict ourselves  to the case $c_1=1$ and $c_2=c\leq \epsilon_0\ll 1$ without loss of generality.
We consider $0<c<c_0$, where $c_0$ small enough so that Sections 2, 3 and 4 apply for any
$0<c'<2c_0$.

\subsection{Nonexistence of pure $2$-soliton solutions. Proof of Theorem~\ref{NONEX4}}

This section is devoted to the proof of Theorem \ref{NONEX4}.
First, we recall the following.

\begin{proposition}
\label{2SOL}
    Let $p=2,3$ or $4$. Let $0<c<c_0$, for $c_0$ small enough.
    \begin{enumerate}
    \item Existence and exponential decay. 
    Let $x_1,x_2\in \mathbb{R}$. There exists a unique solution  
    $U_{c,x_1,x_2}(t,x)=U(t,x)\in C(\mathbb{R},H^1(\mathbb{R}))$
     of \eqref{gkdv} such that 
    \begin{equation}\label{2SOL0}
        \lim_{t\to -\infty} \|U(t)-Q(x-t-x_1)-Q_c(x-ct-x_2)\|_{H^1(\mathbb{R})}=0.
    \end{equation}
    Moreover, $U(t)$  satisfies,
    for all  $t\leq \frac {x_2-x_1}{1-c}-\frac {T_c} {32}$,
        \begin{equation}\label{2SOL1}
            \|U(t)-Q(x-t-x_1)-Q_c(x-ct-x_2)\|_{H^1(\mathbb{R})}
            \leq K e^{\frac 14{\sqrt{c}((1-c)t-(x_2-x_1))}},
        \end{equation}
        \item Uniqueness of the asymptotic $2$-soliton solution at $-\infty$:
        If $u(t)$ is an $H^1$ solution of \eqref{gkdv} satisfying
        \begin{equation}\label{2SOL2}
            \lim_{t\to -\infty}\|u(t)- Q(x-\rho_1(t))-Q_c(x-\rho_2(t))\|_{H^1(\mathbb{R})}=0,
        \end{equation}
        for $\rho_1 $, $\rho_2 :\mathbb{R}\to \mathbb{R}$,  then $u(t)$ 
        satisfies \eqref{2SOL0} for some $x_1$, $x_2$, and so $u(t)\equiv U_{c,x_1,x_2}(t)$.
    \end{enumerate}
\end{proposition}

This result was essentially proved in \cite{Martel}, using tools from \cite{MM1} and \cite{MMT}.
However, some statements in Proposition \ref{2SOL} are slightly more precise than the main result
in \cite{Martel}, so we justify them  in Appendix D.2.

\medskip

The main ingredient of the proof of Theorem \ref{NONEX4} is the following proposition
related to the approximate solution constructed in Section 2.
We keep the notation of Section 3, in particular, $v(t,x)$, $b_{2,0}$ and $V_1$.

\begin{proposition}\label{VDIESE} Let $p=4$.
    Let $\Delta$ and $\Delta_c$ be defined by \eqref{defDeltabis}.
    Let
    \begin{equation}\label{Vdiese}
        v_{\#}(t,x)=v(t,x)+w_{\#}(t,x)\quad \text{where}\quad
        w_{\#}(t,x)=(Q_c^2)'(y_c) b_{2,0} (1+V_1(y)),
    \end{equation}
    and
    \begin{equation}\label{Vdiese1}
        S_{\#}(t,x)=\partial_{t}v_{\#} + \partial_{x} (\partial_{x}^2 v_{\#}- v_{\#} + v_{\#}^4),
    \end{equation}
    where $v$ is the function constructed in Proposition \ref{AVRIL4}.
    Then, for all $0<c<c_0$, for $c_0$ sufficiently small,
    \begin{enumerate}
        \item Approximate solution: for $j=0,1,2$,
        \begin{equation}\label{Vdiese2}
            \forall t\in [-T_c,T_c],\quad
            \|\partial_x^j S_{\#}(t)\|_{L^2} \leq c^{3/2},
        \end{equation}
        \item Closeness to a pure two soliton at $t=-T_c$:
        \begin{equation}\label{Vdiese3}
            \|v_{\#}(-T_c)-Q(.{+}\tfrac{\Delta}2)-Q_c(.{-}(1{-}c)T_c{+}\tfrac {\Delta_c}2)\|_{H^1}\leq K c,
        \end{equation}
        \item Non-matching with a pure two-soliton solution at $t=T_c$:
        \begin{equation}\label{Vdiese4}
            \|v_{\#}(T_c)-Q(.{-}\tfrac{\Delta}2)-Q_c(.{+}(1{-}c)T_c{-}\tfrac {\Delta_c}2)
            -2b_{2,0} (Q_c^2)'(.{+}(1{-}c)T_c{-}\tfrac {\Delta_c}2)\|_{H^1}\leq K c.
        \end{equation}
    \end{enumerate}    
\end{proposition}

\noindent\emph{Remark.} Recall that $\|(Q_c^2)'(y_c)\|_{H^1}> K c^{11/12}$ and $b_{2,0}<0$. Thus at $T_c$, the function
$v_{\#}$ differs from a two-soliton solution of a factor $c^{11/12}$. At $-T_c$ it is close to a two-soliton solution up to a factor $c$ and it is an approximate solution of the gKdV equation in the sense \eqref{Vdiese2}. This will be sufficient to prove 
Theorem \ref{NONEX4} applying Proposition \ref{INTERACT4}.

The function $v_{\#}(t,x)$ is not exactly of the form imposed by Proposition \ref{SYSkl4}.
Indeed, the function $1+V_1$ is even, and thus the function $w_{\#}(t,x)$ has not the required structure.
This will have no consequence in applying Proposition \ref{INTERACT4}, which does not rely on the parity structure.
In contrast, the presence of $w_{\#}$ in $v_{\#}$ is definitely a problem to follow the procedure of Proposition \ref{SYSkl4}.
Indeed, this term creates a new term $F_{5,0}$ which has a nonzero even part, not orthogonal to $Q$, which is a problem to determine a suitable $A_{5,0}$. Thus, we can not improve \eqref{Vdiese2} up to any power. However, \eqref{Vdiese2} is sufficient for
 our purposes, and the function $v_{\#}$ is closer to a $2$-soliton solution at $t=-T_c$ than the function $v$ itself.

It would be interesting to investigate further improvements of the function $v_{\#}$ since it would help understanding the behavior for $t>0$ of solutions which are pure two-soliton solutions at $t\to -\infty$.

\medskip

\noindent\emph{Proof of Proposition \ref{VDIESE}.}
We have
\begin{align*}
    S_{\#}(t,x) & =\partial_{t}v_{\#} + \partial_{x} (\partial_{x}^2 v_{\#}- v_{\#} + v_{\#}^4) \\
                & =S(t,x)+\partial_x((v+w_{\#})^4)-v^4-4Q^3w_{\#})+\partial_t w_{\#}-\partial_x (\overline{\mathcal{L}} w_{\#}),
\end{align*}
where $\overline{\mathcal{L}}$ is defined in \eqref{defLy}.

\medskip

\noindent\emph{1a. Estimate of the linear part in $S_{\#}$}.

We estimate $\partial_t w_{\#}-\partial_x (\overline{\mathcal{L}} w_{\#})$, 
where $w_{\#}(t,x)=(Q_c^2)'(y_c) b_{2,0} (1+V_1(y))$.
Recall that from Claim \ref{ClaimVun}, $\mathcal{L}(1+V_1)=1-4 Q^3 + \mathcal{L}V_1=1$, and
thus $(\mathcal{L}(1+V_1))'=0$. Claim \ref{QckB} gives an explicit expression for
$\partial_t w_{\#}-\partial_x (\overline{\mathcal{L}} w_{\#})$, where the first term in the
second-hand member is zero. For the other terms, we use:
  $$(1+V_1)'=V_1'\in \mathcal{Y}, ~Q^3(1+V_1)\in \mathcal{Y};
  \quad \|(Q_c^2)'\|_{L^\infty}=K c^{7/6}, ~\|\beta\|_{L^\infty}\leq K c^{1/3};$$
   $$\|(Q_c^2)''\|_{L^\infty}=K c^{5/3}, ~\|(Q_c^2)^{(4)}\|_{L^2}= K c^{29/12};$$
so that
$
    \|\partial_t w_{\#}-\partial_x (\overline{\mathcal{L}} w_{\#})\|_{L^2} \leq K c^{3/2}.
$

We 
 obtain, for all $j=0,1,2,$
$
    \|\partial_x^j\left(
    \partial_t w_{\#}-\partial_x (\overline{\mathcal{L}} w_{\#})\right)\|_{L^2} \leq K_j c^{3/2}.
$

\medskip

\noindent\emph{1b. Estimate of the nonlinear part in $S_{\#}$}

Note that
$
    (v+w_{\#})^4-v^4- 4 Q^3 w_{\#} = 4 (v^3-Q^3) w_{\#} + 6 v^2 w_{\#}^2 + 4 v w_{\#}^3
    + w_{\#}^4,
$
so that
\begin{align*}
    \partial_x\left[(v+w_{\#})^4-v^4- 4 Q^3 w_{\#}\right] & = 
    4 \partial_x (v^3 - Q^3) w_{\#} + 4 (v^3-Q^3) \partial_x w_{\#} 
    + 6 \partial_x (v^2w_{\#}^2)\\
    &\quad +4 \partial_x ( v w_{\#}^3 ) + \partial_x (w_{\#}^4).
\end{align*}
Moreover,
\begin{align*}
     \partial_x (v^3 - Q^3) &=
     \partial_x (v-Q) (v^2+ vQ +Q^2) + (v-Q)\partial_x(v^2+vQ+Q^2),
\end{align*}
\begin{equation*}
    \partial_x(v^2+ vQ + Q^2)= \partial_x (3Q^2 + (v^2-Q^2) + (v-Q)Q).
\end{equation*}
Thus,
\begin{equation*}
    \|\partial_x(v^3 - Q^3)\|_{L^2}\leq 
    K ( \|\partial_x (v-Q)\|_{L^2} + \|v- Q\|_{L^\infty}) \leq K c^{1/3}.
\end{equation*}
We also have
$$\|w_{\#}\|_{L^\infty}\leq K c^{7/6},~\|v^3-Q^3\|_{L^2} \leq K c^{1/2}, ~
\|\partial_x w_{\#}\|_{L^\infty}\leq K c^{5/3},~\|w_{\#}^2\|_{L^\infty}\leq K c^{5/3}.$$
Thus,
\begin{equation*}
    \|\partial_x\left[(v+w_{\#})^4-v^4- 4 Q^3 w_{\#}\right]\|_{L^2}
    \leq K c^{3/2}.
\end{equation*}
Similarly, for $j=1,2,$
$
    \|\partial_x^{(j+1)}\left[(v+w_{\#})^4-v^4- 4 Q^3 w_{\#}\right]\|_{L^2}
    \leq K_j c^{3/2}.
$

Taking  $k_0$, $\ell_0$ large enough, so that
$\|\partial_c^{(j)} S\|_{L^2} \leq K c^{3/2}$, 
by Proposition \ref{AVRIL4}, we have proved
$\|\partial_c^{(j)} S_{\#}\|_{L^2} \leq K c^{3/2}$.

\medskip

\noindent\emph{2. Analysis at $t=\pm T_c$.}
By the proof of Proposition \ref{AVRIL4} (see \eqref{3-16}), we have
\begin{equation*} \begin{split}
&    \|v(-T_c)-Q(.+\tfrac \Delta 2)-Q_c(.-(1-c)T_c+\Delta_c/2)+b_{2,0} (Q_c^2)'(y_c)\|_{H^1}\leq K c,\\
 &   \|v(T_c)-Q(.-\tfrac \Delta 2)-Q_c(.+(1-c)T_c-\Delta_c/2)-b_{2,0} (Q_c^2)'(y_c)\|_{H^1}\leq K c.
\end{split}\end{equation*}
Note that by the definition of $v_{\#}$ and Claim \ref{SIX}, we have
\begin{equation*}
    \|v(\pm T_c)-\left(v_{\#}(\pm T_c)-b_{2,0} (Q_c^2)'(y_c)\right)\|_{H^1}
    =\|b_{2,0} (Q_c')^2(y_c)V_1(y)\|_{H^1}\leq K c^{7/6}.
\end{equation*}
Thus,
\begin{equation*}\begin{split}
&    \|v_{\#}(-T_c)-Q(.+\tfrac \Delta 2)-Q_c(.-(1-c)T_c+\Delta_c/2)\|_{H^1}\leq K c\\&
	\|v_{\#}(T_c)-Q(.-\tfrac \Delta 2)-Q_c(.+(1-c)T_c-\Delta_c/2)- 2 b_{2,0} (Q_c^2)'(y_c)\|_{H^1}\leq K c.
	\end{split}
\end{equation*}
By $\|(Q_c^2)'(y_c)-(Q_c^2)'(.+(1-c)T_c-\Delta_c/2)\|\|_{H^1}
=\|(Q_c^2)' -(Q_c^2)'(.-\Delta_c/2)\|\|_{H^1}\leq K c^{\frac 76}$, since $\Delta_c$
is a constant independent of $c$, we obtain the result.

\medskip

\noindent\emph{Proof of Theorem \ref{NONEX4}.}

Step 1. Proof of nonexistence of a pure 2-soliton solution.
First, we claim that if there exists a global $2$-soliton
solution, then the speeds parameters at $+\infty$, $c_1^+<c_2^+$ and at $-\infty$, $c_1^-<c_2^-$ satisfy $c_1^+=c_1^-$ and $c_2^+=c_2^-$.
Indeed, by the conservation of mass and energy, and strong
limit in $H^1(\mathbb{R})$, the following holds ($q=\frac {1}{p-1}-\frac 14$)
\begin{equation*}
         (c_1^+)^{2q}+(c_2^+)^{2q}=(c_1^-)^{2q}+(c_2^-)^{2q},
 \quad
    (c_1^+)^{2q+1}+(c_2^+)^{2q+1}=(c_1^-)^{2q+1}+(c_2^-)^{2q+1}.
\end{equation*}
Set 
$
    \gamma=\frac {2q+1}{2q},$ $ b=\left(\frac {c_1^+}{c_1^-}\right)^{2q},
   $ $ a^+=\frac {c_1^+}{c_2^+}<1,$ $ a^-=\frac {c_1^-}{c_2^-}<1.
$
The first identity yields
$b(1+a^+)=1+a^-$, and the second identity yields
$b^\gamma(1+(a^+)^\gamma)=1+(a^-)^\gamma$. Thus,
\begin{equation*}
    \left(\frac {1+a^+}{1+a^-}\right)^{\gamma}=\frac {1+(a^+)^\gamma}{1+(a^-)^\gamma}
    \quad \text{and}\quad
    \frac {1+(a^+)^\gamma}{\left({1+a^+}\right)^{\gamma}}=
    \frac {1+(a^-)^\gamma}{\left({1+a^-}\right)^{\gamma}}.
\end{equation*}
The function $x\mapsto \frac {1+x^\gamma}{(1+x)^{\gamma}}$
is strictly decreasing on $[0,1]$, thus $a^+=a^-$ and  $b=1$.

\medskip

\noindent\emph{(i) Behavior at $-\infty$.}

Let $u(t)$ be an asymptotic $2$-soliton solution at $-\infty$ with speed parameters $1$ and $c$
($c$ small enough), 
in the sense of Definition \ref{DEF1}. Then, by the uniqueness part  of Proposition \ref{2SOL},
there exists $x_1^-$, $x_2^-\in \mathbb{R}$ such that, for all $t \leq  \frac{x_2^--x_1^-}{(1-c)}-\tfrac 1{32} {T_c}$,
\begin{equation}\label{truc1}
        \|u(t)-Q(x-t-x_2^-)-Q_c(x-ct-x_1^-)\|_{H^1}
    \leq K e^{\frac 14 {\sqrt{c}((1-c)t-(x_2^--x_1^-))} }.
\end{equation}
Let  
$$  T_c^-=T_c+\frac{x_1^- - x_2^-}{1-c}+\frac 12 \frac{\Delta-\Delta_c}{1-c}
\geq -\frac {x_2^--x_1^-}{1-c} +\frac 14 T_c
\quad \text{and}\quad   a=\tfrac \Delta 2 - (T_{c}^--x_{2}^-).$$
Recall from \eqref{defDeltabis} that
$|\Delta|\leq K c^{\frac 16}$ and $\Delta_c$ is a constant independent of $c$.
Then, applying \eqref{truc1} to $t=-T_c^-$, we obtain
\begin{align*}
    & \|u(-T_{c}^-,.+a)-Q(.+\tfrac \Delta 2)-
    Q_c(.-(1-c)T_{c}+\tfrac {\Delta_{c}} 2)\|_{H^1}
     \\& \leq K e^{\frac 14 {\sqrt{c}(-(1-c)T_{c}^- -(x_2^--x_1^-))} } \leq K e^{-\frac 14 {\sqrt{c}((1-c)T_{c}+\frac 12({\Delta_{c}-\Delta}))} }
     \leq K c,
\end{align*}
for $c$ small enough.
By translation in time and space, we can assume $T_{c}^-=T_{c}$ and $a=0$, so that
\begin{equation*}
	\|u(-T_{c})-Q(.+\tfrac \Delta 2)-
    Q_c(.-(1-c)T_{c}+\tfrac {\Delta_{c}} 2)\|_{H^1} \leq K c,
\end{equation*}	
(i.e., we consider $\tilde u(t,x)=u(t-T_c^-+T_c,x+a)$ instead of $u(t,x)$, and we still call it
$u(t)$).

\medskip

\noindent\emph{(ii) Behavior at $t=T_c$.}

By \eqref{Vdiese3}, and the above estimate, we deduce
\begin{equation*}
	\|u(-T_{c})-v_{\#}(-T_{c})\|_{H^1} \leq K c.
\end{equation*}	
Now, we apply Proposition \ref{INTERACT4} for $v_{\#}$ concerning the interaction region, with $\theta=1-\frac 1{100}$ and
$T_{0}=-T_{c}$. Thus, 
\begin{equation*}
	\forall t\in [-T_{c},T_{c}],\quad
	\|u(t)-v_{\#}(t,.-\rho(t))\|_{H^1} \leq K c^{1-\frac 1{100}},
\end{equation*}
for some $\rho(t)$ satisfying $|\rho'(t)|\leq K c^{1-\frac 1{100}}$.
In particular,
$
	\|u(T_{c})-v_{\#}(T_c,.-\rho(T_{c}))\|_{H^1)} \leq K c^{1-\frac 1{100}},
$
and so by Proposition \ref{VDIESE}, we obtain for $a_-$, $b_-\in \mathbb{R}$ such that
$a_- -b_->\frac 12 T_c$,
\begin{equation}\label{atTc2}
	\|u(T_{c})-Q(.-a_-)-Q_c(.-b_-)
	-2b_{2,0} (Q_c^2)'(.-b_-)\|_{H^1 } \leq K c^{1-\frac 1{100}}.
\end{equation}

\noindent\emph{(iii) Behavior as $t\to +\infty$.}

First, since $\|(Q_c^2)'\|_{H^1}\leq K c^{\frac {11}{12}}$,
estimate \eqref{atTc2} implies that for $t=T_{c}$:
\begin{equation}\label{cdixbis}
	\|u(T_{c})-Q(.-a_-)-Q_c(.-b_-)\|_{H^1 } \leq K c^{\frac{11}{12}}.
\end{equation}
We apply Proposition \ref{ASYMPTOTIC} to $u(t)$ (stability of the
$2$-soliton structure after interaction) with $\alpha=K c^{\frac 13}$, so that
for $w(t)=u(t)-Q(.-\rho_1(t))-Q_c(.-\rho_2(t))$
\begin{equation}\label{uneetoile}
	\forall t\geq T_c,\quad
	\sqrt{c} \|w(t)\|_{H^1}\leq  \|w_x(t)\| + \sqrt{c} \|w(t))\|_{L^2} \leq K c^{\frac{11}{12}},
\end{equation}
with $\rho_1(t)$, $\rho_2(t)$ satisfying
\begin{equation}\label{deuxetoile}
	|\rho_1(T_c)-a_-|\leq K c^{\frac {11}{12}},\
	|\rho_2(T_c)-b_-|\leq K c^{\frac 1{3}},\
\quad
	\forall t\geq T_c,\
	\rho_1(t)-\rho_2(t)  \geq \frac {T_c} 2 + \frac 12(t-T_c).
\end{equation}
Assume now that $u(t)$ is also an asymptotic $2$-soliton solution at $+\infty$.
By Proposition \ref{2SOL} (applied to $+\infty$) there exist $x_1^+$, $x_2^+$ such that,
for all $	  t\geq \frac {x_2^+-x_1^+} {1-c} + \frac {T_c}{32}$,
\begin{equation}\label{troisetoile}
	\|u(t)-Q(.-t-x_1^+)-Q_c(.-ct-x_2^+)\|_{H^1 }
	\leq K e^{-\frac {\sqrt{c}}4 ((1-c) t - (x_2^+-x_1^+))}.
\end{equation}
We define 
$$T_c^+=\frac {x_2^+-x_1^+}{1-c} +\frac {T_c}{32}.$$
By \eqref{uneetoile} and \eqref{troisetoile}, we have for all 
$t\geq \max(T_c,T_c^+)$,
\begin{equation}\label{quatreetoile}
	|\rho_1(t)-(t+x_1^+)|\leq K c^{\frac {5}{12}},\quad
	|\rho_2(t)-(ct+x_2^+)|\leq K c^{-\frac 1{6}}.
\end{equation}
This is proved by considering the smallness of the $L^2$ norm of
$Q(.-\rho_1(t))+Q_c(.-\rho_2(t))-Q(.-t-x_1^+)-Q_c(.-ct-x_2^+)$
in the two regions $x> \frac 12 (\rho_1(t)+\rho_2(t))$ and
$x< \frac 12 (\rho_1(t)+\rho_2(t))$) and the fact that for $a$ small
\begin{equation}\label{poura}
|a| \leq K \|Q-Q(.-a)\|_{L^2}\quad 
|a| \leq K c^{-\frac 7{12}} \|Q_c-Q_c(.-a)\|_{L^2}.
\end{equation}

Let us prove that $T_c>T_c^+$. By contradiction, if
$T_c^+\geq T_c$, then by \eqref{quatreetoile} we have
$$
|\rho_1(T_c^+)-\rho_2(T_c^+)|\leq |(1-c)T_c^++x_1^+-x_2^+|+K c^{-\frac 1{6}}
=\frac {T_c}{32} (1-c)+K c^{-\frac 1{6}}
\leq \frac {T_c}{30}.$$
From \eqref{deuxetoile},
$$
|\rho_1(T_c^+)-\rho_2(T_c^+)|\geq \frac 14 T_c +\frac 12 (T_c^+-T_c)\geq \frac 14 T_c.$$
We obtain a contradiction from these two estimates and thus  $T_c>T_c^+$.

\medskip

\noindent\emph{(iv) Conclusion of the proof.}

Let $a_+=T_c+x_1^+$ and $b_+=cT_c+x_2^+$. 
By \eqref{troisetoile}, we know
\begin{equation*}
	\|u(T_c)-Q(.-a_+)-Q_c(.-b_+)\|_{H^1} 
	\leq K e^{-\frac {\sqrt{c}} 4 (1-c)(T_c-T_c^++\frac 1 {32} T_c)}
	\leq K c.
\end{equation*}
Thus, from \eqref{atTc2}-\eqref{cdixbis} and Proposition \ref{ASYMPTOTIC}
\begin{equation}\begin{split}
 & |a_--a_+|\leq |a_--\rho_1(T_c)|+|\rho_1(T_c)-(T_c+x_1^+)| \leq K c^{\frac {11}{12}}, \\
 & |b_--b_+|\leq |b_--\rho_2(T_c)|+|\rho_2(T_c)-(cT_c+x_2^+)| \leq K c^{\frac {1}{3}},
\end{split}\end{equation}
and
\begin{equation*}
	\|(Q(.-a_-)-Q(.-a_+)) + (Q_c(.-b_-)-Q_c (.-b_+))
	+ 2 b_{2,0} (Q_c^2)'(.-b_-)\|_{H^1}\leq K c^{1-\frac 1{100}}.
\end{equation*}
Considering the $L^2$ norm in the region $x<\frac 12 (a_++b_+)$, we
obtain
\begin{equation*}
\|Q_c+2b_{2,0} (Q_c^2)'-Q_{c}(.-(b_+ - b_-))\|_{L^2} \leq K c^{1-\frac 1{100}},
\end{equation*}
where $\min(a_-,a_+)>\max(b_-,b_+)+\frac 14 T_c$.
By scaling, it gives for   $b_c=\sqrt{c}(b_+ - b_-)$,
\begin{equation*}
	\| Q+ 2b_{2,0} c^{5/6} (Q^2)' -Q(.-b_c)\|_{L^2}\leq K c^{\frac {11}{12}-\frac 1{100}},
\end{equation*}
where $|b_c|\leq K c^{\frac {5}{6}}$. Thus,
$Q(x)-Q(x-b_c)=\lambda_c c^{5/6} Q'(x) + c^{5/6} o(1)$, where $|\lambda_c|\leq K$, 
so that 
\begin{equation*}
\|\lambda_c Q'- 2 b_{2,0}(Q^2)'\|_{L^2}=o(1),
\end{equation*}
which is a  contradiction with the fact
$b_{2,0}\neq 0$ (Lemma \ref{RESUME}).

\medskip

\noindent{Step 2. Behavior as  $t\to +\infty$ of $u(t)$.}

As in the previous step, we consider the solution $u(t)$ which is
an asymptotic $2$-soliton solution at $-\infty$  i.e. satisfying \eqref{truc1}.
Recall that we have just proved:

- $u(t)$ is not an asymptotic  $2$-soliton solution at $+\infty$.

- There exist $\rho_1(t)$, $\rho_2(t)$ such that $w(t,x)=u(t,x)-(Q(x-\rho_1(t))+Q_c(x-\rho_2(t)))$,
satisfies \eqref{uneetoile}, \eqref{deuxetoile}, in particular
\begin{equation}\label{upperb}
\forall t\geq T_c,\quad \|w(t)\|_{H^1}\leq c^{-\frac 12}\|w(t)\|_{H^1_c}
\leq K c^{\frac {5}{12}}.
\end{equation}

\medskip

(i) Stability properties of $u(t)$ for $t\geq T_c$. 
First, we claim
\begin{equation}\label{le29nov}
\int_{x>0} x^2 u^2(T_c,x) dx\leq K.
\end{equation}
This follows directly from integration of the following estimate: for all $x_0>0$,
\begin{equation}\label{decayasymp}
\begin{split}
\int_{x>x_0} u^2(T_c,x+\rho_1(T_c)) dx \leq K e^{-\frac 1{16} x_0} + K \exp(-c^{-\frac 1{400}})\,
e^{-\frac 1{16} \sqrt{c} x_0}.
\end{split}
\end{equation}
Let us prove \eqref{decayasymp}.
On the one hand, by monotonicity arguments on $u(t)$ as in Lemma 1 of \cite{Martel}
\begin{equation}\label{le29I}
\int u^2(T_c,x)\psi(x-\rho_1(T_c)-x_0) dx
\leq \int u^2(-T_c,x)\psi(x-\rho_1(-T_c)-x_0-\tfrac {T_c}2) dx
+ K e^{-\frac 1 {16} x_0}.
\end{equation}
On the other hand, using $\mathcal{I}_{\sigma, y_0}$ for
$\sigma=c$, $y_0=\rho_1(-T_c) + x_0+\tfrac {T_c}2$, we get for any $t<-T_c$,
\begin{equation*}
\begin{split}
& \int u^2(-T_c,x)\psi(\sqrt{c}(x-\rho_1(-T_c) -x_0-\tfrac {T_c}2)) dx
\\ & \leq \int u^2(t,x)\psi(\sqrt{c}(x-\rho_1(-T_c)-x_0-\tfrac {T_c}2-\tfrac c 4 t)) dx
+ K e^{-\frac  1{16} {\sqrt{c}} (x_0+\frac 12 T_c)}.
\end{split}\end{equation*}
By \eqref{truc1} and letting $t\to -\infty$, we obtain
\begin{equation}\label{le29II}
\begin{split}
& \int u^2(-T_c,x)\psi(\sqrt{c}(x-\rho_1(-T_c) -x_0-\tfrac {T_c}2)) dx
 \leq   K \exp(-c^{-\frac 1{400}}) \, e^{-\frac  1{16} {\sqrt{c}}  x_0 }.
\end{split}\end{equation}
Therefore, from \eqref{le29I}, \eqref{le29II}, 
$$ \int_{x>x_0} u^2(T_c,x+\rho_1(T_c)) dx \leq \frac 12 
\int u^2(T_c,x)\psi(x-\rho_1(T_c)-x_0) dx
\quad \text{and}\quad \psi(\sqrt{c}y)\geq \frac 12 \psi(y),$$
we obtain \eqref{decayasymp}.

Now, from \eqref{le29nov} and \eqref{cdixbis}, we can apply Proposition \ref{ASYMPTOTIC}
to $u(.+T_c)$, for $t\geq 0$ (with $\alpha=K c^{\frac 13}$).
It follows that there exists $c_1^+,\ c_2^+>0$, $x_1^+, x_2^+\in \mathbb{R}$ such that
$w^+(t)= u(t)-Q_{c_1^+}(.-x_1^+-c_1^+t)
-Q_{c_2^+}(.-x_2^+-c_2^+t)$ satisfies
\begin{equation}\label{asymptoti}
\lim_{t\to +\infty} \| w^+(t)\|_{H^1(x>ct/{10})}=0
\quad \text{with}\quad |c^+-1|\leq K c^{\frac {11}{12}}, \quad 
\left|\frac {c_2^+}{c}-1\right|\leq K c^{\frac 13}.\end{equation}
Note also that  from the stability \eqref{cdixbis} and \eqref{asymptoti}, we obtain the following upper bound on $w^+(t)$ for $t$ large enough:
\begin{equation*}\begin{split}
\|w^+(t)\|_{H^1}\leq \|w^+(t)\|_{H^1(x<\frac 1{10} c t)} + \|w^+(t)\|_{H^1(x>\frac 1{10} c t)}
\leq \|w(t)\|_{H^1(x<\frac 1{10} c t)}+ o(1)\leq K c^{\frac {5}{12}}.
\end{split}\end{equation*}
Therefore, to finish the proof of Theorem \ref{NONEX4}, we only have to prove
the lower bounds on $w^+(t)$, $c_1^+-1$ and $1-\frac {c_2^+}{c}$.

\smallskip

(ii) Lower bounds on the defects.
Let $\eta(t)$, $g(t)$ and $c_j(t)$ ($j=1,2$)
 be defined from $u(t)$ for $t\geq T_c$ as in Claim \ref{pourstab} and
 satisfying
\begin{equation}\label{rolls}
\|\eta(t)\|_{H^1(x\geq \frac c{10} t)}\to 0,\quad 
 c_j(t)\to c_j^+ \quad \text{as $t\to +\infty$ ($j=1,2$)}.
\end{equation}
In particular, it is sufficient
to prove the lower bounds on $\eta(t)$  to obtain lower bounds  
  on $w^+(t)$ for large time. We claim
  \begin{equation}\label{lower}
  \forall t\geq T_c,\quad 
\|\eta(t)\|_{H^1_c}\geq K_1 c^{\frac {17}{12}} \quad (K_1>0).
\end{equation}
Proof of \eqref{lower}. 
To prove this lower bounds  using the defect $(Q_c^2)'$ in \eqref{atTc2}, we need
to apply an argument of stability backwards in time, locally around the soliton $R_2(t)$.
For this, we will use monotonicity type results on $\eta(t)$ as in Claim \ref{NEW}.

First, we claim
\begin{equation}\label{claimeta}
\int_{x\leq \rho_2(T_c)+\frac 14 T_c}
  \eta^2(T_c,x) dx\geq K_0 c^{\frac {11}{6}}\quad (K_0>0). 
\end{equation}
Proof of \eqref{claimeta}.
Let $\varepsilon>0$ to be fixed later and assume for the sake of contradiction
that $\int_{x\leq \rho_2(T_c)+\frac 14 T_c}
   \eta^2(T_c,x)dx \leq \varepsilon^2 c^{\frac {11}{6}}.$ 
Recall from \eqref{atTc2} that
\begin{equation}\label{rappel5.11}
\|u(T_c)-Q(.-a_-)-Q_c(.-b_-)-2 b_{2,0} (Q_c^2)'(.-b_-)\|_{L^2}\leq Kc.
\end{equation}
Thus, as in step 1  (iv), we obtain for $c$ small enough,
$$\|Q_c(.-b_-)+2 b_{2,0} (Q_c^2)'(.-b_-) 
-Q_{c_2(T_c)}(.-b_+)\|_{L^2}\leq K  \varepsilon c^{\frac {11}{12}},
$$
and after scaling
$$\|Q +2 b_{2,0}c^{\frac 5{6}} (Q^2)'  
-Q_{\Lambda}(.-b_c)\|_{L^2}\leq K  \varepsilon c^{\frac {5}{6}},
$$
for $\Lambda=\frac {c_2(T_c)}c$, $b_c=\sqrt{c}(b_+-b_-)$.
From orthogonality of even and odd functions in $L^2$ and parity of $\frac {d^k}{dc^k} Q_c$
for any $k\geq 0$, we obtain
$$\|Q +2 b_{2,0} c^{\frac 5{6}}  (Q^2)'  
-Q(.-b_c)\|_{L^2}\leq K  \varepsilon c^{\frac {5}{6}},
$$
which is a contradiction for $\varepsilon$ small enough, as in Step 1 (iv) ($b_{2,0}\neq 0$).
Thus, \eqref{claimeta} is proved.

\smallskip

Let $\varepsilon>0$ to be fixed later and assume for the sake of contradiction
that for some $  t'\geq T_c$, 
\begin{equation}\label{lowerII}
 \|\eta(t')\|_{H^1_c} \leq \varepsilon c^{\frac {17}{12}}.
 \end{equation}
Let
$\tilde \psi_2(t,x)=1-\psi(\sqrt{c}(x-\rho_2(t)-\tfrac 14 {T_c} -\tfrac 12 (t_0-t))$
where   $\psi$ is defined in \eqref{surphiff} and
$$
\tilde{\mathcal{M}}_2(t)=\int \eta^2t) \tilde \psi_2,\quad
\tilde{\mathcal{E}}_2(t)=\int \Big[\tfrac 12 \eta_x^2 - \tfrac 1{5} ((R_1+R_2+\eta)^5
-5 R_1^4 \eta - 5 R_2^4 \eta - (R_1+R_2)^5)\Big] \tilde \psi_2.
$$
 From \eqref{lowerII} and the properties of $R_1$, $R_2$,  we have
$c\tilde{\mathcal{M}}_2(t')+|\tilde{\mathcal{E}}_1(t')| \leq K \varepsilon^2 c^{\frac {17}6}.$
Thus from Claim \ref{NEW}, integrated on $[T_c,t']$,
we have
\begin{equation*}\begin{split}
&
(c_2^{2q}(T_c)-c_2^{2q}(t')) \int Q^2  \leq - \tilde{\mathcal{M}}_2(T_c) + K \varepsilon^2 c^{\frac {17}6} ,
\\ &
\left( \frac {2q}{2q+1} (c_2^{2q+1 }(T_c) -c_2^{2q+1}(t')) 
-\tfrac c{100} (c_2^{2q}(T_c)-c_2^{2q}(t'))\right)
 \int Q^2  \\& \quad\geq  2 \mathcal{E}_2(T_c) +\tfrac c{100} \mathcal{M}_2(T_c)
  - K \varepsilon^2 c^{\frac {17}6}.
\end{split}\end{equation*}
From this, using the coercive functional of $\eta(t)$:
$\tilde{\mathcal{E}}_2(t)+\frac 12 c_2(t') \tilde{\mathcal{M}}_2(t)$, 
and proceeding as in \cite{MMas2} Appendix B.3, we
obtain successively:
$$
|c_2(T_c)-c_2(t')|\leq  K \int (\eta_x^2+c\eta^2)(T_c) \tilde \psi_2  
+ K\varepsilon^2 c^{\frac {17}6},
$$
$$
\int (\eta_x^2+c\eta^2)(T_c) \tilde \psi_2 
\leq K \varepsilon^2 c^{\frac {17}6} + K |c_2(T_c)-c_2(t')|^2 \leq K' \varepsilon^2 c^{\frac {17}6},
$$
which contradicts \eqref{claimeta} for $\varepsilon$ small enough.

\smallskip

Finally, we prove the lower bounds on $c_1^+-1$ and $1-\frac {c_2^+}{c}$,
using  the two conservation laws, written as $t\to \pm\infty$ and the bounds on $w^+(t)$.
By \eqref{truc1} and \eqref{asymptoti}, we have, for $t$ large,
\begin{equation*}
\begin{split}
& \int u^2(0)= \int Q^2 + \int Q_c^2  = \int Q_{c_1^+}^2 + \int Q_{c_2^+}^2 +\int (w^+)^2(t)+o(1),
\\ & 
E(u(0))=E(Q)+E(Q_c)=E(Q_{c_1^+}) + E(Q_{c_2^+}) + E(w^+(t)) +o(1).
\end{split}
\end{equation*}
By Gagliardo-Nirenberg inequality and the estimate $\|w^+(t)\|_{H^1}\leq K c^{\frac 5{12}}$, we have
$\int (w^+)^5\leq K \|w^+\|_{H^1}^3 \int (w^+)^2 \leq K c^{\frac 54} \int (w^+)^2$
and thus
 for $t$ large enough   
$$\left|E(w^+(t)) - \tfrac 12 \int (w_x^+(t))^2\right|\leq K c^{\frac 54}\int (w^+(t))^2 .$$
Thus, by Claim \ref{LemmaA1},
we obtain for $t$ large
\begin{align}
& \left|(c^{2q} - (c_2^+)^{2q}) + (1-(c_1^+)^{2q}) -
\frac 1{\int Q^2} \int (w^+(t))^2 \right|\leq K c^4, \label{UM}
\\ & 
 \left|(c^{2q+1} - (c_2^+)^{2q+1}) + (1-(c_1^+)^{2q+1}) + \frac 1{2 |E(Q)|}  \int (w^+_x(t))^2
 \right|\leq  K c^{\frac 54}\int (w^+(t))^2 +  K c^4.
\label{DOIS}
\end{align}
Let $a= (c^{2q+1} - (c_2^+)^{2q+1})/(c^{2q+1} - c\,(c_2^+)^{2q})$,
then $\frac 12 \frac {2q+1}{2q} \leq a \leq \frac 32 \frac {2q+1}{2q} $.
Multiplying \eqref{UM} by $c \, a$ and summing  \eqref{DOIS}, we obtain, for $c$ small enough,
$$
 K (c_1^+-1)\geq (c_1^+)^{2q+1}-1\geq K  \int (w^+_x)^2+c (w^+)^2)(t)  - K c^4
\geq K_0 c^{\frac {17}6}.
$$
Similarly, set $b=(1-(c_1^+)^{2q})/(1-(c_1^+)^{2q+1})$, then $\frac 12 \leq b\leq \frac 32$,
and multiplying  \eqref{DOIS} by $-b$ and summing  \eqref{UM}, we obtain, for $c$ small enough,
($q=\frac 1{12}$)
$$
K c^{\frac 16} \left(1-\frac {c_2^+}{c}\right)\geq
c^{2q} - (c_2^+)^{2q}\geq K \int ((w^+_x)^2+(w^+)^2)(t) \geq K c^{\frac {17}6}.$$
This completes the proof of Theorem \ref{NONEX4}.

\medskip

\noindent\emph{Proof of Remark 1.} \quad
The remark is based on the fact that for $p=4$,
$$
\int Q_c = c^{-\frac 16} \int Q.
$$
In the framework of the proof of Theorem \ref{NONEX4}, we consider $u(t)$ the
asymptotic $2$-soliton solution at $-\infty$ with speed parameters $1$ and $c$
($c$ small enough). Let us prove by contradiction that $u(t)$ is not an asymptotic $N$-soliton
solution at $+\infty$.

Assume   that $\|u(t)- \sum_{j=1}^N Q_{c_j^+}(.-x_j^+-c_j^+ t)\|_{H^1}\to 0$ as $t\to +\infty$, where $c_2^+>c_3^+>...>c_N^+$.
Using the methods of \cite{MMT}, \cite{Martel}  and the fact that $u(t)$ 
is an asymptotic $N$-soliton solution both at $\pm \infty$, we have, for some
$T_0>0$ large enough:
$$\forall t\geq T_0,\forall x\in \mathbb{R},\quad
|u(t,x)|\leq K  \sum_{j=1}^N Q_{c_j^+}^{\frac 14} (.-x_j^+-c_j^+ t),$$
 which proves that
 $u(t)\in L^1(\mathbb{R})$, and
in particular $\int u(t)=I_0$ is well-defined and constant in time.
Moreover,
$u(t)- \sum_{j=1}^N Q_{c_j^+}(.-x_j^+-c_j^+ t)\to 0$ as $t\to +\infty$ in $L^1(\mathbb{R})$,
from the $H^1$ convergence. A similar convergence in $L^1$ holds at $-\infty$.

On the one hand, at $-\infty$,
$  I_0 = \lim_{t\to -\infty} \int u(t) = 
\int Q_{c_1} + \int Q_{c_2} = (c_1^{-\frac 16} + c_2^{-\frac 16}) \int Q.
$ 
On the other hand, at $+\infty$,
$ 
I_0=
\sum_{j=1}^N (c_j^+)^{-\frac 16} \int Q.
$ 
Since by Theorem \ref{NONEX4},
$\|w^+(t)\|_{L^2}\leq K c_2^{\frac 7{12}}$, we have
$c_3^+\ll (c_2^+)^4$. Thus,
$I_0\gg (c_2^+)^{-\frac 23} \int Q$, which is a contradiction, for $c_2$ small.

\subsection{Existence of a $2$-soliton-like solution. Proof of Theorem \ref{EXIST4}}
We     consider first the case $c_1=1$ and $c_2=c$,
the general case   following from a scaling argument.
For any $c>0$ small enough,
we consider $u_c(t)$ the global solution of 
$$\partial_t u_c+\partial_x(\partial_x^2 u_c+u_c^4)=0\ ;\quad
u_c(0,x)=v_c(0,x),$$
where $v_c(t)$ is the approximate solution 
constructed in Proposition \ref{AVRIL4}, for $k_0$, $\ell_0$ large enough but fixed.
Recall also that $\Delta$ and $\Delta_c$ are defined in Proposition \ref{AVRIL4}.
By the parity property of $x\to v_c(0,x)$ and since   equation \eqref{gkdv} is invariant
under the transformation $x\to -x$, $t\to -t$, the solution $u_c(t)$ has the following
symmetry: 
\begin{equation}\label{SYMM}
u_c(t,x)=u_c(-t,-x).
\end{equation}
Thus, we shall only study $u_c(t)$ for $t\geq 0$.
We claim the following concerning $u_c(t)$.

\begin{proposition}
\label{DANSTHM12}
There exist $c_0>0$ such that for all $0<c<c_0$,
there exist $c_1^+(c)$, $c_2^+(c)>0$, and $x_1^+(c)$, $x_2^+(c)\in \mathbb{R}$ such that for
$$
w^+(t,x)=u_c(t,x)-Q_{c_1^+(c)}(x-c_1^+(c) t -x_1^+(c))-Q_{c_2^+(c)}(x-c_2^+(c) t -x_2^+(c))
$$
\begin{enumerate}
\item Asymptotic behavior:
\begin{equation}\label{thm12I}
\lim_{t\rightarrow +\infty}\| w^+(t)\|
_{H^1(x> ct/10)}=0.
\end{equation}
\begin{equation}\label{thm12II}\begin{split}
 & |x_1^+(c) -\tfrac 12 \Delta |\leq K c^{\frac 38},
\quad |x_2^+(c)- \tfrac 12 \Delta_c|\leq K c^{\frac 1{12}},\\
& |c_1^+(c) - 1|\leq K c^{\frac {11}{12}}  ,\quad
 \left|\frac {c_2^+(c)}{c} - 1\right|\le K c^{\frac 13},
 \end{split}\end{equation}
\begin{equation}\label{thm12IV}
 \text{for $t$ large},\quad 
 \tfrac 1K c^{\frac {17}{12}} \leq \|w^+(t)\|_{H^1_c}\leq K\min_{\rho_1,\rho_2\in \mathbb{R}}\|w_{\rho_1,\rho_2}(t)\|_{H^1_c} \leq
  K^2 c^{\frac {17}{12}},
\end{equation}
\begin{equation}\label{thm12IVV}
\text{where}\quad w_{\rho_1,\rho_2}(t,x)=u(t,x)-Q_{c_1^+(c)}
(x-\rho_1)-Q_{c_2^+(c)}(x-\rho_2).
\end{equation}
\item   $c\mapsto c^{+}_j(c)$ for $j=1,2$ are continuous.
\end{enumerate}
\end{proposition}

\noindent\emph{Proof of Theorem \ref{EXIST4} assuming Proposition \ref{DANSTHM12}.}
We claim that a rescaled version of $u_{\tilde c}(t)$
for some $\tilde c\sim c$ satisfies the conclusions of Theorem \ref{EXIST4}.

From Proposition \ref{DANSTHM12}, the function $h(c)=\frac {c_2 ^+(c)}{c_1^+(c)}$
is continuous on $(0,c_0]$, moreover
$\frac 12 c\leq h(c)\leq \frac 32 c$. It follows that $h((0,c_0])$
is an interval containing $(0,\frac 12 c_0]$. Thus, for any
$c\in (0,\frac 12 c_0]$,  there exists $\tilde c$ such that 
\begin{equation}\label{ettt2}
\frac 12 c\leq \tilde c\leq 2c,\quad 
h(\tilde c)=c.
\end{equation}
Let 
\begin{equation}\label{ettt}
 \varphi_{1,c}(t,x)=
 \varphi(t,x)=  {c^{-\frac 13}_1(\tilde c)}
u_{\tilde c}\Big( {c_1^{-\frac 32}(\tilde c)} t,
 {c_1^{-\frac 12}(\tilde c)}x\Big). 
\end{equation}
From Proposition \ref{DANSTHM12}, \eqref{ettt}, \eqref{ettt2},
\eqref{SYMM}, it follows that $\varphi$ satisfies \eqref{TH1A}.
Moreover,  \eqref{TH1E} follows from \eqref{thm12IV}.

Let $c_1>0$ and $c_2>0$ such that $c=\frac {c_2}{c_1} <\epsilon_0$ small.
Let 
$$
\varphi_{c_1,c_2}(t,x)=c_1^{\frac 13} \varphi_{1,c}\Big(c_1^{\frac 32} t , c_1^{\frac 12} x\Big),
\quad
\Delta_j=\Delta_j(c_1,c_2)=c_1^{-\frac 12} x_j^+(\tilde c),~j=1,2.
$$
Then $\varphi_{c_1,c_2}$ verifies the conclusion of Theorem \ref{EXIST4}.
Note in particular that \eqref{TH1B} follows from \eqref{thm12II} and \eqref{defDeltabis}, \eqref{AV2}.

\bigskip

\noindent\emph{Proof of Proposition \ref{DANSTHM12}.}
In steps 1 and 2 of this proof, we omit the $c$ dependency.

\smallskip

\noindent\emph{Step 1.}  Control of the modulation of $u(t)$ 
for $t\geq T_c$.

Applying Proposition~\ref{INTERACT4} for $t\in [0,T_c]$,
with $\theta=n_0-\frac 12-\frac 1{100}$,
we obtain, for some $\rho(t)$,
\begin{equation}\label{atc}
\forall t\in [0,T_c],\quad \|u(t)-v(t,.-\rho(t))\|_{H^1} \leq K c^{\theta},
\end{equation}
where $|\rho'(t)-1|\leq K c^{\theta}$, $ \rho(0)=0$
and so   $|\rho(T_c)-T_c|\leq K c^{\theta-
\frac 12 -\frac 1{100}}$ by $T_c= c^{-\frac 12-\frac 1{100}}$.

By \eqref{AV4} and \eqref{atc}, 
and then by
$\|(Q_c^2)'\|_{H^1_c} = K c^{\frac {17}{12}}$ and  
\eqref{3-16}-\eqref{3-15},  we have, for $\theta\geq 2$,
\begin{equation}\label{ennrr}
\|u(T_c)- Q(.-a) - Q_c(.-b)\|_{H^1(x> T_c/4)}\leq K c^{\theta}
\leq K c^{2},
\end{equation}
\begin{equation}\label{atcdeux}
\sqrt{c}\|u(T_c)- Q(.-a) - Q_c(.-b)\|_{H^1} \leq
\|u(T_c)- Q(.-a) - Q_c(.-b)\|_{H^1_c} \leq K c^{\frac {17}{12}},
\end{equation}
for $a=\frac 12 \Delta +\rho(T_c)$, $b=(1-c)T_c + \frac 12 \Delta_c + \rho(T_c)$,
so that $a-b\geq \frac 12 T_c$.

Therefore, from Claim \ref{pourstab} and  Proposition \ref{ASYMPTOTIC},
we have the decomposition 
of $u(t)$ in terms  of $\eta(t)$, $c_j(t)$, $\rho_j(t)$ ($j=1,2$)
defined for all $t\geq T_c$. 
 
\begin{lemma}\label{sureta}
 For all $t\geq T_c$,
$ 
\tfrac 1K c^{\frac{17}{12}} \leq \|\eta(t)\|_{H^1_c}\leq K c^{\frac{17}{12}} .
$ 
\end{lemma}

\noindent\emph{Proof of Lemma \ref{sureta}.}
 \emph{(i) Upper bounds by stability properties.}
 We use Claim \ref{pourstab}, which is a refinement of  Proposition \ref{ASYMPTOTIC}
 (see proof of Proposition 2 in \cite{MMas2}).
 Let $g(t)$ be  defined from $\eta(t)$  by \eqref{defg}. Remark
from 
\eqref{ennrr} and the proof of Claim \ref{pourstab} in \cite{MMas2}, that
$|c_1(t)-1|+|a_- -\rho_1(t)|\leq  K c^{2}$
 and 
\begin{equation}\label{dddnc}
\|\eta(T_c)\|_{H^1(x>T_c/4)}\leq
 K c^{2} .
 \end{equation}
 Similarly, we obtain $\|\eta(T_c)\|_{H^1_c(x<T_c/4)}\leq K c^{\frac {17}{12}}$
 from \eqref{atcdeux}.
Thus, 
  $$
 \sqrt{g(T_c)}\leq K \|\eta(T_c)\|_{H^1_c(x<T_c/4)}
 +\|\eta(T_c)\|_{H^1(x>T_c/4)}\leq
 K c^{\frac {17}{12}}.$$
By Claim \ref{pourstab}, for all $t\geq T_c$,
$\|\eta(t)\|_{H^1_c}\leq 
\sqrt{g(t)}\leq K (\sqrt{g(T_c)}+\exp(-c^{-\frac 1{400}}))\leq K c^{\frac {17}{12}}$.

\smallskip

 \emph{(ii) Lower bounds by backwards stability.} See  proof of \eqref{lower}
 (Theorem~\ref{NONEX4}).
 
\medskip

\noindent\emph{Step 2.}  Proof of asymptotic stability.


From properties of $v$, we claim the following:
 \begin{equation}\label{MOMENT}
\int_{x>\frac {11}{12}|\ln c|} x^2 u^2(T_c,x+T_c+\tfrac 12 \Delta) dx \leq K c^{\frac 54},
\end{equation}
\begin{equation}\label{23nov}
|\rho_1(T_c)-T_c-\tfrac \Delta 2| \leq  K c^{\frac {11}{12}},\quad 
 | \rho_2(T_c)-c T_c-\tfrac {\Delta_c} 2|\leq K c^{\frac 13}.
\end{equation}

See Appendix D.3 for the proof of  \eqref{MOMENT} and \eqref{23nov}.
Note that the proof of \eqref{MOMENT} is based on
monotonicity arguments on $z(t)=u(t)-v(t,.-\rho(t))$ as defined in \eqref{defz} in the proof of Proposition \ref{INTERACT4}.

From \eqref{atcdeux}--\eqref{MOMENT},
we    apply Proposition \ref{ASYMPTOTIC}  
to $u(.+T_c)$ with $\alpha=Kc^{\frac 13}$:
there exist 
$c_1^+$, $c_2^+>0$, $x_1^+$, $x_2^+\in \mathbb{R}$
 such that
\begin{equation}\label{5-40}
c_j(t)\to c_j^+, \quad 
\rho_j(t)-c_j^+ t \to x_j^+,
\quad \text{as $t\to +\infty$,  $j=1,2$,}
\end{equation}
\begin{equation}\label{neuf44}
\lim_{t\rightarrow +\infty}\|u(t)-w^+(t)\|
_{H^1(x> ct/10)}=0,
\end{equation}
where 
$$
w^+(t,x)=Q_{c_1^+}(x-c_1^+ t -x_1^+)-Q_{c_2^+}(x-c_2^+ t -x_2^+),
$$
\begin{equation}\label{sept44}
|c_1^+-1|\leq   K  c^{\frac {11}{12}}  ,\quad
 \left|\frac {c_2^+}{c} - 1\right|\le K c^{\frac 13}.
\end{equation}
 \begin{equation}\label{sept45}
|x_1^++c_1^+ T_c -\rho_1(T_c)|\leq K c^{\frac 58},
\quad |x_2^++c_2^+ T_c -\rho_2(T_c)|\leq K c^{\frac 1{12}}.
\end{equation}
From \eqref{23nov} and \eqref{sept45},
we finish the computation of $x_j^+$.
For $x_1^+$,  inserting  \eqref{23nov} in \eqref{sept45}, we obtain:
$|x_1^+ - (1-c_1^+) T_c - \tfrac 12 \Delta|\leq K c^{\frac 58}$.
Since $|1-c_1^+| T_c\leq K c^{\frac {11}{12}} T_c \leq K c^{\frac 38}$,
we conclude $|x_1^+  - \tfrac 12 \Delta|\leq K c^{\frac 38}$.
Similarly for $x_2^+$, we obtain from \eqref{23nov} and  \eqref{sept45},
$|x_2^+  -  \tfrac 12 \Delta_c|\leq K c^{\frac 1{12}}$.

\medskip
From \eqref{5-40}, $\|\eta(t)-w^+(t)\|_{H^1_c}\to 0$
as $t\to +\infty$ and thus, from Lemma \ref{sureta}, we obtain
$\tfrac 1 K  c^{\frac {17}{12}} \leq \|w^+(t)\|_{H^1_c}\leq
K  c^{\frac {17}{12}}$ for $t$ large.
From \eqref{neuf44}, 
$
 \|w^+(t)\|_{H^1_c}\leq 
\min_{\rho_1,\rho_2} \|w_{\rho_1,\rho_2}(t)\|_{H^1_c}+o(1)$
for $t$ large 
where $w_{\rho_1,\rho_2}(t)$ is defined in \eqref{thm12IVV}
 and thus \eqref{thm12IV} follows.
 This concludes the proof of the first part of Proposition \ref{DANSTHM12}.

\medskip

\noindent\emph{Step 3.} Continuity of $c_1^+(c)$ and $c_2^+(c)$.

Now, we prove that the maps
$c\mapsto c_1^+(c)$ is continuous.
Let us denote by
$\eta_c(t)$,  $c_{c,j}(t)$, $c_j^+(c)$, the parameters  in the decomposition of   $u_c(t)$.
We claim:
\begin{claim}\label{lllff}
For all $t\geq T_c$,
\begin{equation}\label{unifc}\begin{split}
|c_1^+(c)-c_{c,1}(t)| &\leq 
  K_0  \int (\eta_{c,x}^2+\eta_c^2)(t,x) \psi(x - \rho_1(t) +\tfrac {t }4) dx
+  K_0 e^{-\frac 1 {64}Ê\sqrt{c} t }.
\end{split}\end{equation}
\end{claim}

Assuming this claim, let us complete the proof of continuity of $c_1^+(c)$.
Let $0<\bar c<c_0$ and let $\varepsilon>0$.
 Since 
$\|\eta_{\bar c}(t)\|_{H^1(x>ct/10)} \to 0$ as $t\to +\infty$, there exits
$T_\varepsilon>0$ such that 
 \begin{equation*}\begin{split}
     K_0 \int (\eta_{\bar c,x}^2+\eta_{\bar c}^2)(T_\varepsilon,x) \psi(x- \rho_1(T_\varepsilon)
     +\tfrac {T_\varepsilon}4 ) dx
+   K_0  e^{-\frac 1 {64}Ê\sqrt{c} T_\varepsilon}\leq \varepsilon.
\end{split}\end{equation*}
We fix $T_\varepsilon>0$ to such value. Then, by continuous dependence in $H^1$ of $u_c(t)$ solution
of \eqref{gkdv} upon the initial data (see \cite{KPV}), and the fact that
$u_c(0)=v_c(0)$ is continuous upon the parameter $c$, there exists $\delta>0$ such that
if $|c-\bar c|\leq \delta$, then
 \begin{equation*}\begin{split}
     K_0  \int (\eta_{c,x}^2+\eta_c^2)(T_\varepsilon,x) \psi(x- \rho_1(T_\varepsilon)
     +\tfrac {T_\varepsilon}4 ) dx
+   K_0  e^{-\frac 1 {64}Ê\sqrt{c} T_\varepsilon}\leq 2 \varepsilon,
\end{split}\end{equation*}
\begin{equation*}
|c_{\bar c,1}(T_\varepsilon)-c_{c,1}(T_{\varepsilon})|  \leq  \varepsilon.
\end{equation*}
From Claim \ref{lllff},  applied to $\eta_c$, $\eta_{\bar c}$, we have 
$
|c_{1}^+(c)-c_{c,1}(T_\varepsilon)| \leq 2\varepsilon
$ and $|c_{ 1}^+(\bar c)-c_{\bar c,1}(T_\varepsilon)|\leq \varepsilon$. 
Therefore, $
|c_{1}^+(\bar c)-c_{1}^+(c)| \leq 4\varepsilon.
$
Thus, $c\mapsto c_{1}^+(c)$ is continuous. We argue similarly  for $c\mapsto c_{2}^+(c)$ using a claim similar to Claim \ref{lllff}
on $|c_2^{+}(c)-c_{c,2}(t)|$ (related to 
$\mathcal{M}_2(t)$ and $\mathcal{E}_2(t)$) and the previous 
result on $c_1^+(c)$.
This concludes the proofs of Proposition \ref{DANSTHM12} and
of Theorem~\ref{EXIST4}.

\medskip

\noindent\emph{Proof of Claim \ref{lllff}.} The proof follows
closely some arguments in \cite{MMas2}.
For $T_c\leq t_0\leq t$, let $\mathcal{M}_1(t)$ and $\mathcal{E}_1(t)$ be defined 
in Secion 4.3, with $x_0=\frac {t_0}4$.
From  the conclusions of  Claim \ref{NEW} integrated on $[t_0,t]$, 
we obtain
\begin{equation*}\begin{split}
& (c_1^{2q}(t)-c_1^{2q}(t_0))\int Q^2 \leq   
( \mathcal{M}_1(t_0)-\mathcal{M}_1(t))
 + K e^{-\frac 1{64} \sqrt{c} t_0},
\\ & \left(\frac {2q}{2q+1} (c_1^{2q+1}(t)-c_1^{2q+1}(t_0)) -\frac 1{100} (c_1^{2q}(t)-c_1^{2q}(t_0))\right)  \int Q^2 
\\ & \quad \geq 2\mathcal{E}_1(t)- 2 \mathcal{E}_1(t_0) +\frac 1{100} ( \mathcal{M}_1(t)-\mathcal{M}_1(t_0)) - K  e^{-\frac 1{64} \sqrt{c} t_0}.
\end{split}\end{equation*}
Note in particular that 
$
\int_{t_0}^t e^{-\frac 1{16} (t-t_0+x_0)} g_1(t) dt
\leq K e^{-\frac 1 {16} x_0}Ê\leq K e^{-\frac 1 {64} t_0}.
$
Letting $t\to +\infty$, by the asymptotic stability, this gives
\begin{equation*}\begin{split}
 & ((c_1^+)^{2q}-c_1^{2q}(t_0))\int Q^2 \leq   
 \mathcal{M}_1(t_0) 
 + K e^{-\frac 1{64} \sqrt{c} t_0},
\\ & \left(\frac {2q}{2q+1} ((c_1^+)^{2q+1}-c_1^{2q+1}(t_0)) -\frac 1{100} ((c_1^+)^{2q} -c_1^{2q}(t_0))\right)  \int Q^2 
\\ & \quad \geq - 2 \mathcal{E}_1(t_0) -\frac 1{100}   \mathcal{M}_1(t_0) 
- K e^{-\frac 1{64} \sqrt{c} t_0}.
\end{split}\end{equation*}
Thus, we obtain
\begin{equation*}\begin{split}
|c_1^+-c_1(t_0)| &\leq 
K \int (\eta_x^2+\eta^2)(t_0,x) \psi(x - \rho_1(t_0) +\tfrac {t_0}4 ) dx
- K e^{-\frac 1 {64}Ê\sqrt{c} t_0}.
\end{split}\end{equation*}

\subsection{Stability of the $2$-soliton structure.  Proof of Theorem \ref{STAB4}}

Without loss of generality, 
we prove Theorem \ref{STAB4} in the case $c_1=1$ and $c_2=c$.
We assume
\begin{equation*}\label{hypstab}
\|u(0)-\varphi(0)\|_{H^1}\leq K c^{\delta+\frac {7}{12}},
\end{equation*}
for $\delta>0$, where $\varphi$ is the solution constructed  in Theorem \ref{EXIST4}.
Let $\tilde c>0$ small satisfy $\frac {c_2^+(\tilde c)}{c_1^+(\tilde c)}=c$
and $\lambda=1/{c_1^+(\tilde c)}$. Then,
\begin{equation*}
\|\lambda^{\frac 13}u(0,\sqrt{\lambda} x)-\lambda^{\frac 13}\varphi(0,\sqrt{\lambda} x)\|_{H^1}
\leq K c^{\delta+\frac {7}{12}}.
\end{equation*}
By construction of $\varphi(t)$ in Theorem \ref{EXIST4}, $\lambda^{\frac 13}\varphi(0,\sqrt{\lambda} x)
=v(0)$ where $v$ is the approximate solution introduced in Proposition \ref{AVRIL4} corresponding to $\tilde c$
for $k_0$, $\ell_0$
large enough.
Since the solution of \eqref{gkdv} corresponding to $\lambda^{\frac 13}u(0,\sqrt{\lambda} x)$
is $\lambda^{\frac 13} u(\lambda t, \sqrt{\lambda} x)$, it is enough to prove the Theorem 
in the case 
\begin{equation}\label{mlml}
\|u(0)-v(0)\|_{H^1} \leq K c^{\delta + \frac 7{12}}.
\end{equation}
By invariance of \eqref{gkdv} by the transformation $x\to -x$, $t\to -t$, 
it is enough to prove the result for $t\geq 0$.

\medskip

\emph{(i) Estimates on $[0,T_c]$.} By \eqref{mlml} and
Proposition \ref{INTERACT4}, we obtain, for all $t\in [-T_c,T_c]$, for some $\rho(t)$,
\begin{equation*}
\|u(t)-v(t,x-\rho(t))\|_{H^1}  \leq K  c^{\delta+ \frac {7}{12}}.
\end{equation*}
From Proposition \ref{AVRIL4}, we deduce, for some $a$, $b$, with
$a-b\geq \frac 12 T_c$,
\begin{equation}\label{oioi}
\|u(T_c)-Q(.-a)-Q_c(.-b)\|_{H^1}\leq K (c^{\delta+ \frac {7}{12}}+ c^{\frac {11}{12}}).
\end{equation}

\medskip

\emph{(ii) Estimates on $[T_c,+\infty)$.} By \eqref{oioi} and Propositions \ref{ASYMPTOTIC} and
\ref{ASYMPTOTIC},  for all $t\in [T_c,+\infty)$, there exist
$\rho_1(t)$, $\rho_2(t)$
and $c_1^+$, $c_2^+$,  such that (recall that for $p=4$, $q+\frac 12=\frac 7{12}$)
\begin{equation*}\begin{split}
& \|u(t)-Q_{c_1^+}(.-\rho_1(t))-Q_{c_2^+}(.-\rho_2(t))\|_{H^1} \leq K (c^{\delta+ \frac {1}{12}}
+c^{\frac 5{12}}),\\
& \left|{c_1^+}-1 \right| \leq K (c^{\delta+\frac {7}{12}}+c^{\frac {11}{12}})Ê , 
\quad  \left|\frac {c_2^+}c -1 \right| \leq K (c^{\delta}+ c^{\frac 13})Ê .
\end{split}\end{equation*}

\appendix

\section{Appendix -- Proof of Proposition \ref{SYSTEME}}

To prove Proposition \ref{SYSTEME},
we decompose each of the terms $\mathbf{I},$ $\mathbf{II},$ $\mathbf{III}$ and ${\mathbf{IV}}$
obtained in \eqref{decS} in series of $c^\ell Q_c^k$, $c^\ell (Q_c^k)'$.
In this decomposition (for future use in solving the systems $(\Omega_{k,\ell})$), we will separate terms depending on $(k,\ell)$ and terms depending on $(k',\ell')\prec (k,\ell)$. 

\begin{claim}\label{surDD}
       \begin{enumerate}
        \item 
        For  $r>0$,
    $        Q_c^r(y_c) \beta(y_c) = \sum_{\substack{1+r\leq k\leq k_{0}+r \\ 0\leq \ell \leq    \ell_{0}}} 
        c^\ell Q_c^k(y_c) a_{k-r,\ell}.
   $        \item Decomposition of $\beta''$, $\beta^2$, $\beta'\beta$ and $\beta^3$. There exist
   ${a^{1*}_{k,\ell}}$, ${a^{2*}_{k,\ell}}$, ${a^{3*}_{k,\ell}}$ and ${a^{4*}_{k,\ell}}$ depend on   $(a_{k',\ell'})$ for $(k',\ell')\prec (k,\ell)$ such that
    \begin{align*}
        & \beta''(y_c) = \sum_{\substack{1\leq k\leq    k_{0}+p-1 \\ 0\leq \ell \leq    \ell_{0}+1 }}
        c^\ell Q_c^{k}(y_c) {a^{1*}_{k,\ell}},\quad
        \beta^2(y_c) = \sum_{\substack{2\leq k\leq  2 k_{0} \\ 0\leq \ell \leq 2    \ell_{0} }}c^{\ell} Q_c^{k}(y_c)    {a^{2*}_{k,\ell}},\\
        &\beta'(y_c)\beta(y_c) = \sum_{\substack{2\leq k\leq    2 k_{0} \\ 0\leq \ell \leq 2    \ell_{0} }} c^{\ell} (Q_c^{k})'(y_c) {a^{3*}_{k,\ell}}, \quad
        \beta^3(y_c)= \sum_{\substack{3\leq k\leq    3 k_{0} \\ 0\leq \ell \leq 3    \ell_{0} }}c^{\ell} Q_c^{k}(y_c)    {a^{4*}_{k,\ell}}.    \end{align*}
    \end{enumerate}
\end{claim}

\noindent\emph{Proof of Claim \ref{surDD}.}
The first formula follows immediately from the decomposition of $\beta(y_c)$:
\begin{equation}\label{decbeta}
    \beta(y_c)=\sum_{(k,\ell)\in \Sigma_0} a_{k,\ell} \, c^\ell Q_{c}^k(y_c).
\end{equation}
- Decomposition of $\beta''$.
Using Lemma \ref{surQc},
\begin{align*}&
 \beta''(y_c) = \sum_{(k,\ell)\in \Sigma_0} c^\ell (Q_c^{k})''(y_c) a_{k,\ell} 
 = \sum_{\substack{1\leq k\leq k_{0}
\\ 0\leq \ell \leq \ell_{0}}} c^\ell \left(ck^2 Q_c^k(y_c) -    \frac {k(2 k {+}p{-}1)}{p+1} Q_c^{k+p-1}(y_c)\right) a_{k,\ell}
\\
&= \sum_{\substack{1\leq k\leq   k_{0} \\ 1\leq \ell \leq       \ell_{0} +1  }}
c^\ell Q_c^{k}(y_c) k^2 a_{k,\ell-1} +
\sum_{\substack{p\leq k\leq  k_{0}+p-1 \\ 0\leq \ell \leq    \ell_{0} }}
c^\ell Q_c^{k}(y_c) \left(-\frac {(k{-}p{+}1)(2k{-}p{+}1)}{p+1} a_{k-p+1,\ell}\right).
\end{align*}
Thus,
$    \beta''(y_c)
= \sum_{\substack{1\leq k\leq    k_{0}+p-1 \\ 0\leq \ell \leq    \ell_{0}+1 }}
c^\ell Q_c^{k}(y_c) {a^{1*}_{k,\ell}},
$
where ($\mathbf{1}$ denoting the characteristic function)
\begin{align}
    {a^{1*}_{k,\ell}}       & = k^2 a_{k,\ell-1}\mathbf{1}_{\left\{\substack{1\leq k\leq k_0 \\1\leq \ell\leq \ell_0+1}\right\}}
                              +
                            \frac {(k{-}p{+}1)(2k{-}p{+}1)}{p+1} a_{k-p+1,\ell} \mathbf{1}_{\left\{\substack{p\leq k\leq k_0+p-1\\ 0\leq \ell\leq \ell_0}\right\}}
                                                          .\label{addot}
\end{align}
Thus, the coefficient ${a^{1*}_{k,\ell}}$ 
depend on some  $(a_{k',\ell'})$ only for $k',$ $\ell'$ such that $(k',\ell')\prec (k,\ell)$
(more precisely, either $k'\leq k$ and $\ell'\leq \ell-1$ or $k'\leq k-p+1$ and $\ell'\leq \ell$).

- Decomposition of $\beta^2$. By \eqref{decbeta},
\begin{equation*}\beta^2(y_c)   = \sum_{\substack{1\leq k_1,k_2\leq    k_{0} \\ 0\leq \ell_1,\ell_2\leq    \ell_{0} }} 
c^{\ell_1+\ell_2} Q_c^{k_1+k_2}(y_c) a_{k_1,\ell_1} a_{k_2,\ell_2} = \sum_{\substack{2\leq k\leq    2 k_{0} \\ 0\leq \ell \leq 2    \ell_{0} }}c^{\ell} Q_c^{k}(y_c)    {a^{2*}_{k,\ell}}, 
\end{equation*}
where
\begin{equation} \label{astar} 
    {a^{2*}_{k,\ell}}=
    \sum_{\substack{\max(k-k_0,1)\leq k_1 \leq \min(k-1,k_0)\\ \max(\ell-\ell_0,0)\leq \ell_1 \leq  \min(\ell,\ell_{0})}}
\hskip -1cm a_{k_1,\ell_1} a_{k-k_1,\ell-\ell_1}     .
\end{equation}
Note that       the expression of ${a^{2*}_{k,\ell}}$ above involves $a_{k_1,\ell_1}$ with $k_1\leq k-1$ and
$a_{k-k_1,\ell-\ell_1}$ with $k-k_1\leq k-1$ since $k_1\geq 1$.
Thus it is checked that $a_{k,\ell}$ does not appear in the expression of ${a^{2*}_{k,\ell}}$.

- Decomposition of $\beta'(y_c)\beta(y_c)$.
\begin{equation*} 
 \beta'(y_c)\beta(y_c) = \sum_{\substack{1\leq k_1,k_2\leq    k_{0} \\ 0\leq \ell_1,\ell_2\leq    \ell_{0} }} 
c^{\ell_1+\ell_2} (Q_c^{k_1+k_2})'(y_c) \frac {k_1}{k_1{+}k_2} a_{k_1,\ell_1} a_{k_2,\ell_2}
= \sum_{\substack{2\leq k\leq    2 k_{0} \\ 0\leq \ell \leq 2    \ell_{0} }} c^{\ell} (Q_c^{k})'(y_c) {a^{3*}_{k,\ell}}, 
\end{equation*}
where
\begin{equation}\label{astarbar} {a^{3*}_{k,\ell}}=
\sum_{\substack{\max(k-k_0,1)\leq k_1 \leq \min(k-1,k_0)\\ \max(\ell-\ell_0,0)\leq \ell_1 \leq  \min(\ell,\ell_{0})}}
\frac {k_1}ka_{k_1,\ell_1} a_{k-k_1,\ell-\ell_1}    .
\end{equation}

- Decomposition of $\beta^3(y_c)$.
By $\beta^3(y_c)=\beta(y_c)\beta^2(y_c)$ and the decomposition of $\beta^2$,
\begin{align}\beta^3(y_c)
& =\sum_{\substack{1\leq k_1\leq    k_{0} \\ 0\leq \ell_1 \leq  \ell_{0} }}  c^{\ell_1} Q_c^{k_1}(y_c) a_{k_1,\ell_1}       \times \sum_{\substack{2\leq k_2\leq    2 k_{0} \\ 0\leq \ell_2\leq  2 \ell_{0} }} c^{\ell_2} Q_c^{k_2}(y_c) a_{k_2,\ell_2}^*\nonumber
= \sum_{\substack{3\leq k\leq    3 k_{0} \\ 0\leq \ell \leq 3    \ell_{0} }}c^{\ell} Q_c^{k}(y_c)    {a^{4*}_{k,\ell}},
\end{align}
where
\begin{equation} \label{astarstar}
     {a^{4*}_{k,\ell}}=
\sum_{\substack{\max(k-2 k_0,1)\leq k_1 \leq \min(k-2,k_0)\\ \max(\ell-2\ell_0,0)\leq\ell_1 \leq    \min(\ell,\ell_{0})}}
\hskip -1cm a_{k_1,\ell_1} a^{2*}_{k-k_1,\ell-\ell_1}       .
\end{equation}

\subsection{Decomposition of $\mathbf{I}=\partial_t R + \partial_x(\partial_x^2 R -R + R^p)$}

\begin{lemma}[Equation of $R(t)$]\label{PROPI}
    \begin{align}
    \mathbf{I} & = 
    \sum_{(k,\ell)\in \Sigma_0} c^\ell\left( Q_c^k(y_c)  a_{k,\ell} (-3 Q+2 Q^p)'(y)    + (Q_c^k)'(y_c)  a_{k,\ell} (-3 Q'')(y)\right) \label{IF} \\ &
    \quad + \sum_{\substack{1\leq k\leq \max(3 k_{0},k_0+p-1) \\ 0\leq \ell \leq \max(3 \ell_{0},\ell_0+1)}}
    c^\ell\left( Q_c^k(y_c)  F_{k,\ell}^{\mathbf{I}}(y)     +    (Q_c^k)'(y_c) G_{k,\ell}^{\mathbf{I}}(y)\right),
    \label{IG}
    \end{align}
    where $F_{k,\ell}^{\mathbf{I}}$ and $G_{k,\ell}^{\mathbf{I}}$ are functions defined on $\mathbb{R}$ satisfying :\begin{itemize}
   \item[{\rm (i)}]  $F_{k,\ell}^{\mathbf{I}}$, $G_{k,\ell}^{\mathbf{I}}\in \mathcal{Y}$.
       \item[{\rm (ii)}]  $F_{k,\ell}^{\mathbf{I}}$ and $G_{k,\ell}^{\mathbf{I}}$ 
        depend only on $(a_{k',\ell'})$ for $k',$ $\ell'$ such that $(k',\ell') \prec (k,\ell)$; 
                \item[{\rm (iii)}]  $F_{k,\ell}^{\mathbf{I}}$ is odd and $G_{k,\ell}^{\mathbf{I}}$ is even. 
        \end{itemize}  
   Moreover,
    $ 
        F_{1,0}^\mathbf{I}=0, \text{ and for all $\ell\geq 0$, }   G_{1,\ell}^\mathbf{I}=0,
    $ and
    \begin{itemize}
        \item If $p=2$ then
             $
                F_{2,0}^\mathbf{I}=a_{1,0} Q' + 3 a_{1,0}^2 Q^{(3)},\quad G_{2,0}^\mathbf{I}= \frac 32 a_{1,0}^2 Q''.
           $
        \item If   $p=4$ then  
            $
                F_{2,0}^\mathbf{I}= 3 a_{1,0}^2 Q^{(3)},\quad G_{2,0}^\mathbf{I}= \frac 32 a_{1,0}^2 Q''.
            $
        \end{itemize}
\end{lemma}
\begin{claim}\label{forH} 
    Let $h(t,x)=g(y)=g(x-\alpha(y_c))$, where $g$ is a $C^3$ function. Then,
\begin{align*}
 \partial_t h(t,x) &=-(1-c) \beta(y_c) g'(y),\quad 
 \partial_x h(t,x)=(1-\beta(y_c)) g'(y),\\
 \partial_x^2 h(t,x)&=(1-2 \beta(y_c) + \beta^2(y_c)) g''(y) - \beta'(y_c) g'(y),\\
 \partial_x^3 h(t,x)&=(1-3\beta(y_c) + 3 \beta^2(y_c) - \beta^3(y_c)) g^{(3)}(y) \\ 
                    &\quad + (-3 \beta'(y_c) + 3 \beta'(y_c) \beta(y_c)) g''(y) - \beta''(y_c) g'(y).
\end{align*}
\end{claim}

\noindent\emph{Proof of Claim \ref{forH}.} Recall    that       $y_c=x+(1-c) t$ and $\alpha'(s)=\beta(s)$.  Thus,
\begin{equation*}
    \partial_t h(t,x)=-\frac {\partial y_c}{\partial t} \alpha'(y_c) g'(y)=-(1-c) \beta(y_c) g'(y),
\end{equation*}
\begin{equation*}
    \partial_x h(t,x)=\left(1-\frac {\partial y_c}{\partial x} \alpha'(y_c)\right) g'(y)=(1- \beta(y_c)) g'(y).
\end{equation*}
Next, $\partial_x^2 h(t,x)=(1-\beta(y_c))^2 g''(y)-\beta'(y_c) g'(y), $ and so
\begin{align*}
\partial_x^3 h(t,x)&=-2 \beta'(y_c) (1-\beta(y_c)) g''(y) + (1-\beta(y_c))^3 g^{(3)}(y)\\ & \quad-\beta''(y_c) g'(y)
-\beta'(y_c)(1-\beta(y_c))g''(y)\\
&\quad =(1-\beta(y_c))^3 g^{(3)}(y)-3 \beta'(y_c)(1-\beta(y_c))g''(y) - \beta''(y_c)g'(y). 
\end{align*}

\noindent\emph{Proof of Lemma \ref{PROPI}.} \quad 
- Expression of $\mathbf{I}$.
We claim
 \begin{align*}
\mathbf{I}&= \beta(y_c) (-3 Q +2 Q^p)'(y) + \beta'(y_c) (-3 Q'')(y) + c \beta(y_c) Q'(y) + \beta''(y_c) (-Q')(y) 
\\&\quad    + \beta^2(y_c) (3 Q^{(3)})(y) + \beta'(y_c) \beta(y_c) (3 Q'')(y) +\beta^3(y_c) (-Q^{(3)})(y)
\\&= \mathbf{I}_1+\mathbf{I}_2+\mathbf{I}_3+\mathbf{I}_4+\mathbf{I}_5+\mathbf{I}_6+\mathbf{I}_7.
\end{align*}
Indeed, since $R(t,x)=Q(y)$, by Claim \ref{forH}, we have
\begin{align*}
 \partial_t R(t,x) &=-(1-c) \beta(y_c) Q'(y),\\
 \partial_x^3 R(t,x)&=(1-3\beta(y_c) + 3 \beta^2(y_c) - \beta^3(y_c)) Q^{(3)}(y) \\ 
                    &\quad + (-3 \beta'(y_c) + 3 \beta'(y_c) \beta(y_c)) Q''(y) - \beta''(y_c) Q'(y).\\
 -\partial_x R(t,x)&=-(1-\beta(y_c)) Q'(y)  ,\quad
  \partial_x(R^p)=(1-\beta(y_c)) (Q^p)'(y).
\end{align*}
Thus, by arranging terms by increasing order of derivatives and powers of $\beta(y_c)$, we get
\begin{align*}
\mathbf{I}&=\partial_t R + \partial_x(\partial_x^2 R -R + R^p)\\
&= (Q''-Q+Q^p)'(y) + \beta(y_c) (-3Q''-Q^p+cQ)'(y)+ \beta'(y_c) (-3 Q'')(y)\\
&\quad +\beta''(y_c)(-Q')(y)+\beta^2(y_c) (3 Q^{(3)})(y) + \beta'(y_c) \beta(y_c) (3 Q'')(y)+\beta^3(y_c) (-Q^{(3)})(y).
\end{align*}
By the equation of $Q$, i.e. $Q''-Q+Q^p=0$, the claim is proved.

\smallskip

- Decomposition of $\mathbf{I}_1$ and $\mathbf{I}_2$. These two terms give \eqref{IF}.
\begin{equation*} 
    \mathbf{I}_1 = \beta(y_c) (-3 Q + 2Q^p)'(y) = \sum_{(k,\ell)\in \Sigma_0}
    c^\ell Q_c^{k}(y_c)  a_{k,\ell} (-3 Q + 2Q^p)'(y), 
\end{equation*}
\begin{equation*}
    \mathbf{I}_2 = \beta'(y_c) (-3 Q'')(y) = \sum_{(k,\ell)\in \Sigma_0}
     c^\ell (Q_c^{k})'(y_c) a_{k,\ell} (-3 Q'')(y).
\end{equation*}
- Decomposition of $\mathbf{I}_3= c \beta(y_c) Q'(y)$.
\begin{equation*}
    \mathbf{I}_3    = \sum_{(k,\ell)\in \Sigma_0} c^{\ell+1} Q_c^{k}(y_c) a_{k,\ell}Q'(y)
     = \sum_{\substack{1\leq k\leq k_{0} \\ 1\leq \ell \leq \ell_{0}+1}}    c^{\ell} Q_c^{k}(y_c)    F_{k,\ell}^{\mathbf{I}_3}(y),  
    \quad \text{where}\quad F_{k,\ell}^{\mathbf{I}_3}=a_{k,\ell-1} Q'.
\end{equation*}  

\smallskip

- Decomposition of $\mathbf{I}_4$,  $\mathbf{I}_5$, $\mathbf{I}_6$ and  $\mathbf{I}_7$.
For these terms, we use Claim \ref{surDD}.
{\allowdisplaybreaks
\begin{align*}
    &\mathbf{I}_4=\sum_{\substack{1\leq k\leq    k_{0}+p-1 \\ 0\leq \ell \leq    \ell_{0}+1 }}
    c^{\ell} Q_c^{k}(y_c) F_{k,\ell}^{\mathbf{I}_4}(y)
        \quad\text{where}\quad
    F_{k,\ell}^{\mathbf{I}_4}(y) = {a^{1*}_{k,\ell}}(-Q'(y)).
\\&  \mathbf{I}_5= \sum_{\substack{2\leq k\leq  2 k_{0}\\ 0\leq \ell \leq 2  \ell_{0} }}
    c^{\ell} Q_c^{k}(y_c)    F_{k,\ell}^{\mathbf{I}_5}(y)
        \quad\text{where}\quad
    F_{k,\ell}^{\mathbf{I}_5}(y)= {a^{2*}_{k,l}} (3 Q^{(3)})(y).
\\&  \mathbf{I}_6= \sum_{\substack{2\leq k\leq  2 k_{0}\\ 0\leq \ell \leq 2  \ell_{0} }}c^{\ell} (Q_c^{k})'(y_c) 
    G_{k,\ell}^{\mathbf{I}_6}(y)
        \quad\text{where}\quad
    G_{k,\ell}^{\mathbf{I}_6}(y) = {a^{3*}_{k,\ell}} (3 Q'')(y).
\\&  \mathbf{I}_7= \sum_{\substack{3\leq k\leq  3 k_{0}\\ 0\leq \ell \leq 3  \ell_{0} }}c^{\ell} Q_c^{k}(y_c) 
    F_{k,\ell}^{\mathbf{I}_7}(y),
    \quad\text{where}\quad
    F_{k,\ell}^{\mathbf{I}_7}(y)= {a^{4*}_{k,\ell}} (- Q^{(3)})(y).
\end{align*}
}
We check that $F_{k,\ell}^{\mathbf{I}_3}$, $F_{k,\ell}^{\mathbf{I}_4}$, $F_{k,\ell}^{\mathbf{I}_5}$, $G_{k,\ell}^{\mathbf{I}_6}$ and $F_{k,\ell}^{\mathbf{I}_7}$  satisfy properties {\rm (i)}, {\rm (ii)} and {\rm (iii)}.
Set
$F_{k,\ell}^\mathbf{I}=F_{k,\ell}^{\mathbf{I}_3}+F_{k,\ell}^{\mathbf{I}_4}+F_{k,\ell}^{\mathbf{I}_5}+F_{k,\ell}^{\mathbf{I}_7}$ and
$G_{k,\ell}^\mathbf{I}=G_{k,\ell}^{\mathbf{I}_6}$, they satisfy {\rm (i)}, {\rm (ii)} and {\rm (iii)}.

\medskip

\noindent To finish the proof of Lemma \ref{PROPI}, we compute $F_{1,0}^\mathbf{I}$, $G_{1,\ell}^\mathbf{I}$, $F_{2,0}^\mathbf{I}$ and $F_{2,0}^\mathbf{I}$.

- $k=1$: We check that $F_{1,0}^{\mathbf{I}_3}=0$, $F_{1,0}^{\mathbf{I}_4}=a^{1*}_{1,0} (-Q')=0$, $F_{1,0}^{\mathbf{I}_5}=F_{1,0}^{\mathbf{I}_7}=0$, so that
$F_{1,0}^{\mathbf{I}}=0$. Moreover, for any $\ell\geq 0$, we have $G_{1,0}^{\mathbf{I}}=G_{1,0}^{\mathbf{I}_6}=0$.

- $k=2$: We check $F_{2,0}^{\mathbf{I}_3}=0$. The    term $F_{2,0}^{\mathbf{I}_4}=a^{1*}_{2,0}(-Q')$ depends on the value of $p$: from \eqref{addot},
if $p=2$ then $F_{2,0}^{\mathbf{I}_4}=a_{1,0} Q'$, and if $p=3$ or $4$, then $F_{2,0}^{\mathbf{I}_4}=0$.
By \eqref{astar}, we have $F_{2,0}^{\mathbf{I}_5}=3 {a^{2*}_{2,0}} Q^{(3)}= 3 a_{1,0}^2 Q^{(3)}$ and by \eqref{astarstar},
$F_{2,0}^{\mathbf{I}_7}=-a^{4*}_{2,0} Q^{(3)}=0$.
Thus, if $p=2$, we obtain $F_{2,0}^{\mathbf{I}}=a_{1,0} Q' + 3 a_{1,0}^2 Q^{(3)}$ and if $p=3$ or $4$, we obtain
$F_{2,0}^{\mathbf{I}}=  3 a_{1,0}^2 Q^{(3)}$.

Similarly, $G_{2,0}^{\mathbf{I}}=G_{2,0}^{\mathbf{I}_6}=a^{3*}_{2,0} (3 Q'') = \frac 32 a_{1,0}^2 Q''$.

\subsection{Decomposition of ${\mathbf{II}}=\partial_x((R+R_c)^p-R^p-R_c^p).$}

\begin{lemma}[Interaction term between $R$ and $R_c$]\label{PROPII}
    \begin{equation}\label{IIFG}
        \mathbf{II} = 
        \sum_{\substack{1\leq k\leq k_{0}+p-1 \\ 0\leq \ell \leq    \ell_{0}}} c^\ell \left(Q_c^k(y_c)  F_{k,\ell}^{\mathbf{II}}(y)
        + (Q_c^k)'(y_c)  G_{k,\ell}^{\mathbf{II}}(y)\right),
    \end{equation}  
    where    
    for any $k\geq 1$, $\ell\geq 0$, $F_{k,\ell}^{\mathbf{II}}$, $G_{k,\ell}^{\mathbf{II}}$ satisfy properties {\rm (i)}, {\rm (ii)} and {\rm (iii)} as in Lemma \ref{PROPI}.
    
    Moreover, 
    $ F_{1,0}=p (Q^{p-1})',$ $G_{1,0}=p Q^{p-1},  
       $ $F_{1,\ell}^{\mathbf{II}}=G_{1,\ell}^{\mathbf{II}}=0$, for any $\ell\geq 1$,
      \begin{itemize}
        \item If $p=2$ then
        \begin{equation*}
            \text{$F_{k,\ell}^{\mathbf{II}}  = - 2 a_{k-1,\ell} Q'$, for any $k\in \{2,k_0+1\}$, $\ell\in \{ 0,\ell_0 \} $.}
        \end{equation*}
        \item If $p=4$ then 
        \begin{align*}
            & \text{$F_{2,0}^{\mathbf{II}}=(-4 a_{1,0} Q^3+6 Q^2)',$ $G_{2,0}^{\mathbf{II}}=6 Q^2$,
            $G_{2,\ell}^{\mathbf{II}}=0$, for any $\ell\geq 1$,} \\
            & \text{$F_{3,0}^{\mathbf{II}}=(-4a_{2,0}Q^3-6 a_{1,0} Q^2+ 4 Q)',$ $G_{3,0}^{\mathbf{II}}=4 Q$,
            $G_{3,\ell}^{\mathbf{II}}=0$, for any $\ell\geq 1$.}
        \end{align*}
    \end{itemize}
\end{lemma}

\noindent\emph{Proof of Lemma \ref{PROPII}}  

$\bullet$ $p=2$.
Recall that $R(t,x)=Q(y)$ and $R_c(t,x)=Q_c(y_c)$. 
By Claim~\ref{forH} and  Claim \ref{surDD}, we have
\begin{align*}
    {\mathbf{II}} & = 2 \partial_x ( R\, R_c) = 2 (1-\beta(y_c)) Q'(y) Q_c(y_c) + 2Q(y) Q_c'(y_c)\\
    & = Q_c(y_c)    2 Q'(y) + Q_c'(y_c) 2 Q(y)+ Q_c(y_c) \beta(y_c) (-2 Q'(y)) \\
& = Q_c(y_c) 2Q'(y) + Q_c'(y_c) 2 Q(y) +\sum_{\substack{2\leq k\leq k_{0}+1 \\ 0\leq \ell \leq  \ell_{0}}} c^\ell Q_c^k(y_c) a_{k-1,\ell} (-2Q')(y).
\end{align*}

$\bullet$ $p=4$. As before,
{\allowdisplaybreaks
\begin{align*}
{\mathbf{II}} & = \partial_x (4 R^3 R_c + 6R^2 R_c^2+ 4 R R_c^3)\\
& = Q_c(y_c) (4Q^3)'(y)+Q_c'(y_c) (4Q^3)(y)+Q_c(y_c) \beta(y_c)(-4Q^3)'(y)\\
& +Q_c^2(y_c) (6Q^2)'(y)+(Q_c^2)'(y_c) (6Q^2)(y)+Q_c^2(y_c) \beta(y_c)(-6Q^2)'(y) \\
& +Q_c^3(y_c) (4Q)'(y)+(Q_c^3)'(y_c) (4Q)(y)+Q_c^3(y_c) \beta(y_c)(-4Q)'(y)\\
&= Q_c(y_c) (4Q^3)'(y)+Q_c'(y_c) (4Q^3)(y)\\
    & +Q_c^2(y_c) (-4 a_{1,0}Q^3+ 6Q^2)'(y)+(Q_c^2)'(y_c) (6Q^2)(y)
        + \sum_{1\leq \ell\leq \ell_0} c^\ell Q_c^2(y_c) a_{1,\ell} (-4 Q^3)'(y)    \\
    & +Q_c^3(y_c) (-4 a_{2,0}Q^3 - 6 a_{1,0}Q^2+        4Q)'(y)+(Q_c^3)'(y_c) (4Q)(y)\\
    & + \sum_{1\leq \ell\leq \ell_0} c^\ell Q_c^3(y_c)(- 4 a_{2,\ell} Q^3 -6     a_{1,\ell} Q^2)'(y)\\
    & + \sum_{\substack{4 \leq k\leq k_{0}+1 \\ 0\leq \ell \leq  \ell_{0}}} 
        c^\ell  Q_c^k(y_c) (- 4 a_{k-1,\ell} Q^3-6 a_{k-2,\ell} Q^2 -4 a_{k-3,\ell} Q)'(y) \\
    & +\sum_{0\leq \ell \leq \ell_0} c^\ell\left( Q_c^{k_0+2}(y_c) (-6a_{k_0,\ell}Q^2 -4 a_{k_0-1,\ell} Q)'(y)
        +Q_c^{k_0+3}(y_c)(-4 a_{k_0,\ell} Q)'(y)\right).
\end{align*}}

\subsection{Decomposition of $\mathbf{III}=\partial_x W-\partial_x (\overline{\mathcal{L}} W)$}

\begin{lemma}[Linear terms in $W$]\label{PROPIII}
\begin{align}
\mathbf{III} & = 
    \sum_{(k,\ell)\in \Sigma_0} c^\ell \left(Q_c^k(y_c) (-\mathcal{L} A_{k,\ell})'(y) +
    (Q_c^k)'(y_c) \big( 3 A_{k,\ell}'' + pQ^{p-1} A_{k,\ell} -(\mathcal{L} B_{k,\ell})' \big)(y)\right)
    \label{IIIA}\\ & \quad
    + \sum_{\substack{1\leq k\leq 4 k_0+ 2 p - 2    \\ 0\leq \ell \leq  4 \ell_0    + 2}} 
            c^\ell\left( Q_c^k(y_c) F_{k,\ell}^{\mathbf{III}}(y)    + 
             (Q_c^k)'(y_c)  G_{k,\ell}^{\mathbf{III}}(y)\right),\label{IIIB}
\end{align}
where    for any $k\geq 1$, $\ell\geq 0$, $F_{k,\ell}^{\mathbf{III}}$ and $G_{k,\ell}^{\mathbf{III}}$ satisfy:
\begin{itemize}
    \item[{\rm (i)}] Dependence property:
    $F_{k,\ell}^{\mathbf{III}}$ and $G_{k,\ell}^{\mathbf{III}}$ depend only on $(a_{k',\ell'})$ 
    and $(A_{k',\ell'})$, $(B_{k',\ell'})$ for $k',$ $\ell'$ such that $(k',\ell') \prec (k,\ell)$.
    \item[{\rm (ii)}] Parity property: Let $k\in \{1,\ldots,4 k_0+ 2 p - 2\}$, $\ell\in \{0,\ldots,4 \ell_0  + 2\}$. 
    Assume that for  any $(k',\ell')\prec (k,\ell)$,
       $A_{k',\ell'}$ is even and    $B_{k',\ell'}$ is odd, then 
     $F_{k,\ell}^{\mathbf{III}} $ is odd and  $G_{k,\ell}^{\mathbf{III}}$ is    even;
\end{itemize}
Moreover,
$
    F_{1,0}^{\mathbf{III}}=G_{1,0}^{\mathbf{III}}=0;
$
\begin{itemize}
    \item If $p=2$ then
    \begin{align*}
        & F_{2,0}^{\mathbf{III}}=a_{1,0}(-3 A_{1,0}''-2QA_{1,0})'-(3A_{1,0}'+ 3 B_{1,0}''+ 2 Q B_{1,0}), \\
        & G_{2,0}^{\mathbf{III}}=\frac{a_{1,0}}2 (-9 A_{1,0}'-3B_{1,0}''-2QB_{1,0})'-(A_{1,0}+3 B_{1,0}').\
    \end{align*}
    \item If $p=4$ then
    \begin{align*}
        & F_{2,0}^{\mathbf{III}}=a_{1,0}(-3 A_{1,0}''-pQ^{p-1}A_{1,0})', \quad G_{2,0}^{\mathbf{III}}=\frac{a_{1,0}}2 (-9 A_{1,0}'-3B_{1,0}''-pQ^{p-1} B_{1,0})'.
    \end{align*}
\end{itemize}
\end{lemma}

First, we claim two preliminary results concerning  $\mathbf{III}$.

\begin{claim}\label{QckA} 
    Let $k\in \mathbb{N}$ and let $A(x)$ be a class $C^3$ function. Let $w(t,x)=Q_c^k(y_c) A(y)$. Then,
    \begin{align*}
            & \partial_t w - \partial_x (\overline{\mathcal{L}} w) = Q_c^k(y_c) (-\mathcal{L} A)'(y)
            \\ & \quad + Q_c^k(y_c) \beta(y_c) (-3 A'' - p Q^{p-1} A + c A)'(y)+ Q_c^k(y_c) \beta'(y_c) (-3 A'')(y)
            \\ & \quad  + Q_c^k(y_c) \beta''(y_c) (-A')(y)+ Q_c^k(y_c) \beta^2(y_c) (3 A^{(3)})(y)
            \\ & \quad   + Q_c^k(y_c) \beta'(y_c) \beta(y_c) (3 A'')(y)+ Q_c^k(y_c) \beta^3(y_c) (-A^{(3)})(y)
            \\ & \quad +(Q_c^k)'(y_c) (3 A''+pQ^{p-1}A - c A)(y)
            \\ & \quad + (Q_c^k)'(y_c) \beta(y_c) (-6 A'')(y) + (Q_c^k)'(y_c) \beta'(y_c) (-3 A')(y) 
                + (Q_c^k)'(y_c) \beta^2(y_c) (3 A'')(y) 
            \\ & \quad+ (Q_c^k)''(y_c)(3A')(y) + (Q_c^k)''(y_c)\beta(y_c) (-3A')(y)+(Q_c^k)^{(3)}(y_c) A(y).
    \end{align*}
\end{claim}

\begin{claim}\label{QckB}
    Let $k\in \mathbb{N}$ and let $B(x)$ be a class $C^3$ function. Let $w(t,x)=(Q_c^k)'(y_c) B(y)$. Then,
    \begin{align*}
        & \partial_t w - \partial_x (\overline{\mathcal{L}} w) = (Q_c^k)'(y_c) (-\mathcal{L} B)'(y)
        \\ & \quad + (Q_c^k)'(y_c) \beta(y_c) (-3 B'' - p Q^{p-1} B + c B)'(y)+ (Q_c^k)'(y_c) \beta'(y_c) (-3 B'')(y)
        \\ & \quad + (Q_c^k)'(y_c) \beta''(y_c) (-B')(y)+ (Q_c^k)'(y_c) \beta^2(y_c) (3 B^{(3)})(y) 
        \\ & \quad  + (Q_c^k)'(y_c) \beta'(y_c) \beta(y_c) (3 B'')(y) + (Q_c^k)'(y_c) \beta^3(y_c) (-B^{(3)})(y)
        \\ & \quad +(Q_c^k)''(y_c) (3 B''+pQ^{p-1}B - c B)(y)
        \\ & \quad + (Q_c^k)''(y_c) \beta(y_c) (-6 B'')(y)+ (Q_c^k)''(y_c) \beta'(y_c) (-3 B')(y) 
         + (Q_c^k)''(y_c) \beta^2(y_c) (3 B'')(y) 
        \\ & \quad+ (Q_c^k)^{(3)}(y_c)(3B')(y) + (Q_c^k)^{(3)}(y_c)\beta(y_c) (-3B')(y)+(Q_c^k)^{(4)}(y_c) B(y).
    \end{align*}
\end{claim}

\noindent\emph{Proof of Claim \ref{QckA}.} Let $\mathcal{A}(t,x)=A(y)=A(x-\alpha(y_c))$, and $w(t,x)=Q_c^k(y_c)\mathcal{A}(t,x)$.
We first give the expression of $\partial_t w - \partial_x (\overline{\mathcal{L}} w)$ in terms of the partial derivatives of $\mathcal{A}$. First,
\begin{equation*}
    \partial_t w= (1-c) (Q_c^k)'(y_c) \mathcal{A}+ Q_c^k(y_c) \partial_t \mathcal{A}.
\end{equation*}
Since $\overline{\mathcal{L}} (fg)=g \overline{\mathcal{L}} f -2 f' g' - g''$, we have
$
    \overline{\mathcal{L}} w = Q_c^k(y_c) (\overline{\mathcal{L}} \mathcal{A}) - 2 (Q_c^k)'(y_c) \partial_x \mathcal{A} - (Q_c^k)''(y_c) \mathcal{A},
$
and so
\begin{align*}
    -\partial_x (\overline{\mathcal{L}} w) & = - Q_c^k(y_c) \partial_x (\overline{\mathcal{L}} \mathcal{A}) - (Q_c^k)'(y_c) (\overline{\mathcal{L}} \mathcal{A}) + 2 (Q_c^k)''(y_c) \partial_x \mathcal{A} \\ & \quad + 2 (Q_c^k)'(y_c) \partial_x^2 \mathcal{A} + (Q_c^k)^{(3)}(y_c) \mathcal{A} + (Q_c^k)''(y_c) \partial_x \mathcal{A}.
\end{align*}
Thus, by arranging terms by increasing order of derivatives of  $Q_c^k(y_c)$, we get
\begin{align}
 \partial_t w - \partial_x (\overline{\mathcal{L}} w) & = Q_c^k(y_c) (\partial_t \mathcal{A} - \partial_x (\overline{\mathcal{L}} \mathcal{A}))
+ (Q_c^k)'(y_c) ((1-c) \mathcal{A} - (\overline{\mathcal{L}} \mathcal{A}) + 2 \partial_x^2 \mathcal{A}) \nonumber\\ &\quad    + (Q_c^k)''(y_c) (3 \partial_x \mathcal{A})
+(Q_c^k)^{(3)}(y_c) \mathcal{A}.\label{firstDD}
\end{align}

Second, we use Claim \ref{forH} to express the partial derivatives of $\mathcal{A}$ in terms of derivatives of $A$.
We have
\begin{align*}
\partial_t \mathcal{A} - \partial_x (\overline{\mathcal{L}} \mathcal{A}) &= -(1-c)\beta(y_c) A'(y)+ (1-3 \beta(y_c) + 3 \beta^2(y_c)-\beta^3(y_c)) A^{(3)}(y) \\ &\quad +(-3 \beta'(y_c) + 3 \beta'(y_c)\beta(y_c)) A''(y) - \beta''(y_c) A'(y)\\ &\quad +( 1-\beta(y_c)
) (-A+pQ^{p-1}A)'(y)
\\ &= (1-3 \beta(y_c) + 3 \beta^2(y_c)-\beta^3(y_c)) A^{(3)}(y)
 +(-3 \beta'(y_c) + 3 \beta'(y_c)\beta(y_c)) A''(y) \\ &\quad + (1-c\beta(y_c)+\beta''(y_c)) (-A')(y)+(1-\beta(y_c)) (pQ^{p-1}A)'(y).
\end{align*}
Thus, by arranging terms by increasing order of derivatives and powers of $\beta(y_c)$, we get
\begin{align*}
&\partial_t \mathcal{A} - \partial_x (\overline{\mathcal{L}} \mathcal{A}) = (-\mathcal{L} A)'(y)
+ \beta(y_c) (-3 A'' - p Q^{p-1} A + c A)'(y)+  \beta'(y_c) (-3 A'')(y)
\\ & \quad + \beta''(y_c) (-A')(y) +\beta^2(y_c) (3 A^{(3)})(y) + \beta'(y_c) \beta(y_c) (3 A'')(y)
 +  \beta^3(y_c) (-A^{(3)})(y).
\end{align*}
Similarly,
\begin{align*}
& (1-c)\mathcal{A} - (\overline{\mathcal{L}} \mathcal{A})+ 2 \partial_x^2 \mathcal{A}    = - c \mathcal{A} + 3 \partial_x^2 \mathcal{A} + pQ^{p-1}(y) \mathcal{A}
\\ & = -c A(y) + 3(1-2 \beta(y_c) + \beta^2(y_c) ) A''(y)    - 3 \beta'(y_c) A'(y) + p Q^{p-1}(y) A(y)
\\ & =  (3 A''+pQ^{p-1}A - c A)(y) + \beta(y_c) (-6 A'')(y)  + \beta'(y_c) (-3 A')(y) +  \beta^2(y_c) (3 A'')(y),
\end{align*}
and
\begin{equation*}3 \partial_x \mathcal{A} = 3 A'(y) - 3\beta(y_c) A'(y).
\end{equation*}
Inserting all this into \eqref{firstDD}, we obtain  Claim \ref{QckA}.
Proof of Claim \ref{QckB} is the same.

\medskip

\noindent\emph{Proof of Lemma \ref{PROPIII}.}   We recall $
W(t,x)=\sum_{(k,\ell)\in \Sigma_0} c^\ell Q_{c}^k(y_{c}) A_{k,\ell}(y)+  c^\ell (Q_{c}^k)'(y_{c}) B_{k,\ell}(y).
$
To expand $\mathbf{III}$, we use Claims \ref{QckA}-\ref{QckB} on   $W(t,x)$.
We obtain
$
    \mathbf{III}    = \sum_{(k,\ell)\in \Sigma_0} c^\ell \ \mathrm{III}(k,\ell), 
$
where
{\allowdisplaybreaks
\begin{align*}&
    \mathrm{III}(k,\ell)    = 
Q_c^k(y_c) (-\mathcal{L} A_{k,\ell})'(y) + (Q_c^k)'(y_c) (3 A_{k,\ell}''+ p Q^{p-1} A_{k,\ell}-(\mathcal{L} B_{k,\ell})')(y) \\* &
+c (Q_c^k)'(y_c) (- A_{k,\ell})(y)  \\ &
    +\beta(y_c) Q_c^k(y_c) (- 3 A_{k,\ell}''-pQ^{p-1} A_{k,\ell})'(y) + \beta(y_c) (Q_c^k)'(y_c) (-6 A_{k,\ell}'- 3B_{k,\ell}''- p Q^{p-1} B_{k,\ell})'(y)  \\*
    &
    +c \beta(y_c) Q_c^k(y_c) (A_{k,\ell}')(y) + c \beta(y_c) (Q_c^k)'(y_c) (B_{k,\ell}')(y)  \\&
+\beta'(y_c) Q_c^k(y_c) (- 3 A_{k,\ell}'')(y) + \beta'(y_c) (Q_c^k)'(y_c) (-3 A_{k,\ell}'-3  B_{k,\ell}'')(y)    \\&
+\beta''(y_c) Q_c^k(y_c) (-A_{k,\ell}')(y) + \beta''(y_c) (Q_c^k)'(y_c) (-B_{k,\ell}')(y)    \\&
+\beta^2(y_c) Q_c^k(y_c) (3 A_{k,\ell}^{(3)})(y) + \beta^2(y_c) (Q_c^k)'(y_c) (3 A_{k,\ell}''+3 B_{k,\ell}^{(3)})(y)    \\&
+\beta'(y_c)\beta(y_c) Q_c^k(y_c) (3 A_{k,\ell}'')(y) +\beta'(y_c)\beta(y_c) (Q_c^k)'(y_c) (3    B_{k,\ell}'')(y)    \\&
+\beta^3(y_c) Q_c^k(y_c) (- A_{k,\ell}^{(3)})(y) + \beta^3(y_c) (Q_c^k)'(y_c) (- B_{k,\ell}^{(3)})(y)    \\&
+(Q_c^k)''(y_c) (3 A_{k,\ell}'+3 B_{k,\ell}''+pQ^{p-1}B_{k,\ell})(y) + (Q_c^k)^{(3)}(y_c) (A_{k,\ell}+3 B_{k,\ell}')(y)\\* &    +c (Q_c^k)''(y_c) (-    B_{k,\ell})(y) \\&
    +\beta(y_c) (Q_c^k)''(y_c) (- 3 A_{k,\ell}'-6 B_{k,\ell}'')(y) + \beta(y_c) (Q_c^k)^{(3)}(y_c) (-3 B_{k,\ell}')(y)  \\&
    +\beta'(y_c) (Q_c^k)''(y_c) (- 3 B_{k,\ell}')(y) + \beta^2(y_c) (Q_c^k)''(y_c) ( 3B_{k,\ell}'')(y) +(Q_c^k)^{(4)}(y_c) B_{k,\ell}(y)\\&
=\mathrm{III}_1+\mathrm{III}_2+\mathrm{III}_3+\mathrm{III}_4+\mathrm{III}_5+\mathrm{III}_6+\mathrm{III}_7+\mathrm{III}_8+\mathrm{III}_9+\mathrm{III}_{10}+\mathrm{III}_{11}+\mathrm{III}_{12}+\mathrm{III}_{13}.
\end{align*}
}
For   $j\in \{1,\ldots,13\}$, we denote  $\mathbf{III}_j=\sum_{(k,\ell)\in \Sigma_0} \mathrm{III}_j$.
\medskip

- Decomposition of $\mathbf{III}_1$.
This term gives  \eqref{IIIA}:
\begin{equation*}\mathbf{III}_1=\sum_{(k,\ell)\in \Sigma_0} c^\ell\left(Q_c^k(y_c) (-\mathcal{L} A_{k,\ell})'(y) + (Q_c^k)'(y_c) (3 A_{k,\ell}''+ p Q^{p-1} A_{k,\ell}-(\mathcal{L} B_{k,\ell})')(y)\right).
\end{equation*}

For the other terms, by elementary calculations, we obtain:
\begin{equation*}
    \mathbf{III}_2= 
\sum_{\substack{1 \leq k\leq    k_{0} \\ 1\leq \ell\leq  \ell_{0}+1 }}  c^\ell (Q_c^k)'(y_c) G_{k,\ell}^{\mathbf{III}_2}(y)
\quad \text{where}\quad G_{k,\ell}^{\mathbf{III}_2}(y)=(-A_{k,\ell-1})(y).
\end{equation*}
 
\begin{equation*}\mathbf{III}_3   = 
\sum_{\substack{2\leq k\leq  2 k_{0} \\ 0\leq \ell \leq 2    \ell_{0} }}c^{\ell} Q_c^{k}(y_c) F_{k,\ell}^{\mathbf{III}_3}(y)+
\sum_{\substack{2\leq k\leq  2 k_{0} \\ 0\leq \ell \leq 2    \ell_{0} }}c^{\ell} (Q_c^{k})'(y_c) G_{k,\ell}^{\mathbf{III}_3}(y), \quad \text{where}
\end{equation*}
\begin{equation}\label{FIII3}
     F_{k,\ell}^{\mathbf{III}_3}(y)=
\sum_{\substack{\max(k-k_0,1)\leq k_1 \leq \min(k-1,k_0)\\ \max(\ell-\ell_0,0)\leq \ell_1 \leq  \min(\ell,\ell_{0})}}
\hskip -1cm a_{k_1,\ell_1}   (-3 A_{k-k_1,\ell-\ell_1}''-pQ^{p-1} A_{k-k_1,\ell-\ell_1})'(y).
\end{equation}
\begin{align*}G_{k,\ell}^{\mathbf{III}_3}(y) =
&\sum_{\substack{\max(k-k_0,1)\leq k_1 \leq \min(k-1,k_0)
\\ \max(\ell-\ell_0,0)\leq \ell_1 \leq  \min(\ell,\ell_{0})}} \hskip -1cm
a_{k_1,\ell_1} \frac{k-k_1}k    \\&\quad \qquad \times(-6 A_{k-k_1,\ell-\ell_1}'-3 B_{k-k_1,\ell-\ell_1}''-pQ^{p-1} B_{k-k_1,\ell-\ell_1})'(y). \end{align*}
From \eqref{FIII3}, we easily check property {\rm (i)}  since in the sum defining $F_{k,\ell}^{\mathbf{III}_3}$, we have $k_1\leq k-1$ and $k-k_1\leq k-1$; moreover, $0\leq \ell_1\leq \ell$ and $\ell-\ell_1\leq \ell$. The parity statement {\rm (ii)} is also easily checked, as in the rest of this proof.  Thus $F_{k,\ell}^{\mathbf{III}_3}$ satisfies properties {\rm (i)} and {\rm (ii)}.

\begin{equation*}
\mathbf{III}_4=\sum_{\substack{2\leq k\leq  2 k_{0} \\ 1\leq \ell \leq 2    \ell_{0}+1 }} c^{\ell}  \left(Q_c^{k}(y_c) F_{k,\ell}^{\mathbf{III}_4}(y)+  (Q_c^{k})'(y_c) G_{k,\ell}^{\mathbf{III}_4}(y)\right),
\end{equation*}
where
\begin{equation*} 
     F_{k,\ell}^{\mathbf{III}_4}(y)=
\sum_{\substack{\max(k-k_0,1)\leq k_1 \leq \min(k-1,k_0)\\ \max(\ell-\ell_0-1,0)\leq \ell_1 \leq    \min(\ell-1,\ell_{0})}}
\hskip -1cm a_{k_1,\ell_1}   A_{k-k_1,\ell-\ell_1-1}'(y)
\end{equation*}
\begin{equation*}    
    G_{k,\ell}^{\mathbf{III}_4}(y)=
\sum_{\substack{\max(k-k_0,1)\leq k_1 \leq \min(k-1,k_0)\\ \max(\ell-\ell_0-1,0)\leq \ell_1 \leq    \min(\ell-1,\ell_{0})}}
\hskip -1cm a_{k_1,\ell_1}  \frac{k-k_1}k B_{k-k_1,\ell-\ell_1-1}'(y) .
\end{equation*}
 \begin{equation*}
\mathbf{III}_{5}    
  = \sum_{\substack{2\leq k\leq  2 k_{0} \\ 0\leq \ell \leq 2    \ell_{0} }}c^{\ell}
        (Q_c^{k})'(y_c) G_{k,\ell}^{\mathbf{III}_5}(y) 
        +\sum_{\substack{2\leq k\leq    2 k_{0}+p-1 \\ 0\leq \ell \leq 2    \ell_{0}+1 }}c^{\ell} Q_c^{k}(y_c) F_{k,\ell}^{\mathbf
{III}_5}(y) 
       ,\quad \text{where}
\end{equation*}
where
\begin{equation*} G_{k,\ell}^{\mathbf{III}_5}(y)=
\sum_{\substack{\max(k-k_0,1)\leq k_1 \leq \min(k-1,k_0)\\ \max(\ell-\ell_0,0)\leq \ell_1 \leq  \min(\ell,\ell_{0})}}
\frac {k_1}{k} a_{k_1,\ell_1}       (-3A_{k-k_1,\ell-\ell_1}'')(y),
\end{equation*}
\begin{align}    F_{k,\ell}^{\mathbf{III}_5}(y) & =
\sum_{\substack{\max(k-k_0,1)\leq k_1 \leq \min(k-1,k_0)\\ \max(\ell-\ell_0-1,0)\leq \ell_1 \leq    \min(\ell-1,\ell_{0})}}
\hskip -1cm {k_1}(k-k_1) a_{k_1,\ell_1}     (-3 A_{k-k_1,\ell-\ell_1-1}'-3 B_{k-k_1,\ell-\ell_1-1}'')(y)\nonumber \\&\quad
+\hskip -1cm \sum_{\substack{\max(k-k_0-p+1,1)\leq k_1 \leq \min(k-p,k_0)\label{FIII5}
\\ \max(\ell-\ell_0,0)\leq \ell_1 \leq  \min(\ell,\ell_{0})}}
 \hskip -1cm \frac {2k_1(k{-}k_1{-}p{+}1)}{p+1}  a_{k_1,\ell_1} \\ &\hskip 5cm \times(3 A_{k-k_1-p+1,\ell-\ell_1}'+3 B_{k-k_1-p+1,\ell-\ell_1}'')(y)\nonumber.
\end{align}
satisfy properties {\rm (i)} and {\rm (ii)}.
In \eqref{FIII5}, the first sum term has no contribution for $k$ such that
$\max(k-k_0,1)>\min(k-1,k_0)$ (i.e. $k=1$ or $k> 2k_0$), and similarly for the  condition on $\ell$. We will use this notation in all sums appearing in this proof.

For the next terms, we use Claim \ref{surDD}.
\begin{equation*}\mathbf{III}_6=\sum_{\substack{2\leq k \leq 2k_0+p-1\\ 0 \leq \ell \leq 2 \ell_0+1}} c^{\ell}\left( Q_c^{k}(y_c) F_{k,\ell}^{\mathbf{III}_6}(y)    + (Q_c^{k})'(y_c) G_{k,\ell}^{\mathbf{III}_6}(y) \right),
\end{equation*}
where
\begin{equation*} 
    F_{k,\ell}^{\mathbf{III}_6}(y)=
\sum_{\substack{\max(k-k_0,1)\leq k_1\leq    \min(k-1,k_{0}+p-1) \\ \max(\ell-\ell_0,0)\leq \ell_1 \leq \min(\ell,\ell_{0}+1)}}
\hskip -1cm a^{1*}_{k_1,\ell_1} (-A_{k-k_1,\ell-\ell_1}')(y)
\end{equation*}
\begin{equation*} 
    G_{k,\ell}^{\mathbf{III}_6}(y)=
\sum_{\substack{\max(k-k_0,1)\leq k_1\leq    \min(k-1,k_{0}+p-1) \\ \max(\ell-\ell_0,0)\leq \ell_1 \leq \min(\ell,\ell_{0}+1)}}
\frac {k-k_1}{k} a^{1*}_{k_1,\ell_1} (-B_{k-k_1,\ell-\ell_1}')(y).
\end{equation*}

\begin{equation*}
    \mathbf{III}_7=\sum_{\substack{3\leq k \leq 3 k_0\\ 0 \leq \ell \leq 3 \ell_0}}c^{\ell}\left(  Q_c^{k}(y_c) F_{k,\ell}^{\mathbf{III}_7}(y) +(Q_c^{k})'(y_c) G_{k,\ell}^{\mathbf{III}_7}(y) \right),
\end{equation*}
where
\begin{equation*}
    F_{k,\ell}^{\mathbf{III}_7}(y) = \sum_{\substack{\max(k-k_0,2)\leq k_1\leq \min(k-1,2k_0)\\ \max(\ell-\ell_0,0)\leq \ell_1 \leq \min(\ell,2\ell_0)}} {a^{2*}_{k_1,\ell_1}} (3 A_{k-k_1,\ell-\ell_1}^{(3)})(y),
\end{equation*}
\begin{equation*}
    G_{k,\ell}^{\mathbf{III}_7}(y) = \sum_{\substack{\max(k-k_0,2)\leq k_1\leq \min(k-1,2k_0)\\ \max(\ell-\ell_0,0)\leq \ell_1 \leq \min(\ell,2\ell_0)}} \frac {k-k_1}k a^{2*}_{k_1,\ell_1} (3 A_{k-k_1,\ell-\ell_1}''+3 B_{k-k_1,\ell-\ell_1}^{(3)})(y).
\end{equation*}

\begin{equation*}
    \mathbf{III}_8=\sum_{\substack{3\leq k \leq 3 k_0+p-1\\ 0 \leq \ell \leq 3 \ell_0+1}}c^{\ell}\left( Q_c^{k}(y_c) F_{k,\ell}^{\mathbf{III}_8}(y)+ (Q_c^{k})'(y_c) G_{k,\ell}^{\mathbf{III}_8}(y)\right),  
\end{equation*}
where
\begin{equation*}
    G_{k,\ell}^{\mathbf{III}_8}(y) = \sum_{\substack{\max(k-k_0,2)\leq k_1\leq \min(k-1,2k_0)\\ \max(\ell-\ell_0,0)\leq \ell_1 \leq \min(\ell,2\ell_0)}} \frac {k_1}{k} a^{3*}_{k_1,\ell_1} (3 A_{k-k_1,\ell-\ell_1}''(y))
\end{equation*}
\begin{align*}&
    F_{k,\ell}^{\mathbf{III}_8}(y) = \sum_{\substack{\max(k-k_0,2)\leq k_1\leq \min(k-1,2k_0)\\ \max(\ell-\ell_0-1,0)\leq \ell_1 \leq \min(\ell-1,2\ell_0)}} {k_1}(k-k_1) a^{3*}_{k_1,\ell_1} (3 B_{k-k_1,\ell-\ell_1-1}''(y))\\ & + \sum_{\substack{\max(k-k_0
-p+1,2)\leq k_1\leq \min(k-p,2k_0)\\ \max(\ell-\ell_0,0)\leq \ell_1 \leq \min(\ell,2\ell_0)}} -\frac {2{k_1}(k{-}k_1{-}p{+}1)}{p+1} a^{3*}_{k_1,\ell_1} (3 B_{k-k_1-p+1,\ell-\ell_1}''(y)).
\end{align*}

\begin{equation*}
    \mathbf{III}_9=\sum_{\substack{4\leq k \leq 4 k_0\\ 0 \leq \ell \leq 4 \ell_0}}c^{\ell}\left( Q_c^{k}(y_c) F_{k,\ell}^{\mathbf{III}_9}(y)+  (Q_c^{k})'(y_c) G_{k,\ell}^{\mathbf{III}_9}(y)\right),
\end{equation*}
where
\begin{equation*}
    F_{k,\ell}^{\mathbf{III}_9}(y)=
\sum_{\substack{ \max(k-k_0,3) \leq k_1\leq  \min( k-1, 3k_{0}) \\ \max(\ell-\ell_0,0)\leq \ell_1 \leq \min(\ell,3\ell_{0})}}
a^{4*}_{k_1,\ell_1} (-A_{k-k_1, \ell - \ell_1}^{(3)})(y)
\end{equation*}
\begin{equation*}
    G_{k,\ell}^{\mathbf{III}_9}(y)=
\sum_{\substack{ \max(k-k_0,3) \leq k_1\leq  \min(k-1, 3k_{0}) \\ \max(\ell-\ell_0,0)\leq \ell_1 \leq \min(\ell,3\ell_{0})}} \frac {k-k_1}k
a^{4*}_{k_1,\ell_1} (-B_{k-k_1, \ell - \ell_1}^{(3)})(y).
\end{equation*}

Using  from Lemma \ref{surQc} the expression of $(Q_c^k)''$:
\begin{align*}
    \mathbf{III}_{10} & = \sum_{\substack{1\leq k \leq k_0\\ 1 \leq \ell \leq    \ell_0+1 }}\left(c^{\ell} Q_c^{k}(y_c) F_{k,\ell}^{\mathbf{III}_{10}^1}(y)+c^{\ell} (Q_c^{k})'(y_c) G_{k,\ell}^{\mathbf{III}_{10}^1}(y)\right)  \\ &\quad
     +\sum_{\substack{p\leq k \leq k_0+p-1\\ 0 \leq \ell \leq    \ell_0 }}\left(c^{\ell} Q_c^{k}(y_c) F_{k,\ell}^{\mathbf{III}_{10}^2}(y)+c^{\ell} (Q_c^{k})'(y_c) G_{k,\ell}^{\mathbf{III}_{10}^2}(y)\right),
\end{align*}
where
\begin{equation*}
    F_{k,\ell}^{\mathbf{III}_{10}^1}(y) = k^2 (3A_{k,\ell-1}'+3 B_{k,\ell-1}''+pQ^{p-1} B_{k,\ell-1})(y)
,\end{equation*}
\begin{equation*}
    F_{k,\ell}^{\mathbf{III}_{10}^2}(y) = -\frac {(k-p+1)(2k-p+1)}{p+1} (3A_{k-p+1,\ell}'+3 B_{k-p+1,\ell}''+pQ^{p-1} B_{k-p+1,\ell})(y)
,\end{equation*}
\begin{equation*}
    G_{k,\ell}^{\mathbf{III}_{10}^1}(y) = k^2 (A_{k,\ell-1}+3 B_{k,\ell-1}')(y), 
\end{equation*}
\begin{equation*}
G_{k,\ell}^{\mathbf{III}_{10}^2}(y)=    -\frac {(k-p+1)(2k-p+1)}{p+1} (A_{k-p+1,\ell}+3 B_{k-p+1,\ell}')(y).
\end{equation*}
We set 
$F_{k,\ell}^{\mathbf{III}_{10}}=F_{k,\ell}^{\mathbf{III}_{10}^1}+F_{k,\ell}^{\mathbf{III}_{10}^2},\,G_{k,\ell}^{\mathbf{III}_{10}}=G_{k,\ell}^{\mathbf{III}_{10}^1}+G_{k,\ell}^{\mathbf{III}_{10}^2}.$

\begin{equation*}
\mathbf{III}_{11}= \sum_{\substack{1\leq k \leq k_0\\ 2 \leq \ell \leq  \ell_0+2 }}c^{\ell} Q_c^{k}(y_c) F_{k,\ell}^{\mathbf{III}_{11}^1}(y) +\sum_{\substack{p\leq k \leq k_0+p-1\\ 1\leq \ell \leq    \ell_0+1 }}
c^{\ell} Q_c^{k}(y_c) F_{k,\ell}^{\mathbf{III}_{11}^2}(y)
\end{equation*}
where
\begin{equation*} 
F_{k,\ell}^{\mathbf{III}_{11}^1}(y) = k^2 (-B_{k,\ell-2})(y),
\quad 
F_{k,\ell}^{\mathbf{III}_{11}^2}(y) = \frac {(k-p+1)(2k-p+1)}{p+1}  B_{k-p+1,\ell-1}(y).
\end{equation*}
We set
$F_{k,\ell}^{\mathbf{III}_{11}}=F_{k,\ell}^{\mathbf{III}_{11}^1}+F_{k,\ell}^{\mathbf{III}_{11}^2}.
$
\begin{align*}
\mathbf{III}_{12} & = \sum_{\substack{2\leq k \leq 2 k_0\\ 1 \leq \ell \leq  2 \ell_0+1 }}\left(c^{\ell} Q_c^{k}(y_c) F_{k,\ell}^{\mathbf{III}_{12}^1}(y)+c^{\ell} (Q_c^{k})'(y_c) G_{k,\ell}^{\mathbf{III}_{12}^1}(y)\right)    
\\&\quad +\sum_{\substack{p+1\leq k \leq 2 k_0+p-1\\ 0 \leq \ell \leq    2 \ell_0 }}\left(c^{\ell} Q_c^{k}(y_c) F_{k,\ell}^{\mathbf{III}_{12}^2}(y)+c^{\ell} (Q_c^{k})'(y_c) G_{k,\ell}^{\mathbf{III}_{12}^2}(y)\right),
\end{align*}
where
\begin{equation*}
F_{k,\ell}^{\mathbf{III}_{12}^1}(y) = \sum_{\substack{\max(k-k_0,1)\leq k_1 \leq \min(k-1,k_0)\\ \max(\ell-\ell_0-1,0)\leq \ell_1 \leq  \min(\ell-1,\ell_{0})}} \hskip -1cm
a_{k_1,\ell_1}  (k{-}k_1)^2 (-3A_{k-k_1,\ell-\ell_1-1}'-6 B_{k-k_1,\ell-\ell_1-1}'')(y)
,\end{equation*}
\begin{align*}
F_{k,\ell}^{\mathbf{III}_{12}^2}(y) = \sum_{\substack{\max(k-k_0-p+1,1)\leq k_1 \leq \min(k-p,k_0)\\ \max(\ell-\ell_0,0)\leq \ell_1 \leq    \min(\ell,\ell_{0})}}    \hskip -1cm
&    a_{k_1,\ell_1}      \frac {(k{-}k_1{-}p{+}1)(2k{-}2k_1{-}p{+}1)}{p{+}1}    \\& \qquad \times ( 3A_{k-k_1-p+1,\ell-\ell_1}'+6 B_{k-k_1-p+1,\ell-\ell_1}'')(y)
,\end{align*}
\begin{align*}
G_{k,\ell}^{\mathbf{III}_{12}^1}(y) &= \sum_{\substack{\max(k-k_0,1)\leq k_1 \leq \min(k-1,k_0)\\ \max(\ell-\ell_0-1,0)\leq \ell_1 \leq  \min(\ell-1,\ell_{0})}} \hskip -0.5cm \frac { (k{-}k_1)^3}k a_{k_1,\ell_1}(-3 B_{k-k_1,\ell-\ell_1-1}')(y), 
\end{align*}
\begin{align*}
G_{k,\ell}^{\mathbf{III}_{12}^2}(y) = \sum_{\substack{\max(k-k_0-p+1,1)\leq k_1 \leq \min(k-p,k_0)\\ \max(\ell-\ell_0,0)\leq \ell_1 \leq    \min(\ell,\ell_{0})}} & \hskip -0.5cm \frac {(k{-}k_1{-}p{+}1)^2(2k{-}2k_1{-}p{+}1)}{k(p+1)}    a_{k_1,\ell_1}  \\&
\qquad\qquad\times    3 B_{k-k_1-p+1,\ell-\ell_1}' (y).
\end{align*}
We set
$   F_{k,\ell}^{\mathbf{III}_{12}}=F_{k,\ell}^{\mathbf{III}_{12}^1}+F_{k,\ell}^{\mathbf{III}_{12}^2},
\,    G_{k,\ell}^{\mathbf{III}_{12}}=G_{k,\ell}^{\mathbf{III}_{12}^1}+G_{k,\ell}^{\mathbf{III}_{12}^2}.
$

The last term $\mathbf{III}_{13}$.
is the sum of three different terms.

\noindent The contribution of $\beta'(y_c) (Q_c^k)''(y_c) (-3 B'_{k,\ell}(y))$ is
\begin{align*}&
\sum_{\substack{2\leq k \leq 2 k_0\\ 1 \leq \ell \leq    2 \ell_0+1 }} c^{\ell} (Q_c^{k})'(y_c) G_{k,\ell}^{\mathbf{III}_{13}^1}(y)
 +\sum_{\substack{p+1\leq k \leq 2 k_0+p-1\\ 0 \leq \ell \leq    2 \ell_0 }} c^{\ell} (Q_c^{k})'(y_c) G_{k,\ell}^{\mathbf{III}_{13}^2}(y) ,
\end{align*}
where
\begin{equation*}
G_{k,\ell}^{\mathbf{III}_{13}^1}(y)  = \sum_{\substack{\max(k-k_0,1)\leq k_1 \leq \min(k-1,k_0)\\ \max(\ell-\ell_0-1,0)\leq \ell_1 \leq  \min(\ell-1,\ell_{0})}} \hskip -0.5cm \frac {k_1 (k{-}k_1)^2} k a_{k_1,\ell_1} (-3 B_{k-k_1,\ell-\ell_1-1}')(y),
\end{equation*}
and
\begin{align*}
G_{k,\ell}^{\mathbf{III}_{13}^2}(y) = \sum_{\substack{\max(k-k_0-p+1,1)\leq k_1 \leq \min(k-p,k_0)\\ \max(\ell-\ell_0,0)\leq \ell_1 \leq    \min(\ell,\ell_{0})}} & \frac {k_1(k{-}k_1{-}p{+}1)(2k{-}2k_1{-}p{+}1)}{k(p+1)}a_{k_1,\ell_1}    \\&\times  (3 B_{k
-k_1-p+1,\ell-\ell_1}')(y).
\end{align*}
The  contribution of $\beta^2(y_c) (Q_c^k)''(y_c) ( 3 B''_{k,\ell}(y))$ is, using \eqref{astar}
\begin{align*}&
\sum_{\substack{3\leq k \leq 3 k_0\\ 1 \leq \ell \leq    3 \ell_0+1 }} c^{\ell} Q_c^{k}(y_c) F_{k,\ell}^{\mathbf{III}_{13}^1}(y)
 +\sum_{\substack{p+2\leq k \leq 3 k_0+p-1\\ 0 \leq \ell \leq    3 \ell_0 }} c^{\ell} Q_c^{k}(y_c) F_{k,\ell}^{\mathbf{III}_{13}^2}(y) ,
\end{align*}
where
\begin{equation*}
F_{k,\ell}^{\mathbf{III}_{13}^1}(y)  = \sum_{\substack{\max(k-k_0,2)\leq k_1 \leq \min(k-1,2k_0)\\ \max(\ell-\ell_0-1,0)\leq \ell_1 \leq    \min(\ell-1,2 \ell_{0})}}    \hskip -1cm (k{-}k_1)^2 a^{2*}_{k_1,\ell_1}(3 B_{k-k_1,\ell-\ell_1-1}'')(y),
\end{equation*}
\begin{align*}
F_{k,\ell}^{\mathbf{III}_{13}^2}(y) = \sum_{\substack{\max(k-k_0-p+1,2)\leq k_1 \leq \min(k-p,2k_0)\\ \max(\ell-\ell_0,0)\leq \ell_1 \leq    \min(\ell,2\ell_{0})}} & \hskip -1cm \frac{(k{-}k_1{-}p{+}1)(2k{-}2k_1{-}p{+}1)}{p{+}1 }\\&\times a^{2*}_{k_1,\ell
_1}(-3 B_{k-k_1-p+1,\ell-\ell_1}'')(y).
\end{align*}
From Lemma \ref{surQc},
the contribution of $(Q_c^k)^{(4)}(y_c) B(y)$ is
\begin{align*}
& \sum_{\substack{1\leq k \leq k_0\\ 2 \leq \ell \leq    \ell_0+2 }} c^{\ell} Q_c^{k}(y_c) F_{k,\ell}^{\mathbf{III}_{13}^3}(y)
 +\sum_{\substack{p\leq k \leq k_0+p-1\\ 1 \leq \ell \leq    \ell_0+1 }} c^{\ell} Q_c^{k}(y_c) F_{k,\ell}^{\mathbf{III}_{13}^4}(y)
\\ &
 +\sum_{\substack{2p-1\leq k \leq k_0+2p-2\\ 0 \leq \ell \leq    \ell_0 }} c^{\ell} Q_c^{k}(y_c) F_{k,\ell}^{\mathbf{III}_{13}^5}(y),
\end{align*}
where
$
F_{k,\ell}^{\mathbf{III}_{13}^3}(y) = k^4 B_{k,\ell-2} (y)
,$
\begin{align*}
F_{k,\ell}^{\mathbf{III}_{13}^4}(y) &=-  \frac {(k{-}p{+}1)(2k{-}p{+}1)}{p+1}    ((k{-}p{+}1)^2+k^2)    B_{k-p+1,\ell-1} (y),\end{align*}
\begin{equation*}
F_{k,\ell}^{\mathbf{III}_{13}^5}(y) = (k-2p+2)(k-p+1)\frac{(2k{-}3p{+}3)(2k{-}p{+}1)}{(p+1)^2} B_{k-2p+2,\ell}(y).
\end{equation*}
We set 
$  F_{k,\ell}^{\mathbf{III}_{13}}=F_{k,\ell}^{\mathbf{III}_{13}^1}+F_{k,\ell}^{\mathbf{III}_{13}^2}+F_{k,\ell}^{\mathbf{III}_{13}^3}+F_{k,\ell}^{\mathbf{III}_{13}^4}+F_{k,\ell}^{\mathbf{III}_{13}^5},\,    G_{k,\ell}^{\mathbf{III}_{13}}=G_{k,\ell}^{\mathbf{III}_{13}^1}+G_{k,\ell}^{\mathbf{III}_{13}^2},
$
so that 
\begin{equation*}
    \mathbf{III}_{13}=\sum_{\substack{1\leq k\leq \max(3 k_0+p-1, k_0+ 2 p - 2 ) \\ 0\leq \ell \leq  \max(3 \ell_0+1,    \ell_0 + 2)}} 
            c^\ell\left( Q_c^k(y_c) F_{k,\ell}^{\mathbf{III}_{13}}(y)    + 
             (Q_c^k)'(y_c)  G_{k,\ell}^{\mathbf{III}_{13}}(y)\right).
\end{equation*}
Finally, we set $
    F_{k,\ell}^{\mathbf{III}}=\sum_{j=3}^{13} F_{k,\ell}^{\mathbf{III}_j},\,
    G_{k,\ell}^{\mathbf{III}}=\sum_{j=2}^{13} G_{k,\ell}^{\mathbf{III}_j} .
$

\medskip

We now finish the proof of Lemma \ref{PROPIII} by computing explicitely 
$ F_{1,0}^{\mathbf{III}},$ $    G_{1,0}^{\mathbf{III}},$ $ F_{2,0}^{\mathbf{III}}$ and $G_{2,0}^{\mathbf{III}}$.
We first check $F_{1,0}^{\mathbf{III}}=G_{1,0}^{\mathbf{III}}=0$.
For $F_{2,0}^{\mathbf{III}}$, we make the following observations: 
\begin{itemize}
    \item $F_{2,0}^{\mathbf{III}_3}=a_{1,0} (-3 A_{1,0}''-p Q^{p-1}A_{1,0})'$;
     $F_{2,0}^{\mathbf{III}_4}=F_{2,0}^{\mathbf{III}_5}=0$;
      $F_{2,0}^{\mathbf{III}_6}=0$ since $a^{1*}_{1,0}=0$;
    \item $F_{2,0}^{\mathbf{III}_7}=F_{2,0}^{\mathbf{III}_8}=F_{2,0}^{\mathbf{III}_9}=F_{2,0}^{\mathbf{III}_{10}^1}=0$;
    \item For $p=2$, we have $F_{2,0}^{\mathbf{III}_{10}^2}= - (3 A_{1,0}' + 3 B_{1,0}''+ 2 Q B_{1,0} )$, for $p= 4$, we have
    $F_{2,0}^{\mathbf{III}_{10}^2}= 0$.
    \item All the remaining terms in $F_{2,0}^{\mathbf{III}}$ are checked to be zero.
\end{itemize}
Similarly, we check that the only non zero contributions to $G_{2,0}^{\mathbf{III}}$ are
\begin{itemize}
    \item $G_{2,0}^{\mathbf{III}_3}=\frac 12 a_{1,0} (-6 A_{1,0}'- 3 B_{1,0}'' - p Q^{p-1} B_{1,0})'$;
      $G_{2,0}^{\mathbf{III}_5}=\frac 12 a_{1,0} (-3 A_{1,0}')$;
    \item if $p=2$, $G_{2,0}^{\mathbf{III}_9}=-(A_{1,0}+3 B_{1,0}')$, and if $p=4$, $G_{2,0}^{\mathbf{III}_9}=0$.
\end{itemize}
Thus, summing up, Lemma \ref{PROPIII} is proved.

\subsection{Expansion of $\mathbf{IV}=\partial_x\left((R+R_c+W)^p - (R+R_c)^p - p R^{p-1} W\right)$}

\begin{lemma}[Nonlinear terms in $W$]\label{PROPIV}
    \begin{equation*}
        {{\mathbf{IV}}}=\sum_{\substack{2\leq k \leq (p+1)k_0+12 \\ 0 \leq \ell \leq    (p+1)\ell_0+4 }} 
        c^{\ell} \left( Q_c^{k}(y_c) F_{k,\ell}^{{{\mathbf{IV}}}}(y) +  (Q_c^{k})'(y_c) G_{k,\ell}^{{{\mathbf{IV}}}}(y)\right),
    \end{equation*}
    where $F_{k,\ell}^{\mathbf{IV}}$ and $G_{k,\ell}^{\mathbf{IV}}$ are functions defined on $\mathbb{R}$ satisfying
    \begin{itemize}
        \item[{\rm (i)}] Dependence property:
        $F_{k,\ell}^{\mathbf{IV}}$ and $G_{k,\ell}^{\mathbf{IV}}$ depend only on $(a_{k',\ell'})$ 
        and $(A_{k',\ell'})$, $(B_{k',\ell'})$ for $k',$ $\ell'$ such that $(k',\ell') \prec (k,\ell)$.
        \item[{\rm (ii)}] Parity property: Let $k\in \{1,\ldots, (p+1)k_0+12 \}$, $\ell\in \{0,\ldots,(p+1)\ell_0+4\}$. 
        Assume that for any $(k',\ell')$ such that $(k',\ell')\prec (k,\ell)$,
         $A_{k',\ell'}$ is even and    $B_{k',\ell'}$ is odd, then 
         $F_{k,\ell}^{\mathbf{IV}} $ is odd and  $G_{k,\ell}^{\mathbf{IV}}$ is  even.
    \end{itemize}
    Moreover,
    \begin{itemize}
        \item If $p=2$ then
        \begin{align}\label{FIVdeux}
            & F_{2,0}^{{{\mathbf{IV}}}} = \left(2 A_{1,0}+  A_{1,0}^2\right)',\quad 
             G_{2,0}^{{{\mathbf{IV}}}} =    2A_{1,0} + A_{1,0}^2+ \left(B_{1,0}+A_{1,0}B_{1,0}\right)'.
        \end{align}
                 \item If $p=4$ then
        \begin{align}\label{FIVquatre}
            & F_{2,0}^{{{\mathbf{IV}}}} =    \left(12 A_{1,0} Q^2 + 6 A_{1,0}^2 Q^2\right)'\\    
            & G_{2,0}^{{{\mathbf{IV}}}} =    12 A_{1,0} Q^2 + 6 A_{1,0}^2 Q^2 + \left(6 B_{1,0} Q^2 + 6 A_{1,0} B_{1,0} Q^2\right)'.
        \end{align}
    \end{itemize}
\end{lemma}

\noindent\emph{Proof of Lemma \ref{PROPIV}.}
Set
$\mathbf{N}=(R+R_c+W)^p - (R+R_c)^p - p R^{p-1} W.
$
First, we determine $F^{\mathbf{N}}_{k,\ell}$ and $G^{\mathbf{N}}_{k,\ell}$ such that
\begin{equation}\label{IVtildedeux}
    {\mathbf{N}}= \sum_{\substack{k\geq p, ~  \ell \geq 0}}
    c^\ell \left(Q_c^k(y_c)  F_{k,\ell}^{\mathbf{N}}(y)+ (Q_c^k)'(y_c) G_{k,\ell}^{\mathbf{N}}(y)\right).
\end{equation}
Second, we differentiate formula \eqref{IVtildedeux} with respect to $x$ to get the decomposition of ${{\mathbf{IV}}}$. We treat only the case $p=4$, the case $p=2$ is similar and easier.
 
 $\bullet$ $p=4$. 
\begin{align*}
    {\mathbf{N}} 
    & = 4 \left((R+R_c)^3-R^3\right) W + 6 (R+R_c)^2 W^2 + 4 (R+R_c) W^3 + W^4\\
    & = 12 R^2 R_c W + 12 R R_c^2 W + 4 R_c^3 W + 6 R^2 W^2 + 12 R R_c W^2 + 6 R_c^2 W^2 \\
    &\quad +4 R W^3 + 4 R_c W^3 + W^4 \\
    & = {\mathbf{N}}_1+{\mathbf{N}}_2+{\mathbf{N}}_3+{\mathbf{N}}_4+{\mathbf{N}}_5+{\mathbf{N}}_6+{\mathbf{N}}_7+{\mathbf{N}}_8+{\mathbf{N}}_9.
\end{align*}

- Terms ${\mathbf{N}}_1$, ${\mathbf{N}}_2$, ${\mathbf{N}}_3$.
\begin{align*}
    {\mathbf{N}}_1 
    &= Q_c(y_c) (12 W Q^2(y)) \\
    &= \sum_{\substack{2\leq k \leq k_0+1\\ 0 \leq \ell \leq \ell_0}} 
        c^\ell \left(Q_c^k(y_c) ( 12 A_{k-1,\ell}(y) Q^2(y)) + (Q_c^k)'(y_c)    \frac{k{-}1}k
        \left(12 B_{k-1,\ell} Q^2\right)(y)\right);
\end{align*}
\begin{align*}
    {\mathbf{N}}_2 
    &= Q_c^2(y_c) (12 W  Q(y)) \\
    &= \sum_{\substack{3\leq k \leq k_0+2\\ 0 \leq \ell \leq \ell_0}} 
        c^\ell \left(Q_c^k(y_c) ( 12 A_{k-2,\ell}(y) Q(y))
        + (Q_c^k)'(y_c)  \frac{k{-}2}k    \left(12 B_{k-2,\ell} Q\right)(y)\right);
\end{align*}
\begin{align*}
    {\mathbf{N}}_3 
    &= Q_c^3(y_c) (4 W) \\
    &= \sum_{\substack{4\leq k \leq k_0+3\\ 0 \leq \ell \leq \ell_0}} 
    c^\ell \left(Q_c^k(y_c) ( 4 A_{k-3,\ell}(y)) + (Q_c^k)'(y_c)        \frac{k{-}3}k    \left(4B_{k-3,\ell}\right)(y)\right).
\end{align*}

For the next three terms, we first need to expand $W^2$:
\begin{align*}
    W^2 
    & = \sum_{\substack{1\leq k_1 \leq k_0 \\ 0 \leq \ell_1 \leq    \ell_0}}
    c^{\ell_1}\left(Q_c^{k_1}(y_c)A_{k_1,\ell_1}(y)+(Q_c^{k_1})'(y_c)B_{k_1,\ell_1}(y)\right)\\
    &\quad
    \times \sum_{\substack{1\leq k_2 \leq k_0 \\ 0 \leq \ell_2 \leq  \ell_0}}
    c^{\ell_2}\left(Q_c^{k_2}(y_c)A_{k_2,\ell_2}(y)+(Q_c^{k_2})'(y_c)B_{k_2,\ell_2}(y)\right).
\end{align*}
Using Lemma \ref{surQc},
\begin{align*}
    W^2 
    & = \sum_{\substack{2\leq k \leq 2k_0 \\ 0 \leq \ell \leq    2\ell_0}}
    c^{\ell} Q_c^{k}(y_c) \sum_{\substack{\max(k-k_0,1)\leq k_1 \leq \min(k-1,k_0)\\
    \max(\ell-\ell_0,0)\leq \ell_1 \leq  \min(\ell,\ell_{0})}}   \hskip -1cm A_{k_1,\ell_1}(y) A_{k-k_1,\ell-\ell_1}(y) \\
    & +\sum_{\substack{2\leq k \leq 2 k_0 \\ 1 \leq \ell \leq    2 \ell_0+1}}
    c^{\ell} Q_c^{k}(y_c) \sum_{\substack{\max(k-k_0,1)\leq k_1 \leq \min(k-1,k_0)\\
    \max(\ell-\ell_0-1,0)\leq \ell_1 \leq    \min(\ell-1,\ell_{0})}}    \hskip -1cm 
     (k_1(k{-}{k_1}) B_{k_1,\ell_1}B_{k-k_1,\ell-\ell_1-1})(y) \\
    & +\sum_{\substack{5\leq k \leq 2 k_0 +3 \\ 0 \leq \ell \leq    2 \ell_0}}
    c^{\ell} Q_c^{k}(y_c) \sum_{\substack{\max(k-k_0-3,1)\leq k_1 \leq \min(k-4,k_0)\\
    \max(\ell-\ell_0,0)\leq \ell_1 \leq  \min(\ell,\ell_{0})}} 
     \hskip -1cm \left(-\tfrac 25 k_1 (k{-}k_1{-}3) B_{k_1,\ell_1} B_{k-k_1-3,\ell-\ell_1}\right)(y) \\
    & +\sum_{\substack{2\leq k \leq 2 k_0 \\ 0 \leq \ell \leq    2 \ell_0}}
    c^{\ell} (Q_c^{k})'(y_c) \sum_{\substack{\max(k-k_0,1)\leq k_1 \leq \min(k-1,k_0)\\
    \max(\ell-\ell_0,0)\leq \ell_1 \leq  \min(\ell,\ell_{0})}} 
      \hskip -1cm \tfrac {2 (k-k_1)}{k} A_{k_1,\ell_1}(y) B_{k-k_1,\ell-\ell_1}(y).
\end{align*}
Therefore, 
\begin{equation}\label{W2quatre}
    W^2 = \sum_{\substack{2\leq k \leq 2 k_0 + 3\\ 0 \leq \ell \leq 2\ell_0+1}}
    c^\ell \left(Q_c^k(y_c)  A_{k,\ell}^*(y)+ (Q_c^k)'(y_c) B_{k,\ell}^*(y)\right),
\end{equation}
where $A_{k,\ell}^*$ and $B_{k,\ell}^*$ can be extracted from the previous formula.

- Terms ${\mathbf{N}}_4$, ${\mathbf{N}}_5$ and  ${\mathbf{N}}_6$.
\begin{align*}
    {\mathbf{N}}_4 
    &=    6 W^2 Q^2(y) \\
    &= \sum_{\substack{2\leq k \leq 2k_0+3\\ 0 \leq \ell \leq 2\ell_0+1}} 
    c^\ell \left(Q_c^k(y_c) ( 6 A_{k,\ell}^* Q^2)(y) + (Q_c^k)'(y_c)    (6 B_{k,\ell}^* Q^2)(y)\right);
\end{align*}
\begin{align*}
    {\mathbf{N}}_5 
    &= Q_c(y_c) W^2 (12 Q(y)) \\
    &= \sum_{\substack{3\leq k \leq 2 k_0 + 4\\ 0 \leq \ell \leq 2 \ell_0+1 }} 
    c^\ell \left(Q_c^k(y_c) ( 12 A_{k-1,\ell}^* Q)(y) + (Q_c^k)'(y_c)    \frac{k{-}1}k  \left(12 B_{k-1,\ell}^* Q\right)(y)\right);
\end{align*}
\begin{align*}
    {\mathbf{N}}_6 &= Q_c^2(y_c) (6 W^2) \\
    &= \sum_{\substack{4\leq k \leq 2 k_0+5\\ 0 \leq \ell \leq 2\ell_0+1}} 
    c^\ell \left(Q_c^k(y_c) ( 4 A_{k-2,\ell}^*(y)) + (Q_c^k)'(y_c)      \frac{k{-}2}k    \left(6 B_{k-2,\ell}^*\right)(y)\right).
\end{align*}

\medskip

For the next two terms, we expand $W^3= W W^2$ using \eqref{W2quatre}:
We get
 \begin{equation}\label{W3quatre}
    W^3 = \sum_{\substack{3\leq k \leq 3 k_0 + 6\\ 0 \leq \ell \leq 3\ell_0+2}}
    c^\ell \left(Q_c^k(y_c)  A_{k,\ell}^{**}(y)+ (Q_c^k)'(y_c) B_{k,\ell}^{**}(y)\right),
\end{equation}
where $A_{k,\ell}^{**}$ and $B_{k,\ell}^{**}$ are explicit in terms of
$A_{k,\ell} $ and $B_{k,\ell}$,
$A_{k,\ell}^{*}$ and $B_{k,\ell}^{*}$.

- Terms ${\mathbf{N}}_7$ and    ${\mathbf{N}}_8$.
\begin{equation*}
    {\mathbf{N}}_7=4 W^3 Q(y)=\sum_{\substack{3\leq k \leq 3 k_0 + 6\\ 0 \leq \ell \leq 3\ell_0+2}}
    c^\ell \left(Q_c^k(y_c)  (4 Q A_{k,\ell}^{**})(y)+ (Q_c^k)'(y_c) (4 Q B_{k,\ell}^{**})(y)\right),
\end{equation*}
\begin{equation*}
    {\mathbf{N}}_8=4 Q_c(y_c) W^3(y)=\sum_{\substack{4\leq k \leq 3 k_0 + 7\\ 0 \leq \ell \leq 3\ell_0+2}}
    c^\ell \left(Q_c^k(y_c)  (4 A_{k-1,\ell}^{**})(y)+ (Q_c^k)'(y_c) (4 B_{k-1,\ell}^{**})(y)\right).
\end{equation*}

- Term ${\mathbf{N}}_9=W^4$. By using $W^4=W^2 W^2$ and \eqref{W2quatre},   we get
\begin{equation*}
    {\mathbf{N}}_9
    = \sum_{\substack{4\leq k \leq 4k_0+9 \\ 0 \leq \ell \leq  4 \ell_0+3}}
    c^{\ell} \left( Q_c^{k}(y_c) A_{k,\ell}^{***}(y)+ (Q_c^{k})'(y_c) B_{k,\ell}^{***}(y) \right),
\end{equation*}
where $A_{k,\ell}^{***}$ and $B_{k,\ell}^{***}$ are explicit in terms of $A_{k,\ell}^{**}$ and $B_{k,\ell}^{**}$.

Next,
\begin{align*}
    {{\mathbf{IV}}} & =\partial_x ({\mathbf{N}})    =    \sum_{\substack{2\leq k \leq 4 k_0+9\\ 0 \leq \ell \leq 4\ell_0+3}}
    c^\ell \Bigl[(Q_c^k)'(y_c)  F^{\mathbf{N}}_{k,\ell}(y) 
    + Q_c^k(y_c) \left((F^{\mathbf{N}}_{k,\ell})'(y)-\beta(y_c) (F^{\mathbf{N}}_{k,\ell})'(y)\right)\\ 
    & \qquad+ (Q_c^k)''(y_c)G^{\mathbf{N}}_{k,\ell}(y)
    + (Q_c^k)'(y_c)\left((G^{\mathbf{N}}_{k,\ell})'(y)-\beta(y_c) (G^{\mathbf{N}}_{k,\ell})'(y)\right) \Bigr].
\end{align*}
Thus,
\begin{align}
    {{\mathbf{IV}}} & =  
    \sum_{\substack{2\leq k \leq 4 k_0+9\\ 0 \leq \ell \leq 4\ell_0+3}} c^\ell Q_c^k(y_c)    (F^{\mathbf{N}}_{k,\ell})'(y)\nonumber \\
    &+ \sum_{\substack{3\leq k \leq 5 k_0+9\\ 0 \leq \ell \leq 5\ell_0+3}} c^\ell Q_c^k(y_c)    
        \hskip -0.5cm \sum_{\substack{\max(k-4 k_0-9,1)\leq k_1 \leq \min(k-2,k_0)\\\max(\ell-4\ell_0-3,0)\leq \ell_1 \leq  \min(\ell,\ell_{0})}}
        \left(- a_{k_1,\ell_1}(F^{\mathbf{N}}_{k-k_1,\ell-\ell_1})'(y)\right)\nonumber \\
    &+ \sum_{\substack{2\leq k \leq 4 k_0+9\\ 1 \leq \ell \leq 4\ell_0+4}} c^\ell Q_c^k(y_c)    k^2 G^{\mathbf{N}}_{k,\ell-1}(y)
        + \sum_{\substack{5\leq k \leq 4 k_0+12\\ 0 \leq \ell \leq 4\ell_0+3}} c^\ell Q_c^k(y_c)    
        \left(-\tfrac{(k{-}3)(2k{-}3)}{5} G^{\mathbf{N}}_{k-3,\ell}(y)\right) \nonumber\\
    &+ \sum_{\substack{2\leq k \leq 4 k_0+9\\ 0 \leq \ell \leq 4\ell_0+3}} c^\ell (Q_c^k)'(y_c)
        \left(F^{\mathbf{N}}_{k,\ell}(y)+(G^{\mathbf{N}}_{k,\ell})'(y)\right) \nonumber\\
    &+   \sum_{\substack{3\leq k \leq 5 k_0+9\\ 0 \leq \ell \leq 5\ell_0+3}} c^\ell (Q_c^k)'(y_c)    
        \hskip -0.5cm \sum_{\substack{\max(k-4 k_0-9,1)\leq k_1 \leq \min(k-2,k_0)\\\max(\ell-4\ell_0-3,0)\leq \ell_1 \leq  \min(\ell,\ell_{0})}}
        \left(- \tfrac{k-k_1}{k} a_{k_1,\ell_1} (G^{\mathbf{N}}_{k-k_1,\ell-\ell_1})'(y) \right).\label{IVform}
\end{align}
It follows that ${\mathbf{IV}}$ can be written as
\begin{equation}\label{IVFG}
    {{\mathbf{IV}}}=\sum_{\substack{2\leq k \leq \max(5k_0+9,4k_0+12)\\ 0 \leq \ell \leq \max(5\ell_0+3,4\ell_0+4) }} 
    c^{\ell} \left( Q_c^{k}(y_c) F_{k,\ell}^{{{\mathbf{IV}}}}(y) +  (Q_c^{k})'(y_c) G_{k,\ell}^{{{\mathbf{IV}}}}(y)\right),
\end{equation}
where $F_{k,\ell}^{{{\mathbf{IV}}}}$ and $G_{k,\ell}^{{{\mathbf{IV}}}}$ can be extracted from the previous calculations.
 Let us check that $F_{k,\ell}^{{{\mathbf{IV}}}}$ and $G_{k,\ell}^{{{\mathbf{IV}}}}$ satisfy properties {\rm (i)} and {\rm (ii)}.

\smallskip

\noindent \emph{Dependence property {\rm (i)}.} In the decomposition of $\mathbf{N}_1$,
the function in factor of $c^\ell Q_c^k$ is $12 A_{k-1,\ell}Q^2$
and the function in factor of $c^\ell (Q_c^k)'$ is $\frac {k-1}k (12 B_{k-1,\ell}Q^2) $. In the decomposition
of $W^2$, in factor of $c^\ell Q_c^k$,
we have sums where $k_1\leq k-1$ and $k-k_1\leq k-1$ since $k_1\geq 1$; moreover $\ell_1\leq \ell$ and $\ell-\ell_1\leq \ell$. Similar remarks apply to the other terms in $\mathbf{N}$.

Thus, $F^{\mathbf{N}}_{k,\ell}$ and $G^{\mathbf{N}}_{k,\ell}$ contain $(A_{k',\ell'})$ and $(B_{k',\ell'})$
only for $(k',\ell')\prec (k,\ell)$ (in fact $k'\leq k-1$ is always true).
From \eqref{IVform} it is clear that the same is true for $F^{\mathbf{IV}}_{k,\ell}$ and $G^{\mathbf{IV}}_{k,\ell}$.
Note in a similar way that when $(a_{k',\ell'})$ is envolved in some formula for $F^{\mathbf{IV}}_{k,\ell}$ and $G^{\mathbf{IV}}_{k,\ell}$
it is only for $(k',\ell')\prec (k,\ell)$.

\smallskip

\noindent \emph{Parity property {\rm (ii)}.} Assume that all $(A_{k',\ell'})$ are even and all $(B_{k',\ell'})$ are odd.
From the decomposition of the various terms of 
$\mathbf{N}$, it is  easy to observe  that all   $(F^{\mathbf{N}}_{k,\ell})$
are even and all $(G^{\mathbf{N}}_{k,\ell})$ are odd.
Then, formula \eqref{IVform} ensures that all $(F^{\mathbf{IV}}_{k,\ell})$
are odd and all $(G^{\mathbf{IV}}_{k,\ell})$ are even.

\medskip

To complete the proof of Lemma \ref{PROPIV}, we only have to compute $F_{2,0}^{{{\mathbf{IV}}}}$ and $G_{2,0}^{{{\mathbf{IV}}}}$.

\noindent By \eqref{IVform}, we have $F_{2,0}^{{{\mathbf{IV}}}}=(F^{\mathbf{N}}_{2,0})'$, and so we are reduced to compute $F^{\mathbf{N}}_{2,0}$.
We give below the contribution of each ${\mathbf{N}}_j$ for $j=1,\ldots,9$ to $F^{\mathbf{N}}_{2,0}$:
\begin{itemize}
    \item For ${\mathbf{N}}_1$, the contribution is $12 A_{1,0} Q^2$;
    \item The contribution of ${\mathbf{N}}_4$ is $6 A_{2,0}^* Q^2=6 A_{1,0}^2 Q^2 $, by the expression of $W^2$;
    \item The contribution of all the other terms ${\mathbf{N}}_2$, ${\mathbf{N}}_3$ ${\mathbf{N}}_5$, ${\mathbf{N}}_6$, ${\mathbf{N}}_7$, ${\mathbf{N}}_8$ and ${\mathbf{N}}_9$ is zero.
\end{itemize}
Therefore,  
$
    F_{2,0}^{{{\mathbf{IV}}}} =  (F^{\mathbf{N}}_{2,0})' = \left(12 A_{1,0} Q^2 + 6 A_{1,0}^2 Q^2\right)'.
$

\smallskip

By \eqref{IVform}, we have $G_{2,0}^{{{\mathbf{IV}}}}=F^{\mathbf{N}}_{2,0}+(G^{\mathbf{N}}_{2,0})'$. Since $F^{\mathbf{N}}_{2,0}$ was computed above, we are reduced to compute $G^{\mathbf{N}}_{2,0}$.
We give below the contribution of each ${\mathbf{N}}_j$ for $j=1,\ldots,9$ to $G^{\mathbf{N}}_{2,0}$:
\begin{itemize}
    \item For ${\mathbf{N}}_1$, the contribution is $6 B_{1,0} Q^2$;
    \item The contribution of ${\mathbf{N}}_4$ is $6 B_{2,0}^* Q^2=6 A_{1,0}B_{1,0} Q^2 $, by the expression of $W^2$;
    \item The contribution of all the other terms $\mathbf{N}_2$, $\mathbf{N}_3$, ${\mathbf{N}}_5$, ${\mathbf{N}}_6$, ${\mathbf{N}}_7$, ${\mathbf{N}}_8$ and ${\mathbf{N}}_9$ is zero.
\end{itemize}
Therefore, 
$
    G_{2,0}^{{{\mathbf{IV}}}} =  12 A_{1,0} Q^2 + 6 A_{1,0}^2 Q^2 + \left(6 B_{1,0} Q^2 + 6 A_{1,0} B_{1,0} Q^2\right)'.
$

\subsection{End of the proof of Proposition \ref{SYSTEME}}

By  Lemmas \ref{PROPI}--\ref{PROPIV}, we only have to sum the various contributions of $\mathbf{I}$, $\mathbf{II}$, $\mathbf{III}$ and $\mathbf{IV}$ to prove Proposition \ref{SYSTEME}. Setting
\begin{equation*}
F_{k,\ell}=F_{k,\ell}^{\mathbf{I}}+F_{k,\ell}^{\mathbf{II}}+F_{k,\ell}^{\mathbf{III}}+F_{k,\ell}^{\mathbf{IV}},
\quad\text{and}\quad
G_{k,\ell}=G_{k,\ell}^{\mathbf{I}}+G_{k,\ell}^{\mathbf{II}}+G_{k,\ell}^{\mathbf{III}}+G_{k,\ell}^{\mathbf{IV}},
\end{equation*}
we obtain the formula of    Proposition \ref{SYSTEME} for $S(t,x)$.
Properties {\rm (i)} and {\rm (ii)} have been checked on the functions $F_{k,\ell}^{\mathbf{I}},$ $F_{k,\ell}^{\mathbf{II}}$, $F_{k,\ell}^{\mathbf{III}}$, $F_{k,\ell}^{\mathbf{IV}}$
and $G_{k,\ell}^{\mathbf{I}}$, $G_{k,\ell}^{\mathbf{II}}$, $G_{k,\ell}^{\mathbf{III}}$, $G_{k,\ell}^{\mathbf{IV}}$, and so they are also true on $F_{k,\ell}$ and $G_{k,\ell}$.

The expressions  of  $F_{1,0}$, $G_{1,0}$, $F_{2,0}$ and $G_{2,0}$ are obtained from    Lemmas \ref{PROPI}--\ref{PROPIV}. Observe that the only nonzero contribution to $F_{1,0}$ and $G_{1,0}$ comes from $F_{1,0}^{\mathbf{II}}$ and $G_{1,0}^{\mathbf{II}}$,
 we obtain $F_{1,0}=p(Q^{p-1})'$ and $G_{1,0}=pQ^{p-1}$.

\section{Appendix --  Lemma \ref{STRUCT4}}

\begin{lemma}[Structure of $F_{k,\ell}$ and $G_{k,\ell}$]\label{STRUCT4}
    Let $(k,\ell)$ be such that $1\leq k \leq K_0$, $0\leq \ell\leq L_0$, with $(k,\ell)\neq (1,0)$.
    Assume that for all $1\leq k'\leq k_0$, $0\leq \ell'\leq \ell_0$ such that $(k',\ell')\prec (k,\ell)$,
    the functions $A_{k',\ell'}$ and $B_{k',\ell'}$ verify
    \begin{equation}
        A_{k',\ell'}   = \overline{A}_{k',\ell'}  + \widetilde{A}_{k',\ell'}  + \varphi  \widehat{A}_{k',\ell'}, \quad
        B_{k',\ell'}  = \overline{B}_{k',\ell'}  + \widetilde{B}_{k',\ell'}  + \varphi  \widehat{B}_{k',\ell'} ,\label{ABkl2Lemme}
    \end{equation}
    \begin{itemize}
        \item $\overline{A}_{k',\ell'}$, $\overline{B}_{k',\ell'}\in \mathcal{Y}$;
         the function $\overline{A}_{k',\ell'}$ is even and the function $\overline{B}_{k',\ell'}$ is odd;
        \item $\widetilde{A}_{k',\ell'}$ and $\widehat{B}_{k',\ell'}$ are even polynomials;
  $\widehat{A}_{k',\ell'}$ and $\widetilde{B}_{k',\ell'}$ are odd polynomials.
    \end{itemize}
    Then the functions $F_{k,\ell}$ and $G_{k,\ell}$ obtained in Proposition \ref{SYSTEME} from $(a_{k',\ell'})$, $(A_{k',\ell'})$ and
    $(B_{k',\ell'})$ are such that
    \begin{equation*}
        F_{k,\ell}    = \overline{F}_{k,\ell}  + \widetilde{F}_{k,\ell}  + \varphi  \widehat{F}_{k,\ell} ,
        \quad
        G_{k,\ell}    = \overline{G}_{k,\ell}  + \widetilde{G}_{k,\ell}  + \varphi  \widehat{G}_{k,\ell},
    \end{equation*}
 \begin{itemize}
        \item $\overline{F}_{k,\ell}$, $\overline{G}_{k,\ell}\in \mathcal{Y}$;
         the function $\overline{F}_{k,\ell}$ is odd and the function $\overline{G}_{k,\ell}$ is even;
        \item $\widetilde{F}_{k,\ell}$ and $\widehat{G}_{k,\ell}$ are odd polynomials;
       $\widehat{F}_{k,\ell}$ and $\widetilde{G}_{k,\ell}$ are even polynomials.
    \end{itemize}
    Moreover, the following hold.
    \begin{enumerate}
        \item[{\rm (a)}] Let $2\leq k\leq p-1$, $\ell=0$. If for any $1\leq k' <k$,
        \begin{equation*}
            \deg \widetilde A_{k',0}=\deg \widehat A_{k',0}=\deg \widetilde B_{k',0}=\deg \widehat B_{k',0}=0 
            \quad \text{then}\quad F_{k,0}, \ G_{k,0}\in \mathcal{Y}.
        \end{equation*}
        \item[{\rm (b)}] Let $1\leq k \leq k_0$ and  $0\leq \ell \leq \ell_0$ be such that $\frac k{p-1} + \ell\leq 2$.
                         If for any $(k',\ell')\prec (k,\ell)$,
        \begin{equation*}
            \deg \widetilde A_{k',\ell'}=\deg \widehat A_{k',\ell'}=0 \quad\text{and}\quad
            \deg \widetilde B_{k',\ell'}=\deg \widehat B_{k',\ell'}\leq 1 
            \quad \text{then}\quad F_{k,\ell}\in \mathcal{Y}.
        \end{equation*}
        \item[{\rm (c)}] Let 
        \begin{align*}
            d_{AB}(k,\ell) & =\left\{
            \begin{aligned}
                & \max_{\substack{1\leq k'\leq k \\ 0\leq \ell' \leq \ell}}
                \left(\deg \widetilde A_{k',\ell'},\deg \widehat A_{k',\ell'},
                    \deg \widetilde B_{k',\ell'},\deg \widehat B_{k',\ell'}\right) \text{ if $k\geq 1$, $\ell\geq 0$,}\\
                & 0 \text{ otherwise,}
            \end{aligned}\right.
        \\
            d_{\mathbf{N}}(k,\ell) & =\left\{
            \begin{aligned}
                & \max_{\substack{1\leq k_j\leq k_0 \\ 0\leq \ell_j \leq \ell_0}}
                 \biggl( \,\sum_{j=1}^p d_{AB}(k_j,\ell_j)\text{ for } \sum_{j=1}^p k_j\leq k,  \sum_{j=1}^p \ell_j\leq \ell\biggr)
                \text{ if $k\geq p$, $\ell\geq 0$,}  \\
                & 0 \text{ otherwise,}
            \end{aligned}\right.
        \\
            d_{FG}(k,\ell)  & =
                    \max(\deg \widetilde F_{k,\ell},\deg \widehat F_{k,\ell},
                        \deg \widetilde G_{k,\ell},\deg \widehat G_{k,\ell})    \text{ for $1\leq k\leq K_0$, $0\leq \ell\leq L_0$.}
        \end{align*}
        Then, for all $1\leq k\leq K_0$, $0\leq \ell\leq L_0$,
        \begin{equation}\label{estDEGFG}
            d_{FG}(k,\ell)  \leq \max \big( 
            d_{AB}(k{-}1, \ell){-}1, d_{AB}(k {-}p+{1}, \ell), d_{AB}(k, \ell{-}1), d_{\mathbf{N}}(k,\ell)\big).
        \end{equation}
    \end{enumerate}
\end{lemma}

\noindent\emph{Proof of Lemma \ref{STRUCT4}.}
Let $(k,\ell)$ be such that $1\leq k\leq K_0$ and $0\leq \ell\leq L_0$, with $(k,\ell)\neq (1,0)$.
We suppose that for all $1\leq k'\leq k_0$, $0\leq \ell\leq \ell_0$ such that $(k',\ell')\prec (k,\ell)$, $(a_{k',\ell'},A_{k',\ell'},B_{k',\ell'})$ satisfies the assumptions of Lemma \ref{STRUCT4}. 
We consider  $F_{k,\ell}$, $G_{k,\ell}$ defined by Proposition \ref{SYSTEME}
(recall that for given $(k,\ell)$, $F_{k,\ell}$ and $G_{k,\ell}$ depend only on $(a_{k',\ell'},
A_{k',\ell'},B_{k',\ell'})$ for $(k',\ell')\prec (k,\ell)$).
From  the proof of Proposition \ref{SYSTEME} (Appendix A), 
\begin{equation*}
    F_{k,\ell}=F_{k,\ell}^{\mathbf{I}} +F_{k,\ell}^{\mathbf{II}} +F_{k,\ell}^{\mathbf{III}}+F_{k,\ell}^{\mathbf{IV}},
    \quad G_{k,\ell}=G_{k,\ell}^{\mathbf{I}}+G_{k,\ell}^{\mathbf{II}}+G_{k,\ell}^{\mathbf{III}}+G_{k,\ell}^{\mathbf{IV}},
\end{equation*}
where $F_{k,\ell}^{\mathbf{I}}$, $F_{k,\ell}^{\mathbf{II}}$, etc.  are the contributions of $\mathbf{I}$, $\mathbf{II}$, $\mathbf{III}$ and $\mathbf{IV}$ in the decomposition of $S(t,x)$, see \eqref{decS}.

\medskip

\noindent\emph{- Contribution of $\mathbf{I}$ and $\mathbf{II}$.} \quad
From Lemmas \ref{PROPI} and \ref{PROPII},    it follows that $F_{k,\ell}^{\mathbf{I}}$, $F_{k,\ell}^{\mathbf{II}}$, $G_{k,\ell}^{\mathbf{I}}$ and $G_{k,\ell}^{\mathbf{II}}$ belong to $\mathcal{Y}$ and do not depend on $(A_{k',\ell'})$, $(B_{k',\ell'})$ but 
only on the coefficients $(a_{k',\ell'})$.
Moreover,
$F_{k,\ell}^{\mathbf{I}}$ and    $F_{k,\ell}^{\mathbf{II}}$ are odd, and $G_{k,\ell}^{\mathbf{I}}$ and $G_{k,\ell}^{\mathbf{II}}$ are even.
Therefore, they only contribute to $\overline F_{k,\ell}$ and $\overline G_{k,\ell}$, with the desired parity property.

\medskip

\noindent\emph{- Contribution of $\mathbf{III}$.} \quad
We use the notation and calculations of the proof of Lemma \ref{PROPIII}. Note that $\mathbf{III}_1$ does not contribute to $F_{k,\ell}^{\mathbf{III}}$ and $G_{k,\ell}^{\mathbf{III}}$. Observing the other terms, i.e.
$\mathbf{III}_2$, $\mathbf{III}_3^1$, $\mathbf{III}_3^2$, etc. up to $\mathbf{III}_{13}^5$, we note that there are three kinds of terms
depending on the structure of the function of the    variable $y$:
\begin{itemize}
    \item[$T_1$:] Terms depending on $A_{k',\ell'}(y)$ and  $B_{k',\ell'}(y)$ without derivative, for $(k',\ell')\prec (k,\ell)$. A complete list of these terms is given in formula \eqref{listT1} below.
    \item[$T_2$:] Terms depending on derivatives of $A_{k',\ell'}(y)$ and    $B_{k',\ell'}(y)$ (up to order $3$) for $(k',\ell')\prec (k,\ell)$. Examples of such terms are $F_{k,\ell}^{\mathbf{III}_4}$, $G_{k,\ell}^{\mathbf{III}_4}$, a part of $F_{k,\ell}
^{\mathbf{III}_3}$, etc.
    \item[$T_3$:] Terms depending on $(Q^{p-1}A_{k',\ell'})'(y)$ and    $(Q^{p-1}B_{k',\ell'})'(y)$ for $(k',\ell')\prec (k,\ell)$.
    Examples of such terms are a part of $F_{k,\ell}^{\mathbf{III}_3}$, $G_{k,\ell}^{\mathbf{III}_3}$, etc.
\end{itemize}

Terms of type $T_3$ are easily handled. Indeed, since $A_{k',\ell'}$ and $B_{k',\ell'}$ are of the form  \eqref{ABkl2Lemme} and since $Q\in \mathcal{Y}$, it follows that $(Q^{p-1}A_{k',\ell'})'$ and $(Q^{p-1}B_{k',\ell'})'$ belong to $\mathcal{Y}$.
Therefore, this kind of terms only contribute to $\overline F_{k,\ell}$ and $\overline G_{k,\ell}$.
The parity statement for these terms was already checked in the proof of Lemma \ref{PROPIII}.

\smallskip

We now handle terms of type $T_2$. It suffices to remark that when differentiating terms such as $A_{k',\ell'}$ and $B_{k',\ell'}$ of the form  \eqref{ABkl2Lemme}, we obtain terms of the same form, except that the degrees of the polynomial functions decrease by one or more depending on the order of derivation. Indeed, for example, it follows from \eqref{ABkl2Lemme} that:
\begin{equation*}
    A_{k',\ell'}'=(\overline A_{k',\ell'}' + \varphi' \widehat A_{k',\ell'}) +\widetilde A_{k',\ell'}'+\varphi \widehat A_{k',\ell'}',
\end{equation*}
and $\overline A_{k',\ell'}' + \varphi' \widehat A_{k',\ell'}\in \mathcal{Y}$, because of the property $\varphi'\in \mathcal{Y}$.
Thus, for example, we get:
$
    F_{k,\ell}^{\mathbf{III}_4}=\overline F_{k,\ell}^{\mathbf{III}_4} + \widetilde F_{k,\ell}^{\mathbf{III}_4} +\varphi \widehat F_{k,\ell}^{\mathbf{III}_4},
$ where
\begin{align*}
    & \deg \widetilde F_{k,\ell}^{\mathbf{III}_4}\leq \max_{(k',\ell')\prec (k,\ell)}(\deg \widehat A_{k',\ell'}) -1
        \leq \max(d_{AB}(k{-}1,\ell),d_{AB}(k,\ell{-}1)) -1,
    \\ & \deg \widehat F_{k,\ell}^{\mathbf{III}_4}\leq      \max_{(k',\ell')\prec (k,\ell)}(\deg \widetilde A_{k',\ell'}) -1
            \leq \max(d_{AB}(k{-}1,\ell),d_{AB}(k,\ell{-}1)) ) -1 ,
\end{align*}
if $\max(d_{AB}(k{-}1,\ell),d_{AB}(k,\ell{-}1)) )\geq 1$, and 
$\widetilde{F}_{k,\ell}^{\mathbf{III}_4}=\widehat{F}_{k,\ell}^{\mathbf{III}_4}=0$
otherwise.

We obtain similar estimates for all terms of this type. The parity properties are easily checked. For terms of type $T_2$ with higher order derivatives (in fact, only second and third derivative), the argument is the same.

\smallskip

Finally, we look at terms of type $T_1$, i.e. depending on $A_{k',\ell'}$ and    $B_{k',\ell'}$ without derivative:
\begin{equation}\label{listT1}
    G_{k,\ell}^{\mathbf{III}_2}, \quad G_{k,\ell}^{\mathbf{III}_{10}^1}, \quad G_{k,\ell}^{\mathbf{III}_{10}^2},\quad    F_{k,\ell}^{\mathbf{III}_{11}^1}, \quad F_{k,\ell}^{\mathbf{III}_{11}^2},\quad  F_{k,\ell}^{\mathbf{III}_{13}^3}, \quad F_{k,\ell}^{\mathbf{III}_{13}^4},\quad  F_{k,\ell}^{\mathbf{III}_{13}^5}.
\end{equation}
With the assumptions on $A_{k',\ell'}$ and $B_{k',\ell'}$, these terms have the desired structure. We only have to check the estimates on the degrees of the polynomials.

First, we note from the proof of Lemma \ref{PROPIII} that terms 
$
    G_{k,\ell}^{\mathbf{III}_2},$ $G_{k,\ell}^{\mathbf{III}_{10}^1},$ $F_{k,\ell}^{\mathbf{III}_{11}^1}, $ $ F_{k,\ell}^{\mathbf{III}_{11}^2},$ $F_{k,\ell}^{\mathbf{III}_{13}^3},$ $F_{k,\ell}^{\mathbf{III}_{13}^4}
$
depend only on $A_{k',\ell'}$ and $B_{k',\ell'}$ with $k'\leq k$ and $\ell'\leq \ell-1$. Thus, they appear only for $\ell\geq 1$ and contain polynomials with degrees less than or equal to $d_{AB}(k,\ell-1)$.
The other two terms $G_{k,\ell}^{\mathbf{III}_{10}^2}$ and $F_{k,\ell}^{\mathbf{III}_{13}^5}$ depend on $A_{k',\ell'}$ and $B_{k',\ell'}$ with $k'\leq k-p+1$ and $\ell'\leq \ell$. Thus they appear only for $k\geq p$, and    contain polynomials with degrees less than or equal to $d_{AB}(k-p+1,\ell)$.
 
Thus, in conclusion for the term $\mathbf{III}$, we get polynomials of degrees less than
\begin{equation*}
 \max \left( d_{AB}(k{-}1, \ell){-}1, d_{AB}(k {-}p+{1}, \ell), d_{AB}(k, \ell{-}1)\right).    
\end{equation*}
This proves (c) for $d_{FG}^{\mathbf{III}}$.

\medskip

Let us now prove (a) and (b) for $F_{k,\ell}^{\mathbf{III}}$ and $G_{k,\ell}^{\mathbf{III}}$.

Proof of (a).
First, observe that terms of type $T_1$ (see above) do not appear for $k\leq p-1$ and $\ell=0$.
Thus, for such $k$, if we assume $\widetilde A_{k',0} =\widehat A_{k',0} =\widetilde B_{k',0} =0$ and $\widehat B_{k',0}=b_{k',0}\in \mathbb{R}$ then $\widetilde A_{k',0}'=\widehat A_{k',0}'=\widetilde B_{k',0}'=\widehat B_{k',0}'=0$
for all $1\leq k'<k$, and so
$\widetilde{F}_{k,\ell}^{\mathbf{III}}=\widehat{F}_{k,\ell}^{\mathbf{III}}
=\widetilde{G}_{k,\ell}^{\mathbf{III}}=\widehat{G}_{k,\ell}^{\mathbf{III}}=0$, which means $F_{k,\ell}^{\mathbf{III}},G_{k,\ell}^{\mathbf{III}}\in \mathcal{Y}.$ This proves (a) for $F_{k,\ell}^{\mathbf{III}}$ and $G_{k,\ell}^{\mathbf{III}}$.

Proof of (b).
To justify (b) for $F_{k,\ell}^{\mathbf{III}}$, we first observe that for $(k,\ell)$ such that $\frac {k}{p-1}+\ell\leq 2$, there is no term of type $T_1$ contributing to $F_{k,\ell}^{\mathbf{III}}$. Indeed, looking at the expression of all the terms in the list \eqref{listT1} in the proof of Proposition \ref{SYSTEME}, we see that 
$  F_{k,\ell}^{\mathbf{III}_{11}^1},$ $F_{k,\ell}^{\mathbf{III}_{11}^2},$ $F_{k,\ell}^{\mathbf{III}_{13}^3},$ $F_{k,\ell}^{\mathbf{III}_{13}^4},$ $F_{k,\ell}^{\mathbf{III}_{13}^5}$
involves $B_{k',\ell'}$ for $k'\leq k-2(p-1)$ or $\ell'\leq \ell-2$ or simultaneously $k'\leq k-(p-1)$ and $\ell'\leq \ell-1$.
Therefore, these terms do not appear if $\frac {k}{p-1}+\ell\leq 2$. Concerning terms of type $T_2$, we first note that $B_{k',\ell'}$ appear with at least two derivatives, thus any polynomial function    of degree $1$ disappears. Second, $A_{k',\ell'}$ a
re differentiated at least once, and so again any constant term disappears. Thus, there remains no polynomial growth and
$F_{k,\ell}^{\mathbf{III}}\in \mathcal{Y}$ for such $k$, $\ell$.

\medskip

\noindent\emph{- Contribution of $\mathbf{IV}$.}\quad
We focus on the case $p=4$. The other cases, i.e. $p=2$ or $3$ are similar and easier.
We use the notation and calculations of the proof of Lemma \ref{PROPIV}, where we have written
$\mathbf{IV}=\partial_x(\mathbf{N}),$ $\mathbf{N}=(R+R_c+W)^4 - (R+R_c)^4 - 4 R^3 W,$
and where we have decomposed $\mathbf{N}$ into several parts $\mathbf{N}_1,\ldots,\mathbf{N}_9$.
Here, we distinguish two kind of terms: first
$\mathbf{N}_1,\ \mathbf{N}_2,\ \mathbf{N}_4,\ \mathbf{N}_5,\ \mathbf{N}_7,$
which contain the function $Q(y)$, and second,
$    \mathbf{N}_3,\ \mathbf{N}_6,\ \mathbf{N}_8,\   \mathbf{N}_9,$
which depend only on $Q_c$ and $W$.

For the first terms, $\mathbf{N}_1,\ \mathbf{N}_2,\ \mathbf{N}_4,\ \mathbf{N}_5$ and $\mathbf{N}_7,$ since $Q\in \mathcal{Y}$, by the structure of $W$, and the assumptions on $A_{k',\ell'}$ and $B_{k',\ell'}$, the result follows.

For $\mathbf{N}_3,\ \mathbf{N}_6,\ \mathbf{N}_8$ and $\mathbf{N}_9$, we set
\begin{equation*}
    \mathbf{M}=\mathbf{N}_3+ \mathbf{N}_6+\mathbf{N}_8+\mathbf{N}_9=(Q_c+W)^4-Q_c^4.
\end{equation*}
In order to have a simple expression when expanding $(Q_c+W)^4$, it is convenient to set
\begin{equation*}
    \mathbf{A}_{1,0}=1+A_{1,0},\quad 
    \mathbf{A}_{k,\ell}=A_{k,\ell}, \text{ for any $(k,\ell)\neq (1,0)$},\quad 
    \mathbf{B}_{k,\ell}=B_{k,\ell}, \text{ for any $(k,\ell)$}.
\end{equation*}
\begin{equation}\label{DEGbf}
    \deg \mathbf{A}_{k,\ell}=\deg A_{k,\ell},\quad   \deg \mathbf{B}_{k,\ell}=\deg B_{k,\ell}.
\end{equation}
With this notation, we have
$
    Q_c+W= \sum_{(k,\ell)\in \Sigma_0}
    c^\ell\left( Q_c^k(y_c)  \mathbf{A}_{k,\ell}(y)     +    (Q_c^k)'(y_c) \mathbf{B}_{k,\ell}(y)\right).
$
Then, 
{\allowdisplaybreaks
\begin{align} 
     (Q_c+W)^4  = \sum_{\substack{(k_j,\ell_j)\in \Sigma_0   \\ j=1,2,3,4 }} &
    c^{\ell_1+\ell_2+\ell_3+\ell_4} \Biggl\{ Q_c^{k_1+k_2+k_3+k_4}(y_c) \ (\mathbf{A}_{k_1,\ell_1} \mathbf{A}_{k_2,\ell_2} \mathbf{A}_{k_3,\ell_3} \mathbf{A}_{k_4,\ell_4})(y)\nonumber \\
    &   +4 ((Q_c^{k_1})'Q_c^{k_2+k_3+k_4})(y_c)   \ (\mathbf{B}_{k_1,\ell_1} \mathbf{A}_{k_2,\ell_2} \mathbf{A}_{k_3,\ell_3} \mathbf{A}_{k_4,\ell_4})(y)\nonumber\\
    &     +6 ((Q_c^{k_1})'(Q_c^{k_2})'Q_c^{k_3+k_4})(y_c) \   (\mathbf{B}_{k_1,\ell_1} \mathbf{B}_{k_2,\ell_2} \mathbf{A}_{k_3,\ell_3} \mathbf{A}_{k_4,\ell_4})(y)\label{closer}\\
    &     +4 ((Q_c^{k_1})'(Q_c^{k_2})'(Q_c^{k_3})'Q_c^{k_4})(y_c)\    (\mathbf{B}_{k_1,\ell_1} \mathbf{B}_{k_2,\ell_2} \mathbf{B}_{k_3,\ell_3} \mathbf{A}_{k_4,\ell_4})(y)\nonumber\\
    &     +((Q_c^{k_1})'(Q_c^{k_2})'(Q_c^{k_3})'(Q_c^{k_2})')(y_c) \ (\mathbf{B}_{k_1,\ell_1} \mathbf{B}_{k_2,\ell_2} \mathbf{B}_{k_3,\ell_3} \mathbf{B}_{k_4,\ell_4})(y)\Biggr\}.\nonumber
\end{align}
}
Recall that by Lemma \ref{surQc}, we have
\begin{equation*}
     (Q_c^{k_1})'Q_c^{k_2+k_3+k_4}=\tfrac{k_1}{k_1+k_2+k_3+k_4} (Q_c^{k_1+k_2+k_3+k_4})',
\end{equation*}
\begin{equation*}
     (Q_c^{k_1})'(Q_c^{k_2})'Q_c^{k_3+k_4}  =
     {k_1 k_2}  \left(c Q_c^{k_1+k_2+k_3+k_4} -\tfrac {2}{p+1} Q_c^{k_1+k_2+k_3+k_4+3} \right),
\end{equation*}
\begin{equation*}
     (Q_c^{k_1})'(Q_c^{k_2})'(Q_c^{k_3})'Q_c^{k_4} =
     {k_1 k_2 k_3} \left( \frac{c (Q_c^{k_1+k_2+k_3+k_4})'}{k_1{+}k_2{+}k_3{+}k_4}  -\frac {2 (Q_c^{k_1+k_2+k_3+k_4+3})'}{(p+1)(k_1{+}k_2{+}k_3{+}k_4{+}3)}  \right),
\end{equation*}
\begin{equation*}
     (Q_c^{k_1})'(Q_c^{k_2})'(Q_c^{k_3})'(Q_c^{k_4})'  =  k_1 k_2 k_3 k_4  Q_c^{k_1+k_2+k_3+k_4} \left(c^2  -\tfrac {4c}{p+1} Q_c^{ 3} +\tfrac{4}{(p+1)^2}Q_c^{ 6}\right).
\end{equation*}
Therefore, we can write:
\begin{equation}\label{expM}
    \mathbf{M}=\sum_{\substack{4\leq k\leq 4 k_{0} + 6  \\ 0 \leq \ell \leq  4 \ell_{0} +2 }}
    c^\ell\left( Q_c^k(y_c)  F_{k,\ell}^{\mathbf{M}}(y)     +    (Q_c^k)'(y_c)  G_{k,\ell}^{\mathbf{M}}(y)\right),
\end{equation}
where at given $k\geq 4$, $\ell\geq 0$, $F_{k,\ell}^{\mathbf{M}}$ contains only terms of the type
\begin{equation}\label{MFkl}
\mathbf{A}_{k_1,\ell_1} \mathbf{A}_{k_2,\ell_2} \mathbf{A}_{k_3,\ell_3} \mathbf{A}_{k_4,\ell_4},\quad
    \mathbf{B}_{k_1,\ell_1} \mathbf{B}_{k_2,\ell_2} \mathbf{A}_{k_3,\ell_3} \mathbf{A}_{k_4,\ell_4},\quad
    \mathbf{B}_{k_1,\ell_1} \mathbf{B}_{k_2,\ell_2} \mathbf{B}_{k_3,\ell_3} \mathbf{B}_{k_4,\ell_4},
\end{equation}
for $\sum_{j=1}^4 k_j\leq k$ and $\sum_{j=1}^4 \ell_j\leq \ell$, and 
$G_{k,\ell}^{\mathbf{M}}$ contains only terms of the type
\begin{equation}\label{MGkl}
    \mathbf{B}_{k_1,\ell_1} \mathbf{A}_{k_2,\ell_2} \mathbf{A}_{k_3,\ell_3} \mathbf{A}_{k_4,\ell_4},\quad
    \mathbf{B}_{k_1,\ell_1} \mathbf{B}_{k_2,\ell_2} \mathbf{B}_{k_3,\ell_3} \mathbf{A}_{k_4,\ell_4},
\end{equation}
for $\sum_{j=1}^4 k_j\leq k$ and $\sum_{j=1}^4 \ell_j\leq \ell$.

Therefore, we only have to check the structure of the functions in \eqref{MFkl} and \eqref{MGkl}. We check the first term
$\mathbf{A}_{k_1,\ell_1} \mathbf{A}_{k_2,\ell_2} \mathbf{A}_{k_3,\ell_3} \mathbf{A}_{k_4,\ell_4}$, the other terms can be checked similarly.

Recall that $\mathbf{A}_{k_j,\ell_j}=\overline{\mathbf{A}}_{k_j,\ell_j} + \widetilde{\mathbf{A}}_{k_j,\ell_j}+ \varphi  \widehat{\mathbf{A}}_{k_j,\ell_j}$,
where $\overline{\mathbf{A}}_{k_j,\ell_j}\in \mathcal{Y}$, and $\widetilde{\mathbf{A}}_{k_j,\ell_j}$ and $\widehat{\mathbf{A}}_{k_j,\ell_j}$ are polynomials.
In the product $\mathbf{A}_{k_1,\ell_1} \mathbf{A}_{k_2,\ell_2} \mathbf{A}_{k_3,\ell_3} \mathbf{A}_{k_4,\ell_4}$, any term in factor to some $\overline{\mathbf{A}}_{k_j,\ell_j}$ is automatically in $\mathcal{Y}$. The other terms are:
\begin{equation*}
    (\widetilde{\mathbf{A}}_{k_1,\ell_1}+\varphi\widehat{\mathbf{A}}_{k_1,\ell_1}) (\widetilde{\mathbf{A}}_{k_2,\ell_2}+\varphi\widehat{\mathbf{A}}_{k_2,\ell_2})
     (\widetilde{\mathbf{A}}_{k_3,\ell_3}+\varphi\widehat{\mathbf{A}}_{k_3,\ell_3}) (\widetilde{\mathbf{A}}_{k_4,\ell_4}+\varphi\widehat{\mathbf{A}}_{k_4,\ell_4}).
\end{equation*}
In this product we distinguish two  kinds of terms:
\begin{itemize} 
    \item    $\Pi_{j=1}^4 \widetilde{\mathbf{A}}_{k_j,\ell_j}$, $\widetilde{\mathbf{A}}_{k_1,\ell_1}\widetilde{\mathbf{A}}_{k_2,\ell_2} (\varphi^2    \widehat{\mathbf{A}}_{k_3,\ell_3}
    \widehat{\mathbf{A}}_{k_4,\ell_4} )$
    (and similar terms), $\varphi^4 \Pi_{j=1}^4 \widehat{\mathbf{A}}_{k_j,\ell_j}$. 
    Since $1-\varphi^2, 1-\varphi^4\in \mathcal{Y}$, these terms are of the form $\overline F+\widetilde F$, where $\overline F\in \mathcal{Y}$ is even and $\widetilde F$
    is  an even polynomial of degree less than or equal to
    $d_{\mathbf{N}}(k,\ell)$.
    \item    $\Pi_{j=1}^3 \widetilde{\mathbf{A}}_{k_j,\ell_j} (\varphi \widehat{\mathbf{A}}_{k_4,\ell_4} )$,
    $\widetilde{\mathbf{A}}_{k_1,\ell_1} (\varphi^3 \Pi_{j=2}^4 \widehat{\mathbf{A}}_{k_j,\ell_j} )$ (and similar terms). Since $\varphi^3-\varphi\in \mathcal{Y}$, these terms 
    are  of the form $\overline F + \varphi \widehat F $, where $\overline F\in \mathcal{Y}$ and $\widehat F$ is a polynomial function of  degree are less than
$d_{\mathbf{N}}(k,\ell)$.
\end{itemize}

In conclusion, we obtain
\begin{equation*}
    F_{k,\ell}^{\mathbf{M}}  = \overline{F}_{k,\ell}^{\mathbf{M}} + \widetilde{F}_{k,\ell}^{\mathbf{M}}  + \varphi  \widehat{F}_{k,\ell}^{\mathbf{M}} ,\quad 
    G_{k,\ell}^{\mathbf{M}}  = \overline{G}_{k,\ell}^{\mathbf{M}}  + \widetilde{G}_{k,\ell}^{\mathbf{M}}  + \varphi  \widehat{G}_{k,\ell}^{\mathbf{M}} ,
\end{equation*}
\begin{itemize}
    \item $\overline{F}_{k,\ell}^{\mathbf{M}}$, $\overline{G}_{k,\ell}^{\mathbf{M}}\in \mathcal{Y}$;
         $\overline{F}_{k,\ell}^{\mathbf{M}}$ is even and $\overline{G}_{k,\ell}^{\mathbf{M}}$ is odd;
    \item $\widetilde{F}_{k,\ell}^{\mathbf{M}}$ and $\widehat{G}_{k,\ell}^{\mathbf{M}}$ are even polynomials;
$\widehat{F}_{k,\ell}^{\mathbf{M}}$ and $\widetilde{G}_{k,\ell}^{\mathbf{M}}$ are odd polynomials,
satisfying
\end{itemize}
\begin{equation}
\label{estDEGFGM}
    d_{FG}^{\mathbf{M}}(k,\ell)=\max(\deg \widetilde F_{k,\ell}^{\mathbf{M}},\deg \widehat 	F_{k,\ell}^{\mathbf{M}},\deg \widetilde G_{k,\ell}^{\mathbf{M}},\deg \widehat G_{k,\ell}^{\mathbf{M}}),      
 \leq d_{\mathbf{N}}(k,\ell).    
\end{equation}

The last step for    $\mathbf{IV}$ is to use formulas \eqref{IVform} and \eqref{IVFG} to derive the properties of $F_{k,\ell}^{\mathbf{IV}}$ and $G_{k,\ell}^{\mathbf{IV}}$
from the properties of $F_{k,\ell}^{\mathbf{N}}$ and $G_{k,\ell}^{\mathbf{N}}$. We note that $F_{k,\ell}^{\mathbf{IV}}$
involves some $G_{k',\ell'}^{\mathbf{N}}$ and $(F_{k',\ell'}^{\mathbf{N}})'$ for $k'\leq k$ and $\ell'\leq \ell$    and
$G_{k,\ell}^{\mathbf{IV}}$
involves some $(G_{k',\ell'}^{\mathbf{N}})'$ and $F_{k',\ell'}^{\mathbf{N}}$ for $k'\leq k$ and $\ell'\leq \ell$. Thus
${\mathbf{IV}}$ contains polynomials with degrees less than $d_{\mathbf{N}}(k,\ell)$,
and the parity properties are satisfied,
which  proves (c) for $d_{FG}^{\mathbf{IV}}$.

\medskip

Let us now prove (a) and (b) for $F_{k,\ell}^{\mathbf{IV}}$ and $G_{k,\ell}^{\mathbf{IV}}$.

Proof of (a).
Note that from \eqref{closer}--\eqref{expM}, for any $1\leq k\leq p-1=3$, $F_{k,0}^{\mathbf{M}}=G_{k,0}^{\mathbf{M}}=0$.
Thus, $F_{k,0}^{\mathbf{N}}, G_{k,0}^{\mathbf{N}}\in \mathcal{Y}$ for such $k$. From \eqref{IVform} and \eqref{IVFG}
it follows that $F_{k,0}^{\mathbf{IV}}, G_{k,0}^{\mathbf{IV}}\in \mathcal{Y}$ for all $1\leq k\leq p-1$. This proves (a) for
term $\mathbf{IV}$.

Proof of (b).
To prove (b) for $F_{k,\ell}^{\mathbf{IV}}$, we need to give a closer look to \eqref{IVform} and \eqref{IVFG}.
Note that $F_{k,\ell}^{\mathbf{IV}}$ contains only terms of the type $G_{k,\ell-1}^{\mathbf{N}}$, $G_{k-3,\ell}^{\mathbf{N}}$ and    $(F_{k,\ell}^{\mathbf{N}})'$.
For $(k,\ell)$ such that $\frac k 3 +\ell \leq 2$, this provides terms $G_{k',\ell'}^{\mathbf{N}}$ for $\frac {k'}3+\ell\leq 1$.
Since $k'\geq 1$, this condition implies $\ell=0$ and $k'\leq 3=p-1$. 
But, we know from \eqref{closer}--\eqref{expM} that $G_{k',\ell'}^{\mathbf{M}}=0$ for such $(k',\ell')$.
Next, by \eqref{closer}, $F_{k',1}^{M}=0$ for $1\leq k'\leq 3$. Moreover, $F_{k',0}^{M}$ contains only product of $A_{k_j,0}$
for $k'\leq 6$. Indeed, if we look for example at a term of the form 
$
(Q_c^{k_1})'(Q_c^{k_2})'Q_c^{k_3}Q_c^{k_4}B_{1,0}B_{2,0}A_{3,0}A_{4,0},
$
by the formula of $(Q_c^2)'$, it gives a contribution only for $F_{k',\ell'}$, where $k'\geq 7$ or $k'\geq 4$ and $\ell\geq 1$.

Thus, by the assumptions on $A_{k',0}$$, F_{k',0}^{M}$ contains only constant polynomial functions  and so its derivative is in $\mathcal{Y}$.

\section{Appendix -- Identities related to $Q$}

\begin{claim}[Identities for any $p>1$]\label{LemmaA1}
    \begin{equation*}
        \int Q^{p+1} = \frac{2(p+1)}{p+3} \int Q^2, \quad       \int (Q')^2 = \frac{p-1}{p+3} \int Q^2.
    \end{equation*}
    \begin{equation*}
    		\int Q_c^2=c^{2q} \int Q^2,\quad E(Q_c)=c^{2q+1}E(Q)=-\frac {5-p}{2(p+3)} c^{2q+1} \int Q^2.
    \end{equation*}
\end{claim}

\noindent\emph{Proof of Lemma \ref{LemmaA1}.} These are well-known calculations. We have
$Q^p=Q-Q''$ and $\frac 2{p+1} Q^{p+1} = Q^2 -(Q')^2$. Thus, by integration:
$$
\int Q^{p+1} =\int Q^2 +\int (Q')^2,\quad
\frac 2{p+1} \int Q^{p+1} = \int Q^2 - \int (Q')^2.
$$
Therefore, $\int Q^{p+1} = \frac {2(p+1)}{p+3} \int Q^2$ and
  $\int (Q')^2=\int Q^{p+1}-\int Q^2= \frac {p-1}{p+3} \int Q^2$. Moreover,
$E(Q)= \frac 12 \int (Q')^2 -\frac 1{p+1} \int Q^{p+1}= \frac {p-5}{2(p+3)} \int Q^2$.

Since $Q_c(y)=c^{\frac 1{p-1}} Q(\sqrt{c} y)$ and $q= \frac 1{p-1} -\frac 14$, 
we have
$$\int Q_c^2(y) dy=c^{\frac 2{p-1}} \int Q^2(\sqrt{c} y) dy = c^{2q} \int Q^2.$$
Similary, $\int (Q_c')^2=c^{2q+1} \int (Q')^2$
and $\int Q_c^{p+1}= c^{2q+1} \int Q^{p+1}$, and so $E(Q_c)=c^{2q+1} E(Q)$.

\section{Appendix -- Proof of some technical results}
\subsection{Proof of Claim \ref{POSf}}
The proof is based on the  following well-known fact: 
\textit{There exists $\lambda_1>0$ such that
if $v\in H^1(\mathbb{R})$ satisfies $\int  Q v=\int xQ v=0$, then}
\begin{equation}\label{single}
\int v_x^2-p Q^{p-1} v^2 + v^2\ge \lambda_1 \|v\|_{H^1}^2.
\end{equation}
First, we claim from \eqref{single} that if   $\tilde v
\in H^1(\mathbb{R})$ satisfies $\int xQ v=0$, then
\begin{equation}\label{singleprime}
\int \tilde v_x^2-p Q^{p-1} \tilde v^2 + \tilde v^2\ge \lambda_0 \|\tilde v\|_{H^1}^2
-\frac 1{\lambda_0} \left(\int \tilde v Q\right)^2.
\end{equation}
Set $v=\tilde v-\frac {\int \tilde vQ}{\int Q^2}Q$. Then $\int  Q v=\int xQ v=0$ and from 
\eqref{single},  $\int v_x^2-p Q^{p-1} v^2 + v^2 \ge \lambda_1 \|v\|_{H^1}^2$.
Moreover, $\|v\|_{H^1}^2\geq \|\tilde v\|_{H^1}^2 - K (\int \tilde v Q)^2$ and
$\int v_x^2-p Q^{p-1} v^2 + v^2 \le \int \tilde v_x^2-p Q^{p-1} \tilde v^2 + \tilde v^2
+ K (\int \tilde vQ)^2$. Thus \eqref{singleprime} follows.

\medskip

Second, we recall
\begin{equation*}
    \mathcal{F}(t)=\frac 12 \int \left((\partial_x z)^2 + (1+\alpha'(y_c)) z^2)\right)
    -\frac 1{p+1} \int \left((v+z)^{p+1}- v^{p+1}-(p+1) v^p z\right).
\end{equation*}
Since $|\alpha'(s)|\leq K c^{\frac 1{p-1}}$, $\|v-Q\|_{L^\infty}\leq K c^{\frac 1{p-1}}$,
and $\|z\|_{L^\infty}\leq 2 K^* c^\theta$, we have from \eqref{singleprime} and $\int zQ'=0$,
for $c$ small enough,
\begin{equation*}
\begin{split}
    \mathcal{F}(t)&\geq \frac 12 \int \left((\partial_x z)^2 +   z^2 - pQ^{p-1}Êz^2\right)
    -K (c^{\frac 1{p-1}}+K^* c^\theta) \int z^2   \\ &
    \geq \frac {\lambda_0}4 \int \left((\partial_x z)^2 +   z^2 \right)
-\frac 1{2\lambda_0} \left(\int zQ\right)^2.
\end{split}
\end{equation*}

\subsection{Proof of Proposition \ref{2SOL}} 

1. For given $0<c<1$, $x_1$, $x_2\in \mathbb{R}$, the existence
of a solution $U(t)$ satisfying \eqref{2SOL0} is a consequence of Theorem 1 in \cite{Martel}.
Therefore, we only have to check \eqref{2SOL1}, for $c$ small,
 which is a more precise estimate than the one in
\cite{Martel}, giving explicitely the dependency in $c$. This   is obtained by combining
the argument of the proof in \cite{Martel} and estimates depending on $c$ in the proof of Proposition 
\ref{ASYMPTOTIC} of the present paper. 

We work on the time interval $(-\infty, -\tilde T_c]$, for $-\tilde T_c=\frac {x_2-x_1}{(1-c)} -\frac {T_c}{32}$.
Let $R(t,x)=Q(x-t-x_1)+Q_c(x-ct-x_2)$. In the spirit of Proposition 3 in \cite{Martel}, we first claim
the following.

\begin{proposition}\label{B52}
For $c>0$ small enough, 
	if there exists $t^*\leq - \tilde T_c$ such that
	$\forall t\leq t^*$, $\|u(t)-R(t\|_{H^1}\leq \exp(-c^{-\frac q4})$ then
	$\forall t\leq t^*$, $\|u(t)-R(t\|_{H^1}\leq K_0 e^{\frac {\sqrt{c}} 4 ((1-c)t-(x_2-x_1))}$.
\end{proposition}

Assume Proposition \ref{B52}. Since $\lim_{t\to -\infty} \|u(t)-R(t)\|_{H^1}=0$, we can define
$$
t^*=\sup\left\{ t\leq -\tilde T_c \text{ such that } \forall s\leq t^*,\
\|u(s)-R(s)\|_{H^1}\leq \exp(-c^{-r})\right\}.
$$
Since $K_0 e^{\frac {\sqrt{c}} 4 (-(1-c)\tilde T_c-(x_2-x_1))}
\leq K_0 e^{Ð\frac {\sqrt{c}} {128}T_c}\leq \frac 12 \exp(-c^{-r}), $
for $c$ small enough, by a standard continuity argument in $H^1$, we have $t^*=-\tilde T_c$,
and thus the result follows from Proposition \ref{B52} applied on $(-\infty,-\tilde T_c]$.
Therefore, we are reduced to prove Proposition \ref{B52}.

\medskip

\noindent\emph{Sketch of the proof of Proposition \ref{B52}.}
For more details, we refer to the proof of Proposition 3 in \cite{Martel}.
We decompose the solution $u(t)$ on $(-\infty,t^*]$ by Lemma  \ref{pourstab},
with $\alpha=0$ and $\rho_1(t)-\rho_2(t)\leq -\frac t 2 -\frac {T_c} {64}$.
Note that here the two solitons are ordered in a different way, $\rho_2(t)>\rho_1(t)$,
where $\rho_1(t)$ is center of $Q$ and $\rho_2(t)$ is center of $Q_c$.

Then, by \cite{MMas2},
 we have
$|c_1(t)-1|+|c_2(t)-c|\leq K g(t) + K \exp(\frac {\sqrt{c}} 4 ((1-c)t-(x_2-x_1)))$,
where $g(t)$ is defined as in \eqref{defg}. 

Next, similarly as   in \cite{MMas2},
 we use a monotonicity argument,
but since the solitons are ordered in reverse order, we will need the following quantities:
$$
\mathcal{M}(t)=\int u^2(t,x)\psi(x-m(t)) dx,\quad
\tilde{\mathcal{E}}(t)=\int \left(\tfrac 12 u_x^2-\tfrac 1{p+1} u^{p+1} + \tfrac c{100} u^2\right)\psi(x-m(t))dx,
$$
where $m(t)=\frac 12 (\rho_1(t)+\rho_2(t))$.
Similarly as in Lemma 1 of \cite{Martel}, we obtain, for $t'\leq t\leq t^*$,
$$
\mathcal{M}(t)-\mathcal{M}(t')\leq K \exp( \tfrac {1} 4 ((1-c)t-(x_2-x_1))),\quad 
\tilde{\mathcal{E}}(t)-\tilde{\mathcal{E}}(t')\leq K \exp( \tfrac {1} 4 ((1-c)t-(x_2-x_1))).
$$
We set
\begin{equation*}
\begin{split}
\mathcal{F}(t)&=\tfrac 12 \int u^2(t)+E(u(t)) + \left(\tfrac 1{2c} -\tfrac 12\right) \mathcal{M}(t)
+\left(\tfrac 1{c^2} - 1 \right) \int  \left(\tfrac 12 u_x^2-\tfrac 1{p+1} u^{p+1} \right)\psi(x-m(t))
\\ &= \frac 12 \int u^2(t)+E(u(t)) + \left(\tfrac 1{c^2} - 1 \right)\tilde{\mathcal{E}}(t)
+\tfrac 12 \left(\tfrac 1c -1\right) (1-\tfrac 1{100} (\tfrac 1c+1)) \mathcal{M}(t).
\end{split}
\end{equation*}
By the monotonicity results on $\mathcal{M}$ and $\tilde{\mathcal{E}}$, we have for all  $t'\leq t\leq t^*$,
$$
\mathcal{F}(t)-\mathcal{F}(t')\leq K \exp( \tfrac {1} 4 ((1-c)t-(x_2-x_1))),
$$
and using an expansion of $\mathcal{F}(t)$ from \eqref{decompo}, and passing to the limit
$t'\to -\infty$, we obtain the conclusion of Proposition \ref{B52}.

\medskip

\noindent 2. Sharper uniqueness property.

First, we check that for $c$   small enough,   if the solution
$u(t)$ satisfies \eqref{2SOL2} then for $-t$ large, 
$\rho_1(t)-\rho_2(t)\leq -\frac 14 |t|$.
This is a consequence of the asymptotic stability of one soliton.
Indeed, if $c$ is small enough, then for $-t$ large, $u(t)=Q(x-x_1) + \varepsilon(t,x)$,
 and $\varepsilon(t)$ small in $H^1$. Then, by 
stability and asymptotic stability of the soliton (see \ref{MMnonlinearity}),
there exists $\lambda$ such that  $|\lambda-1| \leq \frac 14$ and $\rho(t)$ with $\frac 32 t < \rho(t)< \frac t2$ 
for $-t$ large such that
$\|u(t)-Q_\lambda(x-\rho(t))\|_{H^1(x< t/10)}\to 0$ as $t\to -\infty$.
Thus,
$\|Q(x-\rho_1(t))+Q_c(x-\rho_2(t))-Q_\lambda(x-\rho(t))\|_{H^1(x< t/10)}\to 0$ as $t\to -\infty$.
This clearly implies that $\lambda=1$ and $\rho_2(t)> t/10$ for $-t$ large, and thus
$\rho_1(t)-\rho_2(t)< t/4$.

Using $\rho_1(t)-\rho_2(t)\leq -\frac 14 |t|$, as before by monotonicity arguments,
we have
$\|u(t)-Q(x-\rho_1(t))-Q_c(x-\rho_2(t)\|_{H^1}\leq K \exp( \tfrac {1} 8 ((1-c)t-(x_2-x_1)))$
for $-t$ large.
Therefore,  for this solution $u(t)$, we obtain
$$
|\rho_1'(t)-1|+|\rho_2'(t)-c|\leq  K \exp( \tfrac {1} 8 ((1-c)t-(x_2-x_1))),
$$
which proves  the convergence of $\rho_1(t)-t$ and $\rho_2(t)-ct$ as $t\to -\infty$.
Thus, there exist $x_1$, $x_2$ such that \eqref{2SOL2} holds. We now apply the uniqueness result of \cite{Martel} to conclude.

\subsection{Proofs of   \eqref{MOMENT}--\eqref{23nov}}
\emph{Proof of   \eqref{MOMENT}.}
Consider the decomposition of $u(t,x)$ introduced in the proof of
Proposition \ref{INTERACT4}, i.e. $u(t,x)=v(t,x-\rho(t))+z(t,x-\rho(t))$,
$v(t)$ and $S(t)=\partial_t  v+\partial_x ( \partial_x^2  v -v +v^p)$
satisfy the assumptions of Proposition \ref{INTERACT4}. Recall that  $z(t)$
satisfies equation \eqref{eqz}, and $\sup_{[0,T_c]} \|z(t)\|_{H^1}\leq K c^\theta$,
where $\theta$ is   to be fixed large  (for $\theta \geq \frac 58$).

First, we check that
\begin{equation}\label{MOMz}
\int _{x\geq 0} x^2 z^2(T_c ,x+\tfrac 12 \Delta)dx\leq K c^{2 \theta}
\end{equation}
implies the result.
By the explicit expression of $v(t,x)$ in \eqref{defvthINT},   the decay properties
of $Q$ and $Q_c$, and $\|\alpha\|_{L^\infty}\leq K c^{-\frac 16}$,
($\alpha$ is defined in Section 4.1)
we have the following pointwise estimates:
\begin{equation}\label{dec1}
\forall t\in [0,\tfrac 12 T_c],\ \forall x\geq \tfrac 18 T_c,\
|v(t,x)|+|v_x(t,x)|+|S(t,x)|\leq K\exp(-c^{-r})
 e^{-\frac 12\sqrt{c} x},
\end{equation}
\begin{equation}\label{dec2}
\forall t\in [\tfrac 12 T_c,T_c],\ \forall x\geq \tfrac 12 \Delta,\
|v(t,x)|+|v_x(t,x)|+|S(t,x)|\leq K(e^{-\frac 9{10} (x-\frac 12 \Delta)} + \exp(-c^{-r})
 e^{-\frac 12\sqrt{c} x}).
\end{equation}

By $  u^2(T_c,x ) \leq
2 ( v^2(T_c,x -\rho(T_c)) +   z^2(T_c,x-\rho(T_c)))$,
$|\rho(T_c)-T_c|\leq 1$, and \eqref{dec2} at $t=T_c$, we have
\begin{equation*}\begin{split}
&\int_{x\geq \frac {11}{12} |\ln c|} x^2 u^2(T_c,x+T_c+\tfrac 12 \Delta) dx
\\& \leq K \int_{x\geq \frac {11}{12} |\ln c|} x^2 e^{-\frac 9{5} x} dx
+ \exp(-\tfrac 12 c^{-r})
	+ \int_{x\geq 0} (x+1)^2z^2(T_c ,x+\tfrac 12 \Delta-1)dx
	\\& \leq K c^{\frac 54} + \int_{x\geq 0} x^2 z^2(T_c ,x+\tfrac 12 \Delta)dx.
\end{split}\end{equation*}

Second, we prove \eqref{MOMz},
which  will  finish the proof of   \eqref{MOMENT}.
This is proved by monotonicity arguments on $z(t)$. For $x_0>0$, $t\in [0,T_c]$, let
($\psi$ is defined in \eqref{surphiff})
$$ \mathcal{M}_z(t)=\int z^2(t,x)\psi(x_z) dx, \quad
\text{where} \quad x_z=x-\tfrac 12(T_c-t) -\tfrac 12 \Delta -x_0.$$
Using \eqref{eqz}, we have by direct calculations
\begin{equation*}\begin{split}
&\frac d{dt} \mathcal{M}_z(t) 
 = - 3 \int z_x^2 \psi'(x_z) + \int z^2 \psi'''(x_z) 
-\frac 12 \int z^2 \psi'(x_z) - (\rho'(t)-1) \int z^2 \psi'(x_z)\\
& + \int ((z+v)^4 - v^4 - z^4) (z \psi'(x_z)+ z_x \psi(x_z)) + \frac 4{5} \int z^{5} \psi'(x_z)
+ (\rho'-1) \int v_x z \psi(x_z).
\end{split}\end{equation*}
By \eqref{surphiff}, $\|z(t)\|_{H^1}\leq K c^{\theta}$ small, and then   
\eqref{estxxX}, we obtaiin
\begin{equation*}\begin{split}
\frac d{dt} \mathcal{M}_z(t) 
 & \leq K (\sup_{[0,T_c]} (\|z(t)\|_{L^2}^2 + \|S(t)\|_{L^2}^2)
 (\|v\psi'(x_z)\|_{L^\infty} + \|v_x\psi(x_z)\|_{L^\infty}+ \|v_x\psi(x_z)\|_{L^2}).
\end{split}\end{equation*}
Therefore, by the properties of $\psi$ and \eqref{dec1}-\eqref{dec2}, we obtain
\begin{equation*}
\frac d{dt} \mathcal{M}_z(t) 
 \leq K \exp(- c^{-r}) e^{-\frac 12 \sqrt{c} x_0}  + K c^{2 \theta} e^{-\frac 14 (x_0+
 \frac 12 (T_c-t))}. 
\end{equation*}
Thus, by integating in $t\in [0,T_c]$,
we obtain for all $x_0>0$, 
$\int_{x>x_0} z^2(T_c,x+\frac 12 \Delta) dx \leq 
K \exp(- c^{-r}) e^{-\frac 12 \sqrt{c} x_0}  + K c^{2 \theta} e^{-\frac 14 x_0}$. 
Thus $ \int_{x>0} x^2 z^2 (T_c,x+\frac 12 \Delta) dx \leq K c^{2 \theta}$
and \eqref{MOMz} follows.

\medskip

\emph{Proof of \eqref{23nov}.}
From \eqref{atc},  $|\rho(T_c)-T_c|\leq K c^2$
and $\|v\|_{H^2}\leq K$, 
we have 
$$\|u(T_c)-v(T_c,.- T_c)\|_{H^1}   \leq K c^2.$$
By \eqref{3-16}-\eqref{3-15}, we have
$$
\| v(T_c) - Q(.-\tfrac \Delta 2) - Q_c(.+(1-c)T_c -\Delta_c/2))\|_{H^1_c}
            \leq K c^{\frac {17}{12}},$$
            and thus by the decomposition of $u(T_c)$,
            and \eqref{poura},
            we deduce \eqref{23nov}.



\begin{thebibliography}{10}
 
 \bibitem{Be} T.B. Benjamin, The stability of solitary waves, Proc. Roy. Soc. London A \textbf{328}, (1972) 153--183.
\bibitem{Bo} J.L. Bona, On the stability theory of solitary waves, Proc. Roy. Soc. London A \textbf{349}, (1975) 363--374.
\bibitem{BPS}
J.L. Bona, W.G. Pritchard and L.R. Scott, Solitary-wave interaction, Phys. Fluids 23, \textbf{438}, (1980).
\bibitem{CL} T. Cazenave and P.L. Lions, Orbital stability of
standing waves for some nonlinear Schr\"odinger equations,
Comm. Math. Phys. \textbf{ 85}, (1982) 549--561.
\bibitem{Cohen} A. Cohen, Existence and regularity for solutions
 of the  Korteweg--de Vries equation,
 Arch.    Rat. Mech. Anal. {\bf 71} (1979), 143--175.
\bibitem{Cote1} R. C\^ote, Construction of solutions to the subcritical gKdV equations with a given asymptotic behaviour. Journal of Functional Analysis, to appear.
\bibitem{Cote2} R. C\^ote, Construction of solutions to the $L^2$-critical KdV equation with a given asymptotic behaviour, preprint.
\bibitem{Craig} W. Craig, P. Guyenne, J. Hammack, D. Henderson and C. Sulem, Solitary water wave interactions.  Phys. Fluids  \textbf{18},  (2006),  057106.
\bibitem{EcSc} W. Eckhaus and P. Schuur, The emergence of solutions of the Korteweg--de Vries equation from arbitrary initial conditions, Math. Meth. Appl. Sci., \textbf{5}, (1983) 97--116.
\bibitem{FPU} E. Fermi, J. Pasta and S. Ulam,
Studies of nonlinear problems, I, Los Alamos Report LA1940 (1955); reproduced
in  Nonlinear Wave Motion, A.C. Newell, ed., American Mathematical
Society, Providence, R. I., 1974, pp. 143--156.
\bibitem{HHGY} J. Hammack, D. Henderson, P. Guyenne and Ming Yi, 
Solitary-wave collisions, in \textit{Proceedings of the 23rd ASME Offshore Mechanics and Artic Engineering} (A symposium to honor Theodore Yao-Tsu Wu), Vancouver, Canada, June 2004 (Word Scientific, Singapore, 2004). 
\bibitem{HS}
M. H\u ar\u agu\c s-Courcelle  and D.H. Sattinger,  Inversion of the linearized Korteweg-de Vries equation at the multi-soliton solutions, Z. Angew. Math. Phys. \textbf{49} (1998), 436--469.
\bibitem{HIROTA} R. Hirota, 
Exact solution of the Korteweg-de Vries equation for multiple collisions of solitons, Phys. Rev. Lett., \textbf{27} (1971), 1192--1194.
\bibitem{KB} H. Kalisch and J.L.  Bona,
Models for internal waves in deep water, Discrete and Continuous  Dynamical Systems,
\textbf{6} (2000), 1--20.
\bibitem{KPV} C.E. Kenig, G. Ponce and L. Vega, Well-posedness and scattering results for the generalized Korteweg--de Vries equation via the contraction principle, Comm. Pure Appl. Math. \textbf{46}, (1993) 527--620. 
\bibitem{KRUSKAL} M. D. Kruskal, The Korteweg-de Vries equation and related evolution equations, in  Nonlinear Wave Motion, A.C. Newell, ed.,  American Mathematical
Society, Providence, R. I., 1974,  pp. 61--83. 
\bibitem{LAX1} P. D. Lax, Integrals of nonlinear equations of evolution and solitary waves, Comm. Pure Appl. Math. \textbf{21}, (1968) 467--490.
\bibitem{LS} Yi Li and D.H. Sattinger,  
Soliton collisions in the ion acoustic plasma equations. 
J. Math. Fluid Mech. \textbf{1}  (1999), 117--130. 
\bibitem{Martel} Y. Martel, Asymptotic $N$--soliton--like solutions of the subcritical and critical generalized Korteweg--de Vries equations, Amer. J. Math.  \textbf{127} (2005), 1103-1140.
\bibitem{yvanSIAM} Y. Martel, Linear problems related to asymptotic stability of solitons
of the generalized KdV equations, SIAM J. Math. Anal. \textbf{38} (2006), 759--781.
\bibitem{MM1} Y. Martel and F. Merle, Asymptotic stability of solitons for subcritical generalized KdV equations, Arch. Ration. Mech. Anal. \textbf{157}, (2001) 219--254. 
\bibitem{MM2} Y. Martel and F. Merle,
Stability of blow up profile and lower bounds for blow up rate for the critical generalized KdV equation,
Ann. of Math. \textbf{155}, (2002) 235--280.
\bibitem{MMnonlinearity} Y. Martel and F. Merle, Asymptotic stability of solitons of the subcritical gKdV equations revisited. Nonlinearity \textbf{18} (2005), no. 1, 55--80.
\bibitem{MMas1} Y. Martel and F. Merle, Asymptotic stability of solitons of the gKdV equations with a general nonlinearity, preprint.
\bibitem{MMas2} Y. Martel and F. Merle, Refined asymptotics around solitons for the gKdV equations with a general nonlinearity, preprint.
\bibitem{MMcol2} Y. Martel and F. Merle, Soliton collision for the nonintegrable gKdV
equations with general nonlinearity, preprint
\bibitem{MMb2} Y. Martel and F. Merle, Resolution of coupled linear systems related to the collision  of two solitons for the quartic gKdV equation, preprint.
\bibitem{MMT} Y. Martel, F. Merle and Tai-Peng Tsai,
Stability and asymptotic stability in the energy space of the sum of $N$ solitons
for the subcritical gKdV equations, Commun. Math. Phys. \textbf{231}, (2002) 347--373.
\bibitem{Miura} R.M. Miura, The Korteweg--de Vries equation: a survey of results, SIAM Review \textbf{18}, (1976) 412--459.
\bibitem{Mizu} T. Mizumachi,   Weak interaction between solitary waves of the generalized KdV equations, SIAM J. Math. Anal. \textbf{35} (2003),  1042--1080.
\bibitem{P} G.S. Perelman, Asymptotic stability of multi-soliton solutions for nonlinear Schršdinger equations.  Comm. Partial Differential Equations \textbf{29}, (2004)  1051--1095.
\bibitem{RSS} I. Rodnianski, W. Schlag, A.D. Soffer, Asymptotic stability of $N$-soliton states of NLS, to appear in Comm. Pure. Appl. Math. 
\bibitem{Schuur} P. C. Schuur, Asymptotic analysis of solitons problems, Lecture Notes in Math. {\bf 1232} (1986), Springer-Verlag, Berlin.
\bibitem{SHIH} L.Y. Shih, Soliton--like interaction governed by the generalized Korteweg-de Vries equation,
Wave motion \textbf{2} (1980), 197--206.
\bibitem{Tao} T. Tao, Scattering for the quartic generalised Korteweg-de Vries equation, preprint 2006. 
\bibitem{WT} M. Wadati and M. Toda, The exact $N$--soliton solution of the Korteweg--de Vries equation, J. Phys. Soc. Japan \textbf{32}, (1972) 1403--1411. 
\bibitem{WM} P.D. Weidman and T. Maxworthy, Experiments on strong interactions between
solitary waves,  J. Fluids Mech. \textbf{85}, (1978) 417--431.
\bibitem{We1} M.I. Weinstein, Modulational stability of ground states of nonlinear Schr\"odinger equations, SIAM J. Math. Anal. \textbf{16}, (1985) 472--491. 
\bibitem{We2} M.I. Weinstein, Lyapunov stability of ground states of nonlinear dispersive evolution equations, Comm. Pure Appl. Math. \textbf{39}, (1986) 51--68.
\bibitem{Z} N.J. Zabusky, Solitons and energy transport in nonlinear lattices,
Computer Physics Communications, \textbf{5} (1973), 1--10.
\bibitem{KZ} N.J. Zabusky and M.D. Kruskal, Interaction of ``solitons'' in a collisionless plasma and recurrence of initial states, Phys. Rev. Lett.
\textbf{15} (1965), 240--243.
\end{thebibliography}
\end{document}